\documentclass[11pt,letterpaper,twoside,reqno]{amsart}

\makeatletter{}
\usepackage{subeqnar}
\usepackage{overpic}
\def\sech{\mathop{\rm sech}\nolimits}

\usepackage{fixltx2e}                       
\usepackage[usenames,dvipsnames]{xcolor}
\usepackage{fancyhdr}
\usepackage{amsmath,amsfonts,amsbsy,amsgen,amscd,mathrsfs,amssymb,amsthm}
\usepackage{subfig}
\usepackage{url}

\usepackage{mathtools}

\usepackage[font=small,margin=0.25in,labelfont={sc},labelsep={colon}]{caption}

\usepackage{tikz}
\usepackage{microtype}
\usepackage{enumitem}

\definecolor{dark-gray}{gray}{0.3}
\definecolor{dkgray}{rgb}{.4,.4,.4}
\definecolor{dkblue}{rgb}{0,0,.5}
\definecolor{medblue}{rgb}{0,0,.75}
\definecolor{rust}{rgb}{0.5,0.1,0.1}

\usepackage[colorlinks=true]{hyperref}

\hypersetup{urlcolor=rust}
\hypersetup{citecolor=dkblue}
\hypersetup{linkcolor=dkblue}

\usepackage{setspace}

\usepackage{graphicx}
\usepackage{booktabs,longtable,tabu} 
\usepackage{multirow} 
\usepackage{float}

\usepackage[full]{textcomp}

\usepackage{fourier}
\usepackage{bm} 
\graphicspath{{figures/}}

\theoremstyle{definition}

\floatstyle{plaintop}
\newfloat{recipe}{thp}{lor}
\floatname{recipe}{Recipe}
\numberwithin{recipe}{section}

\renewcommand{\phi}{\varphi}

\title{Advanced Differential Equations:  Asymptotics \& Perturbations}

\author[J. Nathan Kutz]{J.~Nathan~Kutz\\{\em Department of Applied Mathematics, University of Washington, Seattle, WA 98195}}
\thanks{Email: kutz@uw.edu, \url{faculty.washington.edu/kutz}. }

\begin{document}


\begin{abstract}
 Approximation techniques have been historically important for solving differential equations, both as initial value problems and boundary value problems.  The integration of numerical, analytic and perturbation methods and techniques can help produce meaningful approximate solutions for many modern problems in the engineering and physical sciences.   An overview of such methods is given here, focusing on the use of perturbation techniques for revealing many key properties and behaviors exhibited in practice across diverse scientific disciplines. 
\end{abstract}

\maketitle

\tableofcontents

\newpage

\section*{Lecture 1:  Phase-Plane Analysis for Nonlinear Dynamics}

We begin this Chapter by reviewing the primary conclusions of
the last Chapter concerning eigenvalues and eigenvectors.  Thus
we consider the linear system:
\begin{equation}  {\bf x}' = {\bf A} {\bf x} \, . \end{equation}
The equilibrium points of this system are determined
by setting ${\bf x}'=0$.  This then yields
\[  {\bf A} {\bf x} =0 \,\,\,\,\, \rightarrow \,\,\,\,\,
          {\bf x} = 0 \]
since we want to assume that ${\bf A}$ has an inverse.  
The dynamics about the equilibrium point ${\bf x}=0$ can
then be found by solving for the eigenvalues and eigenvectors:
\[  {\bf x} = {\bf v} e^{\lambda t} \,\,\,\,\, \rightarrow \,\,\,\,\,
    ({\bf A} - \lambda {\bf I} ) {\bf v} = 0 \, . \]
Thus once the eigenvalues and eigenvectors are determined, the
phase-plane trajectories can be drawn and the solution behavior      
understood.  There are actually five distinct cases of behavior
which may arise.  They are as follows:\\

\noindent {\bf Case~1:}  eigenvalues -- real, unequal, same sign\\

In this case, once the eigenvalues and eigenvectors are found
we can express the solution as
\[
  {\bf x} = c_1 {\bf v}^{(1)} e^{\lambda_1 t} + c_2 {\bf v}^{(2)}
      e^{\lambda_2 t} 
\]
where we have assumed that $\lambda_1$ and $\lambda_2$ are
different and real.  The prototypical behavior of this case is
depicted in Fig.~11a.  In this figure, we assumed that
$\lambda_1 < \lambda_2 < 0$ so that decay occurs most rapidly
along the eigenvector ${\bf v}^{(1)}$.  This is also called
a {\em node} or {\em nodal sink}.  Note that all trajectories
go to zero so that the critical point is {\em stable}.  If
instead we found that $\lambda_1 > \lambda_2 >0$, then the
direction of the arrow in Fig.~11a would be reversed and
the critical point (called a {\em nodal source}) 
would be {\em unstable}.\\

\noindent {\bf Case~2:}  eigenvalues -- real, opposite sign\\

If the eigenvalues are of opposite sign and real we can
write the general solution
\[
  {\bf x} = c_1 {\bf v}^{(1)} e^{\lambda_1 t} + c_2 {\bf v}^{(2)}
      e^{-\lambda_2 t} 
\]
where we have assumed that $\lambda_1$ and $\lambda_2$ are
real and of the same sign.  The prototypical behavior of this case is
depicted in Fig.~11b.  In this figure, we assumed that
$\lambda_1, \lambda_2 > 0$ so that growth occurs
along the eigenvector ${\bf v}^{(1)}$ and decay along
${\bf v}^{(2)}$.  This is called a {\em saddle point}.  
Note that all trajectories eventually go out to infinity along
the eigenvector ${\bf v}^{(1)}$.  This implies that saddle
points are always unstable.  If
instead we found that $\lambda_1, \lambda_2 < 0$, then the
direction of the arrow in Fig.~11b would be reversed and
solutions would go unstable along ${\bf v}^{(2)}$.\\

\noindent {\bf Case~3:}  eigenvalues -- real and equal (double root)\\

For the case of a double root, two possibilities exist:  either we 
can find two linearly independent eigenvectors so that our solution is
\[
  {\bf x} = c_1 {\bf v}^{(1)} e^{\lambda_1 t} + c_2 {\bf v}^{(2)}
      e^{\lambda_1 t} 
\]
or, there is only one eigenvector, and we must generate a
generalized eigenvector via the methods of the last Chapter
so that our solution takes the form
\[
  {\bf x} = c_1 {\bf v}^{(1)} e^{\lambda_1 t} + c_2 \left[ 
    {\bf v}^{(1)} t  e^{\lambda_1 t} + \boldsymbol{\eta} e^{\lambda_1 t} 
     \right] \, . 
\]
In the first case, a {\em proper node} (or {\em star point}) 
is generated as depicted in Fig.~12a for $\lambda_1>0$.  The
second case generates an {\em improper node} which is depicted
in Fig.~12b, also for the case of $\lambda_1>0$.  Both points
are {\em stable} for $\lambda_1<0$, whereas for $\lambda_1>0$,
the arrows are reversed in Figs.~12a and 12b and the critical
point is {\em unstable}.\\

%
\begin{figure}[t]
\begin{center}
\begin{tabular}[t]{cc}
\includegraphics[width=2.4in]{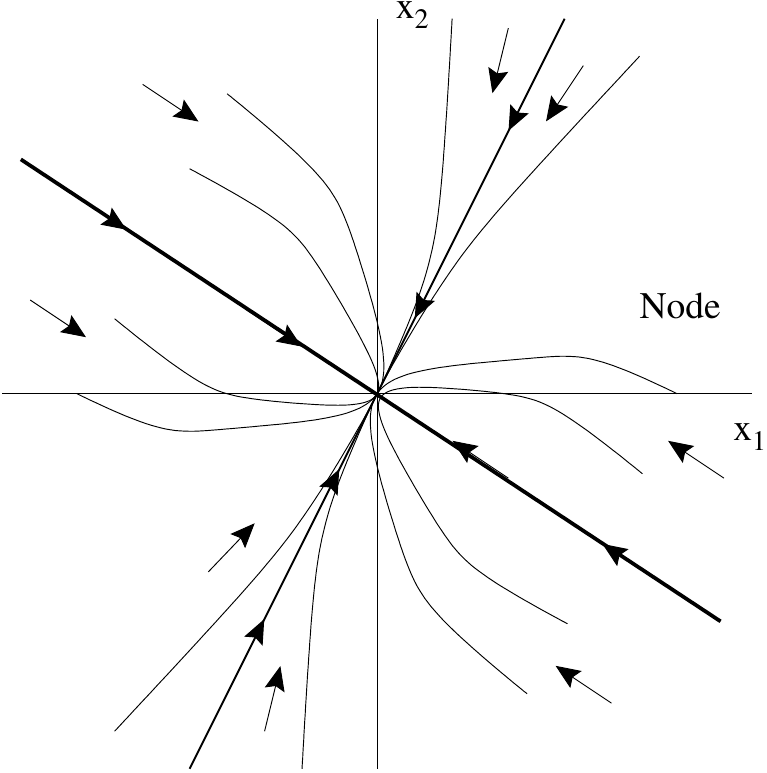}
\includegraphics[width=2.4in]{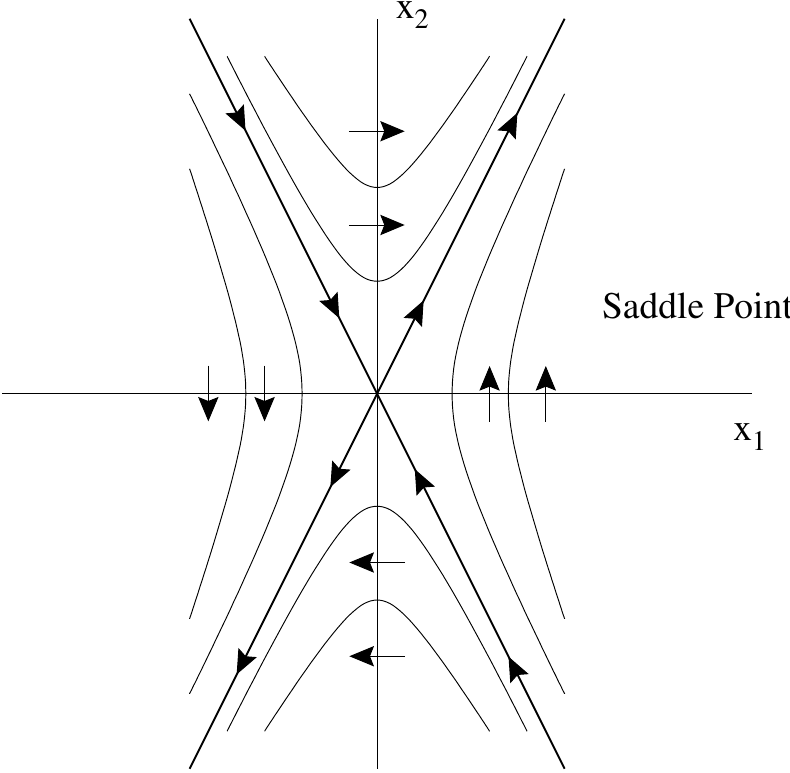}
\end{tabular}
\caption[]{Behavior of Node (a) and Saddle (b) which are
determined by having real distinct eigenvalues which are
of the same sign (a) or of opposite sign (b).  Changing the
signs of the eigenvalues simply changes the directions of
the arrows.  Thus the Node can be stable or unstable whereas
the Saddle is always unstable.}
\end{center}
\end{figure}
%

%
\begin{figure}[t]
\begin{center}
\begin{tabular}[t]{cc}
\includegraphics[width=2.4in]{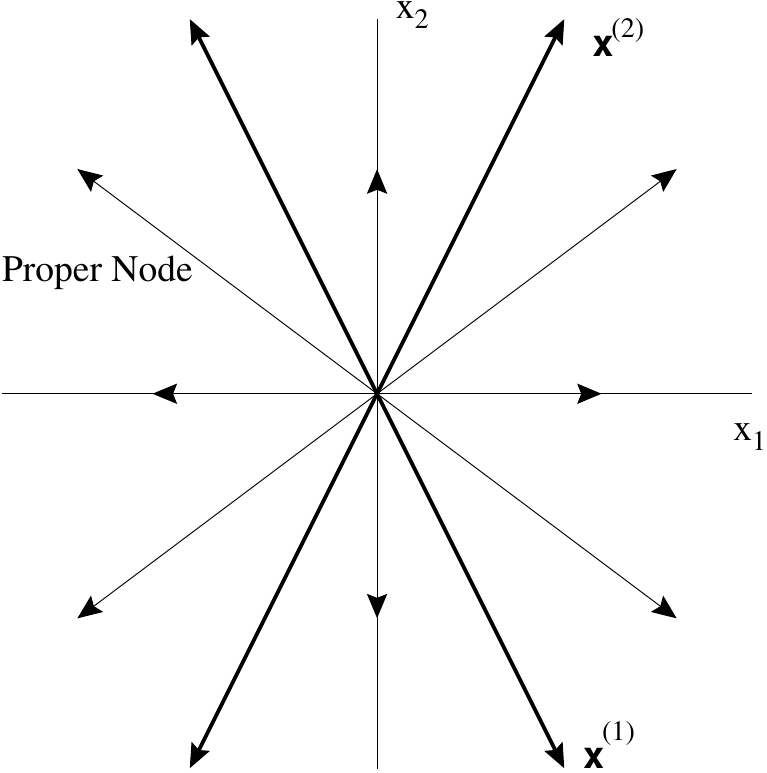}
\includegraphics[width=2.4in]{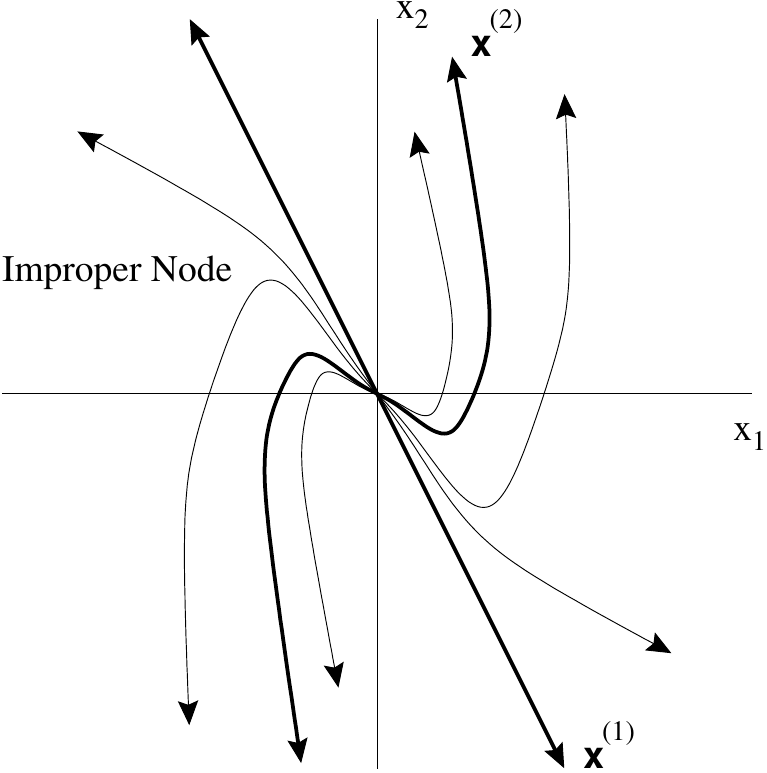}
%
\end{tabular}
\caption[]{Proper (a) and Improper (b) Nodes corresponding to
the double root case with two independent eigenvectors (a) or
one eigenvector and a second generalized eigenvector (b).  The
nodes can be stable or unstable depending on the sign of the
eigenvalue, i.e. the arrows are switched with a switch of the
sign of the eigenvalue.}
\end{center}
\end{figure}
%

%
\begin{figure}
\begin{center}
\begin{tabular}[t]{cc}
\includegraphics[width=2.4in]{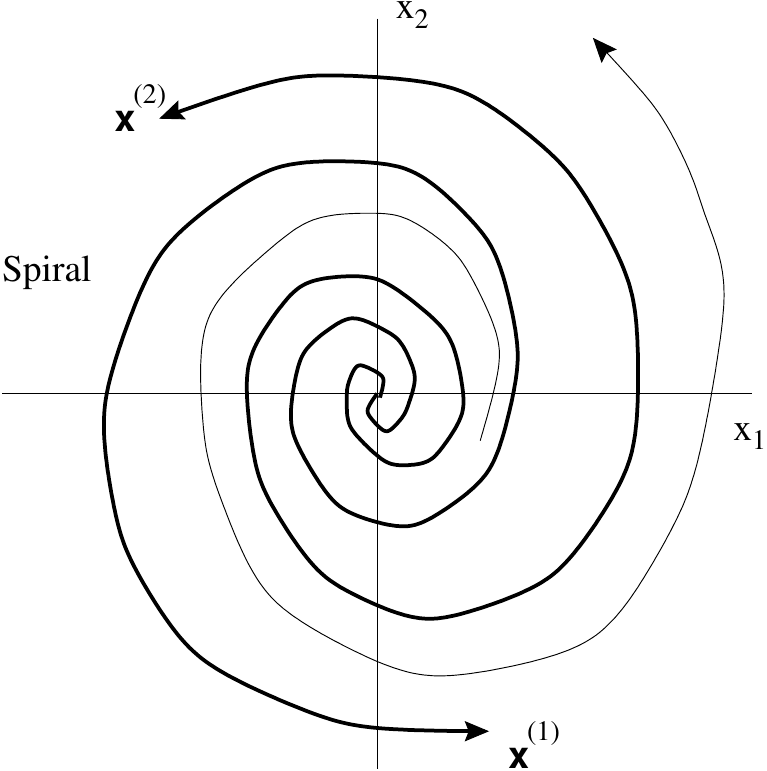}
\includegraphics[width=2.4in]{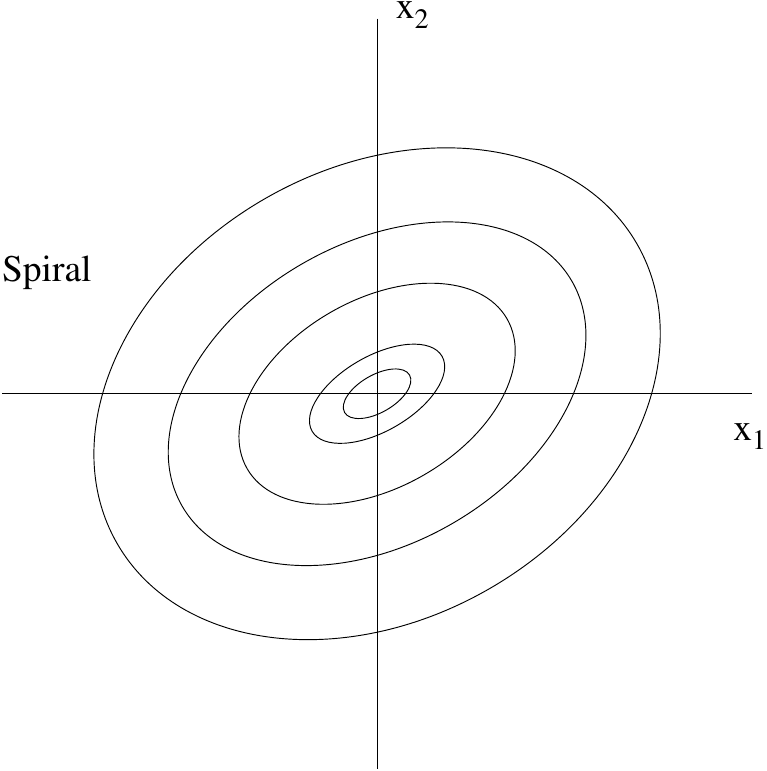}
%
\end{tabular}
\caption[]{Spiral (a) and Center (b) behavior when the eigenvalues
are complex conjugates.  The spiral has a nontrivial real part
which determines whether the trajectories spiral in or out in (a) 
whereas the center has strictly periodic behavior.}
\end{center}
\end{figure}
%

\noindent {\bf Case 4:}  eigenvalues -- complex eigenvalues\\  

For complex roots, we know that the eigenvalues and eigenvectors
come in complex pairs.  In particular, the eigenvalues
are given by
\[  \lambda_\pm = \beta \pm i \mu \, . \]
The resulting behavior is a {\em spiral} where the stability is
strictly determined from the real part $\beta$.  For $\beta>0$
the solutions spiral outward as in Fig.~13a so that the
equilibrium is unstable.  For $\beta<0$, the solutions spiral
inward (reverse the arrows in Fig.~13a) so that the the
equilibrium point is stable.\\

\noindent {\bf Case 5:}  eigenvalues -- purely imaginary\\  

In the case of purely imaginary eigenvalues:
\[  \lambda_\pm = \pm i\mu \]
our solutions are completely oscillatory and we have a {\em neutrally
stable} situation with a {\em center} critical point.  This
behavior is demonstrated in Fig.~13b where generically, the
solution trajectories are {\em ellipses}.  Thus solutions neither
grow or decay:  they simply display periodic motion.\\

Now that the behavior in linear systems is completely categorized,
let's move on to consider a classical problem which is 
actually {\em nonlinear}.  Such is the case of the pendulum.
The schematic of the pendulum is shown in Fig.~14.  By conservation
of angular momentum, we find that the swing of the pendulum
can be described by the equation:
\[  m L^2 \frac{d^2 \Theta}{dt^2} = - cL \frac{d\Theta}{dt} 
    - mg L \sin \Theta \]
where $m$ is the pendulum mass, $L$ is the length, $g$ is
the acceleration due to gravity, and $c$
measures the frictional/damping forces acting on the pendulum.
We can rewrite this equation as
\begin{equation}   \Theta'' + \gamma \Theta' + \omega^2 \sin \Theta = 0 \end{equation}
where $\gamma=c/mL$ and $\omega^2=g/L$.  

\begin{figure}
\centerline{\includegraphics[width=2.8in]{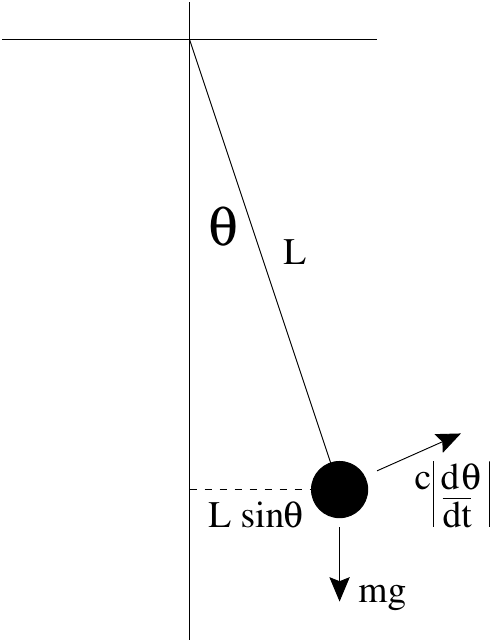}}
\caption{Schematic of pendulum oscillating from a fixed support
subject to forces of gravity and damping.}
\end{figure}

To convert this into a system of equations, we define
\[  x=\Theta \hspace*{.5in} \mbox{and} \hspace*{.5in} y=\frac{d\Theta}{dt} \]
which then results in the {\em nonlinear} system
\[  \left( \begin{array}{c} x' \\
     y' \end{array} \right)
     = \left( \begin{array}{c} y \\
     -\omega^2 \sin x - \gamma y \end{array} \right) 
\, . \]
Equilibrium solutions are found by letting $x'=y'=0$ which
yields
\[    y=0 \hspace*{.5in} \mbox{and} \hspace*{.5in} \sin x = 0 \]
so that the critical points are
\[  y=0 \hspace*{.5in} \mbox{and} \hspace*{.5in} x=\pm n\pi 
     \,\,\,\,\, n=0,1,2,... \, .\]
Thus unlike our linear systems above, we have more than one
equilibrium point.  In fact, we have an infinity of them
which lie at multiples of $\pi$ along the the $x$-axis.
We will learn how to deal with this in the next lecture.\\

\subsection{The Pendulum and Perturbation Theory}

In this section, we consider one of the classical problems of physics:
the pendulum.  Up to now, you have probably understood the pendulum
problem as a case of simple harmonic motion which was introduced
in physics.  However, you may recall that in this behavior was only
{\em approximate}, i.e. it relied on the pendulum swings being rather
small.  Here we consider the fully {\em nonlinear} pendulum dynamics.
Our equations of motion were given in the last lecture as:
\begin{eqnarray*}
 && x' =  y  \\
 && y' =  -\omega^2 \sin x - \gamma y  \, .
\end{eqnarray*}
To simplify the analysis, we begin by considering the case in
which there is no damping so that $\gamma=0$ and our governing
equations are
\begin{eqnarray*}
 && x' =  y  \\
 && y' =  -\omega^2 \sin x  \, .
\end{eqnarray*}
As in the last lecture, the key now is to find the critical points
and their stability.  The critical points (equilibrium) are
determined for $x'=y'=0$ so that
\[  y = 0 \hspace*{.5in} \mbox{and} \hspace*{.5in} x=\pm n\pi  
    \,\,\, n=0,1,2,... \]
which were given in the last lecture.  The idea now is to look
very close to one of the equilibrium points using the ideas
of {\em perturbation theory}.  Therefore we let
\begin{eqnarray*}
   && x = \pm n\pi + \tilde{x} \\
   && y = 0 + \tilde{y}
\end{eqnarray*}
where $\tilde{x}$ and $\tilde{y}$ are both very small.  Thus this
implies we are very near one of the fixed points.  Plugging this
into our governing equations for no damping yields the system:
\begin{eqnarray*}
 && \tilde{x}' =  \tilde{y} \\
 && \tilde{y}' =  -\omega^2 \sin (\pm n\pi + \tilde{x}) \, .
\end{eqnarray*}
To begin, we consider the well known example of the pendulum which
oscillates about the equilibrium $x=0$ (which corresponds to $\Theta=0$).
In this case $n=0$ so that we have
\[ \sin (\pm n\pi + \tilde{x}) =\sin \tilde{x} = \tilde{x} 
   - \frac{\tilde{x}^3}{3!} + \frac{\tilde{x}^5}{5!} + \cdots \approx \tilde{x}
\]
where we have used the Taylor series representation of sine and 
approximated everything by $\tilde{x}$ since all the other terms 
are much smaller (provided, of course, that $\tilde{x}$ is very small).
This is the standard trick that is used in introductory physics in order
to turn the nonlinear system into a linear one.  In particular, if we
plug this result into the above equation we preceding equation we
find.
\begin{eqnarray*}
 && \tilde{x}' =  \tilde{y} \\
 && \tilde{y}' =  -\omega^2 \tilde{x} \, .
\end{eqnarray*}
which results in the {\em linear} system:
\[  {\bf x}'
     = \left( \begin{array}{cc} 0 & 1 \\
     -\omega^2 & 0 \end{array} \right) {\bf x} \]
where ${\bf x} = (\tilde{x} \,\,\, \tilde{y})^T$.  We learned
how to solve this system in the last Chapter and in the first
lecture of this Chapter.  Thus we let ${\bf x} = {\bf v} e^{\lambda t}$
which yields the eigenvalue problem
\[  \left( \begin{array}{cc} -\lambda & 1 \\
     -\omega^2 & -\lambda \end{array} \right) {\bf v} =0 \, . \] 
The eigenvalues are found by taking the determinant of the above
matrix to be zero.  This then yields
\[  \lambda^2 + \omega^2 = 0 \,\,\,\,\, \rightarrow \,\,\,\,\,
      \lambda = \pm i \omega \]
which are purely imaginary eigenvalues.  Thus the equilibrium
point for $n=0$, i.e. $(x,y)=(0,0)$ is a center.  Thus solutions
near the critical point are all elliptic trajectories.

More generally, we can consider all the equilibrium points that
are multiples of $2\pi$ away from the origin $(x,y)=(0,0)$.
Thus we consider the perturbation theory for these points:
\begin{eqnarray*}
   && x = \pm 2 n\pi + \tilde{x} \\
   && y = 0 + \tilde{y}
\end{eqnarray*}
where $\tilde{x}$ and $\tilde{y}$ are both very small.  This
implies we are very near one of the fixed points located at
multiples of $2\pi$ from the origin.  Plugging this
into our governing equations for no damping now yields the system:
\begin{eqnarray*}
 && \tilde{x}' =  \tilde{y} \\
 && \tilde{y}' =  -\omega^2 \sin (\pm 2 n\pi + \tilde{x}) \, .
\end{eqnarray*}
But since
\[  \sin (\pm 2n\pi + \tilde{x}) = \sin \tilde{x} \approx \tilde{x} \, ,\]
we then arrive at exactly the same linearized equations as before.
Therefore, we can conclude that {\em all} the equilibrium points
in multiples of $2\pi$ from the origin are centers with periodic
solutions near each critical points (see Fig.~15).

\begin{figure}
\centerline{\includegraphics[width=5.0in]{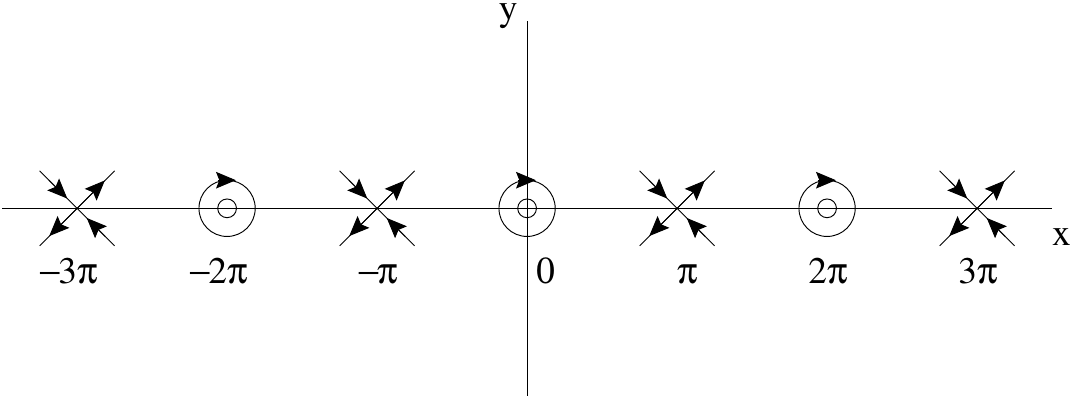}}
\caption{Behavior of solutions near each of the fixed
points which are multiples of $\pi$ from the origin.  Note
that multiples of $2\pi$ produce centers while multiples
of odd $\pi$ are saddles.}
\end{figure}

This is not the case for critical points which are odd multiples
of $\pi$ from the origin.  For these we can perturb around each
critical point by letting
\begin{eqnarray*}
   && x = \pm 2 n\pi + \pi + \tilde{x} \\
   && y = 0 + \tilde{y}
\end{eqnarray*}
where $\tilde{x}$ and $\tilde{y}$ are both very small.  This
implies we are very near one of the fixed points located at
odd multiples of $\pi$ from the origin (i.e. $\pm \pi, \pm 3\pi ,...$).  
Plugging this into our governing equations for no damping now yields 
the system:
\begin{eqnarray*}
 && \tilde{x}' =  \tilde{y} \\
 && \tilde{y}' =  -\omega^2 \sin (\pm 2 n\pi +\pi + \tilde{x}) \, .
\end{eqnarray*}
But since
\[  \sin (\pm 2n\pi + \pi + \tilde{x}) 
   = \sin(\pi+ \tilde{x}) = \sin\pi\cos\tilde{x} + \cos\pi \sin\tilde{x}
   = -\sin\tilde{x} \approx - \tilde{x} \, ,\]
we then arrive at a slightly different set of linearized equations.
In matrix form, this can be written as
\[  {\bf x}'
     = \left( \begin{array}{cc} 0 & 1 \\
     \omega^2 & 0 \end{array} \right) {\bf x} \]
where the only difference now is in the sign of $\omega^2$.  
Letting ${\bf x} = {\bf v} e^{\lambda t}$ yields the eigenvalue problem
\[  \left( \begin{array}{cc} -\lambda & 1 \\
     \omega^2 & -\lambda \end{array} \right) {\bf v} =0 \, . \] 
whose eigenvalues are 
\[  \lambda^2 - \omega^2 = 0 \,\,\,\,\, \rightarrow \,\,\,\,\,
      \lambda = \pm  \omega \]
which are purely real eigenvalues or opposite sign.  Thus the equilibrium
point for odd multiples of $\pi$ are saddles.  The eigenvectors can then
be found:
\begin{eqnarray*}
 &&  \lambda=\omega: \,\,\,\,\,  \left( \begin{array}{cc} -\omega & 1 \\
     \omega^2 & -\omega \end{array} \right) {\bf v} =0  
       \,\,\, \rightarrow \,\,\, -\omega v_1 + v_2 = 0 \,\,\,
       \rightarrow \,\,\, {\bf v}^{(1)} = \left( \begin{array}{c} 
        1\\ \omega \end{array} \right) \\
 &&  \lambda=-\omega: \,\,\,\,\,  \left( \begin{array}{cc} \omega & 1 \\
     \omega^2 & \omega \end{array} \right) {\bf v} =0  
       \,\,\, \rightarrow \,\,\, \omega v_1 + v_2 = 0 \,\,\,
       \rightarrow \,\,\, {\bf v}^{(2)} = \left( \begin{array}{c} 
        1\\ -\omega \end{array} \right) \, . 
\end{eqnarray*}
Thus a complete description of the saddle is given near each
critical point in odd multiples of $\pi$ from the origin.
The resulting dynamics is depicted in Fig.~15 which shows the
results of our perturbation calculations locally near each
of the fixed points.

\begin{figure}
\centerline{\includegraphics[width=5.0in]{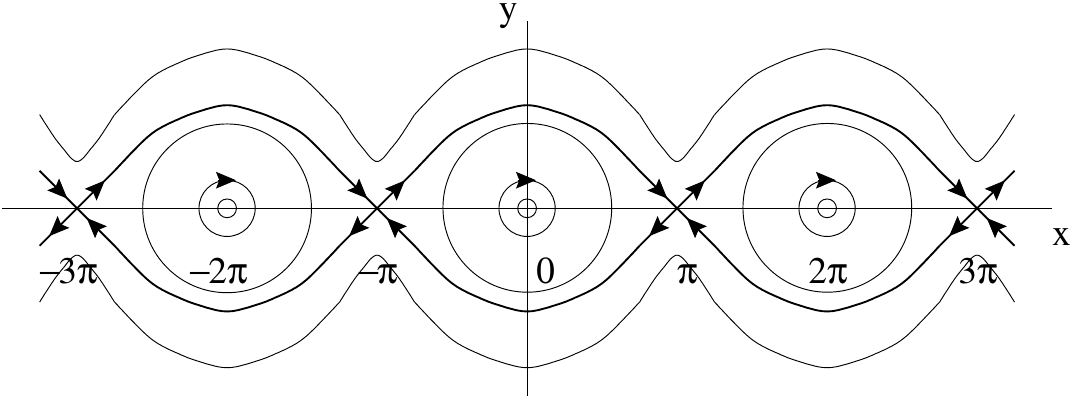}}
\caption{Full nonlinear behavior of solutions near each of the fixed
points and beyond.  This behavior is deducted directly from
Fig.~15 and is the only consistent behavior since all the critical
points are on the line $y=0$.  The bolded lines are the {\em separatrix}
which separate the oscillatory behavior from the behavior above
it corresponding to a pendulum which continuously rotates around.}
\end{figure}

Although the perturbation results are insightful, their
full power is not realized until we generalize our
thinking of Fig.~15.  In particular, since trajectories cannot
intersect, we can utilize the picture in Fig.~15 to develop
a full qualitative understanding of the dynamics.  Specifically,
we can describe the behavior far from the critical points
by their behavior near the critical points.  In Fig.~16, we
develop the full nonlinear qualitative behavior by simply taking
the phase-plane picture and generalizing it in the only way
possible.  This results in a dynamical picture which makes a
great deal of sense.  Note that near the critical points in
multiples of $2\pi$ (which corresponds to the rest position of
the pendulum), the behavior is exactly as expected:  oscillatory.
Whereas for odd $\pi$ values (which corresponds to a pendulum
sticking straight up), the behavior is given by a saddle and is
unstable.  Note that the trajectory separating the oscillatory
behavior from the trajectories above it is called the {\em separatrix}.
The separatrix projects along the unstable eigenvector of one saddle
into the stable eigenvector of a neighboring eigenvector.  The
behavior above this corresponds to the undamped pendulum
swinging around and around its support.  This is the case if we
give it a strong enough initial speed.  And since there is no
damping in this model, it will continue to circle around and around
and will never fall into the oscillatory back and forth motion
predicted near the center equilibrium.\\

\begin{figure}
\centerline{\includegraphics[width=5.0in]{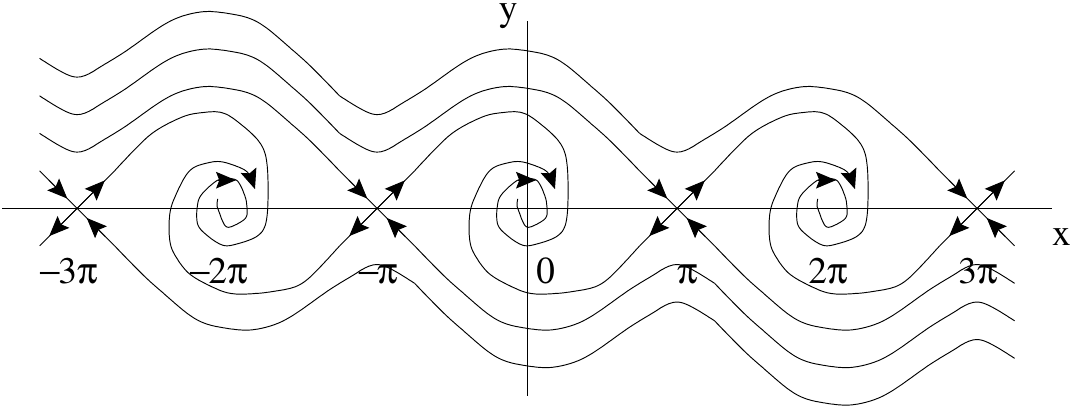}}
\caption{Full nonlinear behavior of solutions near each of the fixed
points when under-damping is applied to the pendulum.  
There are no separatrix in this case.}
\end{figure}

We now generalize our treatment in order to treat the case
of the damped pendulum.  Recall that the system in this case
is given by
\begin{eqnarray*}
 && x' =  y  \\
 && y' =  -\omega^2 \sin x - \gamma y  \, .
\end{eqnarray*}
As in the undamped case, the key now is to find the critical points
and their stability.  The critical points (equilibrium) are
determined for $x'=y'=0$ so that
\[  y = 0 \hspace*{.5in} \mbox{and} \hspace*{.5in} x=\pm n\pi  
    \,\,\, n=0,1,2,... \]
exactly as in the undamped case.  We once again perturb
about these equilibrium points to determine stability.
Therefore we let
\begin{eqnarray*}
   && x = \pm n\pi + \tilde{x} \\
   && y = 0 + \tilde{y}
\end{eqnarray*}
where $\tilde{x}$ and $\tilde{y}$ are both very small.  Plugging this
into our damped equations yields the system:
\begin{eqnarray*}
 && \tilde{x}' =  \tilde{y} \\
 && \tilde{y}' =  -\omega^2 \sin (\pm n\pi + \tilde{x}) - 
     \gamma \tilde{y} \, .
\end{eqnarray*}
As in the undamped pendulum case, there are two interesting cases
to consider.  The first is when the critical point is at the
origin or at multiples of $2\pi$ from it.  Thus we have
\[ \sin (\pm 2 n\pi + \tilde{x}) =\sin \tilde{x} \approx \tilde{x}
\]
where we have again approximated sine by $\tilde{x}$ since it is small.
Plugging this result into the linear damped equation above
results in the system:
\[  {\bf x}'
     = \left( \begin{array}{cc} 0 & 1 \\
     -\omega^2 & -\gamma \end{array} \right) {\bf x} \]
where ${\bf x} = (\tilde{x} \,\,\, \tilde{y})^T$.  
Letting ${\bf x} = {\bf v} e^{\lambda t}$ yields the eigenvalue problem
\[  \left( \begin{array}{cc} -\lambda & 1 \\
     -\omega^2 & -\gamma -\lambda \end{array} \right) {\bf v} =0 \, . \] 
The eigenvalues are found from the determinant to be
\[  -\lambda (-\gamma - \lambda) + \omega^2 =
  \lambda^2 + \gamma \lambda + \omega^2 = 0 \,\,\,\,\, \rightarrow \,\,\,\,\,
      \lambda = -\frac{\gamma}{2} \pm i \sqrt{\omega^2 - \frac{\gamma^2}{4}} 
    \, . \]
So depending on the quantity $\omega^2 - \gamma^2/4$, the equilibrium
is either a spiral ($\omega^2>\gamma^2/4$), an improper node 
($\omega^2=\gamma^2/4$), or a node ($\omega^2<\gamma^2/4$).  
The three different 
cases are referred to as underdamped, critically damped, and 
overdamped respectively.  In any
case, the real part of the eigenvalue is negative so that the
equilibrium point is asymptotically stable.  In what we depict in
Fig.~17, we assume that $\omega^2>\gamma^2/4$ so that the equilibrium
points at multiples of $2\pi$ are all spirals.

When the critical points are at odd multiples of $\pi$
from the origin, we once again have
\[ \sin (\pm 2 n\pi +\pi + \tilde{x}) =\sin(\pi+ \tilde{x})=-\sin \tilde{x} 
  \approx - \tilde{x}
\]
Plugging this result into the linear damped equation results in the system:
\[  {\bf x}'
     = \left( \begin{array}{cc} 0 & 1 \\
     \omega^2 & -\gamma \end{array} \right) {\bf x} \]
where ${\bf x} = (\tilde{x} \,\,\, \tilde{y})^T$.  
Letting ${\bf x} = {\bf v} e^{\lambda t}$ yields the eigenvalue problem
\[  \left( \begin{array}{cc} -\lambda & 1 \\
     \omega^2 & -\gamma -\lambda \end{array} \right) {\bf v} =0 \, . \] 
whose eigenvalues are found from the determinant to be
\[  -\lambda (-\gamma - \lambda) - \omega^2 =
  \lambda^2 + \gamma \lambda - \omega^2 = 0 \,\,\,\,\, \rightarrow \,\,\,\,\,
      \lambda = -\frac{\gamma}{2} \pm  \sqrt{\omega^2 + \frac{\gamma^2}{4}} 
    \, . \]
This yields two real eigenvalues which are of opposite sign.  Thus
a saddle is once again generated for all equilibrium points which
are odd multiples of $\pi$ from the origin.  The complete nonlinear
dynamics is depicted in Fig.~17 where the interaction of the spirals
and saddle nodes is shown.  Note how in this case, the solutions eventually
end up in one of the spiral points.\\

\newpage
\section*{Lecture 2:  Predator-Prey Models}

We now have enough background material to develop a
more general theory and understanding of nonlinear
systems.  There are two key concepts in nonlinear
systems that determine all the resulting dynamics.
The two primary issues are:\\

\begin{itemize}
\item Equilibrium ({\em Critical Points})
\item Stability \\
\end{itemize}

\noindent  Both of these concepts, which were mentioned
at the introduction of this chapter, are rather intuitive
in nature and have been illustrated in the previous two
lectures.

We begin by considering the following general system of
equations
\begin{subeqnarray}
 && x' =  F(x,y,t) \\
 && y' =  G(x,y,t) \, ,
\end{subeqnarray}
where $F(x,y,t)$ and $G(x,y,t)$ are some general functions
of $x$, $y$, and time $t$.  We will simplify this for the
present by considering the {\em autonomous} system:
\begin{eqnarray*}
 && x' =  F(x,y) \\
 && y' =  G(x,y) \, ,
\end{eqnarray*}
where $F$ and $G$ are not explicitly time dependent.  

We begin to analyze this system by considering
the concept of equilibrium.  Equilibrium occurs when
there is no ``motion'' in the system, i.e. when
both $x'=0$ and $y'=0$.  The point at which this
occurs is the equilibrium point $(x_0, y_0)$ which
satisfies:
\begin{eqnarray*}
 && F(x_0,y_0)=0 \\
 && G(x_0,y_0)=0 \, ,
\end{eqnarray*}
since $x'=y'=0$.  This is all there is to equilibrium.  We
simply find the points (there may be more than one) which
satisfy the above equations simultaneously.  Once this is
done, the behavior of the system can be determined entirely
from the stability of each equilibrium (critical) point.

The stability of each critical point may be determined by
looking very near each individual point.  Thus we
assume that
\begin{eqnarray*}
 && x=x_0 + \tilde{x} \\
 && y=y_0 + \tilde{y} \, ,
\end{eqnarray*}
where $\tilde{x}$ and $\tilde{y}$ are both very small
so that they can be considered to be in a very small
neighborhood of the critical point.  Plugging this
into our original equations gives us
\begin{eqnarray*}
 && {\tilde x}' =  F(x_0+\tilde{x},y_0+\tilde{y}) \\
 && {\tilde y}' =  G(x_0+\tilde{x},y_0+\tilde{y}) \, ,
\end{eqnarray*}
where we recall that since $x_0$ and $y_0$ are constants
then $x_0'=y_0'=0$.  The key now is to remember our Taylor
expansion formula from the series chapter.  Thus to
expand about some point, we have
\[  f(x_0 + \tilde{x}) = f(x_0) + \tilde{x} f'(x_0) + \frac{\tilde{x}^2}{2!}
       f''(x_0) + \frac{\tilde{x}^3}{3!} f'''(x_0) + \cdots \, . \]
Keeping only the first few terms in this approximation is good 
provided $\tilde{x}$ is small.

In the full problem, we now have to expand about both
$x_0$ and $y_0$.  Doing so yields the following:
\begin{eqnarray*}
  && \tilde{x}' =  F(x_0,y_0) + \tilde{x} F_x (x_0,y_0) + \tilde{y} 
        F_y (x_0,y_0) + \cdots \\
  && \tilde{y}' =  G(x_0,y_0) + \tilde{x} G_x (x_0,y_0) + \tilde{y} 
        G_y (x_0,y_0) + \cdots 
\end{eqnarray*}
where we have neglected all terms which are smaller than $\tilde{x}^2$,
$\tilde{y}^2$, and $\tilde{x}\tilde{y}$.  In matrix form, this
{\em linearized} system can be written as:
\[  \left( \begin{array}{c} \tilde{x} \\ \tilde{y} 
           \end{array} \right)' =
  \left( \begin{array}{cc} F_x (x_0,y_0) & F_y (x_0,y_0) \\
     G_x (x_0,y_0) & G_y (x_0,y_0) \end{array} \right) 
   \left( \begin{array}{c} \tilde{x} \\ \tilde{y} 
           \end{array} \right)      \, . \] 
We call this a {\em linearized} system since we turned
the original {\em nonlinear} system into a {\em linear} system
near the critical points.  The methods developed in the
previous chapter and reviewed in the introductory
lecture of this chapter are now applicable.  Thus we simply need
to determine the eigenvalues of the above system and the
resulting global, i.e. nonlinear, dynamics can be understood
qualitatively.

To make use of these ideas, we turn to some specific
examples to help illustrate the key ideas.  We begin by
considering what are called {\em predator--prey models}.
These models consider the interaction of two species:
predators and their prey.  It should be obvious that such
species will have significant impact on one another.
In particular, if there is an abundance of prey, then the
predator population will grow due to the surplus of
food.  Alternatively, if the prey population is low, then
the predators may die off due to starvation.  

To model the interaction between these species, we begin
by considering the predators and prey in the absence of
any interaction.  Thus the prey population (denoted by $x(t)$)
is governed by
\[  \frac{dx}{dt} = a x \]
where $a>0$ is a net growth constant.  The solution to this
simple differential equation is $x(t)= x(0) \exp (a t)$
so that the population grows without bound.  We have
assumed here that the food supply is essentially unlimited
for the prey so that the unlimited growth makes sense since
there is nothing to kill off the population.

Likewise, the predators can be modeled in the absence of
their prey.  In this case, the population (denoted by $y(t)$) 
is governed by
\[  \frac{dy}{dt} = -c y \]
where $c>0$ is a net decay constant.  The reason for the decay
is that the population basically starves off since there
is no food (prey) to eat.

We now try to model the interaction.  Essentially, the interaction
must account for the fact the the predators eat the prey.  
Such an interaction term can result in the following system:
\begin{eqnarray*}
  && \frac{dx}{dt} = a x - \alpha x y \\
  && \frac{dy}{dt} = -c x + \alpha x y
\end{eqnarray*}
where $\alpha>0$ is the interaction constant.  Note that $alpha$
acts as a decay to the prey population since the predators will
eat them, and as a growth term to the predators since they
now have a food supply.  These nonlinear and autonomous equations
are known as the {\em Lotka--Volterra equations}.

We rely on the methods introduced in this lecture to 
study this system.  In particular, we consider the equilibrium
points and their associated stability in order to determine
the qualitative dynamics of the Lotka-Volterra equations. 

The critical points are determined by setting $x'=y'=0$ which
gives 
\begin{eqnarray*}
  && ax-\alpha xy = x(a-\alpha y) = 0 \\
  && -cy+\alpha xy = y(\alpha x -c) = 0 \, .
\end{eqnarray*}
This gives two possible fixed points
\begin{eqnarray*}
 && \mbox{I}. \,\,\,  x=0 \,\,\, \mbox{and} \,\,\, y=0 \\
 && \mbox{II}. \,\,\, x=c/\alpha \,\,\, \mbox{and} \,\,\, y=a/\alpha \, .
\end{eqnarray*}
Each of these fixed points needs to be investigated separately in
order to determine the full (qualitative) nonlinear dynamics.

We begin with the critical point I. $(x_0,y_0)=(0,0)$.  Following
the methods outlined above we calculate the following:
\begin{eqnarray*}
 \begin{array}{lclcl}
   F(x,y) = ax - \alpha xy & \longrightarrow 
     &  F_x = a - \alpha y & \longrightarrow & F_x(0,0)=a  \\
     & &  F_y = -\alpha x &   & F_y(0,0)=0 \\[.2in]
   G(x,y) = -cy+\alpha xy  & \longrightarrow  
     &  G_x = \alpha y     & \longrightarrow & G_x(0,0)=0  \\
     & & G_y = -c +\alpha x & & G_y(0,0)=-c \, .
 \end{array}
\end{eqnarray*}
The resulting linearized system is
\[ {\bf w}' = 
   \left( \begin{array}{cc} a & 0 \\
     0 & -c \end{array} \right) 
             {\bf w}      =0 \] 
where ${\bf w}=(\tilde{x} \,\, \tilde{y})^{T}$.  As with
all previous linear systems, we make the substitution
${\bf w}={\bf v}\exp(\lambda t)$ in order to yield the
eigenvalue problem:
\[   
   \left( \begin{array}{cc} a-\lambda & 0 \\
     0 & -c-\lambda \end{array} \right) 
             {\bf v}      =0 \, . \] 
Setting the determinant to zero gives the characteristic
equation
\[   (a-\lambda)(-c-\lambda) =0 \]
whose eigenvalues are
\[  \lambda = a \,\,\, \mbox{and} \,\,\, \lambda=-c \]
Thus the eigenvalues are real and of opposite sign giving
us a saddle at the critical point $(0,0)$.

The eigenvectors can also be easily determined
for this case.  They are as follows:
\[ \lambda=a: \,\,\, 
     \left( \begin{array}{cc} a-a & 0 \\
           0 & -c-a \end{array} \right) 
             {\bf v}      =
     \left( \begin{array}{cc} 0 & 0 \\
           0 & -(c+a) \end{array} \right) 
             {\bf v}=0  \,\,\, \longrightarrow \,\,\,
         {\bf v}^{(1)} = \left( \begin{array}{c} 1 \\ 0 \end{array} 
     \right) \]
and
\[ \lambda=-c: \,\,\, 
     \left( \begin{array}{cc} a-c & 0 \\
           0 & -c+c \end{array} \right) 
             {\bf v}      =
     \left( \begin{array}{cc} a-c & 0 \\
           0 & 0 \end{array} \right) 
             {\bf v}=0  \,\,\, \longrightarrow \,\,\,
         {\bf v}^{(2)} = \left( \begin{array}{c} 0 \\ 1 \end{array} 
     \right) \]
The behavior near the critical point at the origin in thus
completely determined.  Figure~18a depicts the saddle behavior
near the origin.

\begin{figure}
\centerline{\includegraphics[width=5.0in]{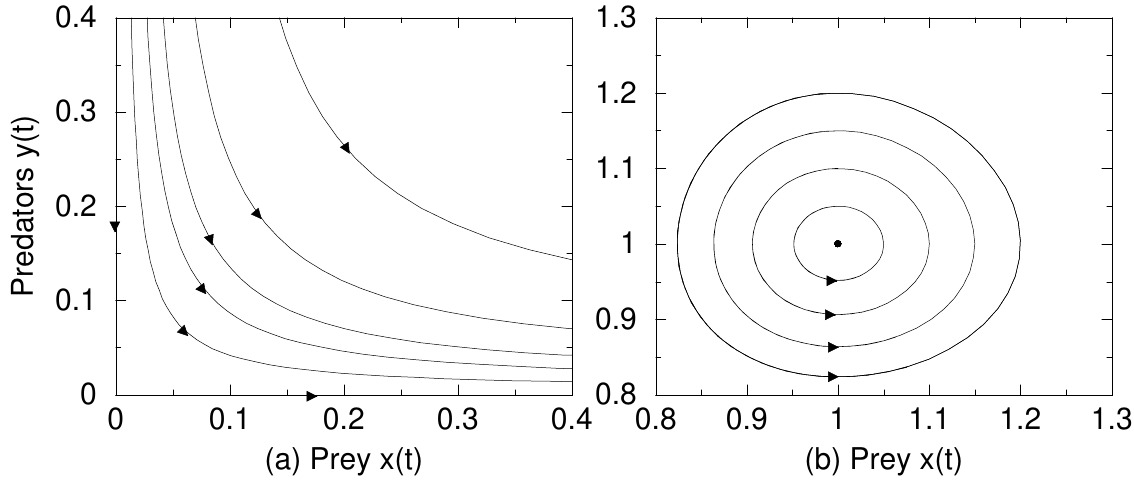}}
\caption{Behavior near each of the fixed points.  Here
     we have assumed that $a=c=\alpha=1$ so that the
     two critical points are the saddle at $(0,0)$ and 
     the center at $(1,1)$.}
\end{figure}

We now consider critical point II. $(x_0,y_0)=(c/\alpha,a/\alpha)$.  
Following the previous calculation we find:
\[  F_x(c/\alpha,a/\alpha)=0, \,\, F_y(c/\alpha,a/\alpha)=-c, \,\,
    G_x(c/\alpha,a/\alpha)=a, \,\, G_y(c/\alpha,a/\alpha)=0 \, . \]
The resulting linearized system is
\[ {\bf w}' = 
   \left( \begin{array}{cc} 0 & -c \\
     a & 0 \end{array} \right) 
             {\bf w}      =0 \, .\] 
Making the substitution
${\bf w}={\bf v}\exp(\lambda t)$ yields the eigenvalue problem:
\[   
   \left( \begin{array}{cc} -\lambda & -c \\
     a & -\lambda \end{array} \right) 
             {\bf v}      =0 \, . \] 
Setting the determinant to zero gives the characteristic
equation
\[   \lambda^2 + ac =0 \]
whose eigenvalues are
\[  \lambda_\pm = \pm i \sqrt{ac} \]
Thus the eigenvalues are imaginary giving a center
at the critical point $(c/\alpha,a/\alpha)$.  Before
calculating the eigenvectors for this case, we
note that the periodic behavior goes counter--clockwise in order
to be consistent with the flow of the critical point I. 
The behavior near the critical point is depicted in
Fig.~18b.  The full nonlinear dynamics is depicted in 
Fig.~19 which shows the saddle behavior near critical point I.
and periodic motion around the critical point II.

\begin{figure}
\centerline{\includegraphics[width=5.0in]{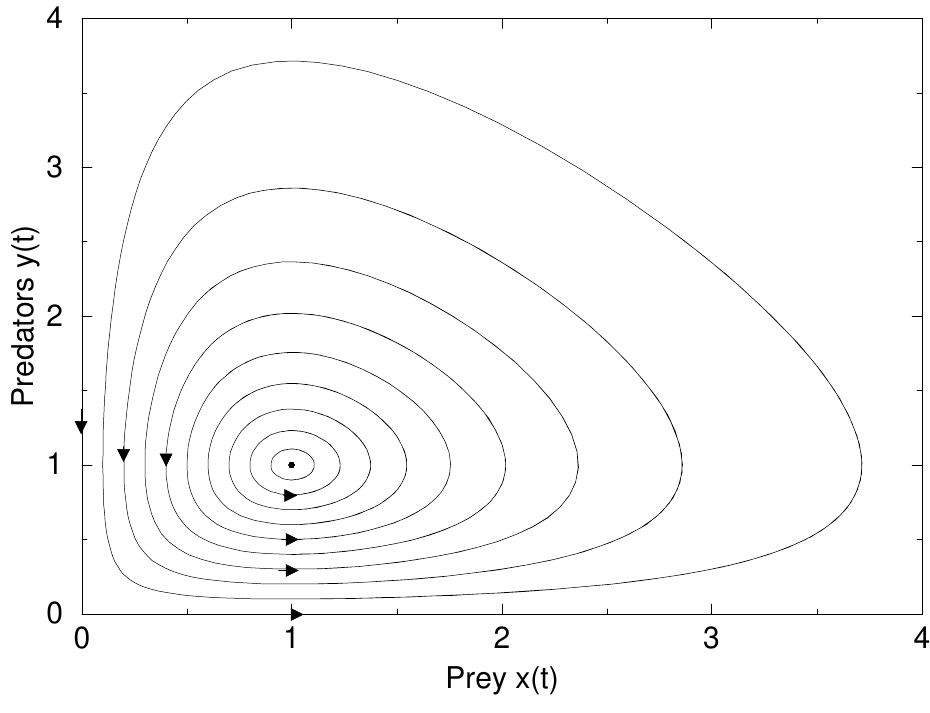}}
\caption{Fully nonlinear behavior of the predator--prey system.
     As with Fig.~18, we have assumed that $a=c=\alpha=1$ so that the
     two critical points are the saddle at $(0,0)$ and 
     the center at $(1,1)$.  Note the periodic behavior between
     these two fixed points.}
\end{figure}

To get a better idea of the periodic motion, we can calculate
the eigenvectors associated with critical point II.
\[ \lambda=i\sqrt{ac}: \,\, 
     \left( \begin{array}{cc} -i\sqrt{ac} & -c \\
           a & -i\sqrt{ac} \end{array} \right) 
             {\bf v}      = 0
        \, \longrightarrow \, -i\sqrt{ac} v_1 - cv_2 = 0
        \, \longrightarrow \,   
         {\bf v}^{(1)} = \left( \!\! 
     \begin{array}{c} \sqrt{ac} \\ -ia \end{array} 
    \!\! \right) . \]
Rewriting the eigenvector then gives
\[  {\bf w}^{(1)} = \left( \begin{array}{c} \sqrt{ac} \\ -ia \end{array} 
     \right) \exp \left( i\sqrt{ac} t \right) 
   = \left( \begin{array}{c} \sqrt{ac} \cos \sqrt{ac}t \\
       a \sin \sqrt{ac} t  \end{array} \right) + i
     \left( \begin{array}{c} \sqrt{ac} \sin \sqrt{ac}t \\
       -a \cos \sqrt{ac} t  \end{array} \right)
      \, .
\]
Rewriting our solution in terms of a purely real solution
then is easily done by combining the real and imaginary parts
to form
\[ {\bf w}  
   = c_1 \left( \begin{array}{c} \sqrt{ac} \cos \sqrt{ac}t \\
       a \sin \sqrt{ac} t  \end{array} \right) + c_2
     \left( \begin{array}{c} \sqrt{ac} \sin \sqrt{ac}t \\
       -a \cos \sqrt{ac} t  \end{array} \right)
      \, . \]
%
\begin{figure}
\centerline{\includegraphics[width=5.0in]{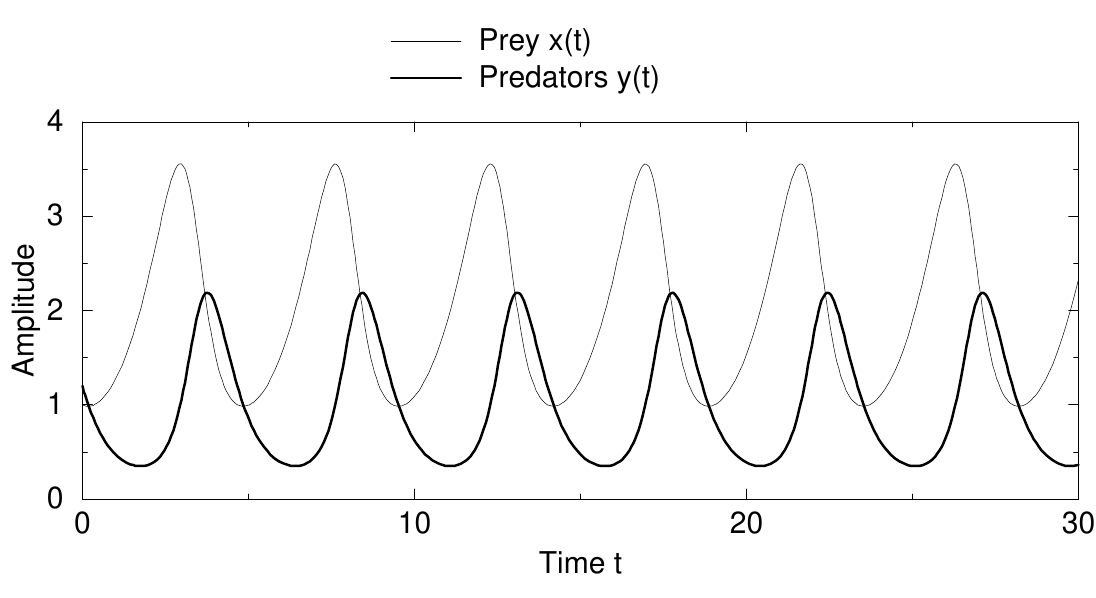}}
\caption{Fully nonlinear behavior of the predator--prey system
     as a function of time.  Here $a=\alpha=1$ and $c=2$.  Note
     the $\pi/2$ lag between the solutions.}
\end{figure}
%
Thus the population of predators ($y(t)$) and prey ($x(t)$)
can be calculated explicitly:
\begin{eqnarray*}
 && x(t) = \frac{c}{\alpha} + c_1 \sqrt{ac} \cos \sqrt{ac}t
                            + c_2 \sqrt{ac} \sin \sqrt{ac}t \\
 && y(t) = \frac{a}{\alpha} + c_1 a \sin \sqrt{ac}t
                            - c_2 a \cos \sqrt{ac}t \, .
\end{eqnarray*}
To simplify this further, we can replace the constants $c_1$
and $c_2$ by the two new constants $K$ and $\phi$ so that
\begin{eqnarray*}
 && x(t) = \frac{c}{\alpha} + \frac{c}{\alpha} K \cos (\sqrt{ac}t+\phi) \\
 && y(t) = \frac{a}{\alpha} + \frac{a}{\alpha} \, \sqrt{\frac{c}{a}}
       K\sin (\sqrt{ac}t+\phi) \, .
\end{eqnarray*}
This representation allows us to see explicitly the fact
that the populations are $\pi/2$ out of phase.  Thus the
reaction to changes of population occur (near the critical point)
one quarter of a cycle out of phase.  This behavior is 
demonstrated in Fig.~20.\\

\newpage
\section*{Lecture 3:  Linear Operators and their Adjoints}

Our goal in this chapter is to solve problems which take the form
\begin{equation}
   Lu=f
   \label{eq:luf_spec}
\end{equation}
where $L$ will be a linear, differential operator on the domain $x\in[0,l]$ and with boundary conditions specified at $x=0$ and $x=l$ which will be specified shortly.  Just like in the matrix version of this problem, i.e. ${\bf A}{\bf x}= {\bf b}$, a solution can be achieved by computing the inverse of the differential operator so that we can produce the solution
\begin{equation}
    u(x)=L^{-1}f .
\end{equation}
The spectral theory developed here will rely on the principle of linear superposition to construct the inverse operator $L^{-1}$ using eigenfunctions and eigenvalues of $L$.

Just as we solved ${\bf A}{\bf x}= {\bf b}$ with an eigenvector expansion, we can do the same for the continuous version of this (\ref{eq:luf_spec}).  To begin, we highlight the corresponding problem for each in Table~\ref{ta:spec}.  This shows that both in a vector and function space, there exists metrics around orthogonality and norms.   It also introduces for us the concept of an inner product in function space which is defined by
\begin{equation}
  \langle u,v \rangle =\int_0^l u v^* dx
\end{equation}
where integration occurs over the domain definied $x\in[0,l]$ and the $*$ denotes complex conjugation.  In function space, this is equivalent to a dot product, or projection, of one function onto another.  If the inner product is zero, then two functions are said to be orthogonal.  Thus just like vectors, there is a concept of orthogonal directions which will play a critical role in using functions as basis expansions to prescribe solutions.  We can also normalize a given function so that $\langle u, u \rangle=1$ which will allow us later to define an orthonormal basis.

Before developing solution techniques, the conditions under which (\ref{eq:luf_spec}) has solutions needs to be considered.
Recall  that there are solvability conditions which must be satisfied for linear systems such as 
${\bf A}{\bf x}= {\bf b}$.  In particular, we showed previously that the {\em Fredholm Alternative Theorem} determines when solutions can be found.  To review, consider the two equations:
\[   {\bf A} {\bf x} = {\bf b} \hspace*{.5in} \mbox{and}
      \hspace*{.5in} {\bf A}^* {\bf y}=0 \]
where we recall that ${\bf A}^*$ is the adjoint of ${\bf A}$.
Taking the inner product of the first equation with respect
to ${\bf y}$ gives:
\[  {\bf A}{\bf x} \cdot {\bf y}={\bf b} \cdot {\bf y} \, . \]
By noting that ${\bf A}{\bf x}\cdot{\bf y}={\bf x}\cdot{\bf A}^* {\bf y}$,
which is the definition of the adjoint,
and that ${\bf A}^* {\bf y}=0$ by definition, we find
\[  {\bf b} \cdot {\bf y}=0 \, .\]
Thus ${\bf y}$ must be orthogonal to ${\bf b}$ in order for the
equation to make sense. 

For $Lu=f$, the Fredholm Alternative theorem applies to function spaces.  In this case, we consider  (\ref{eq:luf_spec}) and the related homogeneous equation
\begin{equation}
   L^\dag v = 0
\end{equation}
where $L^\dag$ is the adjoint of the operator $L$ so that by definition $\langle v, Lu \rangle=\langle L^\dag v, u \rangle$.  In what follows, the construction of the adjoint will be explicitly considered.  For now, assume that $v(x)$ determines the null space of the adjoint operator.   The Fredholm Alternative theorem is established by taking the inner product of (\ref{eq:luf_spec}) on both sides with respect to the adjoint null space function $v(x)$.  Thus we have
\begin{eqnarray*}
 && \langle v, Lu \rangle = \langle v, f \rangle \\
 && \langle L^\dag v, u \rangle = \langle v, f \rangle
\end{eqnarray*}
where since $L^\dag v =0$ we find
\begin{equation}
 \langle v, f \rangle = 0.
\end{equation}
Thus just as in the vector case, the Fredholm Alternative theorem states that the forcing function $f(x)$ in (\ref{eq:luf_spec}) must be orthogonal to the null space of the adjoint operator in order for (\ref{eq:luf_spec}) to be solvable.

As with the vector case, there are important consequences to solvability and especially in regards to the null space. Specifically, the null space plays a critical role in the uniqueness of solutions.  Any null space solution can be simply thought of as eigenvector with eigenvalue zero.  When solving $Lu=f$, the null space solution can be added to the solution and it is still a solution.   For instance, consider the null space vector
\[ L u_0= 0  , \]
which is again called the kernel of the operator $L$.  Thus if we find a solution $\tilde{u}(x)$ to $Lu=f$, then  the solution
\[    u(x) = \tilde{u}(x) + c u_0 (x) \]
also satisfies $Lu=f$ for any value of the constant $c$.  Thus there are an infinite number of solutions.  As noted already, Noether's theorem was the first to highlight the role of zero eigenvalues and their relation to invariances in a linear system.  In this case, the simple zero eigenvalue gives a one-parameter family of solutions due to the null space (invariant) function $u_0(x)$.  Indeed, the null space functions correspond to critically important symmetries or constraints in many physical systems.  Thus null spaces should always be treated with special care as they reveal a great number of important properties of a system.

\begin{table}[t]
\begin{tabular}{lcc}
 &  \mbox{vector space} & \mbox{function space}    \\ \hline
 & &  \\
 \mbox{eigenvalue problem} & ${\bf A}{\bf x}=\lambda {\bf x}$ & $Lu=\lambda u$   \\
 \mbox{inner product \& orthogonality} & ${\bf x}\cdot {\bf y} =0$ &  $\langle u, v\rangle =0$  \\
 \mbox{norm} &  ${\bf x}\cdot {\bf x}$  &  $\langle u, u\rangle = \| u\|^2$
\end{tabular}
\caption{Comparison of the vector versus function space methods for spectral decompositions.  For function spaces, 
$Lu=f$ is equivalent to ${\bf A}{\bf x}= {\bf b}$. \label{ta:spec}}
\end{table}

\subsection*{Computing Adjoint Operators}
Adjoint operators play a critical role for understanding linear operators.  Adjoint operators can be computed for a given operator $L$ and its given boundary conditions.  To start, we consider the second order operator
\begin{equation}
   L = a(x) \frac{d^2}{dx^2} + b(x) \frac{d}{dx} + c(x)
\end{equation}
on the domain $x\in[a,b]$ with boundary conditions
\begin{subeqnarray}
&& \alpha_1 u(a) + \beta_1 \frac{du(a)}{dx} = 0 \\
&& \alpha_2 u(b) + \beta_2 \frac{du(b)}{dx} = 0
\end{subeqnarray}   
which are known as Robin boundary conditions.

By definition $\langle v, Lu \rangle=\langle L^\dag v, u \rangle$.  Thus our objective is to compute the adjoint $L^\dag$ and its associated  boundary conditions so that this definition holds.  We begin by noting 
\begin{eqnarray*}
   \langle v, Lu \rangle &=& \int_a^b  v Lu dx = \int_a^b  v \left( a(x) \frac{d^2u}{dx^2} + b(x) \frac{du}{dx} + c(x)u \right)  dx \\
   &=&  \int_a^b \left( a(x) v \frac {d^2u}{dx^2} + b(x) v \frac{du}{dx} + c(x)vu \right)  dx .
\end{eqnarray*}
Note that technically, we would need the complex conjugate on the second term in the inner product.  However, we will assume for the purposes of illustration that the functions of interest are real valued.
The three components of the integral can be transformed by using integration by parts.  For the first term, this is done twice.  For the second term, this is done once.  The following calculations are noted for the first two terms where use the notation that the subscript denotes differentiation
\begin{eqnarray*}
 && \int_a^b  a(x) v  u_{xx} dx = avu_x |_a^b - \int_a^b (av)_x  u_x dx 
   = avu_x |_a^b - (av)_x u |_a^b + \int_a^b (av)_{xx} u dx\\
 && \int_a^b  b(x) v u_x dx = b uv |_a^b - \int_a^b (bv)_x u dx .
\end{eqnarray*}
This then gives
\begin{eqnarray*}
   \langle v, Lu \rangle &=& \left( avu_x - (av)_x u + uvb \right) |_a^b + \int_a^b
     \left[  (av)_{xx} u - (bv)_x u + cv u \right] dx  \\
   &=& J(u,v) + \int_a^b
     \left[  (av)_{xx} - (bv)_x  + cv  \right] u dx \\
     &=& J(u,v) + \langle L^\dag v, u  \rangle .
\end{eqnarray*}
Thus the {\em formal adjoint} is given by
\begin{equation}
  L^\dag v = \frac{d^2}{dx^2} (a(x) v) - \frac{d}{dx} (b(x) v) + c(x) v  .
\end{equation}
The formal adjoint defines the operator structure of $L^\dag$ but does not include the boundary conditions.  The function $J(u,v)$, which is known as the {\em bilinear concomitant} or {\em conjunct}, contains the information on the boundary conditions of the adjoint operator.  Specifically, we have the important relationship
\begin{equation}
  \langle v, Lu \rangle - \langle L^\dag v,u \rangle = J(u,v) 
\end{equation}
which is known as {\em Green's formula}, or alternatively as {\em Lagrange's Identity}.  The adjoint problem is determined by setting $J(u,v)=0$.  Thus the formal adjoint with boundary conditions making the conjunct zero specifies the adjoint operator.  Setting the conjunct to zero gives
$\langle v, Lu \rangle = \langle L^\dag v,u \rangle$, which was used in our Fredholm Alternative theorem.  Note that if $L^\dag = L$, the operator is {\em formally self-adjoint}.  If in addition, the conjunct is zero so that $\langle v, Lu \rangle = \langle L v,u \rangle$, then the operator is self-adjoint.  Self-adjoint operators have many important properties.  They are equivalent to Hermitian or symmetric matrix operators.

\subsection*{Computing the formal adjoint and the adjoint}
As an example, consider the linear operator
\begin{equation}
  L = \frac{d^2}{dx^2} - \frac{d}{dx}
\end{equation}
on the domain $x\in[0,l]$ with boundary conditions
\begin{subeqnarray}
 && u(0) = 0 \\
 && u_x(l)= 0 .
\end{subeqnarray}
To compute the adjoint operator, we first construct the inner product
\begin{eqnarray*}
   \langle v, Lu \rangle &=& \int_0^l v (u_{xx} - u_x ) dx  \\
   &=& \left[ vu_x - v_x u -v u \right]_0^l + \int (v_{xx}+v_x ) u dx \\
     &=& \langle L^\dag v, u \rangle + \left[ v(l)u_x(l)-v_x(l)u(l) - v(l)u(l) -v(0)u_x(0) + v_x(0)u(0) +v(0) u(0) \right] \\
     &=&\langle L^\dag v, u \rangle - \left[ v_x(l)u(l) - v(l)u(l) - v(0)u_x(0) \right] .
\end{eqnarray*}
In order to remove the conjunct term, we must choose
\begin{subeqnarray}
&& v(0)=0 \\
&& v_x(l) + v(l) = 0.
\end{subeqnarray}
In total then, the formal adjoint is given by
\begin{equation}
   L^\dag =  \frac{d^2}{dx^2} + \frac{d}{dx}
\end{equation}
with the adjoint having the additional boundary conditions given by $v(0)=0$ and $v_x(l)+v(l)=0$.  The operator is not self-adjoint as $L\neq L^\dag$.  The boundary conditions are also different for the adjoint problem.

\newpage
\section*{Lecture 4: Eigenfunction Expansions}

With the concept of the adjoint operator in hand, we will now construct solutions to $Lu=f$ using eigenfunction expansions.  This gives a solution that computes the inverse operator so that $u=L^{-1} f$.  Before proceeding to the construction of eigenfunction expansion solutions, the relationship to such a solution technique for ${\bf A}{\bf x} = {\bf b}$ is considered.  In particular, consider the more familiar
\begin{equation}
  {\bf A} {\bf x}_n = \lambda_n {\bf x}_n \,\,\,\,\, \mbox{for} \,\, n=1, 2, \cdots, N
\end{equation}
where the square matrix  ${\bf A}\in\mathbb{R}^{N\times N}$ is assumed to be Hermitian (non-singular).  Thus the adjoint operator of ${\bf A}^*={\bf A}$.  The eigenvectors and eigenvalues of the matrix ${\bf A}$ can be used to construct a solution for ${\bf x}=b$.  In particular, the solution is assumed to be a sum of the eigenvectors
\begin{equation}
   {\bf x} = \sum_{n=1}^N c_n {\bf x}_n .
\end{equation}
For a non-singular matrix, the eigenvectors span the range of the matrix ${\bf A}$ and can be used as a basis for the $N$-dimensional solution vector ${\bf x}$.  If in addition, the matrix is Hermitian, then the eigenvectors are orthonormal so that 
\begin{equation}
  {\bf x}_n \cdot {\bf x}_m = \delta_{nm} = \left\{ \begin{array}{ll}  1 & n=m \\ 0 & n\neq m \end{array} \right.
\end{equation}
where $\delta_{nm}$ is the Kronecker delta function.

A solution for ${\bf A}{\bf x}={\bf b}$ can be constructed by expanding the solution in the orthonormal basis of the eigenvectors
\begin{equation}
  {\bf A} \sum_{n=1}^N c_n {\bf x}_n = {\bf b} .
\end{equation}
Only the weighting coefficients $c_n$ are unknown, and these are determined by inner products and orthogonality.  Specifically, one can multiply on the right by the eigenvector ${\bf x}_m$ to produce
\begin{eqnarray*}
   {\bf A} \sum_{n=1}^N c_n {\bf x}_n \cdot {\bf x}_m &=& {\bf b} \cdot {\bf x}_m \\
    \sum_{n=1}^N c_n {\bf A}{\bf x}_n \cdot {\bf x}_m &=& {\bf b} \cdot {\bf x}_m \\
    \sum_{n=1}^N c_n \lambda_n {\bf x}_n \cdot {\bf x}_m &=& {\bf b} \cdot {\bf x}_m \\
  c_m \lambda_m &=& {\bf b} \cdot {\bf x}_m \\
  c_m &=& {\bf b} \cdot {\bf x}_m/\lambda_m .
\end{eqnarray*}
The solution is then given by
\begin{equation}
 {\bf x}= \sum_{n=1}^N c_n {\bf x}_n =  \sum_{n=1}^N \frac{{\bf b} \cdot {\bf x}_n}{\lambda_n} {\bf x}_n
\end{equation}
which is a unique representation for a given right-hand side vector ${\bf b}$ satisfies the Fredholm Alternative theorem.

The steps taken for solving for ${\bf A}{\bf x}={\bf b}$ using an eigenvector expansion are identical with those required when solving $Lu=f$ with an eigenfunction expansion.  Specifically, we begin by considering the eigenvalue problem
\begin{equation}
   L u_n = \lambda_n u_n \,\,\,\,\, \mbox{for} \,\, n=1, 2, \cdots, \infty
\end{equation}
where again we assume a self-adjoint operator for convenience so that $L^\dag = L$.  The eigenfunctions for such an operator are orthonormal so that
\begin{equation}
    \langle u_n, u_m \rangle = \delta_{nm} 
\end{equation}
just as in the vector case.  But now, functions are orthogonal in an $\ell_2$ space of functions that are defined on $x\in[0,l]$.

As with the vector case, a solution can be constructed by the eigenfunction expansion 
\begin{equation}
   u(x) = \sum_1^\infty c_n u_n(x)
\end{equation}
where the $c_n$ will be determined by inner products and orthogonality.   A solution for $Lu=f$ can be constructed by expanding the solution in the orthonormal basis of the eigenfunctions
\begin{equation}
  L u = L  \sum_1^\infty c_n u_n = f.
\end{equation}
Only the weighting coefficients $c_n$ are unknown, and these are determined by taking the inner product on both sides with respect to the eigenfunction $u_m$:
\begin{eqnarray*}
   \langle L \sum_{n=1}^\infty c_n u_n , u_m \rangle  &=& \langle f, u_m \rangle \\
     \langle  \sum_{n=1}^\infty c_n L u_n , u_m \rangle  &=& \langle f, u_m \rangle \\
   \langle  \sum_{n=1}^\infty c_n \lambda_n u_n , u_m \rangle &=& \langle f, u_m \rangle \\
   \sum_{n=1}^\infty c_n \lambda_n \langle    u_n , u_m \rangle &=& \langle f, u_m \rangle \\
  c_m \lambda_m &=&  \langle f, u_m \rangle \\
  c_m &=&  \langle f, u_m \rangle/\lambda_m .
\end{eqnarray*}
The solution is then given by
\begin{equation}
 u (x) = \sum_{n=1}^\infty c_n u_n = \sum_{n=1}^\infty  \frac{\langle f, u_n \rangle}{\lambda_n} u_n
 \label{eq:eigenexpansion}
\end{equation}
which is a unique representation for a given right-hand side function $f(x)$ satisfies the Fredholm Alternative theorem.
This simple calculation of taking inner products and exploiting orthogonality results in identical solution techniques for both vector spaces with matrix operators and function spaces with differential operators.  In fact, numerical discretization of the boundary value problem $Lu=f$ results directly in a finite difference approximation which is ${\bf A}{\bf x}={\bf b}$.

\subsection*{Bounded domains and distinct eigenfunctions and eigenvalues}
To demonstrate how eigenfunctions and eigenvalues are determined in practice, we consider a simple illustrative example.  Specifically, consider the boundary value problem
\begin{equation}
    Lu = -u_{xx} = f(x)
\end{equation}
on the domain $x\in[0,l]$ with the boundary conditions
\begin{subeqnarray}
 && u(0)=0 \\
 && u(l)=0 .
\end{subeqnarray}
The eigenvalue problem associated with this problem is given by
\begin{equation}
  - u_{xx} = \lambda u \,\,\,\,\, \rightarrow \,\,\,\,\, u_{xx} + \lambda u = 0
\end{equation}
which has solutions
\begin{equation}
    u(x) = c_1 \sin \lambda x + c_2 \cos \lambda x .
\end{equation}
Imposing the boundary conditions limits the solution form.  Specifically, the boundary condition 
at the left implies
\begin{equation}
   u(0) = c_2 = 0  
\end{equation}
so that the solution is given by
\begin{equation}
    u(x) = c_1 \sin \lambda x .
\end{equation}
The boundary condition on the right gives
\begin{equation}
   u(l) = c_1 \sin \lambda l =0  .
\end{equation}
Since a non-trivial solution is sought ($c_1\neq 0$), then we must have
\begin{equation}
   \sin \lambda l =0  .
\end{equation}
which give
\begin{equation}
   \lambda l = \pm n\pi  \,\,\,\,\,\, n=0, 1, 2, \cdots , \infty .
\end{equation}
This gives the infinite collection of eigenvalues
\begin{equation}
  \lambda_n = \pm \frac{n\pi}{l}
\end{equation}
which has the corresponding eigenfunction
\begin{equation}
  u_n = c_n \sin \left( \frac{n\pi x}{l} \right) .
\end{equation}
This simple boundary value problem shows that there are an infinite number of solutions that satisfy both the governing equation and its associated boundary conditions.  This is fundamentally different than initial value problems where a second order differential equation would have two linearly independent solutions.

How is the general solution to the original $Lu=f$ problem determined?  especially if there are an infinite set of basis functions to construct it with.  The eigenfunction expansion suggests that one simply represent the solution as a linear superposition of the eigenfunctions
\begin{equation}
  u(x) = \sum_{n=1}^\infty c_n \sin \left( \frac{n\pi x}{l} \right) 
\end{equation}
where the weighting of each eigenfunction $c_n$ depends on the specific function $f(x)$ forcing the operator.  Note that the indexing starts at $n=1$ since the negative index values give the same eigenfunctions, i.e. $\sin(-n\pi x/l) = -\sin(n\pi x/l)$.

Before proceeding to the solution, we need to verify that the eigenfunctions themselves are linearly independent so that such a linear superposition is valid.  To compute linear independence, we compute the Wronskian between two eigenfunctions $u_n$ and $u_m$ so that
\begin{eqnarray*}
  W[u_n, u_m] &=& u_n {u_m}_x - {u_n}_x u_m = \sin \left( \frac{n\pi x}{l} \right) \frac{m\pi}{l} \cos \left( \frac{n\pi x}{l} \right) - \frac{n\pi}{l} 
  \cos \left( \frac{n\pi x}{l} \right) \sin \left( \frac{n\pi x}{l} \right) \\
  &=& \frac{\pi}{l}  \left[   \frac{m}{2} \left( \sin \left[ (n+m) \frac{\pi x}{l} \right]  +
  \sin \left[ (n-m) \frac{\pi x}{l} \right] \right)- 
   \frac{n}{2} \left( \sin \left[ (n+m) \frac{\pi x}{l} \right]  -
  \sin \left[ (n-m) \frac{\pi x}{l} \right] \right)  \right] \\
  &=&  \frac{\pi}{2l} \left[ (m-n) \sin \left[ (n+m) \frac{\pi x}{l} \right] + (m+n)
  \left[ (n-m) \frac{\pi x}{l} \right]  \right]   .
\end{eqnarray*}
Thus as long as $n\neq m$, then $W \neq 0$ implying that the $u_n$ and $u_m$ are linearly independent. 

It can be further shown that the eigenfunctions are orthogonal by considering the inner product
\begin{eqnarray*}
 \langle u_n,u_m \rangle &=& \int_0^l \sin \left( \frac{n\pi x}{l} \right) \sin \left( \frac{m\pi x}{l} \right) dx \\
 &=& \frac{1}{2} \int_0^l  \left[ \cos \left[ (n-m) \frac{\pi x}{l} \right] - \cos \left[ (n+m) \frac{\pi x}{l} \right]   \right]  dx \\
 &=& \frac{l}{2\pi}  \left[  \frac{1}{n-m} \sin \left[ (n-m) \frac{\pi x}{l} \right] 
  - \frac{1}{n+m} \sin \left[ (n+m) \frac{\pi x}{l} \right] \right]_0^l \\
  &=&  0
\end{eqnarray*}
since $\sin[(n-m)\pi] = 0$.  This shows that the eigenfunctions are orthogonal for different values of $n$ and $m$.  The eigenfunctions can also be made to be orthonormal by making the norm of the function unity.  Note that
\begin{equation}
\langle u_n,u_n \rangle = \int_0^l \sin \left( \frac{n\pi x}{l} \right) \sin \left( \frac{n\pi x}{l} \right) dx
=\frac{l}{2} .
\end{equation}
Orthonormality requires that $\langle u_n,u_n \rangle =1$.  We can then redefine normalized eigenfunctions so that this is indeed the case.  Thus the normalized eigenfunctions are given by
\begin{equation}
  u_n =  \sqrt{ \frac{2}{l} } \sin \left( \frac{n\pi x}{l} \right) .
\end{equation}
For these eigenfunctions, the desired orthonormality properties hold
\begin{equation}
   \langle u_n,u_m \rangle = \delta_{nm} .
\end{equation}
The solution to the original $Lu=f$ problem can then be determined
by formula (\ref{eq:eigenexpansion}) with the normalized eigenfunctions determined here
and the eigenvalues $\lambda_n= n\pi/l$.

\newpage
\section*{Lecture 5:  Sturm-Liouville Theory}

One of the most important generalizations of boundary value problems is to the Sturm-Liouville formulation.  Many of the most celebrated physics models of the 19th and 20th century take this form, including those named after Schr\"odinger, Laguerre, Hermite, Legendre, Chebychev and Bessel, for instance.  The Sturm-Liouville problem is indeed ubiquitous in mathematical physics and the basis of much of our knowledge and understanding of many canonical models in physics and engineering.  From the heat equation to quantum mechanics and electrodynamics, understanding the eigenvalues and eigenfunctions of these systems provide highly interpretable solutions that can be used to construct arbitrary solutions.

We thus consider the following second-order boundary value problem
\begin{equation}
   {L}u = \mu r(x) u + f(x)  
\end{equation}
on the domain $x\in [a,b]$
with boundary conditions
\begin{subeqnarray}
 && \alpha_1 u(a) + \beta_1  \frac{\partial u(a)}{\partial x} = 0 \\
 && \alpha_2 u(b) + \beta_2  \frac{\partial u(b)}{\partial x} = 0
\end{subeqnarray}
and the operator ${L}$ taking the form
\begin{equation}
   {L}u=- \frac{\partial}{\partial x} \left[ p(x)  \frac{\partial u}{\partial x} \right] + q(x) u .
\end{equation}
where $p(x), q(x)$ and $r(x)>0$ on $x\in[a,b]$.
The Sturm-Liouville form determines many advantageous properties that can be exploited for solutions and their interpretation.

The eigenvalue problem associated with Sturm-Liouville takes the form
\begin{equation}
   L u_n = \lambda_n r(x) u_n
\end{equation}
where the function $r(x)$ is a weighting function that plays a role in establishing orthogonality for different eigenfunctions. Specifically, the inner product is now defined with respect to the weighting function $r(x)$ so that
\begin{equation}
 \langle u, v\rangle_r = \int_a^b r(x) u v^* dx 
\end{equation}
where the subscript is included to explicitly denote the weighting.
With this definition, the eigenfunctions have the following orthogonality properties
\begin{equation}
 \langle u_n, u_m\rangle_r = \delta_{nm} .
\end{equation}
The Sturm-Liouville problem admits pairs of eigenvalues and eigenfunctions such that
\begin{eqnarray*}
  && \mbox{eigenvalues:} \,\,\,\,\, \lambda_1 < \lambda_2 < \lambda_3 < \cdots < \lambda_n < \cdots \\
  && \mbox{eigenfunctions:} \,\,\,\,\, u_1,  u_2, u_3, \cdots  u_n, \cdots
\end{eqnarray*}
where each eigenfunction is normalized so that $\langle u_n, u_n \rangle=1$.  The Sturm-Liouville operator is self-adjoint which implies the following properties:  (i) eigenvalues and eigenfunctions are real, (ii) each eigenvalue is simple (no multiple roots), (iii) eigenfunctions are orthogonal with respect to $r(x)$, (iv) eigenvalues are ordered as above with $\lambda_n\rightarrow \infty$ as $n\rightarrow \infty$, and (v) eigenfunctions form a complete set so that any function can be represented as an appropriately weighted linear sum of the $u_n$.  These facts are left as an exercise for the reader.  Each one of them is proved by simply applying projection and orthogonality of the eigenfunctions.

Since eigenfunctions of the Sturm-Liouville operator form a complete set, any solution can be represented as the linear superposition
\begin{equation}
   u  = \sum_{n=1}^\infty c_n u_n
\end{equation}
where the weighting coefficients $c_n$ are determined by considering the specific forcing $f(x)$ and orthogonality.  In a similar fashion, the forcing function $f(x)$ can be expanded using the complete set of eigenfunctions.  But instead of $f(x)$, we consider the following
\begin{equation}
    \frac{f(x)}{r(x)} =  \sum_{n=1}^\infty b_n u_n .
\end{equation}
The values of $b_n$ can be computed by taking the inner product of both sides with respect to $u_m$:
\begin{eqnarray*}
  \langle f(x)/r(x), u_m \rangle_r &=& \langle \sum_{n=1}^\infty b_n u_n, u_m \rangle_r \\
  \int_a^b r(x) [ f(x)/r(x) ] u_m dx &=&  \int_a^b r(x)  \left( \sum_{n=1}^\infty b_n u_n \right) u_m dx \\
   \int_a^b f(x) u_m dx &=& b_m
\end{eqnarray*}
or more succinctly
\begin{equation}
  b_n = \langle f, u_n \rangle
\end{equation}
and no subscript is included since this is the standard $\ell_2$ norm without weighting by $r(x)$.

The eigenfunction expansions can be used to solve for a solution of the Sturm-Liouville problem.  This then gives the following
\begin{eqnarray*}
   {L}u &=& \mu r(x) u + f(x)  \\
  L \left( \sum_{n=1}^\infty c_n u_n \right) &=& \mu r(x) \left( \sum_{n=1}^\infty c_n u_n \right)
   + f(x) \\
    \left( \sum_{n=1}^\infty c_n L u_n \right) &=& \mu r(x) \left( \sum_{n=1}^\infty c_n u_n \right)
   + f(x) \\
   \left( \sum_{n=1}^\infty c_n \lambda_n r(x) u_n \right) &=& \mu r(x) \left( \sum_{n=1}^\infty c_n u_n \right)
   + f(x)
\end{eqnarray*}
Taking the inner product of both sides with respect to $u_m$ gives
\begin{eqnarray*}
  \sum_{n=1}^\infty c_n \lambda_n \langle r(x) u_n, u_m \rangle &=& 
  \mu \sum_{n=1}^\infty c_n \langle r(x) u_n, u_m \rangle  + \langle f, u_m \rangle \\
  \sum_{n=1}^\infty c_n \lambda_n \langle u_n, u_m \rangle_r &=& 
  \mu \sum_{n=1}^\infty c_n \langle u_n, u_m \rangle_r  + \langle f, u_m \rangle  \\
   c_m \lambda_m &=& \mu c_m + b_m .
\end{eqnarray*}
This then gives the following condition which must be satisfied for each $c_n$
\begin{equation}
   (\lambda_n - \mu ) c_n = b_n \,\,\,\,\,\, n=1, 2, 3, \cdots 
\end{equation}
where $b_n$ is determined from the projection of $f(x)$ onto each eigenfunction $u_n$.  Note further in this computation how the weighted norm $\langle \cdot, \cdot \rangle_r$ was used for orthogonality.

There are three cases of interest which are as follows:
\begin{eqnarray*}
&& \mbox{Case 1:} \,\, \mu\neq \lambda_n \,\,\,\,\,  c_n = \frac{ b_n}{\lambda_n - mu} \\
&& \mbox{Case 2:} \,\, \mu= \lambda_n, b_n \neq 0 \,\,\,\,\, \mbox{no solution} \\
&& \mbox{Case 3:} \,\, \mu= \lambda_n, b_n =0  \,\,\,\,\,  c_n \,\, \mbox{undetermined, no unique solution}
\end{eqnarray*}
The last two cases are simply a restatement of the Fredholm Alternative theorem and solvability.  In the last case, the presence of a non-trivial null space makes the solution not unique as any null space function can be added to the solution.  In contrast, the case in which there exits no solution arises when the forcing function is not in the range of the operator $L$.  The first case if the more typical case which produces the solution
\begin{equation}
  u(x) = \sum_{n=1}^\infty \frac{ \langle f, u_n \rangle}{\lambda_n - \mu} u_n(x) .
\end{equation}
This is the generic form of the solution for the Sturm-Liouville problem when the null space of the linear operator is the trivial solution.

\subsection*{Sturm-Liouville example problem}
To illustrate the solution technique of eigenfunction expansions, consider the specific Sturm-Liouville problem
\begin{equation}
  u_{xx} + 2 u = -x
\end{equation}
on the domain $x\in[0,1]$ with boundary conditions
\begin{subeqnarray}
&& u(0)= 0 \\
&& u(1)+u_x(1) = 0 .
\end{subeqnarray}
The problem is rewritten in Sturm-Liouville form as
\begin{equation}
  -u_{xx} = 2 u + x
\end{equation}
where for this example, we have the following:  $p(x)=1, q(x)=0, r(x)=1, \mu=2$ and $f(x)=x$.
We note that we could have re-arranged and instead taken $\mu=0$ and $q(x)=-2$, but it is more standard to take $q(x)>0$.

The eigenvalue problem associated with this Sturm-Liouville problem ($r(x)=1$) is
\begin{equation}
    -{u_n}_{xx} = \lambda_n u_n
\end{equation}
which has the general solution
\begin{equation}
  u_n = c_1 \sin \sqrt{\lambda_n} x + c_2 \cos \sqrt{\lambda_n}x .
\end{equation}
Satisfying the first boundary condition $u_n(0)=0$ gives $c_2=0$.  The solution is
then given by
\begin{equation}
  u_n = c_1 \sin \sqrt{\lambda_n} x  .
\end{equation}
Application of the second boundary condition $u_n(1)+{u_n}_x (1) = 0$ gives the transcendental
equation for the eigenvalues
\begin{equation}
  \sin \sqrt{\lambda_n} + \sqrt{\lambda_n} \cos \sqrt{\lambda_n} = 0 \,\,\,\,\, \rightarrow \,\,\,\,\, 
  \sqrt{\lambda_n} + \tan \sqrt{\lambda_n} =0.
\end{equation}
Eigenvalues are then typically determined by root finding algorithms.  The normalized eigenfunctions can then be determined by enforcing $\langle u_n, u_n\rangle = 1$.  This gives the eigenfunctions
\begin{equation}
   u_n = \frac{\sqrt{2}}{(1 + \cos^2 \sqrt{\lambda_n})^{1/2}} \sin \sqrt{\lambda_n} x \,\,\,\,\,\,\,\, n=1, 2, 3, \cdots
\end{equation}
where are an orthonormal set that are real-valued with corresponding eigenvalues that are real-valued.

The Sturm-Liouville solution requires the computation of $b_n$ which can be computed with the normalized eigenfunctions
\begin{equation}
   b_n = \langle f, u_n\rangle =  \int_0^1 \frac{\sqrt{2}}{(1 + \cos^2 \sqrt{\lambda_n})^{1/2}} x \sin \sqrt{\lambda_n} x = \frac{2\sqrt{2} \sin\sqrt{\lambda_n} }{\lambda_n (1 + \cos^2 \sqrt{\lambda_n})^{1/2} } .
\end{equation}
This then gives
\begin{equation}
   f(x) = 4 \sum_{n=1}^\infty \frac{ \sin \sqrt{\lambda_n}  \sin \sqrt{\lambda_n} x}{ \lambda_n (1 + \cos^2 \sqrt{\lambda_n})} 
\end{equation}
and the corresponding solution
\begin{equation}
  u(x) = 4 \sum_{n=1}^\infty \frac{ \sin \sqrt{\lambda_n}  \sin \sqrt{\lambda_n} x}{ \lambda_n (\lambda_n - 2) (1 + \cos^2 \sqrt{\lambda_n})} 
\end{equation}
with the distinct, real and simple eigenvalues determined by $\sqrt{\lambda_n} + \tan \sqrt{\lambda_n} =0$.  This represents a typical approach and strategy for computing the various projections onto the eigenfunction space in order to construct the solution.

\newpage
\section*{Lecture 6:  Green's Functions, Delta Functions and Distribution Theory}

The discussion of the Green's function requires that we first consider a class of functions that are not functions in any standard sense.  Specifically, we will want to consider what is called the {\em Dirac delta function}, which is often referred to as a generalized function.  The Dirac delta function has the property that is is zero everywhere, infinity at a single location, and when integrated over produces a value of unity.  Dirac introduced this concept in an effort to consider {\em impulse} responses on a given system.  In our treatment of Laplace transforms, we already introduced the delta function in order to understand how a system responds to a kick of the system.  The usefulness of the Dirac delta function led mathematicians to develop a generalized theory of functions which is known as {\em distribution theory}.  This will be considered briefly here with the goal of developing the fundamental solution for differential equations, what is known as the Green's function.

As with the last chapter, our goal is to consider problems of the form
\begin{equation}
   Lu=f
   \label{eq:greenluf}
\end{equation}
where $L$ will be a differential operator on the domain $x\in[0,l]$ and with boundary conditions specified at $x=0$ and $x=l$. Green's functions will allow us to compute the inverse of the differential operator so that we can produce the solution
\begin{equation}
    u(x)=L^{-1}f .
\end{equation}
The Green's function will rely on the principle of linear superposition to construct the inverse operator $L^{-1}$ and a fundamental solution for (\ref{eq:greenluf}).

To begin our consideration of the concepts of distribution theory, we consider the forcing $f(x)$ in 
(\ref{eq:greenluf}) that takes the form
\begin{equation}
   f(x) = \left\{  \begin{array}{cl} f_0(x) & x_0-\xi<x<x_0+\xi \\0 & \mbox{elsewhere} \end{array} \right.
\end{equation}
Thus the forcing term acts only in a limited spatial range $x\in[x_0-\xi,x_0+\xi]$ of the differential equation centered around $x_0$.

The {\em impulse} of the forcing $f(x)$ is defined as follows
\begin{equation}
   I(\xi) = \int_0^l f(x)dx = \int_{x_0-\xi}^{x_0+\xi} f_0(x) dx .
\end{equation}
The impulse $I(\xi)$ is then seen to be the total integrated value of the function $f(x)$ on the entire domain $x\in[0,l]$.  To consider a specific example of the forcing function, we consider
\begin{equation}
   f(x) = \left\{  \begin{array}{cl} {1}/({2\xi}) & x_0-\xi<x<x_0+\xi \\0 & \mbox{elsewhere} \end{array} \right. .
\label{eq:delta_sequence} 
\end{equation}
This shows the forcing function to be a rectangle of width $2\xi$ and height $1/(2\xi)$.  The impulse is then given by
\begin{equation}
   I(\xi) =  \int_{x_0-\xi}^{x_0+\xi} \frac{1}{2\xi}  dx =  \frac{1}{2\xi}  \int_{x_0-\xi}^{x_0+\xi}  dx  = 1.
\end{equation}
This simple example shows that the impulse is always unity.  

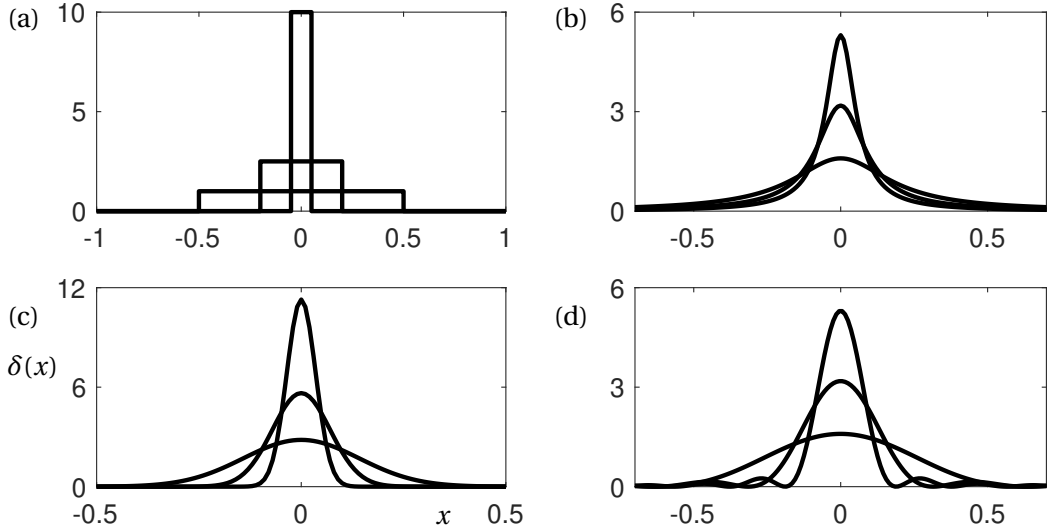
\begin{figure}[t]
\begin{overpic}[width=0.8\textwidth]{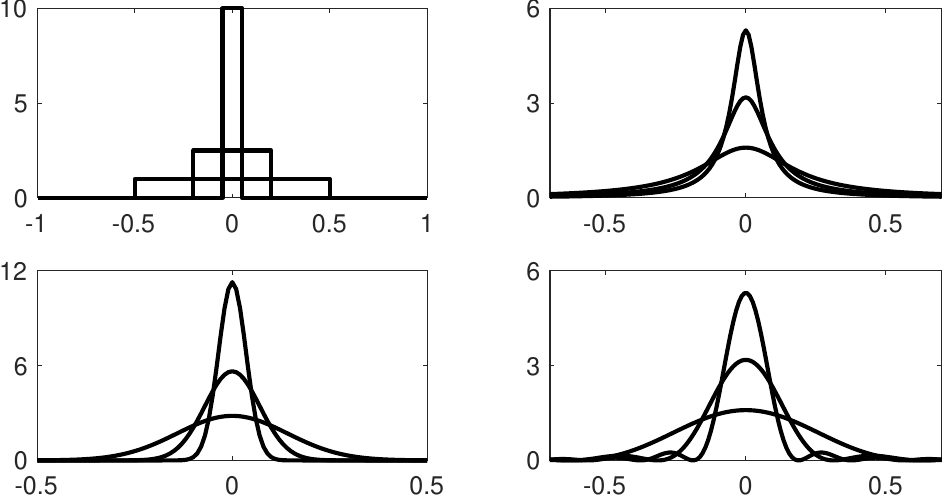}
\put(-5,50){(a)}
\put(50,50){(b)}
\put(-5,20){(c)}
\put(50,20){(d)}
\put(38,0){$x$}
\put(-5,15){$\delta(x)$}
\end{overpic}
\caption{Four potential Delta function sequences given by (\ref{eq:delta_sequence}) and (\ref{eq:delta_sequence2}).   (a)  The step function delta sequence (\ref{eq:delta_sequence})
for $\xi=0.05, 0.2, 0.5$, (b) the algebraic function (\ref{eq:delta_sequence2}a) with
$\xi=0.06, 0.1, 0.2$, (c) the Gaussian function (\ref{eq:delta_sequence2}b) with $\xi=0.05, 0.1, 0.2$, and (d) the sinc function (\ref{eq:delta_sequence2}d) with
$\xi=0.06, 0.1, 0.2$.  Each sequence collapses as $\xi \rightarrow 0$ to a function that is zero everywhere aside from at $x=0$ where it takes on the value of infinity in such as way that integrating of the delta sequence gives an impulse of unity.}
 \label{fig:delta_sequence}
\end{figure}

The example given is an example of how to construct a Dirac delta function.  In fact, there are quite a few {\em Dirac delta function sequences}, all of which produce the desired behavior of 
being zero everywhere, infinity at a single location, and of impulse unity.  The delta function is actually produced by considering the limit as $\xi$ goes to zero so that
\begin{equation}
   \delta(x-x_0) = \lim_{\xi\rightarrow 0} \left\{  \begin{array}{cl} {1}/({2\xi}) & x_0-\xi<x<x_0+\xi \\0 & \mbox{elsewhere} \end{array} \right. 
\end{equation}
where $\delta(x-x_0)$ is the delta function.  The limit shows that the solution is zero everywhere except at $x_0$, where it goes to infinity.  However, it goes to infinity in such as way that the integral defining the impulse is finite, specifically unity.  Thus the Dirac delta function is not a traditional function at all, but rather a generalized function that only makes sense when integrated against.   Figure~\ref{fig:delta_sequence} shows the delta function sequence considered here, showing how the function changes as a function of $\xi$.

Distribution theory generalizes the concept of functions to those functions such as the Dirac delta function that only make sense when integrated.  Another example function that fits within distribution theory is the Heaviside function $H(\xi)$ which is useful in many applications.  In fact, the Heaviside function is the derivative of the delta function as we will show later, i.e. 
$H'(x-x_0) = \delta(x-x_0)$.   The following delta function sequences could have also been used instead of 
(\ref{eq:delta_sequence})
\begin{subeqnarray}
  && \delta(x) = \lim_{\xi\rightarrow 0} \frac{1}{\pi\xi} \frac{1}{1+x^2/\xi^2} \\
&& \delta(x) = \lim_{\xi\rightarrow 0} \frac{1}{\sqrt{\pi} \xi} \exp (-x^2/\xi^2)   \\
&& \delta(x) = \lim_{\xi\rightarrow 0} \frac{\xi}{\pi} \frac{\sin^2 x/\xi}{x^2}
\label{eq:delta_sequence2} 
\end{subeqnarray}
where all these functions have an impulse of unity and where they take on a value of zero everywhere except $x=0$ where they are infinity.  The advantage of these alternative Dirac delta function sequences is that they can be explicitly differentiated in order to evaluate quantities such as $\delta' (x)$ or $\delta'' (x)$ for instance.

\subsection*{The Sifting Property}
The Dirac delta function is fundamental for constructing the Green's function.  Specifically, its sifting property make it perfect for producing a solution representation in terms of inner products.  The sifting property follows from (\ref{eq:delta_sequence}) and the following inner product
\begin{equation}
 \langle  f(x), \delta(x-x_0)   \rangle = \int_0^l f(x) \delta(x-x_0) dx = \lim_{\xi\rightarrow 0}  \int_{x_0+\xi}^{x_0+\xi} \frac{f(x)}{2\xi} dx
    = \lim_{\xi\rightarrow 0}  \frac{1}{2\xi} \int_{x_0+\xi}^{x_0+\xi} f(x) dx =  \lim_{\xi\rightarrow 0}  \frac{1}{2\xi} f(\bar{x}) = f(x_0)
\end{equation}
where $\bar{x}$ is from the mean-value theorem of integral calculus.  Note that as $\xi\rightarrow 0$, the integral value is squeezed to the value $f(x_0)$.  Thus the important result is achieved
\begin{equation}
 \langle  f(x), \delta(x-x_0)   \rangle = f(x_0)
\end{equation}
which is known as the sifting property.  The sifting property of the delta function selects the value of the function $f(x)$ wherever the delta function is applied.  This makes sense since the delta function is zero everywhere aside from a single point $x_0$.  At that point, the unit impulse response selects the value of the function at that location.  

\subsection*{The Green's function} 
The sifting property will be used to construct the Green's function solution of interest.  To do so, we consider the two fundamental problems
\begin{subeqnarray}
 && Lu=f   \\
 && L^\dag G=\delta(x-\xi)
 \label{eq:Lgreen}
\end{subeqnarray}
where $x, \xi\in [0,l]$ and $L^\dag$ is the adjoint operator associated with $L$ and its boundary conditions.  The function $G=G(x,\xi)$ will denote the Green's function.   Taking the inner product of both sides of (\ref{eq:Lgreen}a) with respect to $G(x,\xi)$ gives
\begin{subeqnarray*}
 &&   \langle  Lu, G   \rangle = \langle f, G   \rangle \\
 &&   \langle  u,L^\dag G   \rangle = \langle f, G   \rangle \\
 &&   \langle  u, \delta(x-\xi)   \rangle = \langle f, G   \rangle \\
  &&   u(\xi) = \langle f, G   \rangle .
\end{subeqnarray*}
This gives the solution $u(\xi)=\int_0^l f(x) G(x,\xi) dx$.  These small number of manipulations are critical to being fluent with the method.  But since $x$ and $\xi$ are simply dummy variables, we can switch them to give the solution
\begin{equation}
   u(x) = \int_0^l f(\xi) G(\xi,x) d\xi
   \label{eq:green_solution}
\end{equation}
and thus the inverse operator $L^{-1}[f]  =  \int_0^l f(\xi) G(\xi,x) d\xi$ is computed.  So once the Green's function is found, it is integrated against the forcing $f(x)$ to produce the unique solution.  Note that the Green's function satisfies the impulsively forced (\ref{eq:Lgreen}b).
Thus once the solution is achieved for a given $\xi$, the more general forcing $f(x)$ is a superposition via integration of all possible impulse forcings.  The solution technique relies fundamentally on the linearity of the governing equations (\ref{eq:Lgreen}).  

\subsection*{Solving for the Green's Function}
To illustrate how to solve for the Green's function, consider the example boundary value problem
\begin{equation}
   u_{xx} = f(x)
\end{equation}
on the domain $x\in [ 0 ,l]$ with the boundary conditions
\begin{subeqnarray}
  && u(0)=0 \\
  && u_x (l)=0 .
\end{subeqnarray}
This problem is self-adjoint so that $L^\dag=L = d^2/dx^2$.  The Green's functions then satisfies.
\begin{equation}
   G_{xx} = \delta(x-\xi)
\end{equation}
on the domain $x, \xi \in [ 0 ,l]$ with the boundary conditions
\begin{subeqnarray}
  && G(0)=0 \\
  && G_x (l)=0 .
\end{subeqnarray}
In order to solve this problem, we first integrate the governing equation for the Green's function across the impulse at $x=\xi$.  Since the delta function is zero everywhere except $x=\xi$, we then localize the integration from $\xi^-$ to $\xi^+$ which is the value of $x$ immediately to the left and right of $x=\xi$.  Thus we have
\begin{equation}
  \int_{\xi^-}^{\xi^+} G_{xx} dx = G_x |_{\xi^-}^{\xi^+}  = \left[ G_x \right]_\xi = \int_{\xi^-}^{\xi^+} \delta(x-\xi) dx =1
\end{equation}
where $\left[ G_x \right]_\xi$ denotes the value of the function evaluated on the right $\xi^+$ minus the value on left $\xi^-$. In this case, the integration yields the jump in the first derivative at the impulse $x=\xi$ of 
\begin{equation}
  \left[ G_x \right]_\xi = 1 .
\end{equation}
In contrast, the Green's function is assumed to be a continuous function so that
\begin{equation}
  \left[ G \right]_\xi = 0.
\end{equation}
Thus the impulse creates a discontinuous, but quantifiable, jump in the first derivative while maintaining a continuous solution.

The Green's function solution is then achieved in two parts, a solution for $x<\xi$ and a solution for $x>\xi$.  For both these cases, the boundary value problem is homogenous, making it much easier to solve.  Specifically, we have the following:\\

\noindent {\em I.  Solution for $x<\xi$:}  For this case, the governing equation reduces to the homogeneous form $G_{xx}=0$ with the left boundary condition $G(0)=0$.  This gives the solution in this domain as 
\begin{equation}
   G= Ax + B
\end{equation}
which upon applying the boundary condition gives $B=0$ so that
\begin{equation}
   G= Ax 
\end{equation}
for $x <\xi$.\\

\noindent {\em II.  Solution for $x>\xi$:}  For this case, the governing equation  again reduces to the homogeneous form $G_{xx}=0$ with the right boundary condition $G_x(l)=0$.  This gives the solution in this domain as 
\begin{equation}
   G= Cx + D
\end{equation}
which upon applying the boundary condition gives $C=0$ so that
\begin{equation}
   G= D 
\end{equation}
for $x >\xi$.\\

The left and right solution must now satisfy the two extra constraints that the solution is continuous at $x=\xi$ and there is a jump in the derivative at $x=\xi$.    These conditions applied give
\begin{subeqnarray}
  &&   \left[ G \right]_\xi =0  \,\,\,\, \rightarrow \,\,\,\, G(\xi^+)-G(\xi^-) = D-A\xi = 0 \\
   &&   \left[ G_x \right]_\xi =1  \,\,\,\, \rightarrow \,\,\,\, G_x(\xi^+)-G_x(\xi^-) =0 - A = 1
\end{subeqnarray}
which results in $A=-1$ and $D=A\xi=-\xi$.  Note that to the left and right of $x=\xi$, there were two linearly independent solutions.  Thus there were four unknown constants to determine.  The left and right boundary conditions, along with the jump conditions in the solution and its derivate provide four constraints to uniquely determine the Green's function.

The Green's function can then be constructed from the two solutions as
\begin{equation}
   G(x,\xi)= \left\{ \begin{array}{cc} -x & x<\xi \\ -\xi & x>\xi\end{array}   \right.
   \label{eq:green_ex1}
\end{equation}
which is the fundamental solution to this problem for any $x, \xi$ pair.  The overall solution to the boundary value problem can then be computed from (\ref{eq:green_solution}).

For the specific example where $f(x)=x$, the Green's function can be explicitly computed to give
\begin{eqnarray}
 u(x) && \hspace*{-.2in}= \int_0^x \xi (-\xi) d\xi + \int_x^l \xi (-x) d\xi \nonumber \\
  && \hspace*{-.2in} = \int_0^x - \xi^2 d\xi - \int_x^l \xi x d\xi \nonumber \\
  &&  \hspace*{-.2in}= - \frac{\xi^3}{3} |_0^x - x \frac{\xi^2}{2} |_x^l \nonumber \\
  &&  \hspace*{-.2in}= \left( \frac{x^2}{2} - \frac{l^2}{2} \right) x - \frac{x^3}{3} 
\end{eqnarray}
which is a compact representation of the solution satisfying the boundary conditions and the equation.

As a second example of a Green's function, consider the following
\begin{equation}
   G_{xx} + k^2 G= \delta(x-\xi)
\end{equation}
on the domain $x, \xi \in [ 0 ,l]$ with the boundary conditions
\begin{subeqnarray}
  && G_x (0)=0 \\
  && G_x (l)=0 .
\end{subeqnarray}
The Green's function solution is again achieved in two parts, a solution for $x<\xi$ and a solution for $x>\xi$.  For both these cases, the boundary value problem is homogenous, making it much easier to solve.  Specifically, we have the following:\\

\noindent {\em I.  Solution for $x<\xi$:}  For this case, the governing equation reduces to the homogeneous form $G_{xx}+k^2 G=0$ with the left boundary condition $G_x (0)=0$.  This gives the solution in this domain as 
\begin{equation}
   G= A \cos kx 
\end{equation}
for $x <\xi$.\\

\begin{figure}[t]
\begin{overpic}[width=0.6\textwidth]{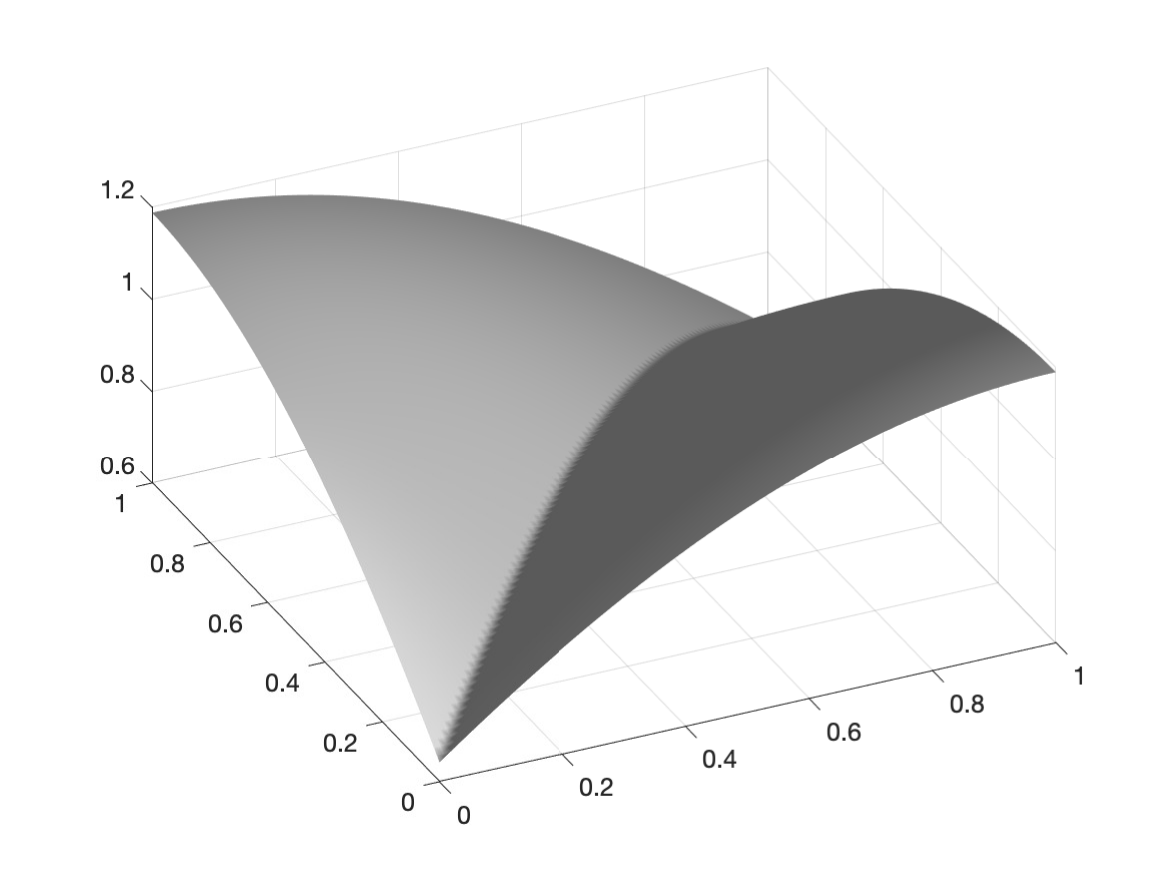}
\put(-3,50){$G(x,\xi)$}
\put(10,20){$\xi$}
\put(78,9){$x$}
\end{overpic}
\caption{Green's function given by (\ref{eq:green1}) with $l=1$ and $k=1$.  Note the solution surface is continuous with a first derivative jump at $x=\xi$.  Integrating this surface against $f(x)$ prescribes the solution $u(x)$ for $Lu=f$.}
 \label{fig:green1}
\end{figure}

\noindent {\em II.  Solution for $x>\xi$:}  For this case, the governing equation  again reduces to the homogeneous form $G_{xx}+k^2 G=0$ with the right boundary condition $G_x(l)=0$.  This gives the solution in this domain as 
\begin{equation}
   G= B \cos k(x-l)
\end{equation}
for $x >\xi$.\\

The left and right solution must now satisfy the two extra constraints that the solution is continuous at $x=\xi$ and there is a jump in the derivative at $x=\xi$.    The conditions are $\left[ G \right]_\xi =0$ and $  \left[ G_x \right]_\xi =1 $.  Applying the jump conditions gives the Green's function solution
\begin{equation}
    G(x,\xi) = \left\{ \begin{array}{cc}   {\cos kx \cos k(\xi-l)}/{(k \sin kl)} & x< \xi \\
      { \cos k(x-l) \cos k\xi }/{(k \sin kl)} & x> \xi \end{array}    \right.  .
      \label{eq:green1}
\end{equation}
This Green's function can be integrated against any forcing $f(x)$ to produce a solution and the inverse of the operator $L$.

\newpage
\section*{Lecture 7:  Green's Function for Sturm-Liouville Problems}

The Green's function provides a principled method to construct the inverse of an operator $L^{-1}$.  Using distribution theory, and specifically the delta function, allows us to solve two simpler homogenous equations that are patched together by constraints on continuity and derivative jumps.  To broaden the framework of the theory, we consider the canonical boundary value problem arising from Sturm-Liouville theory. The Sturm-Liouville operator is given by
\begin{equation}
    L u = - \left[ p(x) u_x \right]_x + q(x) u = f(x)
\end{equation}
on the interval $x\in[0,l]$.  The boundary conditions for the Sturm-Liouville operator are given by
\begin{subeqnarray}
 &&  \alpha_1 u(0) + \beta_1 u_x(0)=0 \\
 &&  \alpha_2 u(l) + \beta_2 u_x(l)=0 .
\end{subeqnarray}
Since the Sturm-Liouville operator is self-adjoint, the Green's function $G(x,\xi)$ then satisfies
\begin{equation}
    L G = - \left[ p(x) G_x \right]_x + q(x) G = \delta(x-\xi)
    \label{eq:LSG}
\end{equation}
on the interval $x, \xi \in[0,l]$.  The boundary conditions for the Sturm-Liouville operator are given by
\begin{subeqnarray}
 &&  \alpha_1 G(0) + \beta_1 G_x(0)=0 \\
 &&  \alpha_2 G(l) + \beta_2 G_x(l)=0 .
\end{subeqnarray}
What remains is to impose conditions on the Green's function and its derivative at the delta function location $x=\xi$.  

The first constraint on the Green's function is that it is continuous.  Thus there are no jumps in the solution even where the delta function is applied.  This gives the condition
\begin{equation}
   [G(x,\xi)]_\xi = G({\xi^+},\xi) - G({\xi^-},\xi) = 0 
\end{equation}
which enforces continuity of the solution.  The jump in the derivative is more difficult to evaluate, but it can be computed by integrating (\ref{eq:LSG}) near $x=\xi$.  This then gives
\begin{eqnarray}
   & \int_{\xi^-}^{\xi^+} \left(  - \left[ p(x) G_x \right]_x + q(x) G \right) dx =
     \int_{\xi^-}^{\xi^+} \delta(x-\xi) dx  \nonumber \\
    &      - \left[ p(x) G_x \right]_{\xi^-}^{\xi^+}    + \int_{\xi^-}^{\xi^+} q(x) G  dx =
     1 \nonumber \\
    & \left[ p(x) G_x \right]_\xi = -1 .
\end{eqnarray}
where we have assumed as usual in the Sturm-Liouville problem that $p(x)$ and $q(x)$ are continuous.  This gives in total the jump condition of the derivative
\begin{equation}
     \left[ G_x (x,\xi) \right]_\xi = - \frac{1}{p(\xi)} .
\end{equation}

The Green's function for the Sturm-Liouville problem can now be calculated.  Specifically, we solve the homogeneous problem for the two regimes $x<\xi$ and $x>\xi$ and enforce continuity of the solution and the jump in the derivate above.  \\

\noindent {\em I.  Solution for $x<\xi$:}  For this case, the governing equation reduces to the homogeneous form $L G=0$ with the left boundary conditions enforced.  This gives the solution in this domain as 
\begin{equation}
   G= A y_1(x)
\end{equation}
where $y_1(x)$ is the homogenous solution satisfying the left boundary conditions.  The parameter $A$ is an unknown constant that is yet to be determined.  Recall that for this second order operator, there are in general two linearly independent solutions.  Once the boundary condition is imposed, there is effectively a single linearly independent solution denoted by $y_1(x)$.\\

\noindent {\em II.  Solution for $x>\xi$:}  For this case, the governing equation  again reduces to the homogeneous form $LG=0$ with the right boundary conditions enforced.  This gives the solution in this domain as 
\begin{equation}
   G= B y_2(x) 
\end{equation}
where $y_2(x)$ is the homogenous solution satisfying the right boundary conditions.  The parameter $B$ is an unknown constant that is yet to be determined.  Again recall that for this second order operator, there are in general two linearly independent solutions.  Once the boundary condition is imposed, there is effectively a single linearly independent solution denoted by $y_2(x)$.\\

The left and right solution must now satisfy the two extra constraints that the solution is continuous at $x=\xi$ and there is a jump in the derivative at $x=\xi$.  Imposing
$\left[ G (x,\xi) \right] = 0$ and $\left[ G_x (x,\xi) \right]_\xi = - {1}/{p(\xi)}$ then gives
the Green's function solution
\begin{equation}
   G(x,\xi) = \left\{ \begin{array}{cc}   {y_1(x)y_2(\xi)}/{(p(\xi)W(\xi))} & x< \xi \\
      {y_1(\xi)y_2(x)}/{(p(\xi)W(\xi))} & x> \xi \end{array}    \right.  .
      \label{eq:green2}
\end{equation}
where $W(\xi)$ is the Wronskian between $y_1(x)$ and $y_2(x)$ evaluated at $x=\xi$.
This is the general form of the solution for all Sturm-Liouville problems.

\subsection*{Example Green's functions for Sturm-Liouville Problems}
To illustrate the computation of the Green's function, we consider a number of examples.  The first is to consider
\begin{equation}
   G_{xx} = \delta (x-\xi)
\end{equation}
with $G(0)=0$ and $G_x(l)=0$.  The solution for $x<\xi$ is given by $y_1(x)=x$ while the solution for $x>\xi$ is given by $y_2(x)=1$.    The Wronskian is computed to be $W[y_1,y_2]=-1$ and $p(\xi)=1$.  This gives in total the Green's function already found
in (\ref{eq:green_ex1}).

As a second example, we again can consider 
\begin{equation}
   G_{xx} +k^2 G= \delta (x-\xi)
\end{equation}
with $G_x(0)=0$ and $G_x(l)=0$.  The solution for $x<\xi$ is given by $y_1(x)=\cos kx$ while the solution for $x>\xi$ is given by $y_2(x)=\cos k(x-l)$.    The Wronskian is computed to be $W[y_1,y_2]=k\sin kl$ and $p(\xi)=1$.  This gives in total the Green's function already found in (\ref{eq:green1}).

As a third example, we consider the Green's function in radial coordinates
\begin{equation}
  \frac{d}{dr} \left[ r \frac{dG}{dr} \right]  = \delta (r-\rho)
\end{equation}
where $r, \rho\in [0,l]$ and with $G(l)=0$ and $G(0)=$finite.  This is a singular problem since when $r=0$, the derivative terms disappear.  The homogenous solution to this problem is constructed from the two linearly independent solutions $G=\{ 1, \ln r\}$.  
For $r<\rho$, the only solution which can satisfy the boundary condition gives $y_1=1$.  For $r>\rho$, the only solution which can satisfy the boundary condition gives $y_2=\ln r-\ln l = \ln (r/l)$.  The Wronskian can be computed to give $W=1/r$ and $pW=r (1/r)=1$.
The Green's function is then 
\begin{equation}
    G(r,\rho) = \left\{ \begin{array}{cc}   \ln (\rho/l) & r< \rho \\
      \ln (r/l) & r> \rho \end{array}    \right.  .
\end{equation}

\subsection*{Example solution using Green's functions}
As a final example, we consider not just construction of the Green's function, but of a solution to a $Lu=f$ problem.  For this example, we consider
\begin{equation}
   u_{xx} +2 u= -x
\end{equation}
with $u(0)=0$ and $u(1)+u_x(1)=0$.  This is a Sturm-Liouville problem with $p(x)=1$, $q(x)=2$ and $f(x)=-x$.  The Green's function satisfies
\begin{equation}
   G_{xx} +2 G= \delta (x-\xi)
\end{equation}
with $G(0)=0$ and $G(1)+G_x(1)=0$.  For $x<\xi$, the solution satisfying the left boundary value is given by $y_1(x)=\sin \sqrt{2} x$.  For $x>\xi$, the solution satisfying the right boundary condition is given by  $y_2(x)=\sin \sqrt{2}(x-1)-\sqrt{2} \cos\sqrt{2}(x-1)$.  The Wronskian can be computed to be $W=\sqrt{2}(\sin\sqrt{2}+\sqrt{2}\cos\sqrt{2})$.  This gives the Green's function
\begin{equation}
   G(x,\xi) = \left\{ \begin{array}{cc}   {\sin\sqrt{2}x \left[ \sin\sqrt{2}(\xi-1) - \sqrt{2} \cos\sqrt{2} (\xi-1)   \right] }/{\left[ \sqrt{2}(\sin\sqrt{2}+\sqrt{2}\cos\sqrt{2}) \right]} & x< \xi \\
      {\sin\sqrt{2}\xi \left[ \sin\sqrt{2}(x-1) - \sqrt{2} \cos\sqrt{2} (x-1)   \right] }/{ \left[ \sqrt{2}(\sin\sqrt{2}+\sqrt{2}\cos\sqrt{2}) \right]  } & x> \xi \end{array}    \right.  .
      \label{eq:green3}
\end{equation}
The solution is then found from performing the integration
\begin{equation}
   u(x)=\int_0^1 f(\xi) G(\xi,x) d\xi
\end{equation}
where it is important to note the arguments of the Green's function have been switching from what is presented in (\ref{eq:green3}).  The integral is done in two parts:  an integration from $\xi\in[0,x]$ and an integration from $\xi\in[x,1]$.  It is quite a long calculation, but it eventually leads to the solution
\begin{equation}
   u(x) = \frac{ \sin \sqrt{2} x }{\sin\sqrt{2}+\sqrt{2}\cos\sqrt{2}} - \frac{x}{2} .
\end{equation}
Thus the forcing $f(x)=-x$ needs to be integrated against the Green's function to produce the solution $u(x)$.

\newpage
\section*{Lecture 8:  Modified Green's Functions}

We return to the concept of solvability and null spaces in relation to the Green's function.  We again consider the boundary value problem
\begin{equation}
 Lu=f   
\end{equation}
with boundary conditions
\begin{subeqnarray}
  && B_1(u(0))=0 \\
  && B_2(u(l))=0
\end{subeqnarray}
given the domain $x, \xi\in [0,l]$ with prescribed boundary conditions given by $B_j(u)=0$.  The Fredholm-Alternative theorem requires that the forcing function $f(x)$ by orthogonal to the null space of the adjoint operator.  Consider a self-adjoint operator so that $L=L^\dag$.  Then the adjoint null space is given by solutions
\begin{equation}
 Lv=0  . 
\end{equation}
If the null space is given by the trivial solution, $v=0$, then the standard Green's function approach given in the last section yields a unique Green's function solution.  However, if the null space is non-trivial, then the Green's function approach must be modified in order to achieve a fundamental solution.   To illustrate this, let $v(x)$ be a normalized and nontrivial null space solution so that $\|v\|=1$.  Then when solving $Lu=f$, if one finds a solutions $u_0(x)$, then another solution is simply 
\begin{equation}
  u(x) = u_0(x) + c v(x)
\end{equation}
where $c$ is an arbitrary constant.  Thus an infinite number of solutions exists since the solution spanning the null space can always be added to a solution $u(x)$.  This will create problems for the Green's function unless the solvability condition
\begin{equation}
   \langle f,v \rangle = 0
\end{equation}
is satisfied.

Consider the implications of a non-trivial null space on the Green's function solution.  Again, this is illustrated with a self-adjoint problem so that 
\begin{equation}
 LG=\delta(x-\xi)   
\end{equation}
with boundary conditions
\begin{subeqnarray}
  && B_1(G(0))=0 \\
  && B_2(G(l))=0
\end{subeqnarray}
given the domain $x, \xi\in [0,l]$.  If the operator $L$ has a non-trivial null-space $v(x)$, then solvability gives
\begin{equation}
   \langle f,v \rangle = \langle \delta(x-\xi), v(x) \rangle = v(\xi) \neq 0
\end{equation}
for any $\xi\in(0,l)$.  Thus solvability is generically violated, even though for a specific value of $\xi$ one might have $v(\xi)=0$.

The modified Green's function overcomes this difficulty by reformulating the Green's function boundary value problem so that
\begin{equation}
 LG=\delta(x-\xi)  - v(x)v(\xi)  = F(x,\xi)
\end{equation}
with boundary conditions
\begin{subeqnarray}
  && B_1(G(0))=0 \\
  && B_2(G(l))=0
\end{subeqnarray}
given the domain $x, \xi\in [0,l]$.  Note the reformulation has the null space solution, posited in a symmetric form between $x$ and $\xi$, as part of the right-hand side forcing function.  With this modified right-hand side, the solvability condition is given by
\begin{equation}
   \langle F,v \rangle = \langle \delta(x-\xi) - v(x)v(\xi), v(x) \rangle = v(\xi) - v(\xi) \langle v,v\rangle = v(\xi) - v(\xi) =0
\end{equation}
where it should be recalled that the null-space solution is normalized to unity without loss of generality. 

The {\em modified Green's function} thus considers the two fundamental problems
\begin{subeqnarray}
 && Lu=f   \\
 && L^\dag G_m=\delta(x-\xi) - v(x)v(\xi)
\end{subeqnarray}
where $x, \xi\in [0,l]$ and $L^\dag$ is the adjoint operator associated with $L$ and its boundary conditions.  The function $G_m=G_m(x,\xi)$ will denote the modified Green's function.   Taking the inner product of both sides of $Lu=f$ with respect to $G_m(x,\xi)$ gives
\begin{subeqnarray*}
 &&   \langle  Lu, G_m   \rangle = \langle f, G_m   \rangle \\
 &&   \langle  u,L^\dag G_m   \rangle = \langle f, G_m   \rangle \\
 &&   \langle  u, \delta(x-\xi) - v(x)v(\xi)   \rangle = \langle f, G_m   \rangle \\
  &&   u(\xi) - v(\xi) \langle u,v \rangle = \langle f, G   \rangle .
\end{subeqnarray*}
This gives the solution
\begin{equation}
    u(\xi)=\int_0^l f(x) G_m(x,\xi) dx + c v(\xi)
\end{equation}  
where the constant $c=\langle u,v \rangle$.  But since $x$ and $\xi$ are simply dummy variables, we can switch them to give the solution
\begin{equation}
   u(x) = \int_0^l f(\xi) G_m(\xi,x) d\xi + c v(x)
   \label{eq:mgreen_solution}
\end{equation}
where $c$ is actually an arbitrary constant and reflects the fact that the non-trivial null space gives this arbitrary contribution to the solution.

\subsection*{Solving for the modified Green's Function}
To illustrate how to solve for the modified Green's function, consider the example boundary value problem which is very similar to the first example of the last section
\begin{equation}
   u_{xx} = f(x)
\end{equation}
on the domain $x\in [ 0 ,l]$ with the boundary conditions
\begin{subeqnarray}
  && u_x(0)=0 \\
  && u_x (l)=0 .
\end{subeqnarray}
The boundary conditions now, both zero flux, is what allows there to be a non-trivial null space solution.  Specifically, the adjoint null space
satisfies
\begin{equation}
   v_{xx} = 0
\end{equation}
on the domain $x\in [ 0 ,l]$ with the boundary conditions
\begin{subeqnarray}
  && v_x(0)=0 \\
  && v_x (l)=0 .
\end{subeqnarray}
Thus the solution $v(x)=$constant is the non-trivial null space solution.  In normalized form, $v(x)=1$.

This problem is self-adjoint so that $L^\dag=L = d^2/dx^2$.  The standard Green's functions would then satisfy
\begin{equation}
   G_{xx} = \delta(x-\xi)
\end{equation}
on the domain $x, \xi \in [ 0 ,l]$ with the boundary conditions
\begin{subeqnarray}
  && G_x (0)=0 \\
  && G_x (l)=0 .
\end{subeqnarray}
Note that the standard approach for the Green's function does not satisfy solvability since
\begin{equation}
  \langle f,v \rangle = \langle \delta(x-\xi), 1 \rangle = 1\neq 0.
\end{equation}

To show where problems arise, one can immediately attempt to solve for the Green's function in the standard way
where $\left[ G \right]_\xi = 0$ and $ \left[ G_x \right]_\xi=1$ must be satisfied.  As before, the Green's function solution is then achieved in two parts, a solution for $x<\xi$ and a solution for $x>\xi$.  For both these cases, the boundary value problem is homogenous, making it much easier to solve.  Specifically, we have the following:\\

\noindent {\em I.  Solution for $x<\xi$:}  For this case, the governing equation reduces to the homogeneous form $G_{xx}=0$ with the left boundary condition $G_x(0)=0$.  This gives the solution in this domain as 
\begin{equation}
   G= Ax + B
\end{equation}
which upon applying the boundary condition gives $A=0$ so that
\begin{equation}
   G= B 
\end{equation}
for $x <\xi$.\\

\noindent {\em II.  Solution for $x>\xi$:}  For this case, the governing equation  again reduces to the homogeneous form $G_{xx}=0$ with the right boundary condition $G_x(l)=0$.  This gives the solution in this domain as 
\begin{equation}
   G= Cx + D
\end{equation}
which upon applying the boundary condition gives $C=0$ so that
\begin{equation}
   G= D 
\end{equation}
for $x >\xi$.\\

The left and right solution must now satisfy the two extra constraints that the solution is continuous at $x=\xi$ and there is a jump in the derivative at $x=\xi$.    These conditions applied give $B=D$ for continuity, but the jump in the derivative is not satisfied since $G_x=0$ for both solutions.  Thus we fail to satisfy the appropriate Green's function conditions.

The modified Green's function accounts for the non-trivial null space solution by modifying the Green's function to the following
\begin{equation}
   {G_m}_{xx} = \delta(x-\xi) - \frac{1}{l} = F(x,\xi)
\end{equation}
on the domain $x, \xi \in [ 0 ,l]$ with the boundary conditions
\begin{subeqnarray}
  && {G_m}_x (0)=0 \\
  && {G_m}_x (l)=0 .
\end{subeqnarray}
The constant $1/l$, which is a multiple of the constant null solution, is added to that the solvability condition holds
\begin{equation}
   \langle F,v \rangle = \langle F,1 \rangle =  \int_0^l \left(  \delta(x-\xi) - \frac{1}{l} \right)  dx = 1 - l \frac{1}{l}  = 0 .
\end{equation}
The modified Green's function can now be solved in the two regions.

\noindent {\em I.  Solution for $x<\xi$:}  For this case, the governing equation reduces to the non-homogeneous form ${G_m}_{xx}=-1/l$ with the left boundary condition ${G_m}_x(0)=0$.  This gives the solution in this domain as 
\begin{equation}
   G_m= Ax + B + \frac{x^2}{2l} 
\end{equation}
which upon applying the boundary condition gives $A=0$ so that
\begin{equation}
   G_m= B + \frac{x^2}{2l} 
\end{equation}
for $x <\xi$.\\

\noindent {\em II.  Solution for $x>\xi$:}  For this case, the governing equation  again reduces to the non-homogeneous form ${G_m}_{xx}=-1/l$ with the right boundary condition ${G_m}_x(l)=0$.  This gives the solution in this domain as 
\begin{equation}
   G_m= Cx + D + \frac{x^2}{2l}
\end{equation}
which upon applying the boundary condition gives $C=0$ so that
\begin{equation}
   G_m= D + \xi - x + \frac{x^2}{2l}
\end{equation}
for $x >\xi$.\\

After applying the conditions of continuity and a discrete jump in the derivative, this yields the solution
\begin{equation}
   G_m(x,\xi)= \left\{ \begin{array}{cc}  c+\frac{x^2}{2l} & x<\xi \\ c+\xi - x + \frac{x^2}{2l} & x>\xi\end{array}   \right.
   \label{eq:green_m1}
\end{equation}
for some constant $c$.  This is the fundamental solution to this problem for any $x, \xi$ pair.  Often the constant $c$ is chosen to make the modified Green's function symmetric.  Specifically, without loss of generality it can be chose so that
\begin{equation}
   G_m(x,\xi)= \left\{ \begin{array}{cc}  \frac{l}{3} - \xi +\frac{x^2+\xi^2}{2l} & x<\xi \\ \frac{l}{3} - x +\frac{x^2+\xi^2}{2l} & x>\xi\end{array}   \right.  .
   \label{eq:green_m2}
\end{equation}
The overall solution to the boundary value problem can then be computed from the modified Green's function as
\begin{equation}
  u(x) = \int_0^l G_m(\xi,x) f(\xi) d\xi + \mbox{constant}
\end{equation}
where the constant reflects the fact that null space is spanned by $v(x)=$constant.  Thus the null space and solvability play an important role in practical solutions of Green's functions.  Indeed, the null space must be explicitly handled in order to achieve a fundamental solution.

\newpage
\section*{Lecture 9:  Regular Perturbation Theory}

Asymptotics and perturbation theory have played a foundational role in the development of our understanding of physical phenomenon and nonlinear systems.  This is largely due to the fact that, with few exceptions, we do not have general methods for solving nonlinear system of equations, differential equations or partial differential equations.  Perturbation theory, however, allows us to consider solutions to nearby problems which are analytically tractable.  The simplest, and one of the oldest physical examples to consider, is the pendulum.  Although we often approximate the pendulum by a linear, second order differential equation, this is actually an asymptotic approximation for the pendulum's equations of motion for small amplitude oscillations. 

Perturbation theory can be formulated for algebraic system of equations, time-dependent differential equations, space-time PDEs, and boundary value problems.  In the treatment considered here, we will consider one of two forms of the the perturbation theory
\begin{subeqnarray}
 &&  Lu = \epsilon F(u, u_t, \cdots, t)  \,\,\,\,\,\,\,  t\in[0,\infty]  \\
 &&  Lu = \epsilon F(u, u_x, u_{xx}, \cdots , x)  \,\,\,\,\,\,\,  x\in[0,l]
 \label{eq:pert_eq}
\end{subeqnarray}
with the corresponding constraints (initial conditions for the former, boundary conditions for the later)
\begin{subeqnarray}
  && u(0)= A, u_t (0)=B \\
  &&  \alpha_1 u(0) + \beta_1 u_x(0) = 0, \,\,\,\, \alpha_2 u(l) + \beta_2 u_x(l) = 0.
  \label{eq:pert_bc}
\end{subeqnarray}
Here $L$ is assumed to be a linear operator and $\epsilon\ll 1$.  Moreover, it is assumed we don't know how to solve these problems when $\epsilon\neq 0$.   When $\epsilon=0$, however, we can solve the resulting problem $Lu=0$.  Thus perturbation theory assumes we can solve a nearby problem since $\epsilon$ is small.

Peturbation theory recasts a nonlinear problem as a sequence, or hierarchy, of forced linear problems.  Thus the difficulty of not having analytic solution techniques for handling the nonlinearity is replaced by solving multiple linear problems with forcings.  To demonstrate how perturbation theory accomplishes this, consider the time-dependent problem (\ref{eq:pert_eq}a) and (\ref{eq:pert_bc}a).  To approximate the solution, we expand the solution in a perturbation series
\begin{equation}
  u = u_0 + \epsilon u_1 + \epsilon^2 u_2 + \cdots
  \label{eq:pert_expand}
\end{equation}
where each term in the series is assumed to get smaller and smaller due to the $\epsilon^n$ term.  Provided each individual contribution $u_n\sim O(1)$, then this will occur.

Inserting the perturbation expansion into (\ref{eq:pert_eq}a) gives 
\begin{equation}
  L (u_0 + \epsilon u_1 +  \cdots) = \epsilon F(u_0 + \epsilon u_1 +  \cdots, {u_0}_t + \epsilon {u_1}_t +  \cdots,  t) 
  \label{eq:Lexpand}
\end{equation}
with initial conditions
\begin{subeqnarray}
 &&  u(0)=u_0(0) + \epsilon u_1(0) + \epsilon^2 u_2(0) + \cdots  = A \\
&&  u_t(0)={u_0}_t (0) + \epsilon {u_1}_t (0) + \epsilon^2 {u_2}_t (0) + \cdots = B .
\end{subeqnarray}

After the formal expansion, the goal is to collect terms at each order of $\epsilon$.  A Taylor expansion of 
the governing equations (\ref{eq:Lexpand}) is required to do this.  This yields the hierarchy of linear equations
\begin{subeqnarray}
&&  O(1) \hspace*{0.5in}  L u_0 = 0 \\
&&  O(\epsilon) \hspace*{0.5in} L u_1 = F_1 (u_0, {u_0}_t, \cdots , t) \\
&&  O(\epsilon^2) \hspace*{0.5in} L u_2 = F_2 (u_0, u_1, {u_0}_t, {u_1}_t, \cdots , t) \\
&&  \,\,\,\,\,\,\,\, \vdots \nonumber \\
&& O(\epsilon^n) \hspace*{0.5in} L u_n = F_n (u_0, u_1, \cdots, u_{n-1}, {u_0}_t, {u_1}_t, \cdots, {u_{n-1}}_t, \cdots , t) 
\label{eq:pert_hier}
\end{subeqnarray}
with the associated initial conditions
\begin{subeqnarray}
   && u_0(0)=A, \,\,\,\,\, u_n (0) = 0 \,\,\, (n\geq 2) \\
  && {u_0}_t(0)=B, \,\,\,\,\, {u_n}_t (0) = 0 \,\,\, (n\geq 2) .
\end{subeqnarray}
The perturbation expansion thus replaces the nonlinear evolution equation with a hierarchy of linear equations whose forcing function is progressively more complex.  Regardless of the complexity, it is a linear problem and can be solved through a variety of techniques already outlined in a previous chapter on differential equations.

The solution of the problem can then be formally constructed as
\begin{equation}
   u= \sum_{n=0}^{\infty} \epsilon^n u_n .
\end{equation}
In practice, however, only the first couple of terms are kept to approximate the solution so that
\begin{equation}
   u= u_0 + \epsilon u_1 + O(\epsilon^2)
\end{equation}
which is a valid approximation up to typically a time $T\sim O(1/\epsilon)$~\cite{kevorkian2013perturbation,murdock1999perturbations}.  Each correction term added to the approximation improves the quality of the approximation, but typically does not extend the range of its validity in time.

The general goal of perturbation theory is then to do the following:  (i) Identify a part of the problem that can be solved analytically, (ii) Determine which terms in the problem are small and which can be ordered by a small parameter $\epsilon$, and (iii) Calculate an approximate solution to a desired accuracy by a hierarchical approach which refines the approximation at each step.  More colloquially, to solve the problem you are incapable of solving, solve one nearby and perturb it towards the actual problem of interest.

\subsection*{Solvability and the Fredholm alternative theorem}

The above asymptotics allows us a mathematical framework whereby we can approximate solutions to governing equations in a principled manner from nearby solutions or equations.  The underlying mathematical strategy will rely on the {\em Fredholm alternative theorem}.  Recall from the chapter on linear operator theory that we can consider a linear operator ${L}$ acting in an infinite dimensional function space
\begin{equation}
  {L} u = f
  \label{eq:luf}
\end{equation}
where $u(x)$ and $f(x)$ are functions on a given domain (bounded or unbounded).  The Fredholm alternative theorem considers under what conditions (\ref{eq:luf}) admits solutions.  This can be made explicit by considering an associated problem which involves computing the null space of the adjoint operator
\begin{equation}
  {L}^\dag v = 0
  \label{eq:lv0}
\end{equation}
where $L^\dag$ is the adjoint operator and $v(x)$ is the null space of the adjoint operator.  

Taking the inner product on both sides of (\ref{eq:luf}) with respect to the adjoint null space $v(x)$, we find
\begin{eqnarray}
  &&  \langle f,v \rangle = \langle L u,v \rangle  \nonumber \\
  &&  \langle f,v \rangle = \langle u,L^\dag v \rangle  \nonumber \\
  &&  \langle f,v \rangle = 0 .
\end{eqnarray}
Thus the right hand side forcing $f(x)$ must be orthogonal to the null space of the adjoint operator $L^\dag$.  Under these conditions, (\ref{eq:luf}) can be solved.   This simple theorem will have profound implications for the development of our asymptotic methods and the hierarchy of equations it generates through (\ref{eq:pert_hier}).

\subsection*{Regular perturbation expansions}

The perturbation method introduced is known as a regular perturbation expansion.  More 
generally, we can consider a PDE of the form
\begin{equation}
   \frac{\partial u}{\partial t} = N(u, u_x, u_{xx}, \cdots, \mu) +\epsilon  G(u,x,t)
   \label{eq:pert_general}
\end{equation}
where $N(\cdot)$ is nonlinear representation of the PDE, $G(u,x,t)$ is a perturbation to
the underlying PDE, and $\mu$ is parameter (bifurcation).   A regular perturbation expansion
takes the form
\begin{equation}
u=u_{0}(x,t)+\epsilon u_{1}(x,t) + \epsilon^2 u_2 (x,t) + \cdots .
\end{equation}
Collecting terms at each order of the perturbation theory gives
\begin{subeqnarray}
 & O(1): \hspace{.5in}  &\frac{\partial u_0}{\partial t} = N (u_0, {u_0}_x, {u_0}_{xx}, \cdots, \mu)  \\
 &  O(\epsilon): \hspace{.5in} & \frac{\partial u_1}{\partial t} = L (u_0) u_1 + F_1 (u_0)  \\
 & O(\epsilon^2): \hspace{.5in} & \frac{\partial u_2}{\partial t} = L (u_0) u_2 + F_2 (u_0, u_1)
\end{subeqnarray}

At leading order, we find the solution
\begin{equation}
  u(x,t)= u_0 (x,t,A,B, \cdots)
\end{equation}
where the parameters $A, B, \cdots$ are constants that parametrize the leading order-solution.

At higher orders in the perturbation expansion, Fredholm alternative requires that the forcing terms $F_j$ be orthogonal to the null space of the adjoint operator $L$ so that
\begin{equation}
  \langle v, F_j \rangle = 0 .
  \label{eq:solvability}
\end{equation}
If the Fredholm alternative is not satisfied, then secular growth terms arise and the correction terms $u_j(x,t)$ grow without bounds.  This establishes a formal procedure for approximating solutions of a given PDE of the form (\ref{eq:pert_general}).  The three examples below demonstrate the use of the regular perturbation theory framework.

\subsection*{Example:  An initial value problem}

To illustrate the mechanics of the regular perturbation expansion framework, consider the following initial value problem
\begin{equation}
    u_{tt} - u = \epsilon u^2
\end{equation}
with $u(0)=A$ and $u_t(0)=B$.  Using the expansion (\ref{eq:pert_expand}) and collecting terms in powers of $\epsilon$ gives the hierarchy of equations to $O(\epsilon^2)$
\begin{subeqnarray}
  O(1) \hspace*{0.5in}  & L u_0 = {u_0}_{tt} - u_0  = 0    & \hspace*{0.5in}   u_0(0)=A, \,\, {u_0}_t (0) =B\\
  O(\epsilon) \hspace*{0.5in} & L u_1 = {u_1}_{tt} - u_1  = {u_0}^2  & \hspace*{0.5in}   u_1(0)=0, \,\, {u_1}_t (0) =0 \\
  O(\epsilon^2) \hspace*{0.5in} & L u_2 = {u_2}_{tt} - u_2  = 2 {u_0} u_1  &  \hspace*{0.5in}  u_2(0)=0, \,\, {u_2}_t (0) =0 .
\end{subeqnarray}
The leading order solution to this problem can be easily found to be
\begin{equation}
   u_0=\frac{A+B}{2} \exp(t) + \frac{A-B}{2} \exp(-t)  .
\end{equation}
The square of this solution is then the nonhomogenous forcing term in the equation for $u_1$ which yields
\begin{equation}
   u_1=\frac{A^2-2AB-2B^2}{6} \exp(t) + \frac{A^2+2AB-2B^2}{6} \exp(-t) 
    + \frac{(A+B)^2}{12} \exp(2t) + \frac{(A-B)^2}{12} \exp(-2t) - \frac{A^2-B^2}{2} .
\end{equation}
Thus the solution is determined up to $O(\epsilon^2)$ by the solution $u\approx u_0+\epsilon u_1$.  This example shows how the hierarchical nature of the problem very quickly gets complicated due to the nonhomogenous forcing term of the differential equation.

\subsection*{Example:  An boundary value problem}

To illustrate the mechanics of the regular perturbation expansion framework on boundary value problems, consider the following 
\begin{equation}
    u_{xx} - u = \epsilon u^2
\end{equation}
with $u(0)=A$ and $u(1)=B$. 
Using the expansion (\ref{eq:pert_expand}) and collecting terms in powers of $\epsilon$ gives the hierarchy of equations to $O(\epsilon^2)$
\begin{subeqnarray}
  O(1) \hspace*{0.5in}  & L u_0 = {u_0}_{xx} - u_0  = 0    & \hspace*{0.5in}   u_0(0)=A, \,\, {u_0}(1) =B\\
  O(\epsilon) \hspace*{0.5in} & L u_1 = {u_1}_{xx} - u_1  = {u_0}^2  & \hspace*{0.5in}   u_1(0)=0, \,\, {u_1}(1) =0 \\
  O(\epsilon^2) \hspace*{0.5in} & L u_2 = {u_2}_{xx} - u_2  = 2 {u_0} u_1  &  \hspace*{0.5in}  u_2(0)=0, \,\, {u_2}(1) =0 .
\end{subeqnarray}
The leading order solution to this problem can be found to be
\begin{equation}
   u_0=\frac{A-Be}{1-e^2} \exp(x) - \frac{eA-B}{1-e^2} \exp(-x)  .
\end{equation}
The nonhomogenous equation for $u_1$ can be determined by a number of methods, including the method of undetermined coefficients or a Green's function solution.  This yields the solution form
\begin{equation}
  u_1 = C_1 + C_2 \exp(x) + C_3 \exp(-x) + C_4 \exp(2x) + C_5 \exp(-2x) 
\end{equation}
where the coefficients $C_j$ are determined after much algebraic manipulation.  The solution is determined up to $O(\epsilon^2)$ by the solution $u\approx u_0+\epsilon u_1$.  Again, this example shows how even for a simple problem the hierarchical nature of the problem very quickly gets complicated due to the nonhomogenous forcing term of the differential equation.

\subsection*{Example:  An eigenvalue problem}

Finally, we consider the perturbed eigenvalue problem
\begin{equation}
   u_{xx} + \lambda u = \epsilon f(u)
\end{equation}
with $u(0)=u(1)=0$.  Such problems are often considered in quantum mechanics as time-independent perturbation theory.  Specifically, the interest in quantum mechanics stems from understanding how perturbations in potentials shift energy levels in atoms.  For eigenvalue problems, our expansion is not only in the solution itself, but also in the eigenvalues:
\begin{subeqnarray}
 &&  u = u_0 + \epsilon u_1 + \epsilon^2 u_2 + \cdots \\
 && \lambda= \lambda_0 + \epsilon \lambda_1 + \epsilon^2 \lambda_2 + \cdots .
  \label{eq:pert_evexpand}
\end{subeqnarray}
Collecting terms at each order of the expansion gives the hierarchy of equations
\begin{subeqnarray}
  O(1) \hspace*{0.5in}  &  {u_0}_{xx} - \lambda_0 u_0  = 0    & \hspace*{0.5in}   u_0(0)=0, \,\, {u_0}(1) =0\\
  O(\epsilon) \hspace*{0.5in} & {u_1}_{xx} - \lambda_0 u_1  = f(u_0) - \lambda_1 u_0  & \hspace*{0.5in}   u_1(0)=0, \,\, {u_1}(1) =0 \\
 \end{subeqnarray}
where now determined the values of the $\lambda_n$ is part of the solution architecture.  The leading order solution to this problem gives
\begin{equation}
  u_0 = A\sin n\pi x
\end{equation}
where $\lambda_0 = n^2 \pi^2$, i.e. there are an infinite number of solutions corresponding to the eigenfunctions and eigenvalues.  At next order this gives
\begin{equation}
   {u_1}_{xx} + n^2 \pi^2 u_1 = f(A \sin n \pi x) - \lambda_1 A \sin n \pi x .
\end{equation}
To solve this problem, we first represent the function $f(\cdot)$ in terms of the eigenfunctions of the leading order problem.  This gives
\begin{equation}
    f = \sum_{m=1}^{\infty} a_m \sin m\pi x .
\end{equation}
The equation for the next order solution is then
\begin{equation}
   {u_1}_{xx} + n^2 \pi^2 u_1 = (a_n- \lambda_1 A) \sin n \pi x 
   +\sum_{\substack{m=1 \\ m\neq n}}^{\infty} a_m \sin m \pi x .
\end{equation}
And so we come to a situation where the Fredholm-Alternative theorem must be applied.  In particular, the nonhomogenous term $\sin n\pi x$ is not orthogonal to the null space of the adjoint operator.  Recall that the operator here is self-adjoint and its null space is spanned by $\sin n\pi x$.  This will lead to the generation of a secular term that does not satisfy the boundary conditions.  To rectify this situation, the value of $\lambda_1$ must be appropriate chosen to remove the secular growth term.  Thus one can choose
\begin{equation}
   \lambda_1 = \frac{a_n}{A}
\end{equation}
and the eigenfunction solution of the next order correction term is then
\begin{equation}
  u_1(x) = C \sin n\pi x  +\sum_{\substack{m=1 \\ m\neq n}}^{\infty} \frac{a_m }{(n^2-m^2)\pi^2}\sin m \pi x .
\end{equation}
This example explicitly shows that the solvability conditions highlighted in the section on linear operators will play a dominant and important role in perturbation theory in general.

\newpage
\section*{Lecture 10:  The Poincare-Lindsted Method}

Although regular perturbation expansions are useful in acquiring approximate solutions to many problems, the theory has a number of significant limitations.  This can be illustrated by a simple example problem.  Consider the linear, perturbed oscillator
\begin{equation}
   u_{tt} + (1 + \epsilon)^2 u  = 0 
\end{equation}
with initial conditions $u(0)=A$ and $u_t (0) = 0$.  This problem is trivial to solve analytically.  The solution is expressed as
\begin{equation}
    u(t) = A \cos (1+\epsilon) t 
\end{equation}
which shows that solution has simple oscillations with a frequency $\omega=1+\epsilon$ and a period of $T=2\pi/(1+\epsilon)$.   Thus the perturbation simply shifts the oscillation frequency by amount $\epsilon$.  It is informative to evaluate this simple problem from a regular perturbation expansion.  Inserting the expansion (\ref{eq:pert_expand}) into this problem gives the hierarchy of linear equations
\begin{subeqnarray}
  O(1) \hspace*{0.5in}  & {u_0}_{tt} + u_0  = 0    & \hspace*{0.5in}   u_0(0)=A, \,\, {u_0}_t(0) =0\\
  O(\epsilon) \hspace*{0.5in} &  {u_1}_{tt} + u_1  = -2 {u_0}  & \hspace*{0.5in}   u_1(0)=0, \,\, {u_1}_t(0) =0 \\
  O(\epsilon^2) \hspace*{0.5in} & {u_2}_{tt} + u_2  = -2 u_1 - u_0   &  \hspace*{0.5in}  u_2(0)=0, \,\, 
   {u_2}_t(0) =0 .
\end{subeqnarray}
The leading order solution is easily solved to yield
\begin{equation}
   u_0 = A\cos t
\end{equation}
which then gives at $O(\epsilon)$
\begin{equation}
    {u_1}_{tt} + u_1  = -2 A \cos t
\end{equation}
with $ u_1(0)=0$ and ${u_1}_t(0) =0$.  Immediately, one can see there will be a problem when considering the Fredholm-Alternative theorem. Specifically, the right hand side $\cos t$ is in the null space of the adjoint operator, and not orthogonal to it.  As a consequence, the solution takes the form
\begin{equation}
  u_1 = A t \sin t
\end{equation}
which gives the perturbative solution to $O(\epsilon^2)$ as
\begin{equation}
   y = A\cos t - \epsilon A t \sin t + O(\epsilon^2).
\end{equation}
This is a secular grown term that shows the perturbation is the same size as the leading order solution at $t\sim O(1/\epsilon)$.  Indeed, if one computes the solution at $O(\epsilon^2)$, there is a $\{t^2 \sin t, t^2 \cos t\}$ component which also becomes $O(1)$ at the same time as the first correction.  In fact, since the nonhomogeneous term is in the null space, all the perturbation terms become $O(1)$ simultaneously and the regular expansion fails to hold.  This is a significant failing for the regular perturbation expansion technique since it cannot even handle the simplest of perturbations.

The problem with the regular perturbation expansion is clear:  it is not constructed to handle perturbative shifts in the frequency, thus leading to secular growth terms.  Unlike the regular perturbation expansion, the Poincare-Lindstedt method gives additional flexibility to the asymptotic expansion by defining a new time scale which is a function of the expansion parameter $\epsilon$.  We again consider the general framework (\ref{eq:pert_general}). The A Poincare-Lindstedt perturbation expansion takes the form
\begin{equation}
u=u_{0}(x,\tau)+\epsilon u_{1}(x,\tau) + \epsilon^2 u_2 (x,\tau) + \cdots
\end{equation}
where
\begin{equation}
  \tau = (\omega_0 +\epsilon \omega_1 +\epsilon^2 \omega_2 + \cdots ) t .
\end{equation}
Thus in addition to determining the solution and its corrections $u_j(x,t)$, one must also determine corrections to the time scale through the $\omega_j$.  The ability to account for frequency shifts was critical for the use of this method in planetary orbit calculations.

Collecting terms at each order of the perturbation theory gives now
\begin{subeqnarray}
 & O(1): \hspace{.5in}  &\frac{\partial u_0}{\partial \tau} = N (u_0, {u_0}_x, {u_0}_{xx}, \cdots, \mu)  \\
 &  O(\epsilon): \hspace{.5in} & \frac{\partial u_1}{\partial \tau} = L (u_0,\omega_0) u_1 + F_1 (u_0,\omega_1,\omega_0)  \\
 & O(\epsilon^2): \hspace{.5in} & \frac{\partial u_2}{\partial \tau} = L (u_0,\omega_0) u_2 + F_2 (u_0, u_1,\omega_2,\omega_1,\omega_0) \, .
\end{subeqnarray}
Applying the solvability conditions (\ref{eq:solvability}) allows us to determine values of $\omega_j$ which prevent the solutions from growing without bound.  Specifically, we choose the values of $\omega_j$ so that the terms in the null space are removed or set to zero.  This overcomes some of the fundamental limitations of the regular perturbation theory expansion.

\subsection*{Example:  A perturbed oscillator}

Nonlinearities fundamentally drive two outcomes:  They generate harmonics and they induce frequency shifts.  Accounting for such frequency shifts is what made the Poincare-Lindstedt method so attractive for its earliest application of planetary motion approximation.  To demonstrate this  basic principle, consider the simple perturbed oscillatory system
\begin{equation}
   u_{tt} + k^2 u = \epsilon f (u)
   \label{eq:pert_osc}
\end{equation}
where $\epsilon\ll1$, $u(0)=0$ and $u_t (0) =A$. 
The Poincare-Lindsted method defines a new time variable 
\begin{equation}
   \tau=\omega(\epsilon) t
\end{equation}
in order to allow time to stretch and accommodate changes in frequency that can arise from the perturbation $f(u)$ in (\ref{eq:pert_osc}).  The change of variables from $t$ to $\tau$ ($\partial t \rightarrow \omega \partial \tau$) gives the governing dynamics 
\begin{equation}
   \omega^2 u_{\tau\tau} + u = \epsilon f (u) .
\end{equation}
Our objective now is to perform an asymptotic expansion that not only solves for the state space $u$, but for the new stretched time $\tau$.  Thus we expand each so that
\begin{subeqnarray}
   && u= u_0 + \epsilon u_1 + \epsilon^2 u_2 + \cdots \\
   && \omega= \omega_0 +\epsilon \omega_1 +\epsilon^2 \omega_2 + \cdots  .
\end{subeqnarray}
Inserting this formal asymptotic expansion into (\ref{eq:pert_osc}) gives the following equations at each order of the $\epsilon$ expansion
\begin{subeqnarray}
\label{eq:pert_osc1}
   && O(1):  \hspace*{.5in}  {u_0}_{\tau\tau} +  u_0 =0 \\
   && O(\epsilon):  \hspace*{.5in}   {u_1}_{\tau\tau} +  u_1 = f(u_0) - 2\omega_1 {u_0}_{\tau\tau} \\
   && O(\epsilon^2):  \hspace*{.5in}  {u_2}_{\tau\tau} +  u_2 = u_0 \frac{\partial  f (u_0)}{\partial u}  -2\omega_1 {u_1}_{\tau\tau}   
    - \left(2\omega_2 + \omega_1^2\right) {u_0}_{\tau\tau}  \\
   &&  \,\,\,\,\,\,\,\, \vdots \nonumber
\end{subeqnarray}
with 
\begin{subeqnarray}
\label{eq:pert_osc2}
   && u_0(0)=0, \,\, {u_0}_\tau (0) =A \\
   && u_1(0)=0, \,\, {u_1}_\tau (0)=0 \\
  && u_2(0)=0, \,\, {u_2}_\tau (0)=0  \\
   &&  \,\,\,\,\,\,\,\, \vdots  \nonumber
\end{subeqnarray}

The Fredholm-Alternative theorem given by (\ref{eq:solvability}) imposes a restriction on the right hand side forcings of (\ref{eq:pert_osc1}).    The left side operator is self-adjoint and the leading order solution is given by 
\begin{equation}
   u_0(\tau) = A \sin \tau. 
\end{equation}
This results in the solvability conditions
\begin{subeqnarray}
\label{eq:pert_osc3}
   && O(\epsilon):  \hspace*{.5in}   \left\langle f(u_0) - 2\omega_1  {u_0}_{\tau\tau}, \sin \tau \right\rangle =0 \\
   && O(\epsilon^2):  \hspace*{.5in}  \left\langle u_0 \frac{\partial  f (u_0)}{\partial u}  -2\omega_1   {u_1}_{\tau\tau}   - \left(2\omega_2 + \omega_1^2\right)  {u_0}_{\tau\tau}, \sin \tau \right\rangle=0  \\
   &&  \,\,\,\,\,\,\,\, \vdots \nonumber
\end{subeqnarray}
Once a specific $f(u)$ is given, then the values of $\omega_j$  and solutions $u_j (\tau)$ can be determined at each order.

\subsection*{The Duffing Equation}

\begin{figure}[t]
\begin{overpic}[width=0.8\textwidth]{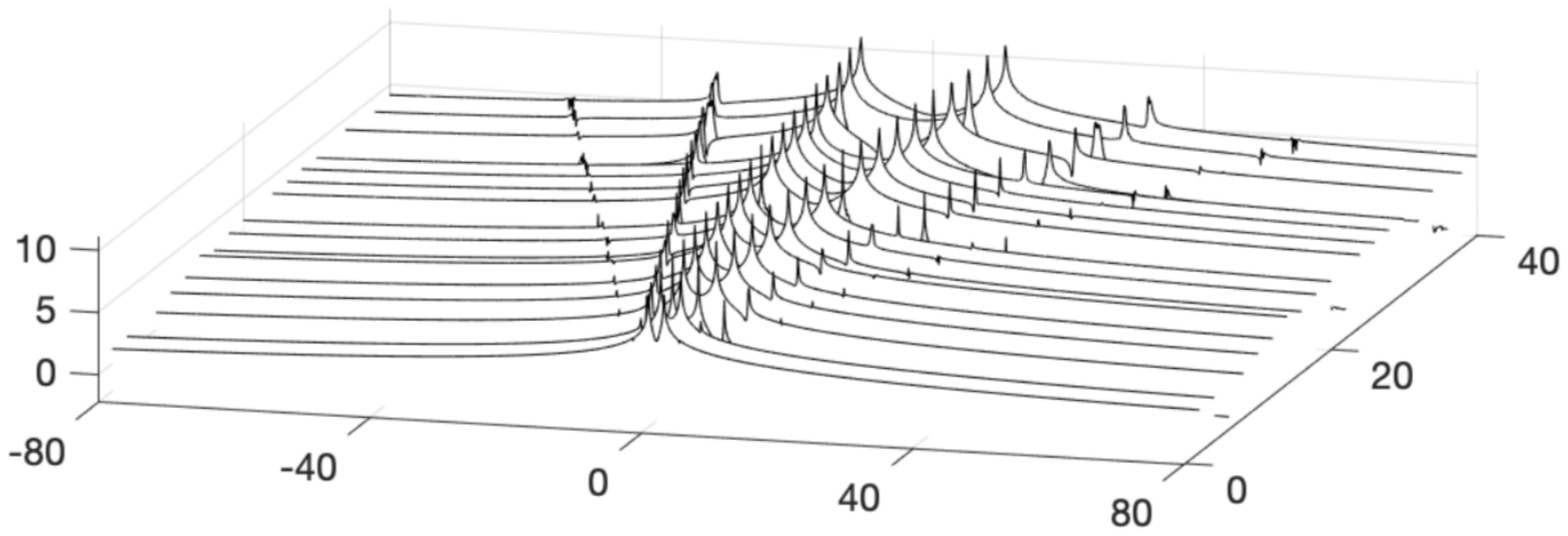}
\put(-5,45){(a)}
\put(-5,10){(b)}
\put(-7,36){$\log (|\hat{u}|)$}
\put(83,27){Amplitude $A$}
\put(35,20){frequency $\omega$}
\put(100,5){Amplitude $A$}
\put(43,-20){frequency $\omega$}
\put(63,54){$\omega_0$}
\put(43,52){$-3\omega_0$}
\put(33,51){$-5\omega_0$}
\put(80,48){$5\omega_0$}
\put(52,55){$-\omega_0$}
\put(70,52){$3\omega_0$}
\end{overpic}\\[-2in]
\begin{overpic}[width=0.8\textwidth]{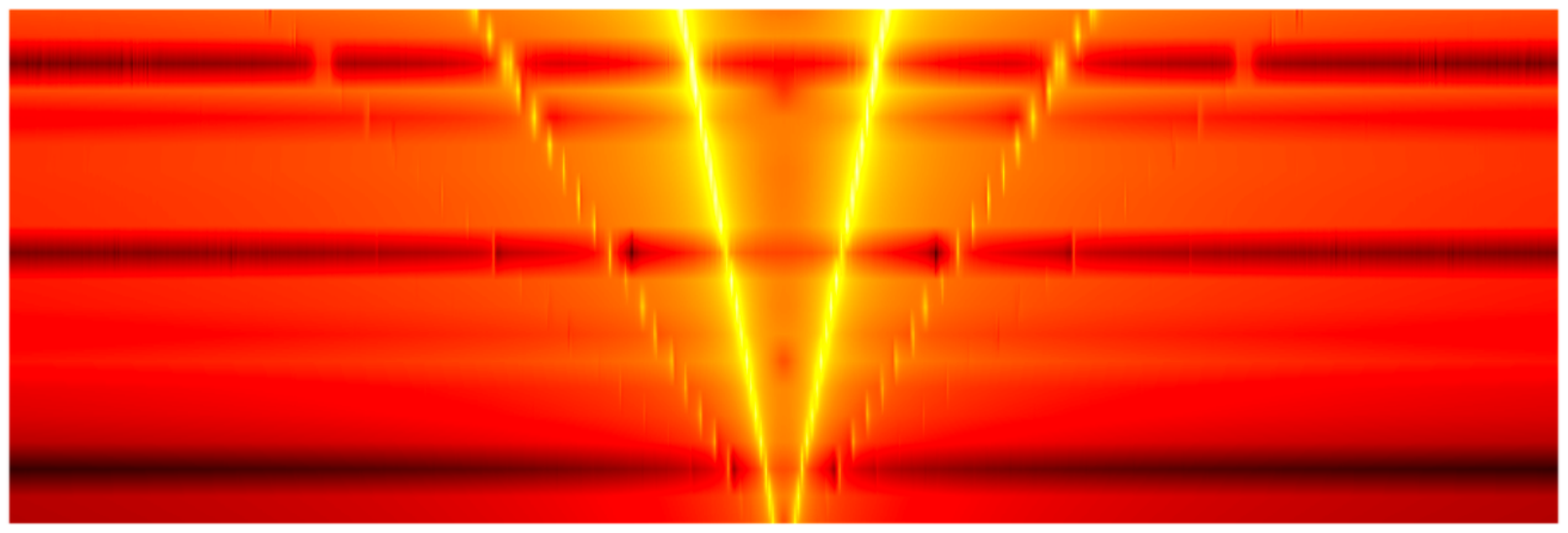}
\end{overpic}
\vspace*{-1in}
\caption{(a) Spectrum and (b) spectrogram of the Duffing oscillator as a function of amplitude $u(0)=A$ with $u_t(0)=0$.  The Fourier transform $\hat{u}(\omega)$ of the signal $u(t)$ is plotted on a logarithmic scale.  The frequency $\omega_0$ is the shifted frequency.  Note the generation of sidebands at $\pm 3\omega_0$ and $\pm 5\omega_0$ as is expected from the perturbation theory.}
 \label{fig:duffing1}
\end{figure}

The Duffing equation is given by
\begin{equation}
      u_{tt} + u + \epsilon u^3 = 0
   \label{eq:duffing}
\end{equation}
so that $f(u)=-u^3$.  The solvability (\ref{eq:pert_osc3}) at $O(\epsilon)$ gives
\begin{equation}
    \left\langle A^3 \sin^3 \tau - 2\omega_1 A \sin \tau, \sin \tau \right\rangle =0 .
\end{equation}
Using the relationship $\sin^3 \tau= (3 \sin \tau -\sin 3\tau)/4$ reduces the solvability to
\begin{equation}
    \left\langle (2\omega_1 - 3A^2/4)  \sin \tau + A^2/4 \sin 3\tau, \sin \tau \right\rangle =0 .
\end{equation}
Noting that the Fourier terms $\sin \tau$ and $\sin 3\tau$ are orthogonal over the periodic interval of integration $\tau\in[0,2\pi]$ gives 
\begin{equation}
  \omega_1 = \frac{3}{8} A^2 .
\end{equation}
This value can then be inserted into (\ref{eq:pert_osc1}b) so that the non-homogenous differential equation at $O(\epsilon)$ can be solved to produce
\begin{equation}
   u_1(\tau)= \frac{3A^3}{32} \sin\tau - \frac{A^3}{32} \sin 3\tau .
\end{equation}
This gives the perturbed solution
\begin{equation}
  u(t) =A\sin\left[ (1 + 3 \epsilon A^2/8)t \right] + \frac{\epsilon A^3}{32} \left\{ 3\sin \left[(1+3\epsilon A^2/8) t\right] - \sin\left[ 3(1+3\epsilon A^2/8) t \right]  \right\}
\end{equation}
where the frequency shift and harmonic generation are clearly illustrated from the asymptotic expansion.  This is shown more clearly in Fig.~\ref{fig:duffing1}.  Note that if the nonlinearity had been an even function like $f(u)=-u^2$, then second harmonics $\sin 2\tau$ would be generated instead of $\sin 3\tau$.  Thus the symmetry of the nonlinearity is in the types of harmonics that are generated.

\subsection*{The Van der Pol Oscillator}

The Van der Pol oscillator is given by the equation
\begin{equation}
   u_{tt} + \epsilon \left(  u^2 - 1 \right) u_t + u = 0
\end{equation}
with initial conditions imposed on the amplitude and its derivative $u(0)$ and $u_t(0)$.  The Van der Pol oscillator exhibits a limit cycle attractor in its dynamics.  To understand this qualitatively, consider the second (damping) term in the governing equations.  If $(u^2-1)<1$ then the second term makes the damping term such that it causes growth of the solution.  If $(u^2-1)>1$, then the damping term is acts like a standard damping and forces the solution to decay.   Thus small amplitudes grow and large amplitudes decay.  This gives rise to a limit cycle.  Moreover, the nonlinear damping term shifts frequencies and generates harmonics as is expected from a cubic perturbation term.

As before, we introduce the new strained time variable $\tau=\omega(\epsilon) t$ which re-expresses the Van der Pol as
\begin{equation}
   \omega^2 u_{\tau\tau} + \epsilon \left(  u^2 - 1 \right) \omega u_\tau + u = 0 .
\end{equation}
Expanding terms
\begin{eqnarray}
   u=u_0 + \epsilon u_1 + \epsilon^2 u_2 + \cdots \\
   \omega=\omega_0 + \epsilon \omega _1 + \epsilon^2 \omega_2 + \cdots
\end{eqnarray}
gives us the ability to evaluate a perturbative solution.  But in addition to expanding the solution and frequency, we will also need to determine the initial conditions that allow for the existence of the limit cycle.  Thus we also have the expansion
\begin{eqnarray}
  && u(0) = A_0 + \epsilon A_1 + \epsilon^2 A_2 + \cdots \\
  && u_t (0) = 0 .
\end{eqnarray}
Note that we could also, or alternatively, expand the initial derivative, but the expansion used simply fixes a specific instance in time that is convenient to work with.  

The leading order solution satisfies
\begin{equation}
   \omega_0^2 {u_0}_{\tau\tau} + u_0 = 0
\end{equation}
with $u_0(0)=A_0$ and ${u_0}_\tau (0) = 0$.  The solution can be found to be $u_0=A_0 \cos (\tau/\omega_0)$.  The value of $\omega_0$ can be picked arbitrarily.  However, we will choose to look for solutions which are $2\pi$-periodic so that $\omega_0=1$ and 
\begin{equation}
   u_0=A_0 \cos \tau .
\end{equation}
Note that the value of $A_0$ will be chosen from solvability conditions applied at higher order.

The next order evolution equation at $O(\epsilon)$ is given by 
\begin{eqnarray}
   {u_1}_{\tau\tau} + u_1 &=& -2\omega_1 {u_0}_{\tau\tau} - (u_0^2 -1) {u_0}_{\tau} \nonumber \\
   &=& 2\omega_1 A_0 \cos \tau + (A_0^2 \cos^2 \tau - 1) A_0 \sin \tau \nonumber \\
   &=& 2 a_0 \omega_1 \cos \tau + A_0 (A_0^2/4 -1)\sin \tau + (A_0^3/4) \sin 3\tau = F_1 .
\end{eqnarray}
The right hand side $F_1$ has two terms in the null space of the adjoint, the $\cos\tau$ and $\sin \tau$ terms.  They can be made to vanish by enforcing
\begin{subeqnarray}
  &&  2 A_0 \omega_1 = 0 \,\,\, \rightarrow \,\,\,  \omega_1 =0 \\
  && A_0  ( A_0^2/4 -1) = 0 \,\,\, \rightarrow \,\,\, A_0=2
\end{subeqnarray}
which is equivalent to enforcing the Fredholm-Alternative which states that the forcing must be orthogonal to the null space of the adjoint operator so that $\langle F_1, \cos \tau \rangle=0$ and $\langle F_1, \sin\tau \rangle =0$.

At the next order $O(\epsilon^2)$, the evolution equation is given by
\begin{eqnarray}
   {u_2}_{\tau\tau} + u_2 &=& -2\omega_2 {u_0}_{\tau\tau} -\omega_1^2 {u_1}_{\tau\tau} - (u_0^2 -1) {u_1}_{\tau}  -\omega_1(u_0^2 -1) {u_0}_{\tau} - 2 u_0 {u_0}_\tau u_1   \nonumber \\
   &=& \left( 4\omega_2 - \frac{7}{4} \right) \cos \tau - 2 A_1 \sin \tau + \cos 3\tau - A_1 \sin 3\tau -\frac{5}{4} \cos 5\tau = F_2 .
\end{eqnarray}
Removal of the secular terms with $\langle F_2, \cos \tau \rangle=0$ and $\langle F_2, \sin\tau \rangle =0$
gives
\begin{subeqnarray}
  &&  4\omega_2 - \frac{7}{4}  = 0 \,\,\, \rightarrow \,\,\,  \omega_2 =\frac{7}{16} \\
  && 2 A_1 = 0 \,\,\, \rightarrow \,\,\, A_1=0 .
\end{subeqnarray}
This gives $u(0)= 2 + O(\epsilon^2)$ and $\tau=(1+7\epsilon^2/16 + O(\epsilon^4))t$ and the solution
\begin{equation}
  u(t)= 2 \cos \left[ (1+7\epsilon^2/16) t \right] + \epsilon \left[ \frac{3}{4} \sin \left[ (1+7\epsilon^2/16) t \right]
  -\frac{1}{4} \sin \left[ 3(1+7\epsilon^2/16) t \right]   \right] + O(\epsilon^2) .
\end{equation}
This solution explicitly captures the frequency shift and harmonic generation characteristic of nonlinear systems up to $O(\epsilon^2)$.  One could make a similar plot to Fig.~\ref{fig:duffing1} showing the shifts of the frequencies and the generation of the harmonics.

\newpage
\section*{Lecture 11:  The Forced Duffing Oscillator}

In the examples that follow, we will consider the origins of the Duffing equation and its dynamics under forcing.  Some of the most interesting problems in physics and engineering arise in such nonlinear dynamcal systems which are damped and driven.  The forced Duffing oscillator highlights some of the underlying dynamics that can occur in such a scenario.  The beginning of our study is the nonlinear pendulum equation
\begin{equation}
   u_{tt} + \sin u = 0
\end{equation}
which is the platonic model of the pendulum without a small angle approximation and/or damping or forcing terms.  If the system includes both damping and a forcing, then the model becomes
\begin{equation}
   u_{tt} + \delta u_t + \sin u = \gamma \cos \omega t
\end{equation}
which is now a damped-driven system parametrized by three parameters $\delta, \gamma$ and $\omega$.  The driving term is a simple harmonic driving at a constant frequency $\omega$.  In addition to the equation of motion, there are the following initial conditions $u(0)=\alpha$ and $u_t(0) =\beta$. 

A standard reduction of this system occurs for small amplitude oscillations where $u\ll 1$ so that $\sin u\approx u$.  If the oscillations get bigger, the linearizing approximation of the sinusoidal dependency no longer holds.  In this case, the next correction gives
\begin{equation}
   \sin u \approx u - \frac{u^3}{3!} +  O(u^5) 
\end{equation}
which shows the dominant nonlinearity to be a cubic by a Taylor series expansion of the sine term.\\

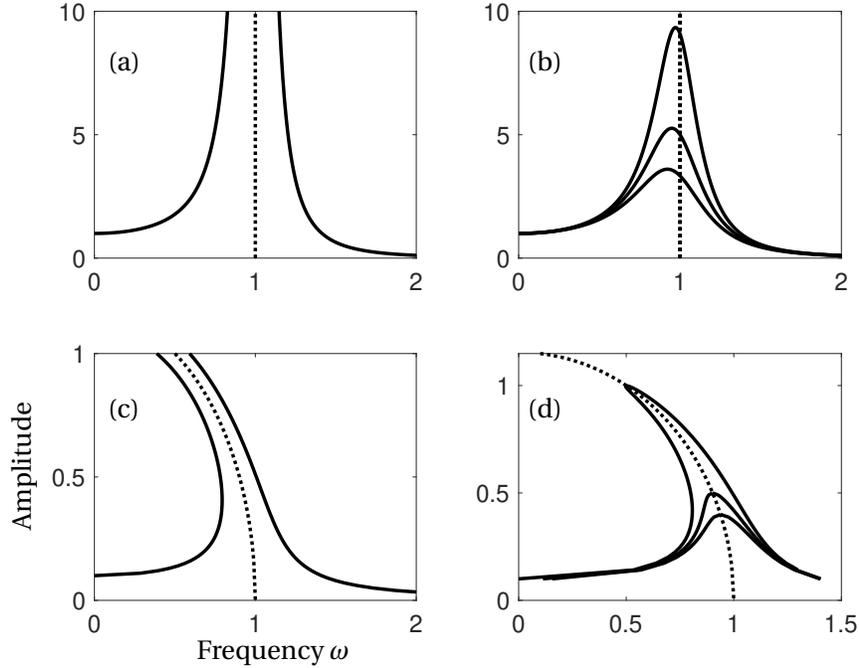
\begin{figure}[t]
\begin{overpic}[width=0.65\textwidth]{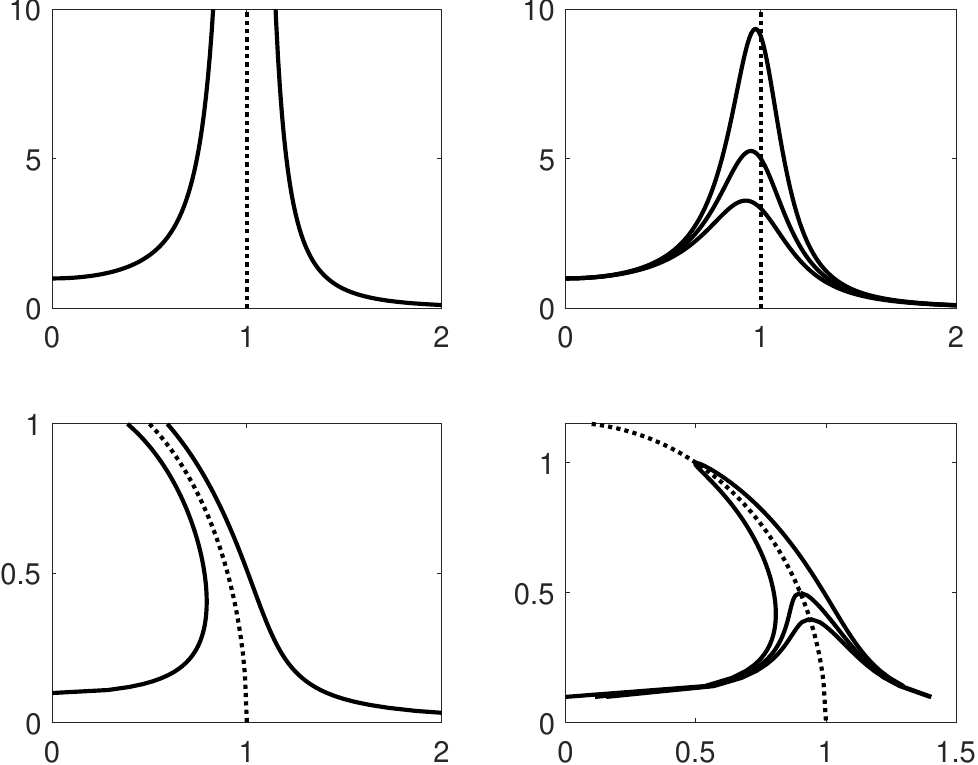}
\put(7,70){(a)}
\put(59,70){(b)}
\put(7,27){(c)}
\put(59,27){(d)}
\put(18,-3){Frequency $\omega$}
\put(-5,13){\rotatebox{90}{Amplitude}}
\end{overpic}
\vspace*{.1in}
\caption{Responce curves for (a) linear oscillator without damping and nonlinearity, (b) linear oscillator with damping and no nonlinearity, (c) nonlinear oscillator with no damping, and (d) nonlinear oscillator with damping.  The damped-driven system in (d) exhibits a cusp catastrophe whereby the response curve becomes multivalued for a frequencies $\omega<1$.}
 \label{fig:forced_duff}
\end{figure}

\noindent {\em Linear Theory:  No Damping --}
We begin the analysis by considering the well known linear theory of forced oscillations with no damping.  In this case, the model becomes
\begin{equation}
   u_{tt} + u = \gamma \cos \omega t .
\end{equation}
The dynamics for this system has  already been detailed.  To recap those findings, the solution to this system
is 
\begin{equation}
   u=\frac{\gamma}{2} t \sin t
\end{equation}
with the initial conditions $u(0)=u_t (0) =0$.  This shows the solution grows linearly with time and the amplitude is unbounded.  At the time originally considered, this growth was articulated as a result of a resonance forcing of the system.  However, since learning about linear operators and the Fredholm-Alternative theorem, we can see that the unbouded growth is mathematically a result of the forcing function being in the null space of the adjoint.  For this case, there is no mechanism for removing the secular growth that occurs due to the resonant forcing.  The response curve is shown in Fig.~\ref{fig:forced_duff}(a).\\   

\noindent {\em Linear Theory:  With Damping --}
The addition of damping generates solutions that no longer have unbounded growth.  This is because the damping term ensures that the forcing is no longer in the null space of the adjoint operator.  For damping, the model is given by  
\begin{equation}
   u_{tt} + \delta u_t + u = \gamma \cos \omega t .
\end{equation}
The details of the solution technique for this new damped oscillator model was already considered.  Of particular importance are the homogenous solutions ($\gamma=0$) that are given by
\begin{equation}
   u= c_1\exp (-\delta t/2) \cos \sqrt{1-\delta^2/4} t + c_2\exp (-\delta t/2) \sin \sqrt{1-\delta^2/4} t ,
\end{equation}
which show that the damping forces a change in frequency.  Moreover, the harmonic forcing is no longer part of the null space of the adjoint operator.  Thus unbounded growth will no longer occur.  

Solutions to the governing dynamics can now be constructed by finding a particular solution to add to the homogeneous solutions.   Since the term $\cos\omega t$ is on the right hand side, a particular solution of the form $u_p=A\cos\omega t + B\sin\omega t$ can be used.  Solving as previously, the persistent solution (as $t\rightarrow \infty$) is given by
\begin{equation}
  u=\frac{\gamma}{(1-\omega^2)^2 + \delta^2 \omega^2} \left( (1-\omega^2)\cos\omega t
  + \omega \delta \sin \omega t \right)
\end{equation}
where we note the homogeneous solutions decay to zero as $t\rightarrow\infty$ due to the $
\exp(-\delta t/2)$ term.  Importantly, solutions are now bounded since the forcing is no longer in the null space of the adjoint operator, thus satisfying the Fredholm-Alternative theorem.  The response curve is shown in Fig.~\ref{fig:forced_duff}(b).\\

\noindent {\em Nonlinear Theory:  No Damping --}
The role of nonlinearity is now considered.  First in the context of no damping, then in the context of a damped system.  Including nonlinearity gives the model 
\begin{equation}
   u_{tt} + u + \kappa u^3 = \gamma \cos \omega t .
\end{equation}
where $\kappa=-1/3!$ for the expansion of the sinusoidal term of the nonlinear pendulum.  The Poincare-Lindsted expansion introduces a new time scale $\tau=\omega t$ so that the dynamics is governed by
\begin{equation}
   u_{\tau\tau} + \frac{1}{\omega^2} u + \frac{\kappa}{\omega^2} u^3 = \frac{\gamma}{\omega^2} \cos \tau.
\end{equation}
This can be rewritten in the form
\begin{equation}
   u_{\tau\tau} +  u = \epsilon \left(  \Gamma \cos \tau -\beta u + u^3 \right) 
\end{equation}
where we have only retained terms to $O(\epsilon)$ and $\gamma/\epsilon^2=\epsilon \Gamma$, 
$1/\omega^2=1+\epsilon \beta$ and $\epsilon=-\kappa/\omega^2$.  Expanding the solution in the expansion $u=u_0 + \epsilon u_1 + \cdots$ gives at leading order
\begin{equation}
  {u_0}_{\tau\tau} + u_0 = 0
\end{equation}
with the leading order solution $u_0=A \cos \tau + B \sin \tau$.  At the next order, the governing equations are given by 
\begin{eqnarray}
   {u_1}_{\tau\tau} + u_1 &=& \Gamma \cos \tau -\beta u_0 + {u_0}^3 \nonumber \\
   &=& \left\{ \Gamma - \beta A + \frac{3}{4} A(A^2+B^2) \right\} \cos\tau  
      + \left\{ - \beta B + \frac{3}{4} B(A^2+B^2) \right\} \sin\tau  \nonumber \\
   & & \hspace*{.4in} +\frac{A}{4} (A^2-3B^2) \cos 3\tau + \frac{B}{4} (3A^2 - B^2) \sin 3\tau .
\end{eqnarray}
Solvability conditions require us to remove the secular terms.  This then gives
$\Gamma=A [ \beta- 3(A^2+B^2)/4]$ and $B [ \beta- 3(A^2+B^2)/4] = 0$.  This requires that
$B=0$ and $\Gamma=A(\beta-3A^2/4)$.  In terms of the original variables, this is then
\begin{equation}
   \omega^2 = 1 + \frac{3}{4} \kappa A^2 - \frac{\gamma}{A} .
\end{equation}
When the forcing strength is zero ($\gamma=0$), this gives $A=0$ or $A=\pm \sqrt{4(\omega^2 -1)/3\kappa }$.  Thus the solution can still blow up like the resonantly forced, undamped linear pendulum.  However, the unbounded growth is only along the line $\gamma=0$.  A finite amplitude solution is found for any $\omega$ otherwise.  Thus nonlinearity inhibits the growth of the oscillations.  Moreover, as the amplitude of the oscillations grow, the resonant forcing frequency changes.  Any kid knows this:  if you are on a swing, then as your swing amplitude gets bigger, you have to change your swing pumping frequency to go even higher.  The response curve is shown in Fig.~\ref{fig:forced_duff}(d).\\

\noindent {\em Nonlinear Theory:  With Damping --}
Damping is now included in the analysis along with the nonlinearity and driving.  This is now the full damped-driven system subject to nonlinear interactions.  In this case, the dynamics is given by
\begin{equation}
   u_{tt} +\delta u_t + u + \kappa u^3 = \gamma \cos \omega t .
\end{equation}
where $\kappa=-1/3!$ for the expansion of the sinusoidal term of the nonlinear pendulum.  The Poincare-Lindsted expansion introduces a new time scale $\tau=\omega t$ so that the dynamics is governed by
\begin{equation}
   u_{\tau\tau} +  u = \epsilon \left(  \Gamma \cos \tau -\beta u + u^3 -\Delta u_\tau \right) 
\end{equation}
where we have only retained terms to $O(\epsilon)$ and as before $\gamma/\epsilon^2=\epsilon \Gamma$, $1/\omega^2=1+\epsilon \beta$ and $\epsilon=-\kappa/\omega^2$ with the addition now that $\Delta=\delta/\omega$.  Expanding the solution in the expansion $u=u_0 + \epsilon u_1 + \cdots$ gives the same leading order solution as before with $u_0=A \cos \tau + B \sin \tau$.
It can be computed that at next order the secular growth terms can be removed provided the following equations are satisfied
\begin{subeqnarray}
  && \Delta B + A (\beta - 3(A^2 + B^2)/4) = \Gamma \\
  && \Delta A - B (\beta- 3(A^2 + B^2)/4) = 0 .
\end{subeqnarray}
By defining the total amplitude to be $R=(A^2 + B^2)^{1/2}$, these equations combined give
\begin{equation}
    R^2 \left\{ \Delta^2 + \left( \beta- \frac{3}{4} R^2 \right)^2 \right\} = \Gamma^2 
\end{equation}
which dictates the relationship between the amplitude $R$ and the frequency $\beta=\omega^2-1$.
In this case, the solutions no longer blowup.  However, the so-called cusp catastrophe occurs which is exhibited in Fig.~\ref{fig:forced_duff}(d). \\

There is more to be said about the Duffing equation.  In what has been considered, we fundamentally are driving an oscillator at leading order.  By flipping the sign of the linear term in the leading order governing equations to $u_{tt}-u=0$, the dynamics takes on a significantly different dynamics.  In particular, the onset of chaos and in such a damped-driven system are quite remarkable.  This is considered in the chapter on dynamical systems.

\newpage
\section*{Lecture 12:  Multiple-Scale Expansions}

Although the Poincare-Lindsted method is a powerful technique for characterizing the effects of perturbations on nonlinear dynamical systems, it is quite limited in scope.  Specifically, it only can accommodate frequency shifts by defining a stretched time variable $\tau=\omega t$.  In practice, this only then has application in problems with periodic solutions.  Thus for the Van der Pol oscillator, it can help discovery the correct frequency of the limit cycle, but it is incapable of characterizing the transient dynamics leading to the limit cycle.

Multiple scale expansions circumvent the fundamental limitations of the Poincare-Lindsted method by defining a second slow time-scale $\tau=\epsilon t$ and letting solutions vary on both fast and slow scales so that $u(x,t)\rightarrow u(x,t,\tau)$.  For spatial-temporal systems, a multiple-scale perturbation expansion takes the form
\begin{equation}
u=u_{0}(x,t)+\epsilon u_{1}(x,t,\tau) + \epsilon^2 u_2 (x,t,\tau) + \cdots
\end{equation}
where $\tau=\epsilon t$ is the slow variable dependence.  The new slow variable $\tau$ is treated as an independent variable in the expansion procedure.

Inserting this expansion into (\ref{eq:pert_general}) and
collecting terms at each order of the perturbation theory then gives
\begin{subeqnarray}
 &&  \frac{\partial u_0}{\partial t} = N (u_0, {u_0}_x, {u_0}_{xx}, \cdots, \mu)  \\
 &&  \frac{\partial u_1}{\partial t} = L (u_0) u_1 + F_1 (u_0, {u_0}_\tau)  \\
 &&  \frac{\partial u_2}{\partial t} = L (u_0) u_2 + F_2 (u_0, u_1, {u_0}_\tau, {u_1}_\tau) .
\end{subeqnarray}
In this case, the leading order solution is given by
\begin{equation}
  u(x,t)= u_0 (x,t,A,B, \cdots)
\end{equation}
where the parameters $A(\tau), B(\tau) \cdots$ parametrize the leading order-solution and its dependence on the slow time scale $\tau$, i.e. they are no longer constants, rather they depend on the slow time scale $\tau$ explicitly.

The Fredholm-Alternative theorem will result in equations for the slow evolution of the parameters so that
\begin{subeqnarray}
 &&  A_{\tau} = f_1(A,B,\cdots)  \\
 &&  B_{\tau} = f_2 (A,B,\cdots)  \\
 &&  \vdots .
\end{subeqnarray}
The slow dependence of the variables, like the Poincare-Lindstedt method, allows one to satisfy the solvability conditions in order to remove secular growth terms.  The multiple-scale expansion is even more flexible and versatile than the Poincare-Lindstedt method as will be shown in the forthcoming examples.

\subsection*{Perturbed, linear damped oscillator}

To see how the multiple scale expansion works, consider the weakly damped linear oscillator
\begin{equation}
   u_{tt} + 2 \epsilon u_t + (1+\epsilon) u = 0
   \label{eq:ms_oscillator}
\end{equation}
with $u(0)=\alpha$ and $u_t(0)=0$.  This is a linear problem that can be solved exactly to yield the solution
\begin{equation}
   u(t) = A \exp (-\epsilon t) \cos \sqrt{1+\epsilon - \epsilon^2} t + \frac{\epsilon A}{\sqrt{1+\epsilon - \epsilon^2}} \exp (-\epsilon t) \sin \sqrt{1+\epsilon - \epsilon^2} t .
   \label{eq:slow_exact}
\end{equation}
There are three key observations about this solution:  (i) the damping occurs on the scale $\epsilon t$, (ii) the leading-order frequency shift also occurs on time $\epsilon t$, and (iii) the leading order oscillations occurs on a time scale $t$.  Given that all the critical dynamics happen on timescales $t$ and $\epsilon t$, we define the slow scale
\begin{equation}
   \tau=\epsilon t
\end{equation}
and define our solution to be dependent upon both slow and fast times scales so that
\begin{equation}
   u= u(t,\tau) .
\end{equation}
These two time scales are considered as independent variables.  

The multiple scales assumption requires an update to the derivatives via the chain rule:
\begin{subeqnarray}
 &&  u_t \rightarrow u_t + u_\tau \tau_t = u_t + \epsilon u_\tau \\
 &&  u_{tt} = u_{tt} + 2 \epsilon u_{t\tau} + \epsilon^2 u_{\tau\tau} . 
\end{subeqnarray}
This then finally allows us to make the perturbation expansion
\begin{equation}
u=u_{0}(t,\tau)+\epsilon u_{1}(t,\tau) + \epsilon^2 u_2 (t,\tau) + \cdots
\label{eq:ms_expand}
\end{equation}
which can be inserted into the governing equation (\ref{eq:ms_oscillator}).   The various powers of $\epsilon$ can then collected.

The hierarchy of linear models generated by the perturbation expansion result in
\begin{subeqnarray}
 &  \hspace{-.30in} O(1) & {u_0}_{tt} + \! u_0 \!=\! 0 \,\, \mbox{with} \,\, u_0(0,0)\!=\!\alpha, \, {u_0}_t (0,0)\!=\!0 \\
 & \hspace{-.30in} O(\epsilon) \hspace{.0in} &  {u_1}_{tt} + \!u_1 \!=\! - \left( u_0 + 2 {u_0}_t + 2 {u_0}_{t\tau}  \right)  
    \,\, \mbox{with} \,\, u_1(0,0)\!=\!0, \, {u_1}_t (0,0)\!=\!-{u_0}_\tau (0,0) \\
 & \hspace{-.30in} O(\epsilon^2) \hspace{.0in} &  {u_2}_{tt} +\! u_2 \!=\! - \left( u_1 + 2 {u_1}_t + 2 {u_1}_{t\tau} +2 {u_0}_\tau + {u_0}_{\tau\tau} \right) \,\, \mbox{with} \,\, u_2(0,0)\!=\!0, \, {u_2}_t (0,0)\!=\!-{u_1}_\tau (0,0) .
\end{subeqnarray}
The solution for this linear hierarchy of problems shows how the multiple scale decomposition can be leveraged to approximate the solution.

At leading order, the solution can be easily obtained
\begin{equation}
  u_0= A(\tau) \cos t + B(\tau) \sin t
  \label{eq:ab_slow}
\end{equation}
where $A(0)=\alpha$ and $B(0)=0$.   Thus the coefficients in front of the sinusoidal terms are no longer constants, but rather functions of the slow time variable $\tau$.  At $O(\epsilon)$, the perturbation theory yields
\begin{eqnarray}
    {u_1}_{tt} + u_1 &=&  - \left( u_0 + 2 {u_0}_t + 2 {u_0}_{t\tau}  \right) \\
     &=& - (A+2B+2B_\tau) \cos t + (-B+2A+2A_\tau) \sin t .
\end{eqnarray}
The terms on the right hand side are in the null space of the adjoint operator.  To remove them, this then requires that
\begin{subeqnarray}
  && A_\tau = \frac{1}{2}B - A \\
 && B_\tau = - B -\frac{1}{2} A
 \label{eq:multi_AB}
\end{subeqnarray}
with $A(0)=\alpha$ and $B(0)=0$.  

Multiplying the top equation by $2A$ and the bottom equation by $2B$ gives
\begin{subeqnarray}
  && 2AA_\tau = AB - 2A^2 \\
 && 2BB_\tau = -2 B^2 - AB
\end{subeqnarray}
which upon adding the equations together yields
\begin{equation}
  (A^2+B^2)_\tau = -2 (A^2 + B^2) .
\end{equation}
The solution for this first order differential equation is given by
\begin{equation}
  A^2 +B^2 = C\exp (-2\tau) 
\end{equation}
with the initial condition $A^2(0) + B^2(0)=\alpha^2 = C$.   A convenient representation of the solution is given by
\begin{subeqnarray}
  && A(\tau) = \alpha \exp(-\tau) \cos \xi(\tau) \\
  && B(\tau) = \alpha \exp(-\tau) \sin \xi(\tau)
\end{subeqnarray}
so that we only need to determine the $\xi(\tau)$.  This can be achieved by plugging in this form of solution into (\ref{eq:multi_AB}a) which gives $\alpha(2\xi_tau +1)\sin \xi =0$.  This then gives $\xi_\tau=-1/2$ and $\xi=\tau/2$.  Using these results in (\ref{eq:ab_slow}) gives
\begin{equation}
  u_0(t,\tau) = \alpha \exp(-\tau) \left( \cos\frac{\tau}{2} \cos t - \sin\frac{\tau}{2} \sin t \right) = \alpha \exp(-\tau) \cos (t+\tau/2) .
\end{equation}
This is the leading order solution to the problem which accounts for both the frequency shift and the slow damping dynamics. 

To make comparison of this perturbative solution to the exact solution (\ref{eq:slow_exact}), consider that the second term in the exact solution is of $O(\epsilon)$.  Thus the leading order solution to the exact problem is given by
\begin{equation}
   u\approx \alpha \exp(-\epsilon t) \cos \sqrt{1+\epsilon - \epsilon^2} t + O(\epsilon).
\end{equation}
The Taylor series expansion for the frequency is given by
\begin{equation}
   \sqrt{1+\epsilon - \epsilon^2} = 1+ \frac{\epsilon}{2} - \frac{5\epsilon^2}{8} + \cdots
\end{equation}
so that 
\begin{equation}
   u\approx \alpha \exp(-\epsilon t) \cos (1+\epsilon/2) t + O(\epsilon).
\end{equation}
By replacing $\tau=\epsilon t$, the approximation of the exact solution matches that derived for $u_0(t,\tau)$ of the multiple scale perturbation theory.

\newpage
\section*{Lecture 13: The Van der Pol Oscillator}
We now apply the multiple scales expansion method to the Van der Pol oscillator.  This problem was previously considered from the point of view of Poincare-Lindsted.  Unfortunately, Poincare-Lindsted was only able to model the limit cycle.  Specifically, the method can only accommodate frequency shifts that occur due to the nonlinearity.  The transient behavior in the system is completely missed by the Poincare-Lindsted method.  Multiple scales overcomes these limitations as it provides a more generic framework for handling perturbations.

We can consider the Van der Pol oscillator which is modeled by
\begin{equation}
  u_{tt}  + \epsilon (u^2-1) u_t + u = 0
\end{equation}
with $u(0)=\alpha$ and $u_t(0)=0$.  We again define the slow time scale $\tau=\epsilon t$ and let $u(t)\rightarrow u(t,\tau)$.  The chain rule now gives the following equation for the Van der Pol oscillator in the two time scales
\begin{equation}
  \left(  u_{tt}  +2\epsilon u_t +\epsilon^2 u_{\tau\tau}  \right)  + \epsilon (u^2-1) \left( u_t + \epsilon u_\tau \right) + u = 0 .
\end{equation}
The expansion (\ref{eq:ms_expand}) can now be used.  Note that in addition to expanding the governing equations, the initial conditions are also expanded.  Thus the initial derivative is modified to $u_t(0)=0\rightarrow 
u_t(0)+\epsilon u_\tau(0) =0$.  This gives the hierarchy of linear problems
\begin{subeqnarray}
 &  \hspace{-.10in} O(1) & {u_0}_{tt} + \! u_0 \!=\! 0 \hspace*{.3in} \mbox{with} \,\, u_0(0,0)\!=\!\alpha, \, {u_0}_t (0,0)\!=\!0 \\
 & \hspace{-.10in} O(\epsilon) \hspace{.0in} &  {u_1}_{tt} + \!u_1 \!=\! - \left( 2 {u_0}_{t\tau} +(u_0^2-1){u_0}_t  \right)  
    \hspace*{.3in} \mbox{with} \,\, u_1(0,0)\!=\!0, \, {u_1}_t (0,0)\!=\!-{u_0}_\tau (0,0) \\
 & \hspace{-.10in} O(\epsilon^2) \hspace{.0in} &  {u_2}_{tt} +\! u_2 \!=\! - \left(   2 {u_1}_{t\tau} + (u_0^2-1){u_1}_t
    +2 u_0 u_1 {u_0}_t +(u_0^2-1){u_0}_\tau + {u_0}_{\tau\tau} \right) \nonumber \\
    &&  \hspace*{1in} \mbox{with} \,\, u_2(0,0)\!=\!0, \, {u_2}_t (0,0)\!=\!-{u_1}_\tau (0,0) .
\end{subeqnarray}
The leading order solution produces the solution
\begin{equation}
   u_0 (t,\tau)=A(\tau) \cos t + B(\tau) \sin t
\end{equation}
where $A(0)=\alpha$ and $B(0)=0$.

At the next order, the Fredholm-Alternative theorem must be applied.  Thus it is critical to determine which terms are in the null space of the adjoint operator.  The right hand side can be found to give
\begin{eqnarray}
  2 {u_0}_{t\tau} +(u_0^2-1){u_0}_t &=& ( -2A_\tau - A^3 + 2AB^2 + A) \sin t + (A^3 -3AB^2) \sin^3 t \nonumber \\
   & & (2B_\tau - 2A^2  B + B^3 - B ) \cos t + (3 A^2 B - B^3)   \nonumber \\
   &=& \left( -2A_\tau -  \frac{A^3}{4} - \frac{AB^2}{4} + A \right) \sin t -\frac{1}{4} \left( A^3 - 3AB^2 \right) \sin 3 t
           \nonumber \\
     & &  \left( 2B_\tau + \frac{B^3}{4} + \frac{A^2 B}{4} -B \right) \cos t -\frac{1}{4} \left( B^3 - 3A^2 B \right) \cos 3 t
\end{eqnarray}
where the trigonometric relations $\sin^3 x = (3/4)\sin x-(1/4)\sin 3x$ and $\cos^3 x = (3/4)\cos x+(1/4)\cos 3x$ have been used to appropriately extract the $\sin t$ and $\cos t$ terms which are in the null space of the adjoint operator.  

Eliminating the terms that produce secular growth gives the slow evolution equations for the leading order amplitudes $A$ and $B$:
\begin{subeqnarray}
  && 2A_\tau +  \frac{A^3}{4} + \frac{AB^2}{4} - A =0 \\
  && 2B_\tau + \frac{B^3}{4} + \frac{A^2 B}{4} -B  =0 .
  \label{eq:ABslow}
\end{subeqnarray}
Multiplying the first equation by $A$ and the second by $B$ and adding gives the evolution dynamics
\begin{equation}
   \rho_\tau + \frac{1}{4} \rho^2 - \rho = 0
\end{equation}
where $\rho=A^2 + B^2$.  The is a first order equation for $\rho$ that can be solved using standard integrating factor methods outlined at the beginning of the differential equations lectures.  This gives
\begin{equation}
   \rho(\tau)= \frac{4\alpha^2}{\alpha^2 + (4-\alpha^2) \exp(-\tau)}
\end{equation}
where $\rho(0)=A^2(0)+B^2(0) = \alpha^2$. To determine the evolution of $A(\tau)$ and $B(\tau)$, the solution $\rho(\tau)$ can be used in either of the equations (\ref{eq:ABslow}).   Indeed, inserting this into (\ref{eq:ABslow}a) leads to a first order differential equation which can be readily solved for $A(\tau)$.  This can then be used in  
(\ref{eq:ABslow}b).   This gives in total
\begin{subeqnarray}
  && A(\tau) = \frac{2\alpha}{(\alpha^2 + (4-\alpha^2) \exp(-\tau) )^{1/2}} \\
  && B(\tau)= 0 .
\end{subeqnarray}
Thus the slow evolution of the parameters of the leading order solution gives the leading order solution
\begin{equation}
   u= \frac{2\alpha}{(\alpha^2 + (4-\alpha^2) \exp(-\tau) )^{1/2}} \cos t + O(\epsilon) .
\end{equation}
We note that the solution captures the transient, slow dynamics that occurs from the initial conditions $u(0)=\alpha$ to the limit cycle $u\rightarrow 2 \cos t$  when $t\rightarrow \infty$.  At higher-orders, frequency corrections can also be acquired using the multiple scales expansion.

\newpage
\section*{Lecture 14:  Boundary Layer Theory}

Asymptotics provides a principled mathematical architecture for turning nonlinear problems into a hierarchy of linear problems.  In what has been considered so far, which are time-dependent differential equations, Poincare-Lindsted and multiple scale techniques enable characterization of the changing temporal dynamics on different time scales.  We now turn our attention to boundary value problems.  And specifically, those that are {\em singular} in nature.  

To illustrate the singular nature of a problem, consider the boundary value problem
\begin{equation}
  \epsilon u_{xx} + u_x + u = 0
\end{equation}
with $u(0)=0$ and $u(l)=A$.  Although this looks very much like our previous perturbed problems, it is strikingly different due to which term is multiplied by the small parameter $\epsilon$.  
To study this problem, we first attempt a regular perturbation expansion
\begin{equation}
   u(x) = u_0(x) + \epsilon u_1(x) + \epsilon^2 u_2(x) + \cdots .
   \label{eq:bl_regular}
\end{equation}
The resulting hierarchy of linear models is given by
\begin{subeqnarray}
   && O(1)  \hspace*{.5in}  {u_0}_{x} +  u_0 =0 \\
   && O(\epsilon)  \hspace*{.5in}   {u_1}_{x} +  u_1 = - {u_0}_{xx} 
\end{subeqnarray}
with 
\begin{subeqnarray}
   && u_0(0)=0, \,\, {u_0}(l) =A \\
   && u_1(0)=0, \,\, {u_1} (l)=0 .
\end{subeqnarray}
The leading order problem is a first order differential equations whose solution is given by
\begin{equation}
   u_0 = C \exp (-x) .
\end{equation}
This solution has a single free parameter $C$, yet must satisfy two boundary conditions given by 
$u_0(0)=0$ and ${u_0}(l) =A$.   This then is an overdetermined situation since the first boundary condition requires $C=0$ while the second requires $C=A\exp(l)$.

The leading order solution immediately leads to a contradiction, which is a hallmark feature of singular perturbation problems.   The failure of the regular perturbation expansion becomes obvious upon reflection.  Specifically, as $\epsilon\rightarrow 0$, the highest derivative vanishes.  But the number of linearly independent solutions is determined by the highest derivative.  So at leading order, the first order differential equation only admits a single linearly independent solution while the boundary conditions that need to be satisfied remains unchanged.   The loss of the highest derivative at leading order is what constitutes a singular perturbation problem.  Boundary layer theory will help us resolve this issue and generate an approximate solution to the problem.

\begin{figure}[t]
\begin{overpic}[width=0.65\textwidth]{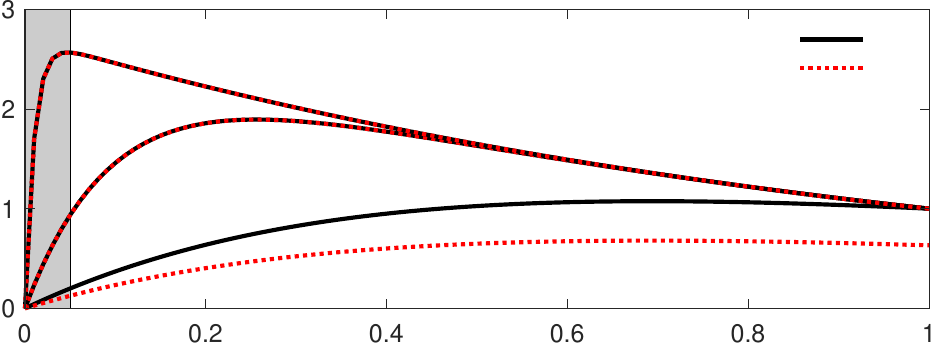}
\put(-6,30){$u(x)$}
\put(46,-3){space $x$}
\put(72,32){exact}
\put(65,28.5){approximate}
\put(9,32){$\epsilon=0.01$}
\put(12,17){$\epsilon=0.1$}
\put(37,15){$\epsilon=0.5$}
\end{overpic}
\caption{Singularly perturbed boundary value problem with a boundary layer at $x=0$.  The exact solution (black) is plotted against the boundary layer approximation (red dotted).  Note that as $\epsilon\rightarrow 0$, a transition region (shaded area) is formed near $x=0$ which is well characterized by the perturbation analysis.}
 \label{fig:bl_example1}
\end{figure}
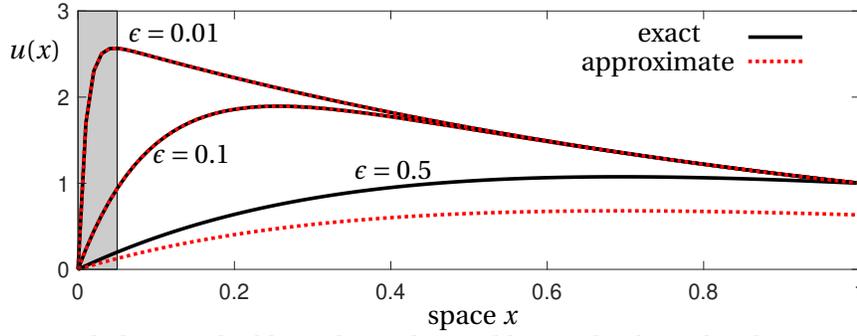

To illustrate the boundary layer method, consider the linear system
\begin{equation}
  \epsilon u_{xx} + (1+\epsilon) u_x + u = 0
\end{equation}
with $u(0)=0$ and $u(1)=1$.  This has the exact solution
\begin{equation}
   u(x) = \frac{\exp(-x) - \exp(-x/\epsilon)}{\exp(-1) - \exp(-1\epsilon)}.
\end{equation}
Figure~\ref{fig:bl_example1} shows the solution as a function of decreasing values of $\epsilon$.  The figure is suggestive:  there appear to be two distinct regions of the solution.  The two regions are as follows:  (i) the {\em outer region} which takes up most of the interval $x\in[0,1]$ and shows simple exponential decay, and (ii) the {\em inner region} which dominates a small part of the interval near $x=0$.  The region near $x=0$ is also known as the boundary layer or transition region.\\

\noindent {\em I.   The Outer Problem ($\delta \ll x \leq 1$) -- }
We again attempt our perturbation expansion (\ref{eq:bl_regular}) for the outer region which results in the hierarchy
\begin{subeqnarray}
   && O(1)  \hspace*{.5in}  {u_0}_{x} +  u_0 =0 \\
   && O(\epsilon)  \hspace*{.5in}   {u_1}_{x} +  u_1 = - {u_0}_{xx} - {u_0}_x
\end{subeqnarray}
with 
\begin{subeqnarray}
   && {u_0}(1) =1 \\
   && {u_1} (l)=0 .
\end{subeqnarray}
Note that in the outer region, only the boundary condition at $x=1$ is applied to the solution.  Given that only the right side boundary conditions are applied, the solution is no longer ill-posed (overdetermined).  The leading order solution can then be found to be
\begin{equation}
  u_0(x) = \exp (1-x) .
\end{equation}
Moreover, we can also find that the higher-order terms vanish so that $u_1(x) = u_2(x) = \cdots = 0$.\\

\noindent {\em II.   The Inner Problem ($0\leq x < \delta \ll 1$) -- } Since the outer problem characterizes most of the right side of the spatial interval, the inner problem must characterize the rapid transition region on the left.  To do this, we introduce a {\em stretching} variable which provides a local variable allowing the boundary region to be of $O(1)$ width.  Thus we introduce
\begin{equation}
  \xi = \frac{x}{\epsilon} .
\end{equation}
This coordinate transformation acts like a microscope to more effectively zoom into the boundary layer of interest at the left side of the domain.  The change of variable gives by following by chain rule
\begin{subeqnarray}
  &&   u_x = u_\xi \xi_x = \frac{1}{\epsilon} u_\xi \\
  &&   u_{xx} = \frac{1}{\epsilon^2} u_{\xi\xi}.
\end{subeqnarray}
In the new inner variables, the governing equation is now
\begin{equation}
   u_{\xi\xi} + (1+\epsilon) u_\xi + \epsilon u = 0
\end{equation}
with the left boundary condition $u(x=0)=u(\epsilon \xi =0) = u(0)=0$.

A perturbation expansion is now performed in the inner region so that
\begin{equation}
   u(\xi) = u_0(\xi) + \epsilon u_1(\xi) + \epsilon^2 u_2(\xi) + \cdots .
\end{equation}
This gives the hierarchy of equations in the inner region to be
\begin{subeqnarray}
   && O(1)  \hspace*{.5in}  {u_0}_{\xi\xi} +  {u_0}_\xi =0 \\
   && O(\epsilon)  \hspace*{.5in}   {u_1}_{\xi\xi} +  {u_1}_\xi = - {u_0}_\xi - u_0
\end{subeqnarray}
with 
\begin{subeqnarray}
   && {u_0}(0) =0 \\
   && {u_1} (0)=0 .
\end{subeqnarray}
The leading order problem is a first order differential equation in ${u_0}_\xi$.  The leading order solution can be found to be
\begin{equation}
  u_0(\xi) = A (1 - \exp(-\xi) )
\end{equation}
where $A$ is an undetermined constant.  Interestingly, since the leading order equation in the inner region is second order, two linearly independent solutions exist.  But there is only one boundary constraint at $x=0$, requiring us to determine the value of $A$ by other means.\\

\noindent {\em III.   Matched Asymptotics and the Uniform Solution -- } 
The final step in the boundary layer theory is to merge the inner and outer solution in some consistent fashion.  Thus far, the solutions found are
\begin{subeqnarray}
  && \mbox{outer solution:}  \hspace*{.2in}  u_{\mbox{\small out}} = \exp(1-x)   \hspace*{.5in} \delta \ll x \leq 1 \\
  && \mbox{inner solution:}  \hspace*{.2in} u_{\mbox{\small in}} =A (1-\exp(-\xi)) \hspace*{.5in} 0\leq x < \delta \ll 1 .
\end{subeqnarray}
Matching these two regions requires that
\begin{equation}
   \lim_{x\rightarrow 0} u_{\mbox{\small out}} = \lim_{\xi\rightarrow\infty} u_{\mbox{\small in}} = u_{\mbox{\small match}}.
\end{equation}
Thus as the two solutions encroach on each other's regions, a consistent solution would have them match in a self-consistent way.

The {\em uniform solution} to the boundary layer problem is defined in the following way
\begin{equation}
  u_{\mbox{\small unif}} =  u_{\mbox{\small in}} -  u_{\mbox{\small out}} - u_{\mbox{\small match}} .
\end{equation}
The matching conditions for this specific example can be determined as follows
\begin{subeqnarray}
  &&  \lim_{x\rightarrow 0} u_{\mbox{\small out}} =  \lim_{x\rightarrow 0}  \exp(1-x) = e \\
  &&  \lim_{\xi\rightarrow\infty} u_{\mbox{\small in}} = \lim_{\xi\rightarrow\infty} A(1-\exp(-\xi)) = A
\end{subeqnarray}
so then
\begin{equation}
   A=e
\end{equation}
and the uniform solution is given by
\begin{equation}
  u_{\mbox{\small unif}} =    \exp(1-x) - \exp(1-x/\epsilon) .
\end{equation}
This is the leading order boundary layer solution.  Note that it matches the exact analytic solution to leading order, which is an encouraging sign that the method produces a self-consistent result.
Figure~\ref{fig:bl_example1} shows the approximate boundary layer solution as a function of decreasing values of $\epsilon$.  The approximation does well in agreeing with the exact solution.

\newpage
\section*{Lecture 15:  Dominant Balance, Distinguishing Limits and Matched Asymptotic Expansions}

Our goal now is to broaden our understanding of boundary layer theory.  In particular, where does the boundary layer appear?  And what is its asymptotic width?  Both these issues were assumed to have been known in the last example.  In general, the concept of balancing the dominant terms in physical systems is important for understanding the underlying physics of a system~\cite{callaham2021learning}.  In what follows, we determine a principled way to determine both of these important issues.  We consider the singularly perturbed boundary value problem
\begin{equation}
  \epsilon u_{xx} + b(x) u_x + c(x) u = 0 
\end{equation}
on the domain $x\in[0,1]$ and with the boundary conditions $u(0)=A$ and $u(1)$=B.
We will further assume $b(x)\neq 0$ on the interval $x\in[0,1]$.  As before, we will separate the problem solution into three portions (i) the outer solution, (ii) the inner solution, and (iii) the matching region.\\

\noindent {\em I.   The Outer Problem  -- }
We attempt our regular perturbation expansion (\ref{eq:bl_regular}) for the outer region which results in the hierarchy
\begin{subeqnarray}
   && O(1)  \hspace*{.5in}  b(x) {u_0}_{x} + c(x) u_0 =0 \\
   && O(\epsilon)  \hspace*{.5in}  b(x) {u_1}_{x} +  c(x) u_1 = - {u_0}_{xx}  .
\end{subeqnarray}
At this point, we do not know where the boundary layer is.  This will determine wether to impose the boundary condition $u_0(0)=A$ or $u_0(1)=B$.  Only one of these can be imposed for a well-posed solution of the outer problem.  For now, we will leave the position of the boundary layer an open question.  At the matching stage, we will be able to determine its location.

Regardless of wether we know the boundary layer location, we can easily solve the leading order problem since it is a simple linear first order differential equation with non-constant coefficients
\begin{equation}
  {u_0}_x + \frac{c(x)}{b(x)} u_0 = 0
\end{equation}
whose solution is given by
\begin{equation}
  u_0(x) = C \exp \left[ \int \frac{c(x)}{b(x)} dx  \right] 
\end{equation}
where $C$ is an undetermined constant of integration.  This is the outer solution.  Once we are given the specific form of $b(x)$ and $c(x)$, the integral can be evaluated.  Additionally, once the location of the boundary layer is known, then $C$ can be evaluated.\\

\noindent {\em II.   The Inner Problem  -- } The inner problem must characterize the rapid transition region of the boundary layer which is located (for this case when $b(x)\neq 0$) on either the left ($x=0$) or right ($x=1$) of the domain.  To do this, we introduce a {\em stretching} variable which provides a local variable allowing the boundary region to be of $O(1)$ width.  Thus we introduce
\begin{equation}
  \xi = \frac{x}{\delta} .
\end{equation}
This coordinate transformation acts like a microscope to more effectively zoom into the boundary layer of interest.  Note that we have at this point assumed we don't know the relative scaling of $\delta\ll1$ to $\epsilon\ll 1$.  This will need to be determined.  Moreover, if the boundary layer is on the right, we would need to expand so that  $\xi = \frac{1-x}{\delta}$.  This would then create a microscopic view around the boundary at $x=1$.

The change of variable gives by following by chain rule
\begin{subeqnarray}
  &&   u_x = u_\xi \xi_x = \frac{1}{\delta} u_\xi \\
  &&   u_{xx} = \frac{1}{\delta^2} u_{\xi\xi}.
\end{subeqnarray}
In the new inner variables, the governing equation is now
\begin{equation}
  \frac{\epsilon}{\delta^2}  u_{\xi\xi} + \frac{b(\delta \xi)}{\delta} u_\xi + c(\delta \xi) u = 0
  \label{eq:bl_dominant}
\end{equation}
where we need to determine $\delta(\epsilon)$.  This would determine the width of the boundary layer.  This is done by a {\em dominant balance}, or {\em distinguished limits}, analysis.  

Given that there are three terms in (\ref{eq:bl_dominant}), there are three potential dominant balance regimes.\\

\noindent {\em (i) Distinguished Limit:  $\epsilon \gg \delta$}\\

In this limit, the leading order solution to (\ref{eq:bl_dominant}) is given by
\begin{equation}
  u_{\xi\xi}=0 
\end{equation}
which upon integrating twice yields the solution $u_{\mbox{\small in}}=D\xi + E$.
The inner solution must be matched to the outer solution in the limits of $\xi\rightarrow\infty$ so then
\begin{equation}
  \lim_{\xi\rightarrow \infty} u_{\mbox{\small in}} =  \lim_{\xi\rightarrow \infty}  D\xi + E = \infty .
\end{equation}
Thus solutions blow up in this distinguished limit and we can rule it out of consideration.\\

\noindent {\em (ii) Distinguished Limit:  $\epsilon \ll \delta$}\\

In this limit, the leading order solution to (\ref{eq:bl_dominant}) is given by
\begin{equation}
  b(\delta\xi) u_{\xi}=0 
\end{equation}
which yields the solution $u_{\mbox{\small in}}=D$.  In this case, there is no problem with the solution being bounded for matching purposes.  However, the constant must satisfy both a boundary condition and the matching conditions, leading to an overdetermined system.  This allows us to rule out this distinguished limit.\\

\noindent {\em (iii) Distinguished Limit:  $\epsilon \sim \delta$}\\

In this limit, we can set $\delta=\epsilon$ for convenience.  
The leading order solution to (\ref{eq:bl_dominant}) is given by
\begin{equation}
  u_{\xi\xi} + b(\epsilon \xi) u_\xi + \epsilon c(\epsilon \xi) u = 0 .
\end{equation}
An perturbation expansion (\ref{eq:bl_regular}) in the new stretched variable gives the hierarchy
\begin{subeqnarray}
   && O(1)  \hspace*{.5in}  {u_0}_{\xi\xi} +  b(\epsilon\xi) {u_0}_\xi =0 \\
   && O(\epsilon)  \hspace*{.5in}   {u_1}_{\xi\xi} + b(\epsilon\xi) {u_1}_\xi = - c(\epsilon \xi) u_0
\end{subeqnarray}
where the boundary conditions are yet to be determined.

The leading order solution will determine the boundary layer location.  Consider the coefficient $b(\epsilon\xi)$ in the leading order equation.  Suppose we define
\begin{equation}
   b_0 = b(\epsilon \xi) = \left\{  \begin{array}{l} b(0) + O(\epsilon) \\
   b(1) + O(\epsilon) \end{array}  \right.
\end{equation}
where $b_0$ simply takes on the value of the function $b(x)$ at the left or right of the domain.  In this case, the distinguished limit takes the asymptotic approximation
\begin{equation}
  {u_0}_{\xi\xi} +  b_0 {u_0}_\xi =0
\end{equation}
which has the solution
\begin{equation}
  u_0 = D \exp(-b_0 \xi) + E .
\end{equation}
Note that this solution has two free constants:  one which can be used to satisfy the boundary condition and a second to satisfy the matching conditions.    The sign of $b_0$ will ultimately determine the location of the boundary layer as will be shown in the matching.\\

\noindent {\em III.   Matching  -- } The final part of the boundary layer theory is constructing the matching region.  The method of matched asymptotic expansions is critical for this task.  Specifically, the inner and outer solutions must match asymptotically for a self-consistent solution to exist.  This self-consistency will determine the position of the boundary layer, which at least in this problem is of thickness $\delta=\epsilon$ by our distinguished limit computation.
 
Recall that the inner and outer solutions are given by 
\begin{subeqnarray}
  && \mbox{outer solution:}  \hspace*{.2in}  u_{\mbox{\small out}} = C \exp \left[ \int \frac{c(x)}{b(x)} dx  \right]  \\
  && \mbox{inner solution:}  \hspace*{.2in} u_{\mbox{\small in}} = D \exp(-b_0 \xi) + E .
\end{subeqnarray}
Two possible scenarios exist for matching
\begin{subeqnarray}
 && \mbox{layer at} x=0:  \hspace*{.2in}  \lim_{\xi \rightarrow \infty} u_{\mbox{\small in}} =
   \lim_{x\rightarrow 0} u_{\mbox{\small out}} \\
 && \mbox{layer at} x=1:  \hspace*{.2in}  \lim_{\xi \rightarrow -\infty} u_{\mbox{\small in}} =
   \lim_{x\rightarrow 1} u_{\mbox{\small out}}.
\end{subeqnarray}
The matching conditions will depend largely on the sign of $b_0$.  Specifically, the possible cases are the following:\\

\noindent {\em Case 1:  $b(x)>0$}\\

In this case $b_0>0$ and  the boundary layer must be at $x=0$ since otherwise the inner solution $u_{\mbox{\small in}} = D \exp(-b_0 \xi) + E$ blows up otherwise.  Specifically, the following holds if the boundary is located at the two key locations:
\begin{subeqnarray}
 &&  x=0:  \hspace*{.2in}  \lim_{\xi \rightarrow \infty} u_{\mbox{\small in}} =
     \lim_{\xi \rightarrow \infty} D\exp(-b_0\xi) + E = E \\
 && x=1:  \hspace*{.2in}  \lim_{\xi \rightarrow -\infty} u_{\mbox{\small in}} =
    \lim_{\xi \rightarrow -\infty} D\exp(-b_0\xi) + E = \infty .
\end{subeqnarray}
Thus a bounded solution is only possible at $x=0$ for $b(x)>0$.\\

\noindent {\em Case 2:  $b(x)<0$}\\

In this case $b_0<0$ and the boundary layer must be at $x=1$ since otherwise the inner solution $u_{\mbox{\small in}} = D \exp(-b_0 \xi) + E$ blows up otherwise.  Specifically, the following holds if the boundary is located at the two key locations:
\begin{subeqnarray}
 &&  x=0:  \hspace*{.2in}  \lim_{\xi \rightarrow \infty} u_{\mbox{\small in}} =
     \lim_{\xi \rightarrow \infty} D\exp(b_0\xi) + E = \infty \\
 && x=1:  \hspace*{.2in}  \lim_{\xi \rightarrow -\infty} u_{\mbox{\small in}} =
    \lim_{\xi \rightarrow -\infty} D\exp(b_0\xi) + E = E .
\end{subeqnarray}
Thus a bounded solution is only possible at $x=1$ for $b(x)<0$.\\

The boundedness of the solutions thus determines the location of the boundary layer and we can finalize our uniform solution.   Letting $b(x)>0$, for example, allows us to satisfy the boundary conditions since $u_{\mbox{\small out}} (1)=B$ and $u_{\mbox{\small in}}(0)=A$.  The matching constant is $E$ from the asymptotic limits above, which gives the uniform solution
\begin{equation}
u_{\mbox{\small unif}} =  D\exp(-b_0 \xi) +  C \exp \left[ \int \frac{c(x)}{b(x)} dx  \right] .
\end{equation}
This is the leading order solution to the singularly perturbed boundary value problem.  Once $b(x)$ and $c(x)$ are prescribed, the constants $D$ and $C$ can be evaluated.  It is constructed from two spatial scales, $x$ and $\xi$, much like our multiple scales differential equations considered earlier.    Interestingly, when $b(x)=0$ in the interval $x\in[0,1]$, internal boundary layers can form as will be shown.

\subsection*{Multiple boundary Layers.}
This example demonstrates that multiple boundary layers can form with different characteristic widths.  Thus distinguished limits must be found for each layer, and some evidence suggesting that a boundary layer will form should be provided.  Consider the following
\begin{equation}
   \epsilon u_{xx} - x^2 u_x - u = 0
\end{equation}
with $u(0)=u(1)=1$.  This is a singular problem since the $\epsilon$ is in front of the highest derivative.  Further $b(x)=-x^2$ and $c(x)=-1$.  Since $b(x)<0$ everywhere except at zero, we would expect a boundary layer to form at $x=1$.  But in this case, there may be a boundary layer at $x=0$ since $b(0)=0$.\\

\noindent {\em Distinguished limit near $x=0$}\\

Near $x=0$ we define the stretching variable
\begin{equation}
   \xi = \frac{x}{\delta}
\end{equation}
which gives the new governing equation
\begin{equation}
   \epsilon u_{\xi\xi} - \delta^3 \xi^2 u_{\xi} - \delta^2 u = 0 .
\end{equation}
The second term in this expansion is much smaller than the third terms since each is $O(\delta^3)$ and $O(\delta^2)$ respectively.  This gives the dominant balance as
\begin{equation}
   \epsilon u_{\xi\xi} - \delta^2 u = 0 .
\end{equation}
The only way to balance these two terms is to select
\begin{equation}
  \delta = \epsilon^{1/2} .
\end{equation}
Thus we can expect that there is potentially a boundary layer at $x=0$ with characteristic width $O(\epsilon^{1/2})$. \\ 

\noindent {\em Distinguished limit near $x=1$}\\

Near $x=1$ we define the stretching variable
\begin{equation}
   \xi = \frac{1-x}{\delta}
\end{equation}
 which gives the governing equation
 \begin{equation}
  \epsilon u_{\xi\xi} + \delta (1- \delta \xi)^2 u_\xi - \delta^2 u = 0
\end{equation}
In this case, the last term is much smaller than the second term, which gives the leading order equation
 \begin{equation}
  \epsilon u_{\xi\xi} + \delta  u_\xi = 0 .
\end{equation}
In order for the two terms to balance, the distinguished limit requires
\begin{equation}
  \delta = \epsilon .
\end{equation}
Thus we can expect a boundary layer at $x=1$ with characteristic width $O(\epsilon)$. \\ 

As before, we will separate the problem solution into three portions (i) the outer solution, (ii) the inner solution (actually there are two of them), and (iii) the matching region (again, two matching regions are found).\\

\noindent {\em III.   Outer Solution  -- }  The outer solution is found by 
simply expanding the solution $u=u_0 + \epsilon u_1 + \cdots$.  Inserting the regular perturbation expansion into the governing equation yields the leading order equation
\begin{equation}
   - x^2 {u_0}_x - u_0 = 0 
\end{equation}
whose solution is given by
\begin{equation}
   u_0 = C \exp (1/x) .
\end{equation}
At this point, we cannot determine the constant $C$ because there are boundary layers at both $x=0$ and $x=1$.  Regardless, we can see immediately there will be problems with this solution when we try to match it in the region $x\rightarrow 0$.

\noindent {\em II.   Inner Solution -- }  There are two boundary layers for this problem.  Thus each is treated separately as they are also of different characteristic widths.

We consider first the inner solution near $x=0$ and with $\xi=x/\epsilon^{1/2}$.  As already shown in our dominant balance arguments, the leading order solution for this inner solution is given by
\begin{equation}
   {u_0}_{\xi\xi} - u_0 = 0
\end{equation}
with $u_0(0)=1$.  This yields the solution $u_0(\xi)=A\exp(-\xi) + B\exp(\xi)$.  For matching, $\xi\rightarrow\infty$ which will require that $B=0$ for bounded solutions.  This leaves us to determine the value of $A$ by imposing the boundary condition $u_0(0)=1$.  Thus the leading order solution is $u_0 (\xi)=\exp(-\xi)$.

Near $x=1$ the stretching variable takes the form $\xi=(1-x)/\epsilon$.
A regular expansion yields the leading order solution
\begin{equation}
   {u_0}_{\xi\xi} + {u_0}_\xi = 0
\end{equation}
with $u_0(1)=1$.  This is a first order differential equation for ${u_0}_\xi$ which can be solved to yield $u_0(\xi)=A \exp(-\xi) + (1-A)$.  This is our second inner solution required for producing the uniform solution.

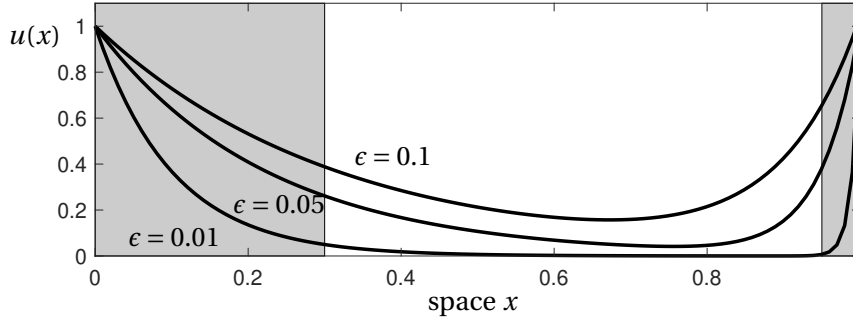
\begin{figure}[t]
\begin{overpic}[width=0.65\textwidth]{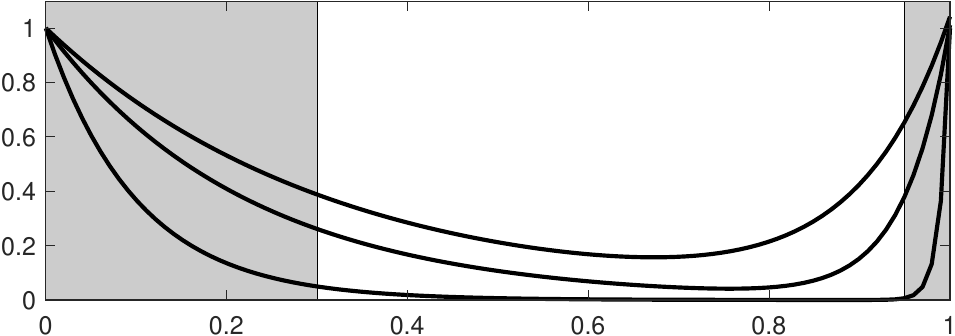}
\put(-6,30){$u(x)$}
\put(46,-3){space $x$}
\put(9,5){$\epsilon=0.01$}
\put(22,9){$\epsilon=0.05$}
\put(37,15){$\epsilon=0.1$}
\end{overpic}
\caption{Singularly perturbed boundary value problem with a boundary layer at $x=0$ and $x=1$.  Each boundary layer region (shaded) has a different distinguished limit.  The boundary layer at $x=0$ has width $O(\epsilon^{1/2})$ while the boundary layer at $x=1$ has width $O(\epsilon)$.  The outer solution reduces to $u_{\mbox{\small out}}=0$ at leading order. }
 \label{fig:bl_double}
\end{figure}

\noindent {\em III.   Matching  -- }  The matching must now be done at both the left and right of the domain.  Thus we have the following conditions that must be met asymptotically
\begin{eqnarray}
  &&  \mbox{at} x=0: \,\,\, \lim_{x\rightarrow 0} u_{\mbox{\small out}} =
                \lim_{x\rightarrow \infty} u_{\mbox{\small in (left)}} \\
 &&  \mbox{at} x=1: \,\,\, \lim_{x\rightarrow 1} u_{\mbox{\small out}} =
                \lim_{x\rightarrow \infty} u_{\mbox{\small in (right)}} 
  \end{eqnarray}
The first limit, $\lim_{x\rightarrow 0} u_{\mbox{\small out}}$, is bounded provided $C=0$ in the outer solution.  This then gives an outer solution of  $u_{\mbox{\small out}}=0$.  The second limit at $x=1$ then requires that 
$\lim_{x\rightarrow \infty} u_{\mbox{\small in (right)}}=0$, giving $A=1$.  This the uniquely determines all constants associated with the inner and outer solutions.
The uniform solution is constructed as
\begin{equation}
     u_{\mbox{\small unif}} = u_{\mbox{\small out}} + u_{\mbox{\small in (left)}} + u_{\mbox{\small in (right)}} - u_{\mbox{\small match (left)}} - u_{\mbox{\small match (right)}} .
\end{equation}
This gives in total the uniform solution at leading order
\begin{equation}
  u_{\mbox{\small unif}} = \exp \left(- \frac{x}{\epsilon^{1/2}} \right) + \exp \left(-\frac{1-x}{\epsilon} \right) .
\end{equation}
Figure~\ref{fig:bl_double} shows the asymptotic solution to this problem as a function of the parameter $\epsilon$.  Note the boundary layer on both sides and their different characteristic widths.

\subsection*{Internal Boundary Layers.}

Not only can multiple boundary layers develop, but they can also develop in the interior of a domain.  In this case, unlike the last example, two outer solutions are required to match one internal layer.  To demonstrate this, consider the 
singularly perturbed boundary value problem
\begin{equation}
  \epsilon u_{xx}  + 2x u_x + (1+x^2) u = 0
\end{equation}
on the domain $x\in[-1,1]$ and with boundary conditions $u(-1)=2$ and $u(1)=1$.
In this case, the function $b(x)=2x$ and $c(x)=(1+x^2)$.  Importantly, the function has a simple zero at $b(0)=0$ which is responsible for creating an internal layer in the domain.

Since $b(x)>0$ for $x>0$ and $b(x)<0$ for $x<0$, using our previous considerations, we can show that no boundary layers form at $x=\pm 1$.  Note however, that if the function $b(x)=2x$, then an internal layer would form at $x=0$ along with two boundary layers at $x=\pm 1$.  

\noindent {\em I.   Outer Solution  -- } We begin by constructing the outer solutions.  In this case, there is an outer solution for $x<0$ and one for $x>0$.  A regular perturbation expansion gives at leading order
\begin{equation}
   2x {u_0}_x + (1+x^2) {u_0} = 0 .
\end{equation}
This is a first order differential equation whose solution is given by
\begin{equation}
   {u_0} (x) = C x^{-1/2} \exp (-x^2/4) .
\end{equation}
Applying the boundary condition on the left ($u(-1)=2$) and right ($u(1)=1$) gives the two outer solutions 
\begin{subeqnarray}
  && u_{\mbox{\small out (left)}} = 2(-x)^{-1/2} \exp ((1-x^2)/4) \hspace*{.3in} x\in[-1,-\delta] \\
  && u_{\mbox{\small out (right)}} = x^{-1/2} \exp ((1-x^2)/4)  \hspace*{.3in} x\in[\delta,1].
\end{subeqnarray}
where the intervals to the left and right of $x=0$ are specified for their validity.  The parameter $\delta\ll 1$ must be determined by a distinguished limit in the inner layer.\\

\noindent {\em II.   Inner Solution -- }  The inner expansion is determined first by the stretching variable $\xi=x/\delta$. This gives
\begin{equation}
  \frac{\epsilon}{\delta^2} u_{\xi\xi}  + 2\xi u_\xi + (1+(\delta \xi)^2) u = 0 .
\end{equation}
The second derivative term cannot be eliminated in the boundary layer since this will simply produce the outer solution once again.  Thus in order to bring the dominant balance terms together so that the second derivative is not eliminated, we must choose $\delta=\epsilon^{1/2}$.  This gives at leading order the equation
\begin{equation}
   {u_0}_{\xi\xi}  + 2\xi {u_0}_\xi +  {u_0} = 0 .
\end{equation}
This is simplified by the transformation $u_0=w(\xi) \exp(-\xi^2/2)$ so that
\begin{equation}
   w_{\xi\xi} - \xi^2 w =0 .
\end{equation}
The solution of this differential equation can be expressed in terms of special functions.  In particular, the solution of the above equation is given by the {\em parabolic cylinder functions}.   This gives the inner solution
\begin{equation}
   w(\xi) = A D_{-1/2} (\xi) + B D_{-1/2} (-\xi) .
\end{equation}
Details of these functions can be found in Abramowitz and Stegun~\cite{abramowitz1948handbook}.  The inner and outer solutions can now be matched, but this is no longer an easy task given the parabolic cylinder function solutions.\\

\begin{figure}[t]
\begin{overpic}[width=0.65\textwidth]{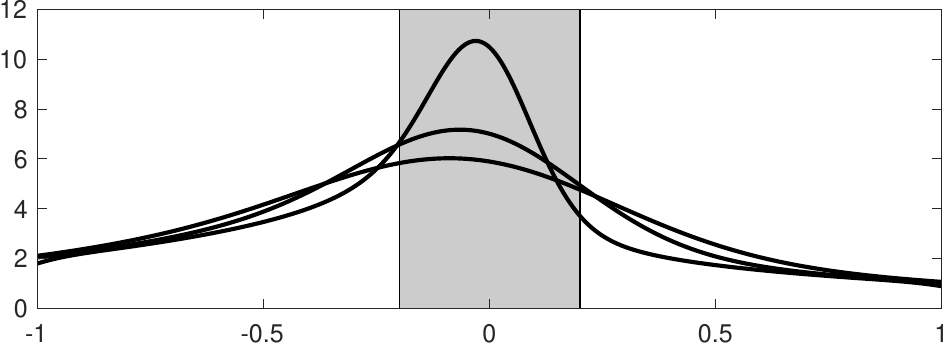}
\put(-8,30){$u(x)$}
\put(46,-3){space $x$}
\put(49,33){$\epsilon=0.02$}
\put(22,9){$\epsilon=0.2$}
\put(31,22){$\epsilon=0.1$}
\end{overpic}
\caption{Singularly perturbed boundary value problem with a boundary layer in the middle of the domain.  The boundary layer region (shaded) at $x=0$ has width $O(\epsilon^{1/2})$.  There are two outer solutions to the left and right of the boundary layer. }
 \label{fig:bl_middle}
\end{figure}
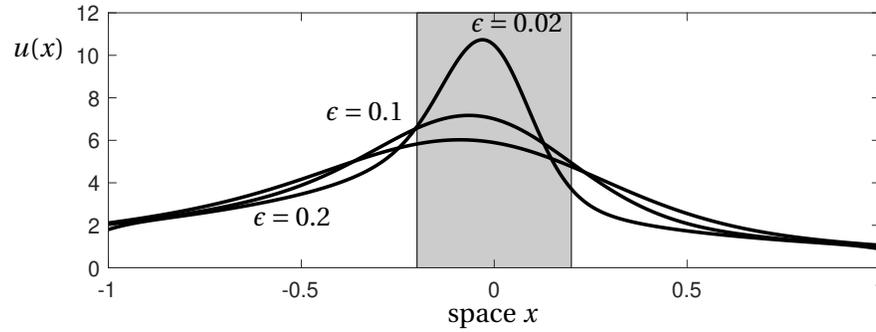

\noindent {\em III.   Matching  -- }
The matching of the inner and outer solution requires a detailed knowledge of the asymptotic solutions of the parabolic cylinder function as $\xi\rightarrow \pm\infty$.  Bender and Orzag provide details of this matching process for this specific problem and more general problems of this form~\cite{bender2013advanced}.  The matching gives the uniform solution
\begin{equation}
  u_{\mbox{\small unif}}  = \exp (-x^2/2\epsilon) \left( \frac{\exp(1-x^2)}{2\epsilon}  \right) 
  \left[  2 D_{-1/2} (x\sqrt{2/\epsilon}) + D_{-1/2} (-x\sqrt{2/\epsilon})
  \right] .
\end{equation}
This can be evaluated using various software packages that include parabolic cylinder functions~\cite{cojocaru2009parabolic}.  Figure~\ref{fig:bl_middle} shows the solution to the boundary value problem as a function of decreasing $\epsilon$.  Note the clear boundary structure that develops in the middle of the domain.

\newpage
\section*{Lecture 16:  Initial  Layers and Limit Cycles}

Boundary layers have been demonstrated to occur in boundary value problems.  The mathematical framework of the perturbation expansion can be easily ported to nonlinear problems where the distinguished limits lead to more tractable dominant balance physics models.  The theoretical ideas can also be carried over to time-dependent problems where initial layers, or even repeated regions of rapid transitions in time are known to occur.  This is a generic manifestation singularly perturbed temporal problems. 

To demonstrate the existence of rapid transition regions, we consider the  Rayleigh oscillator
\begin{equation}
   \epsilon u_{tt} + \left[ -u_t + \frac{1}{3}  (u_t)^3 \right] + u = 0
\end{equation}
with initial conditions $u(0)=A$ and $u_t(0)=B$.  The Rayleigh equation manifests strongly nonlinear limit cycle behaviors. The parameter $\epsilon\ll 1 $ is in front to the highest derivative, leading to temporal boundary layer phenomena as shown for boundary value problems.  As before, we will need to break this problem up into calculations of the outer regions, inner regions, and matching regions.

To study this equation in more detail, we make a variable transformation by letting
\begin{equation}
    v = \frac{du}{dt}
\end{equation}
so that the nonlinear oscillator can be written in the form
\begin{subeqnarray}
     \frac{du}{dt} &=& v \\
    \epsilon \frac{dv}{dt} &=& v - \frac{1}{3} v^3 - u .
\end{subeqnarray}
In parametric form, where the explicit time dependency has been removed from the equations, the coupled equations are represented as
\begin{equation}
  \epsilon \frac{dv}{du} = \frac{v - v^3/3 - u}{v} 
\end{equation}
which is derived by dividing the top equation by the bottom equation.  This then gives $v$ as a function of $u$.\\

\begin{figure}[t]
\begin{overpic}[width=0.65\textwidth]{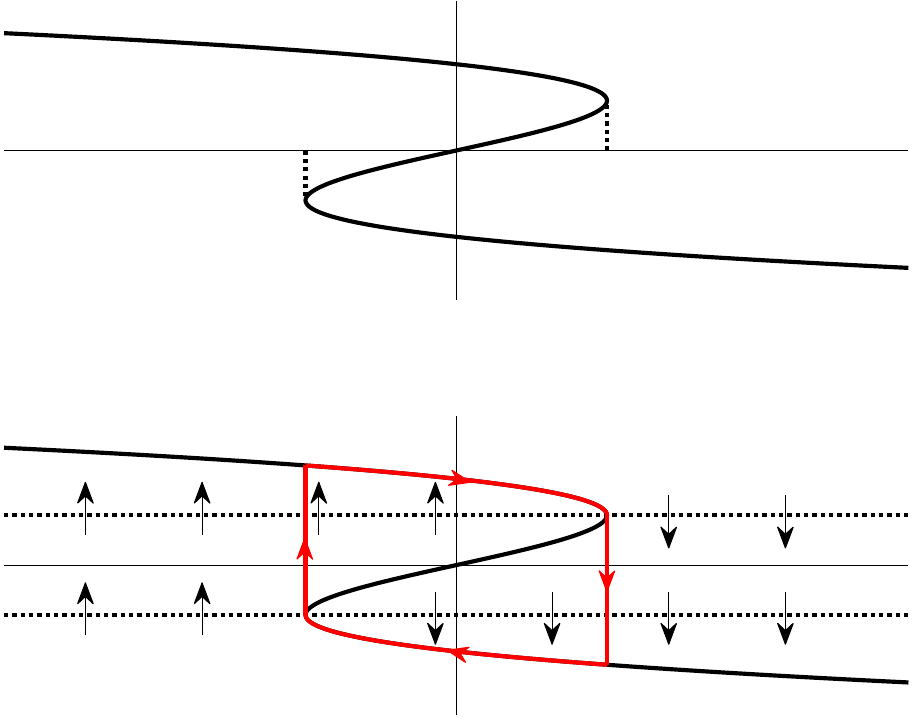}
\put(0,70){(a)}
\put(0,25){(b)}
\put(70,68){$(2/3,1)$}
\put(17,55){$(-2/3,-1)$}
\put(52,78){$v$}
\put(98,59){$u$}
\put(52,32.5){$v$}
\put(98,13.5){$u$}
\put(65,45){$u=v_0 - \frac{1}{3} {v_0}^3$}
\end{overpic}
\caption{(a) The leading order outer solution branch which satisfies the cubic equation  $u=v_0 - \frac{1}{3} {v_0}^3$.  (b) The flow dynamics of the leading order system.  At the lines $v=\pm1$, the dynamics can be shown to flow towards one of the leading order branches.  The resulting limit cycle is shown in red. Note the rapid transitions that occur when the solution jumps from one branch to another, this is the inner solution (transition layer).}
 \label{fig:il_ray}
\end{figure}
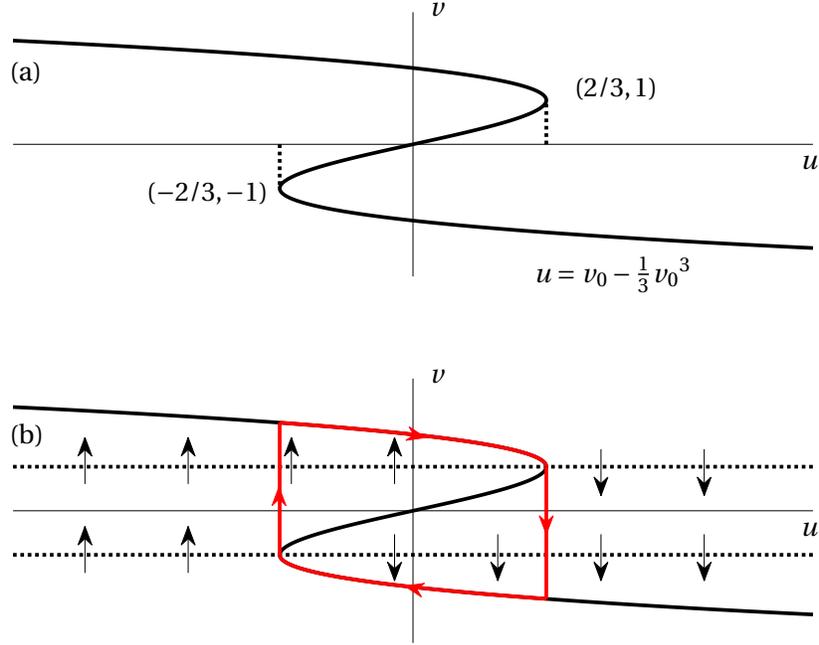

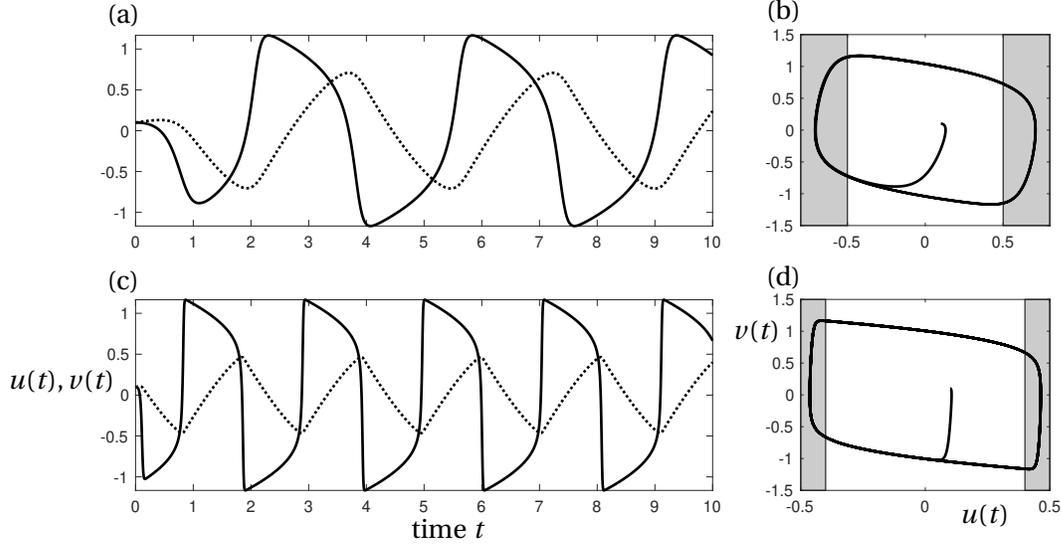
\begin{figure}[t]
\begin{overpic}[width=0.5\textwidth]{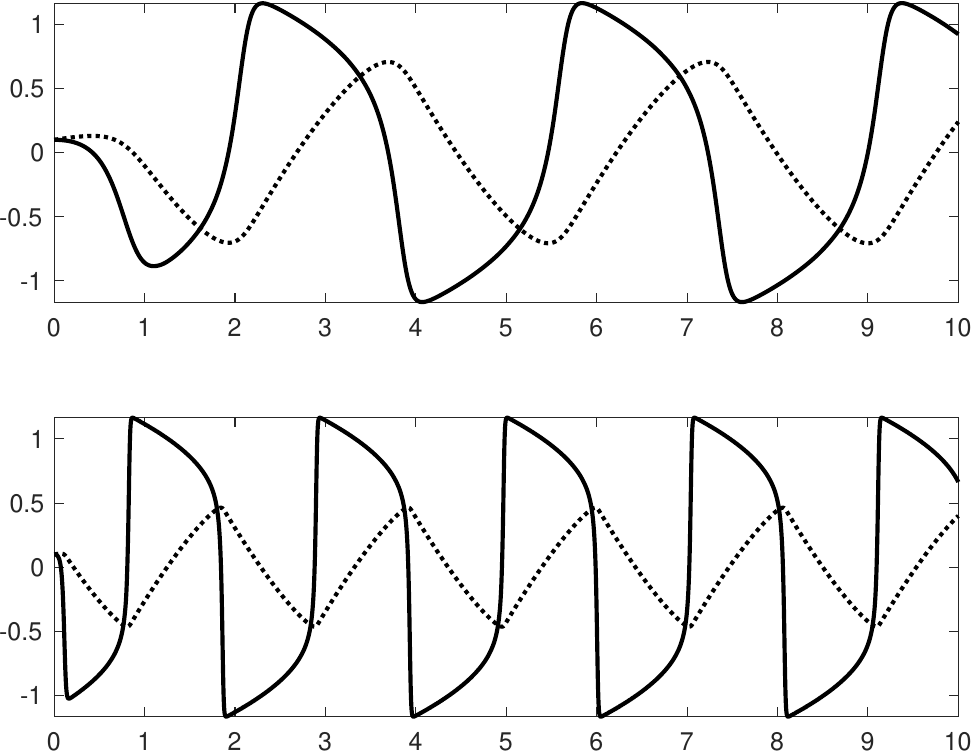}
\put(1,79){(a)}
\put(1,36){(c)}
\put(-15,20){$u(t),  v(t)$}
\put(50,-4){time $t$}
\end{overpic}\hspace*{.2in}
\begin{overpic}[width=0.237\textwidth]{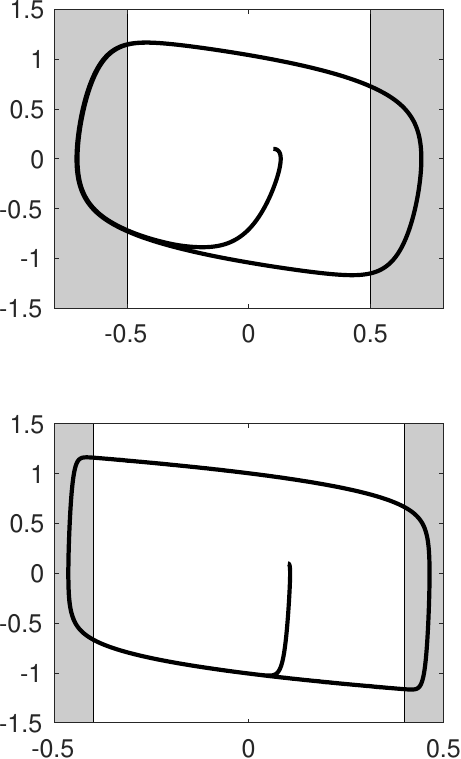}
\put(0,102){(b)}
\put(0,47){(d)}
\put(-7,36){$v(t)$}
\put(40,-2){$u(t)$}
\end{overpic}
\caption{Evolution dynamics for the Rayleigh equation for $\epsilon=0.2$ (top panels) and $0.02$ (bottom panels).  The left panels show the evolution dynamics of $u(t)$ (dotted line) and $v(t)$ (solid line).  The right panels show the phase plane dynamics along with the corresponding {\em fast time scale} dynamics region in the shaded box.}
 \label{fig:rayleighdyn}
\end{figure}

\noindent {\em I.   Outer Solution  -- }
The outer solution can be found by the regular expansion $v=v_0 + \epsilon v_1 + \cdots$.  This gives at leading order
\begin{equation}
    \frac{v_0 - {v_0}^3/3 - u}{v_0} = 0
\end{equation}
which results in
\begin{equation}
   u= v_0 - \frac{1}{3} {v_0}^3 .
\end{equation}
Figure~\ref{fig:il_ray}(a) shows this leading order solution branch.  Note that for  $|u|<2/3$, there are three solution branches possible for $v_0$.  Figure~\ref{fig:rayleighdyn} shows the full dynamics from simulations of the Rayleigh equation.

We can also evaluate the dynamics if we are not on the leading order curve.  For instance, we can look at the dynamics across the line $v=\pm 1$ to see the flow field dynamics.  In particular, we have 
\begin{subeqnarray}
  && \mbox{at} \,\, v=+1:  \,\,\,  \epsilon \frac{dv}{dt} =+\frac{2}{3} - u = \left\{
     \begin{array}{ccc} >0 & \mbox{if} & y<2/3 \\
     <0 & \mbox{if} & y>2/3 \end{array} \right. \nonumber \\
  && \mbox{at} \,\, v=-1:  \,\,\,  \epsilon \frac{dv}{dt} = -\frac{2}{3} - u = \left\{
     \begin{array}{ccc} >0 & \mbox{if} & y<-2/3 \\
     <0 & \mbox{if} & y>-2/3 \end{array} \right. \nonumber .
\end{subeqnarray}
Thus the dynamics of the flow across the lines $v=\pm 1$ are easily determined as shown in Fig.~\ref{fig:il_ray}(b).  Also depicted in this figure is the overall flow field showing how the limit cycle (red) is formed from the three branches of solutions.\\

\noindent {\em II.   Inner Solution  -- }
The transition region in Fig.~\ref{fig:il_ray}(b) suggests the following inner layer expansion
\begin{equation}
  \xi = \frac{u- 2/3}{\epsilon} 
\end{equation}
where a distinguished limit of $\delta=\epsilon$ is correct for this inner layer problem.  This  then gives the following change of variable
\begin{subeqnarray}
  &&  v_u = v_\xi \xi_u = \frac{1}{\epsilon} v_\xi \\
  &&  v_{uu} = \frac{1}{\epsilon^2} v_{\xi\xi} \\
  &&  u= \frac{2}{3} + \epsilon \xi .
\end{subeqnarray}
Plugging this into our governing equations gives
\begin{equation}
  \frac{dv}{d\xi} = \frac{ v - v^3/3 - 2/3 - \epsilon \xi}{v} .
\end{equation}
A perturbation expansion $v=v_0 (\xi) + v_1 (\xi) + \cdots$ then yields the leading order inner solution
\begin{equation}
  \frac{dv_0}{d\xi} = \frac{ v_0 - {v_0}^3/3 - 2/3}{v_0} .
\end{equation}
Separation of variables and integration yields the transcendental equation
\begin{equation}
    -\frac{2}{9} \ln |v_0+2| +  \frac{2}{9} \ln |v_0-1| - \frac{1}{3(v_0-1)} = -\frac{1}{3} (\xi+C)
\end{equation}
where $C$ is an unknown constant of integration.  This is the inner solution characterizing the rapid transition region on the right side highlighted in Fig.~\ref{fig:il_ray}(b).  One can similarly expand  $\xi = {(u+2/3)}/{\epsilon}$ to determine the inner solution on the left side.\\

\noindent {\em III.   Matching  -- }
For the inner layer on the right side of the limit cycle, matching conditions apply at both the bottom and top of where the inner solution connects to the outer solution.  
For the outer solution, the solution approaches the following at the top-right corner of the matching region
\begin{equation}
   \mbox{outer:} \,\, (u,v_0) = (2/3^-, 1) 
\end{equation}
where $2/3^-$ denotes the approach is from values less than $2/3$.  The inner solution already matches this solution since
\begin{equation}
   \mbox{inner:} \,\, \lim_{\xi\rightarrow -\infty} (\xi,v_0) = (-\infty,1)
\end{equation}
since $-1/(v_0-1)\rightarrow -\infty$ or $\ln |v_0-1| \rightarrow -\infty$ as $\xi\rightarrow -\infty$ with $z=1$.

At the bottom-right corner of the matching region
\begin{equation}
   \mbox{outer:} \,\, (u,v_0) = (2/3^-, -2) .
\end{equation}
The inner solution matches this solution since
\begin{equation}
   \mbox{inner:} \,\, \lim_{\xi\rightarrow -\infty} (\xi,v_0) = (-\infty,-2)
\end{equation}
since $\ln |v_0+2| \rightarrow -\infty$ as $\xi\rightarrow -\infty$ with $z=-2$.

As already mentioned previously, a similar inner layer expansion and matching procedure can be carried out on the left side transition region.  In this case, 
the inner expansion  $\xi = {(u+ 2/3)}/{\epsilon}$ determine the inner solution on the left side.  This procedure shows that the different dynamical regimes of the inner and outer variables are produced as a consequence of the singular nature of the perturbation problem in time.  Thus both slow and fast time scale dynamics are important for understanding the manifest dynamics.

\newpage
\section*{Lecture 17:  WKB Theory}

Two things make analytically solving differential equations difficult:  nonlinearity and non-constant coefficients.  The perturbation methods considered to this point are constructed to handle either of these cases.  A final method highlighted here is the so-called WKB method (Wentzel-Kramers-Brilloun), or alternatively the WKBJ method (add Jeffries).  It is a method particularly well suited for singularly perturbed problems.

WKB formulates a strategy by recognizing that a large number of problems in the physical and engineering sciences  are characterized by exponential behavior, either real (growth/decay), imaginary (oscillations) or both.  WKB then seeks a solution of the form
\begin{equation}
   u(x)\sim A(x) \exp \left[  \frac{S(x)}{\delta} \right] 
\end{equation}
with $\delta \ll 1$.  The goal is to determine the two quantities of interest using perturbation methods.  Specifically, we determine
\begin{subeqnarray}
 &&  S(x)  - \mbox{phase} \\
 && A(x) - \mbox{amplitude}
\end{subeqnarray}
which then suggests that the WKB is nothing more than an amplitude-phase decomposition.

A canonical problem considered by WKB theory is the following singularly perturbed differential equation
\begin{equation}
   \epsilon^2 u_{xx} + Q(x) u =0
   \label{eq:wkb_quantum}
\end{equation}
where $Q(x)\neq 0$.  This problem is motivated by both quantum mechanics (particles in Schr\"odinger potentials) and electrodynamics (waveguides).

A formal WKB expansion is a perturbation expansion that takes the form
\begin{equation}
   u(x)\sim \exp  \left[  \frac{S(x)}{\delta} \right]  =\exp \left[ 
   \frac{S_0(x) + \delta S_1(x) + \delta^2 S_2(x) + \cdots }{\delta} \right] .
\end{equation}
To use this trial solution to solve the differential equation, we compute the following quantities
\begin{subeqnarray}
  &&  \hspace*{-.5in} u_x(x)= \left( \frac{{S_0}_x(x) + \delta {S_1}_x(x) + \cdots }{\delta} \right) 
    \exp \left[  \frac{S_0(x) + \delta S_1(x) + \cdots }{\delta} \right]  \\
    &&  \hspace*{-.5in} u_{xx}(x)= \left(  \frac{({S_0}_x(x) + \delta {S_1}_x(x) + \cdots)^2 }{\delta}
    + \frac{{S_0}_{xx}(x) + \delta {S_1}_{xx}(x) + \cdots }{\delta}  \right) 
    \exp \left[  \frac{S_0(x) + \delta S_1(x) + \cdots }{\delta} \right] .
\end{subeqnarray}
Plugging this into the governing differential equation yields
\begin{equation}
  \epsilon^2 \left[  \frac{ {S_0}_x^2(x) }{\delta^2} +  \frac{ 2 {S_0}_x (x) {S_1}_x (x) }{\delta} + \cdots +  \frac{ {S_0}_{xx} }{\delta} + \cdots \right] = Q(x) .
\end{equation}
To establish a dominant balance, the dominant term on the left must balance the $Q(x)\sim O(1)$ term on the right.  This gives the condition
\begin{equation}
   \frac{\epsilon^2}{\delta^2} {S_0}_x^2 (x) \sim O(1)
\end{equation}
which requires $\delta=\epsilon$ as a distinguished limit.  This then give the asymptotic expansion
\begin{equation}
  {S_0}_x^2(x) + 2 \epsilon {S_0}_x (x) {S_1}_x (x) + \cdots + \epsilon {S_0}_{xx} + \cdots = Q(x).
\end{equation}
Collecting power of epsilon gives the hierarchy of equations
\begin{subeqnarray}
  && O(1) \,\,\,\,\,  {S_0}_x^2(x)= Q(x) \,\,\,\,\,\,\,\,\, \mbox{(eikonal equation)} \\
  && O(\epsilon) \,\,\,\,\, {S_0}_{xx} + 2 \epsilon {S_0}_x (x) {S_1}_x (x) =0 
  \,\,\,\,\,\,\,\,\, \mbox{(transport equation)} \\
  && \,\,\,\,\, \vdots \\
  && O(\epsilon^n) \,\,\,\,\, {S_{n-1}}_{xx} + 2 \epsilon {S_{n-1}}_x (x) {S_n}_x (x) 
     + \sum_{j=1}^{n-1} {S_j}_x {S_{j-1}}_x = 0  \,\,\,\,\, n\geq 2 .
\end{subeqnarray}
The leading order solution is easily solved for to give
\begin{equation}
    {S_0}(x) = \pm \int_x \sqrt{ Q(x) } dx .
\end{equation}
This solution is achieved by simple integration of the leading order equation.  At the next order, the solution $S_0(x)$ can be used to solve for $S_1(x)$:
\begin{equation}
   {S_1}_x (x) =  - \frac{{S_0}_{xx}}{2 {S_0}_x} = -\frac{1}{2} \frac{d}{dx}
   \left[  \ln \left(  {S_0}_x \right)  \right] .
\end{equation}
Integrating both sides of this equation gives
\begin{equation}
   S_1(x) = - \frac{1}{4} \ln Q(x) .
\end{equation}
Thus at this point we have both the leading order solution and its first correction specified.  This gives the solution
\begin{eqnarray}
   u(x) &=&   \exp \left[  \frac{S_0(x) + \delta S_1(x) + \cdots }{\epsilon} \right] \nonumber \\
   &=&   \exp \left[  \frac{ \pm  \int_x \sqrt{ Q(\xi) } d\xi  - (\epsilon/4) \ln Q(x) + \cdots }{\epsilon} \right]  \nonumber \\
   &=& \exp \left[ \pm \frac{1}{\epsilon} \int_x \sqrt{ Q(\xi) } d\xi \right]
   \exp \left[ - \frac{\epsilon}{4} \ln Q(x)  \right] \nonumber \\
   &=& Q(x)^{-1/4} \exp \left[ \pm \frac{1}{\epsilon} \int_x \sqrt{ Q(\xi) } d\xi \right] .
\end{eqnarray}
The two signs of solutions gives us the overall WKB approximation
\begin{equation}
  u(x) = C_1 Q(x)^{-1/4} \exp \left[ \frac{1}{\epsilon} \int_x \sqrt{ Q(\xi) } d\xi \right]
  + C_2 Q(x)^{-1/4} \exp \left[ - \frac{1}{\epsilon} \int_x \sqrt{ Q(\xi) } d\xi \right] .
\end{equation}
This is the leading order approximation of the amplitude and phase decomposition which is accurate up to $O(\epsilon)$.  If one continues to higher order in the approximation, one can find the overall hierarchy of solutions
\begin{subeqnarray}
  &&  S_0(x) = \pm \int_x \sqrt{ Q(\xi) } d\xi \\
  &&  S_1(x) = - \frac{1}{4} \ln Q(x) \\
  &&  S_2(x) = \pm \int_x \left[ \frac{Q_{xx}}{8Q^{3/2}} - \frac{5 (Q_x)^2}{32 Q^{5/2}}  \right] d\xi \\
  &&  S_3(x) = - \frac{Q_{xx}}{16Q^2} + \frac{5(Q_x)^2}{64 Q^3} \\
  &&  \,\,\,\, \vdots .
  \label{eq:wkb_s}
\end{subeqnarray}
Thus each even order produces a closed form expression while the odd orders require integrals to be evaluated with two branches ($\pm$) of solutions.

\subsection*{Boundary Layer Theory Revisited}  We can use the WKB method to reconsider our boundary layer theory calculations.  In so doing, we again consider the canonical model
\begin{equation}
  \epsilon u_{xx} + b(x) u_x + c(x) u = 0
\end{equation}
on the interval $x\in[0,1]$ and with $u(0)=A$ and $u(1)=B$.    The WKB expansion assumes the solution takes the form
\begin{equation}
   u(x)  =\exp \left[ 
   \frac{S_0(x) + \epsilon S_1(x) + \epsilon^2 S_2(x) + \cdots }{\epsilon} \right] .
\end{equation}
Inserting this ansatz into our governing boundary value problem gives 
\begin{equation}
   \epsilon^2 \left[  \frac{ {S_0}_x^2(x) }{\epsilon^2} +  \frac{ 2 {S_0}_x (x) {S_1}_x (x) }{\epsilon} + \cdots +  \frac{ {S_0}_{xx} }{\epsilon} + \cdots \right]
   + b(x) \left[  \frac{ {S_0}_x (x) }{\epsilon} +  {S_1}_x (x) + \cdots \right] + c(x) =0 .
\end{equation}
Collecting powers of epsilon gives the hierarchy of equations
\begin{subeqnarray}
  && O(1/\epsilon) \,\,\,\,\,\,\,  {S_0}_x^2(x) + b(x) {S_0}_x (x) \\
  && O(1) \,\,\,\,\,\,\,\,\, {S_0}_{xx} + 2 {S_0}_x (x) {S_1}_x (x) + b(x) {S_1}_x (x) + c(x)=0   .
\end{subeqnarray}
The leading order problem can be rewritten as
\begin{equation}
    {S_0}_x \left[ {S_0}_x + b(x) \right]
\end{equation}
which gives the two solutions ${S_0}_x=0$ or ${S_0}_x =- b(x)$.   

For the case ${S_0}_x=0$ so that $S_0$ is a constant.  At next order this gives
${S_1}_x b(x) + c(x) =0$ so that
\begin{equation}
  S_1 (x) = - \int_{0}^x \frac{c(\xi)}{b(\xi)} d\xi .
\end{equation}
The WKB solution in this case is given by
\begin{equation}
   u(x)  =\exp \left[ \frac{S_0(x)}{\epsilon} +  S_1(x) \right]   =C_1
    \exp \left[  - \int_{0}^x \frac{c(\xi)}{b(\xi)} d\xi  \right]
\end{equation}
which is exactly the outer solution we found from boundary layer theory.

For the case where ${S_0}_x =- b(x)$, then $S_0$ is determined by integrating over $b(x)$.  At the next order, this gives
\begin{equation}
   b(x) {S_1}_x (x) + b_x (x) = c(x) .
\end{equation}
Rearranging gives the equation for $S_1(x)$
\begin{equation}
 {S_1}_x (x)  = \frac{c(x)}{b(x)} - \frac{d}{dx} \left[   \ln b(x) \right]
\end{equation}
which can be integrated to give
\begin{equation}
   S_1(x) = \int_0^x \frac{c(\xi)}{b(\xi)} d\xi  - \ln b(x) .
\end{equation}
Plugging this into the WKB expansion gives
\begin{equation}
   u(x)  =\exp \left[ \frac{S_0(x)}{\epsilon} +  S_1(x) \right]   =C_2
   \frac{1}{b(x)} \exp \left[  \int_{0}^x \frac{c(\xi)}{b(\xi)} d\xi  -\frac{1}{\epsilon}
   \int_{0}^x b(\xi) d\xi \right]
\end{equation}
which was the inner solution we found previously.  This gives in total the solution
\begin{equation}
  u(x) =   C_1
    \exp \left[  - \int_{0}^x \frac{c(\xi)}{b(\xi)} d\xi  \right]  +C_2
   \frac{1}{b(x)} \exp \left[  \int_{0}^x \frac{c(\xi)}{b(\xi)} d\xi  -\frac{1}{\epsilon}
   \int_{0}^x b(\xi) d\xi \right]
\end{equation}
thus reproducing the boundary layer theory results and connecting the WKB method to multiscale perturbation techniques for singularly perturbed problems.

\subsection*{Validity of WKB method}
It is important to assess under what conditions the formal WKB expansion holds.
Thus the validity of the expansion
\begin{equation}
  u(x) \sim \exp \left[  \frac{1}{\delta} \sum_{n=0}^\infty \delta^n S_n(x) \right]
\end{equation}
when $\delta\rightarrow 0$ is evaluated.  Specifically, for the WKB expansion to be valid, it is necessary that the sum be an asymptotic series as $\delta\rightarrow 0$ so that
\begin{subeqnarray}
  && \delta S_1 \ll S_0 \\
  && \delta S_2 \ll S_1 \\
  && \,\,\,\, \vdots \\
  && \delta S_{n+1} \ll S_n
\end{subeqnarray}
for $\delta\rightarrow 0$.  If the series is truncated at  the term $\delta^{n-1} S_n(x)$ the the next term must be small compared to one.  Thus
\begin{equation}
   \delta^n S_{n+1}(x) \ll 1
\end{equation}
for $\delta\rightarrow 0$.  Thus the leading order solution $u\sim \exp(S_0/\delta)$ is never a good approximation.  In fact, we require at least 
$u\sim \exp(S_0/\delta + S_1)$ to produce a reasonable approximation.  This is known as the {\em physical optics} approximation.  This language is borrowed from the Helmholtz equation $u_{xx} + k^2 u=0$, which is a reduction of Maxwell's equations of electromagnetism when assuming plane wave solutions.  Specifically, for optical frequencies, the Helmholtz equation is solved for $k^2\gg 1$, which makes it ideal for WKB theory.

\newpage
\section*{Lecture 18:  WKB Theory and Turning Points}

In what follows, we consider a number of important examples that typically arise from either quantum mechanics or electrodynamics.  The first to consider is the Airy equation
\begin{equation}
  u_{xx} = x u
\end{equation}
which is known as a turning point problem since the exact Airy function solution exponentially decays on one side of $x=0$ and oscillates on the other side of $x=0$.  Thus it models the behavior near the edge of a potential well.  

For this case, we can connect it back to the results for (\ref{eq:wkb_quantum}) with $\epsilon=1$ and $Q(x)=x$.   Evaluating the expressions in (\ref{eq:wkb_s}) gives
\begin{subeqnarray}
 && S_0(x) = \pm \frac{2}{3} x^{3/2} \\
 && S_1(x) = - \frac{1}{4} \ln x \\
 && S_2 (x) =\pm \frac{5}{48} x^{-3/2} .
\end{subeqnarray}
The WKB theory can be used to evaluate the behavior of the solution as $x\rightarrow \infty$ since  $\epsilon S_2 \ll S_1 \ll S_0/\epsilon$ and
$\epsilon S_2 \\ 1$ for $\epsilon=1$ and $x\rightarrow \infty$.  This then
gives the behavior of the solution
\begin{equation}
    u(x)\sim C_{\pm} x^{-1/4} \exp (\pm 2x^{3/2}/3 )
\end{equation}
as $x\rightarrow \infty$.  This identifies the decay rate of the solution and accurately represents its asymptotic behavior.

A second example to consider is that of the parabolic cylinder functions which was already considered in the example of internal boundary layers.  The generic form of the equation is given by
\begin{equation}
  u_{xx} = \left(  \frac{1}{4} x^2 - \nu -\frac{1}{2} \right) u
\end{equation}
so that $\epsilon=1$ and $Q(x)=x^2/4 - \nu - 1/2$.  Solutions for both $S_0$ and $S_1$ can be found in the $x\rightarrow \infty$ limit to yield
\begin{subeqnarray}
  &&  S_0 = \pm \left[ \frac{x^2}{4} + \ln x^{-(\nu+1/2)}  \right] \\
  && S_1 = \ln \left( \frac{1}{4} x^{-1/2} \right) 
\end{subeqnarray}
which gives the solution
\begin{equation}
    u(x)= C_{+} x^{-\nu-1} \exp (x^{2}/4 ) + C_{-} x^{\nu} \exp (-x^{2}/4 ) 
\end{equation}
in the limit as $x\rightarrow \infty$.  The behavior of special functions at extremes like $x\rightarrow \pm \infty$ is critical for understanding many classical problems in mathematical physics.

We consider yet another example where the physical optics approximation is insufficient to give a good approximation.  Consider the governing equations
\begin{equation}
  u_{xx} = \left(  \frac{\ln x}{x} \right)^2 u
\end{equation}
so that $\epsilon=1$ and $Q(x)=(\ln x/x)^2$.  This allows us to compute the following
\begin{subeqnarray}
  &&  S_0 = \pm \frac{1}{2} (\ln x)^2 \\
  && S_1 = \frac{1}{2} \ln x - \frac{1}{2} \ln (\ln x)  \\
  && S_2 = \pm \frac{1}{8} \ln (\ln x) \pm \frac{3}{16} (\ln x)^{-2}
\end{subeqnarray}
where $S_2 \ll S_1 \ll S_0$ for $x\rightarrow \infty$. But if we truncate at $S_2$, then we find that $S_2\gg 1$ for $x\rightarrow \infty$.  Thus the asymptotic series would not provide a good approximation.  At the next order, one can compute
\begin{equation}
  S_3 = \frac{3}{16} (\ln x)^{-4} - \frac{1}{16} (\ln x)^{-2} 
\end{equation}
where $S_3 \ll S_2$ and $S_3 \ll 1$ for $x\rightarrow \infty$.  This then determines where the asymptotic expansion should be truncated.

As a final example, we consider an eigenvalue problem from quantum mechanics.  
\begin{equation}
  u_{xx} + E Q(x) u = 0
\end{equation}
with $Q(x)>0$ and $u(0)=u(\pi)=0$.  The goal is to approximate the eigenfunctions and eigenvalues for large eigenvalues $E\rightarrow \infty$, i.e. $\epsilon = 1/E$ so that this again takes the form $\epsilon u_{xx}+Q(x) u=0$.   Plugging in the specific form given to our calculations gives
\begin{subeqnarray}
  &&  S_0 = \pm \int_x \sqrt{-EQ(\xi)} d\xi = \pm i\sqrt{E} \int_x \sqrt{Q(\xi)} d\xi \\
  &&  S_1 = - \frac{1}{4} \ln Q(x) 
\end{subeqnarray}
so that
\begin{equation}
  u= C_1 Q^{-1//4} \sin \left[ \sqrt{E} \int_x \sqrt{Q(\xi} d\xi \right] 
  + C_2 Q^{-1//4} \cos \left[ \sqrt{E} \int_x \sqrt{Q(\xi} d\xi \right] .
\end{equation}
Applying the boundary conditions then gives the eigenvalues and eigenfunctions
\begin{eqnarray}
   &&  E_n = \left[  \frac{ n\pi}{\int_0^\pi \sqrt{Q(\xi)}d\xi}  \right]^2  \\
   &&  u(x) = A Q^{-1//4} \sin \left[ \sqrt{E} \int_x \sqrt{Q(\xi} d\xi \right] 
\end{eqnarray}
for $n\rightarrow \infty$.  If the specific example of $Q(x)=(x+\pi)^4$ is used, then we find the approximation solutions for the high-energy quantum states
\begin{eqnarray}
   &&  E_n = \frac{9n^2}{49\pi^4}  \\
   &&  u_n(x) = \sqrt{\frac{6}{7\pi^3}} \frac{ \sin \left[ n (x^3+3x^2\pi + 3\pi^2 x) \right]}{\pi+x}  
\end{eqnarray}
with $n\rightarrow \infty$.  One can see that such analytic expressions were important before the advent of computational tools to evaluate such limits.  And regardless of our computing capabilities, the formula reveals important trends in the behavior which is shown analytically.

\newpage
\section*{Lecture 19:  Bifurcation Theory}

One of the interesting things that can happen in dynamical systems is that solutions and behaviors can change as a function of parameter in the governing equations.  Thus we can consider a dynamical system of the form
\begin{equation}
  \frac{d{\bf x}}{dt} = f({\bf x},\mu)
\end{equation}
where ${\bf x}$ is the collection of state-space  variables,  $f(\cdot)$ specifies the governing equations and $\mu$ is a parameter that can be modified to fundamentally alter the dynamics of the system, i.e. it is a bifurcation parameter.  A bifurcation refers the division of something into parts or branches.  In the context of dynamical systems, it refers to the branching of solutions of the governing equations.  Which solution is relevant for the dynamics depends on the stability of the solution itself.   Specifically, Bifurcation theory is the mathematical study of changes in the qualitative or topological structure of a given family, such as the integral curves of a family of vector fields, and the solutions of a family of differential equations.

\subsection*{Saddle node bifurcation}
\begin{figure}[t]
\begin{overpic}[width=0.65\textwidth]{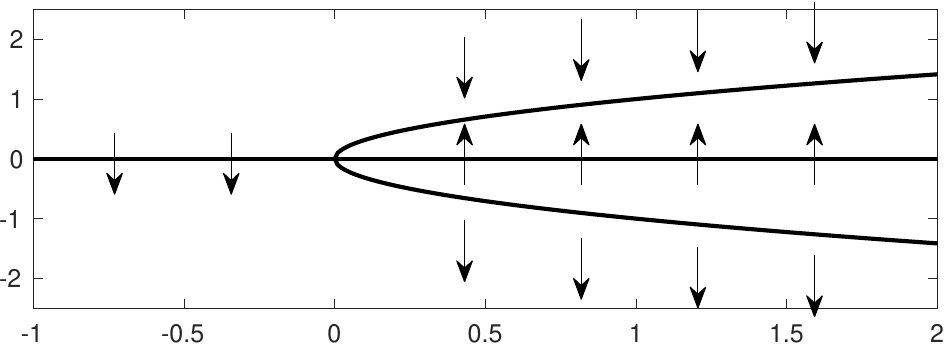}
\put(102,28){$x_0=\sqrt{\mu}$}
\put(102,10){$x_0=-\sqrt{\mu}$}
\put(102,19){$x_0=0$}
\put(86,-3){$\mu$}
\put(-3,25){$x$}
\end{overpic}
\caption{Bifurcation structure of the governing equations $dx/dt=\mu-x^2$.  There are two equilibrium branches of solutions for $\mu>0$.  The top branch $x_0=\sqrt{\mu}$ is stable while the bottom branch $x_0=-\sqrt{\mu}$ is unstable.  There are no equilibria for $\mu<0$.}
 \label{fig:saddle_node}
\end{figure}
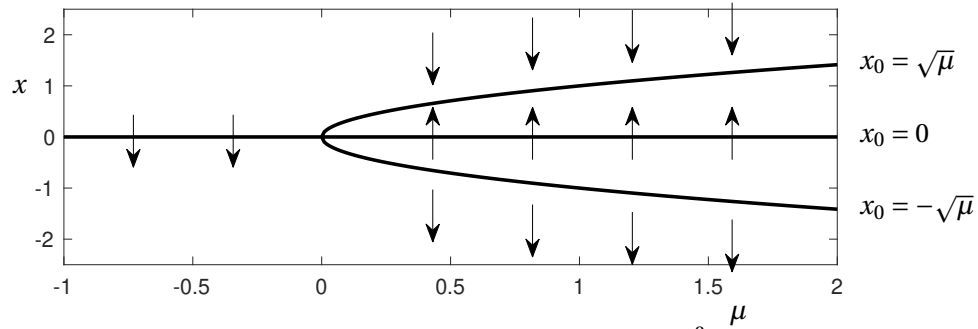

To begin thinking about bifurcation theory, consider the following canonical equation
\begin{equation}
    \frac{dx}{dt}= \mu- x^2
\end{equation}
which is a simple differential equation parametrized by $\mu$.  A closed form solution can easily be constructed for this equation by separation of variables.  However, we instead consider a qualitative approach to the problem.  In particular, we first consider equilibrium solutions where $dx/dt=0$ so that
\begin{equation}
   \mbox{equilibrium:} \,\,\,\,\,  \mu-x^2=0
\end{equation}
which gives
\begin{equation}
    x_0 = \pm \sqrt{\mu} .
\end{equation}
Figure~\ref{fig:saddle_node} shows these equilibrium solutions along with the dynamical trajectories in the $x-\mu$ plane.  Note that the top branch of solutions $x_0=\sqrt{\mu}$ is stable while $x_0=-\sqrt{\mu}$ is unstable.

The stability of the solutions can be explicitly determined by linearization.  Thus we can consider the perturbative solution
\begin{equation}
  x= x_0 + \tilde{x}
\end{equation}
where $\tilde{x}\ll 1$.  Plugging this into the governing equation and keeping terms to $O(\tilde{u})$ gives
\begin{equation}
    \frac{d\tilde{x}}{dt}= -2 x_0 \tilde{x} 
\end{equation}
whose solution is given by
\begin{equation}
   \tilde{x}(t) = \tilde{x}(0) \exp(-2x_0 t) .
\end{equation}
For $x_0=\sqrt{\mu}$, then this solution is $ \tilde{x}(t) = \tilde{x}(0) \exp(-2\sqrt{\mu} t) \rightarrow 0$ as $t\rightarrow 0$.  Thus perturbations decay.  On the contrary, for $x_0=-\sqrt{\mu}$, then this solution is $\tilde{x}(t) = \tilde{x}(0) \exp(2\sqrt{\mu} t) \rightarrow \infty$ as $t\rightarrow 0$.  The arrows in Fig.~\ref{fig:saddle_node} shows these equilibrium solutions along with the dynamical trajectories and stability in the $x-\mu$ plane.
This stability structure is known as a {\em saddle node} or {\em limit cycle} bifurcation.  It is characterized by a turning point in the solution where the stability changes.

\subsection*{Transcritical bifurcation}
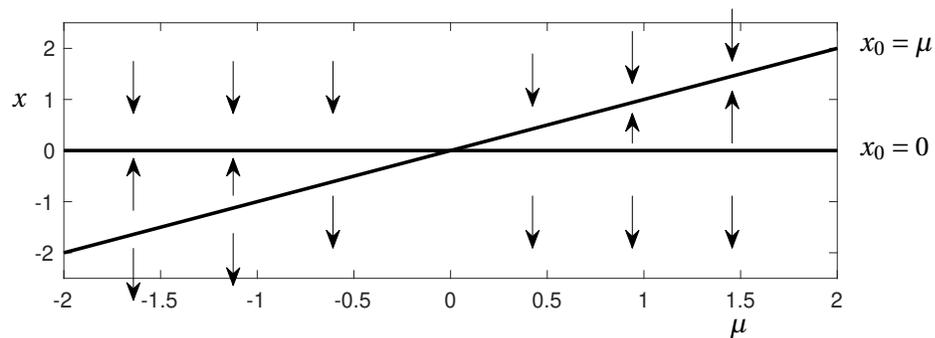
\begin{figure}[t]
\begin{overpic}[width=0.65\textwidth]{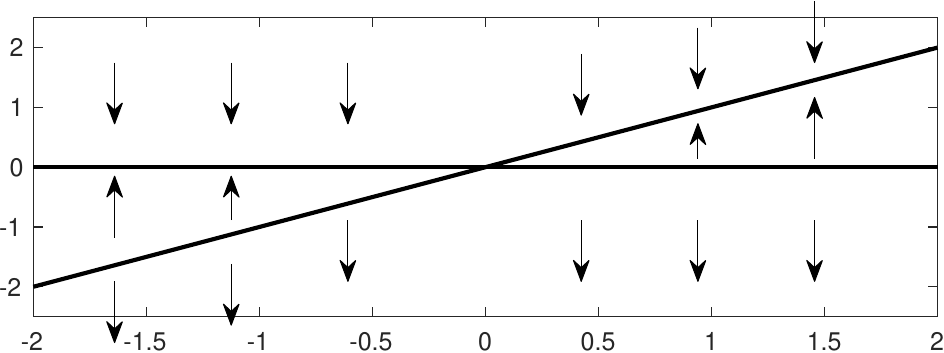}
\put(102,32){$x_0={\mu}$}
\put(102,19){$x_0=0$}
\put(86,-3){$\mu$}
\put(-3,25){$x$}
\end{overpic}
\caption{Bifurcation structure of the governing equations $dx/dt=\mu x-x^2$.  There are two equilibrium branches of solutions:  $x_0=0, \mu$.  At $\mu=0$, the two solutions trade stability.}
 \label{fig:trans_bif}
\end{figure}

Consider the second canonical equation
\begin{equation}
    \frac{dx}{dt}= \mu x- x^2
\end{equation}
which is a simple differential equation again parametrized by $\mu$. 
Equilibrium solutions where $dx/dt=0$ are found to be
\begin{equation}
   \mbox{equilibrium:} \,\,\,\,\,  \mu x-x^2 =x (\mu-x)=0
\end{equation}
which gives
\begin{equation}
    x_0 = 0, \mu .
\end{equation}
Figure~\ref{fig:trans_bif} shows these equilibrium solutions along with the dynamical trajectories in the $x-\mu$ plane.  Note the switch in stability of both solution branches at $\mu=0$.

The stability of the solutions can be explicitly determined by linearization.  Thus we can consider the perturbative solution
\begin{equation}
  x= x_0 + \tilde{x}
\end{equation}
where $\tilde{x}\ll 1$.  Plugging this into the governing equation and keeping terms to $O(\tilde{u})$ gives
\begin{equation}
    \frac{d\tilde{x}}{dt}= \mu \tilde{x} - 2 x_0 \tilde{x} 
\end{equation}
whose solution is given by
\begin{equation}
   \tilde{x}(t) = \tilde{x}(0) \exp [ (\mu-2x_0) t] .
\end{equation}
For $x_0=0$, then this solution is $ \tilde{x}(t) = \tilde{x}(0) \exp(\mu t)$ which is stable for $\mu<0$ and unstable for $\mu>0$.  For $x_0=\mu$, then the solution is $ \tilde{x}(t) = \tilde{x}(0) \exp(-\mu t)$ which is stable for $\mu>0$ and unstable for $\mu<0$.  This stability structure is known as a {\em transcritical} bifurcation.  It is characterized by two linear branches of solutions intersecting and exchanging stability.

\subsection*{Pitchfork bifurcation} 
\begin{figure}[t]
\begin{overpic}[width=0.65\textwidth]{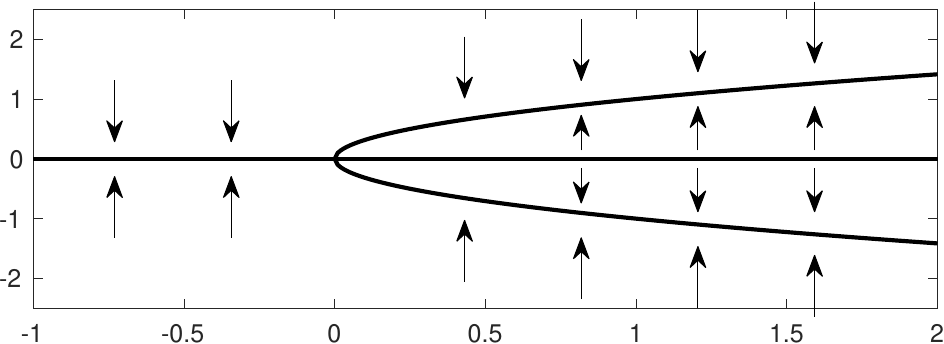}
\put(102,28){$x_0=\sqrt{\mu}$}
\put(102,10){$x_0=-\sqrt{\mu}$}
\put(102,19){$x_0=0$}
\end{overpic}
\caption{Bifurcation structure of the governing equations $dx/dt=\mu x-x^3$.  There are up to three equilibrium branches of solutions.  The stability of the solutions changes at $\mu=0$.}
 \label{fig:pitch_bif}
\end{figure}

Another canonical equation to consider is given by
\begin{equation}
    \frac{dx}{dt}= \mu x- x^3
\end{equation}
which is again parametrized by $\mu$.  We consider equilibrium solutions where $dx/dt=0$ so that
\begin{equation}
   \mbox{equilibrium:} \,\,\,\,\,  \mu x-x^3 = x (\mu-x^2)=0
\end{equation}
which gives
\begin{equation}
    x_0 = 0, \pm \sqrt{\mu} .
\end{equation}
Figure~\ref{fig:pitch_bif} shows these equilibrium solutions along with the dynamical trajectories in the $x-\mu$ plane.  Note that the top and bottom branch of solutions $x_0=\pm\sqrt{\mu}$ are stable while $x_0=0$ becomes unstable at $\mu=0$.

The stability of the solutions can be explicitly determined by linearization.  Thus we can consider the perturbative solution
\begin{equation}
  x= x_0 + \tilde{x}
\end{equation}
where $\tilde{x}\ll 1$.  Plugging this into the governing equation and keeping terms to $O(\tilde{u})$ gives
\begin{equation}
    \frac{d\tilde{x}}{dt}= (\mu -3 x_0^2) \tilde{x} 
\end{equation}
whose solution is given by
\begin{equation}
   \tilde{x}(t) = \tilde{x}(0) \exp [(\mu -3 x_0^2) t ] .
\end{equation}
For $x_0=\pm\sqrt{\mu}$, then this solution is $ \tilde{x}(t) = \tilde{x}(0) \exp(-2\mu t) \rightarrow 0$ as $t\rightarrow 0$.  Thus perturbations decay.  On the contrary, for $x_0=0$, then this solution is $\tilde{x}(t) = \tilde{x}(0) \exp(\mu t) \rightarrow \infty$ as $t\rightarrow 0$ for $\mu>0$.  The arrows in Fig.~\ref{fig:pitch_bif} shows these equilibrium solutions along with the dynamical trajectories and stability in the $x-\mu$ plane.
This stability structure is known as a {\em pitchfork} bifurcation.

\begin{figure}[t]
\begin{overpic}[width=0.65\textwidth]{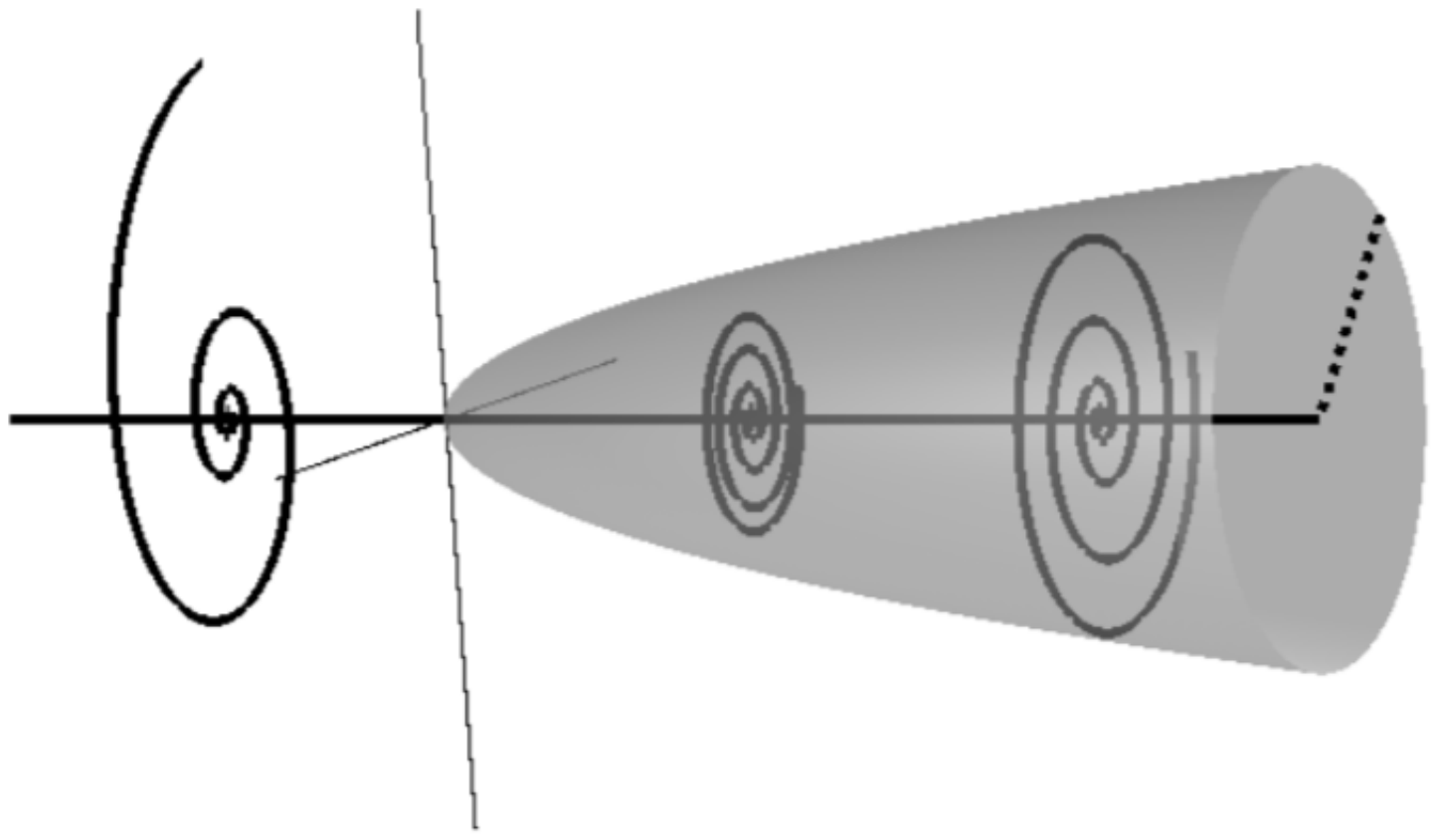}
\put(20,31){$x$}
\put(28,63){$y$}
\put(85,34){$\mu$}
\put(91,42){$r=\mu$}
\end{overpic}
\vspace*{-.5in}
\caption{Bifurcation structure of the Hopf bifurcation.  For this case, two complex conjugate pairs of imaginary eigenvalues cross from the left-half plane to the right-half plane, resulting in the onset of an oscillatory dynamics.  The shaded structure represents the surface on which the oscillatory trajectories reside after bifurcation for $\mu>0$.  In fact, all solutions for $\mu>0$ are attracted to the shown limit cycle.  For $\mu<0$, the only stable solution is the trivial solution. }
 \label{fig:hopf_bif}
\end{figure}

\begin{figure}[t]
\begin{overpic}[width=0.8\textwidth]{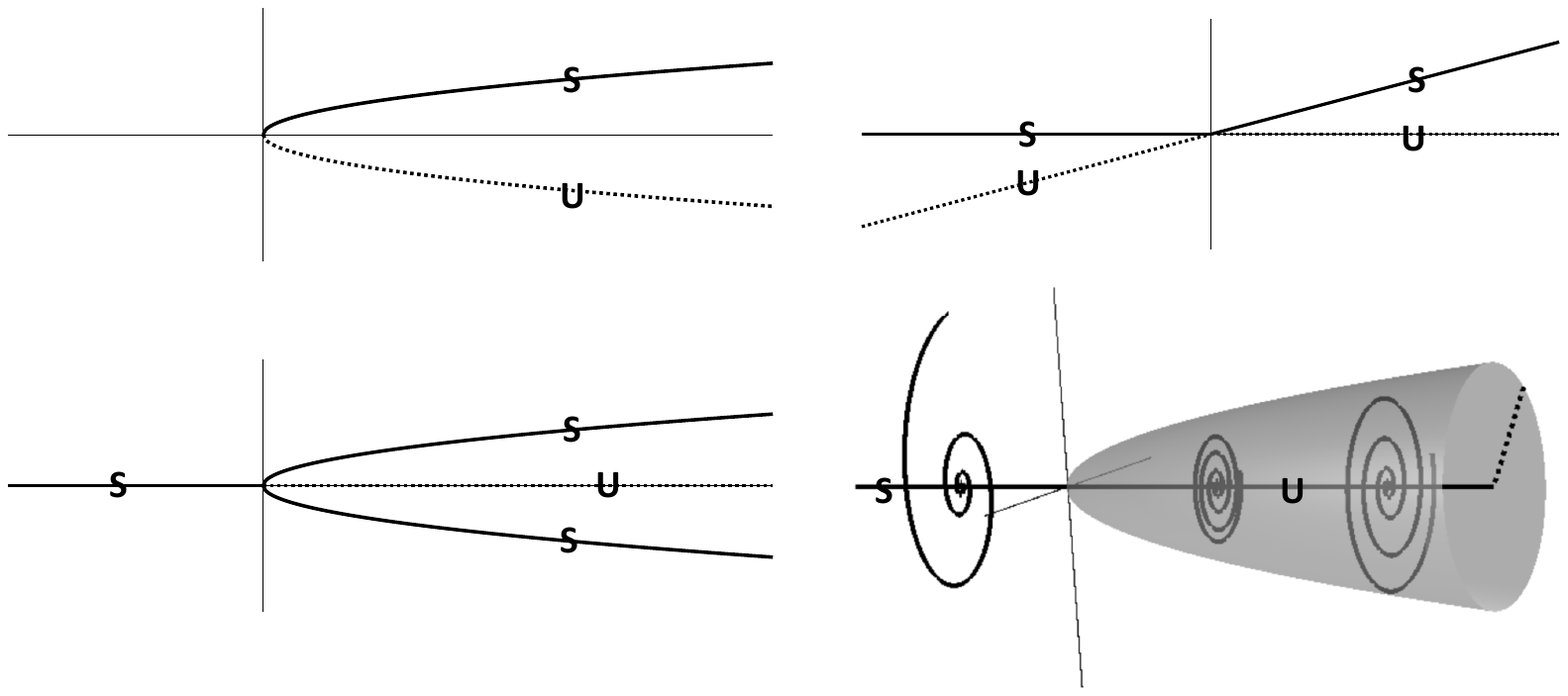}
\put(2,58){(a)}
\put(2,36){(c)}
\put(54,58){(b)}
\put(54,36){(d)}
\put(16,36){$x$}
\put(48,27.5){$\mu$}
\end{overpic}
\vspace*{-.7in}
\caption{The normal form bifurcations typically exhibited at a bifurcation point $\mu=\mu_c$.  These are the canonical forms for switching from one branch of solutions to another. The ${\bf S}$ denotes the stable branch (solid line) of solutions while the ${\bf U}$ denotes the unstable branch (dotted line).  The bifurcation types are (a) saddle node bifurcation, (b) transcritical bifurcation, (c) pitchfork bifurcation, and (d) Hopf bifurcation.}
 \label{fig:normal_bif}
\end{figure}

\subsection*{Hopf bifurcation} 
The final bifurcation to consider is the Hopf bifurcation which is characterized by a pair of complex eigenvalues that cross into the right-half plane from the left-half plane.  As a result, oscillatory behavior is observed at instability.  The Hopf does not have a canonical form.  However, we an consider an exemplar system
\begin{subeqnarray}
 &&  \frac{dx}{dt} = -y + (\mu - x^2 -y^2 ) x \\
 &&  \frac{dy}{dt} = x + (\mu - x^2 -y^2 ) y
\end{subeqnarray}
where $\mu$ once again plays the role of the bifurcation parameter.  Steady-state solutions are found now for both $dx/dt=0$ and $dy/dt=0$ simultaneously.  This gives the algebraic equations
\begin{subeqnarray}
  &&  - y+ (\mu - x^2 -y^2 ) x =0 \\
  && x + (\mu - x^2 -y^2 ) y = 0 .
\end{subeqnarray}
In this case, the origin $(x_0,y_0)=(0,0)$ is a fixed point.  In fact, it is the only one.  Linearization gives
\begin{subeqnarray}
 &&   x= 0 + \tilde{x} \\
  && y= 0 + \tilde{y}
\end{subeqnarray}
so that
\begin{subeqnarray}
  &&  \frac{d\tilde{x}}{dt} = -\tilde{y} + \mu \tilde{x} \\
 &&  \frac{d\tilde{y}}{dt} = \tilde{x} + \mu \tilde{y} .
\end{subeqnarray}
In matrix form, this is given by
\begin{equation}
   \frac{d{\bf x}}{dt} = \left[ \begin{array}{cc} \mu & -1 \\ 1 & \mu \end{array}\right] {\bf x} .
\end{equation}
This an eigenvalue problem which can be solved with the solution ansatz ${\bf x}={\bf v}\exp(\lambda t)$ which gives the eigenvalues
\begin{equation}
   \lambda_\pm = \mu \pm i .
\end{equation}
For $\mu>0$, the eigenvalues characterize growing oscillations.  For $\mu<0$, the behavior is damped oscillations.  This can be more readily seen by the transformation $x=r\cos \theta$ and $y=r\sin \theta$ which gives the alternative coupled equations
\begin{subeqnarray}
 &&  \frac{dr}{dt}= r (\mu-r^2) \\
 && \frac{d\theta}{dt} = 1
\end{subeqnarray}
whose solutions are given by
\begin{subeqnarray}
  &&  \theta=t+\theta_0 \\
  &&  r^2 = \left\{  \begin{array}{lc}  {\mu r_0^2}/ \left[ {r_0^2 + (\mu-r_0^2) \exp(-2\mu t)} \right] &  \mu\neq 0 \\
            {r_0^2}/
      \left[ {1+ 2r_0^2 t} \right] &  \mu=0 \end{array} \right.
\end{subeqnarray}
So although there are no other fixed points aside from the origin, the solution evolves to a limit cycle for $\mu>0$.  There is not limit cycle for $\mu<0$.  In fact, for $\mu<0$, all solutions collapse to the trivial solution.  Figure~\ref{fig:hopf_bif} shows the behavior of the limit cycle dynamics and the resulting Hopf bifurcation.

\clearpage
\newpage
\section*{Lecture 20:  Normal Forms and Imperfections}
We can summarize the preceding analysis:  there are four canonical and common types of instabilities that are fundamentally represented by (a) saddle node bifurcation, (b) transcritical bifurcation, (c) pitchfork bifurcation, and (d) Hopf bifurcation.  These are often called the {\em normal forms} of the bifurcation and are given by
\begin{subeqnarray}
&&  \frac{dx}{dt} = \mu - x^2  \hspace*{.5in} \mbox{saddle-node}\\
&&  \frac{dx}{dt} = \mu x- x^2 \hspace*{.5in} \mbox{transcritical}\\
&&  \frac{dx}{dt} = \mu x- x^3 \hspace*{.5in} \mbox{pitchfork}\\
&&  \frac{dr}{dt} = \mu r- r^3, \,\, \frac{d\theta}{dt}=1 \hspace*{.5in} \mbox{Hopf} 
\end{subeqnarray}
where $\mu$ is the bifurcation parameter.  Specifically, these have been written so that a change of stability occurs for the critical value $\mu_c=0$.  It should be noted that the (canonical) form of the Hopf written is not unique as there are other ways to express the limit cycle dynamics.    Figure~\ref{fig:normal_bif} highlights these canonical representations of instability.  They are at the core of pattern forming instabilities covered later in the book.  They are also often referred to as the {\em universal unfolding of the bifurcation}.

\subsection*{Universal bifurcation unfolding example}
To see the role of normal forms in practice, consider the following example
\begin{equation}
  \frac{dx}{dt} = - x (x^2 - 2 x - \mu) 
\end{equation}
which is a fairly simple dynamical system with up to a cubic nonlinearity.  The form of the equation is not in any of the normal forms considered.

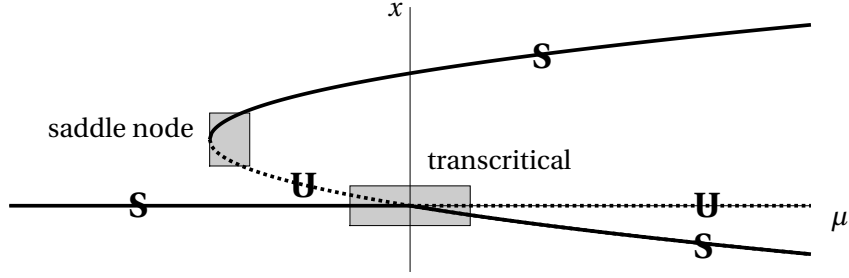
\begin{figure}[t]
\begin{overpic}[width=0.65\textwidth]{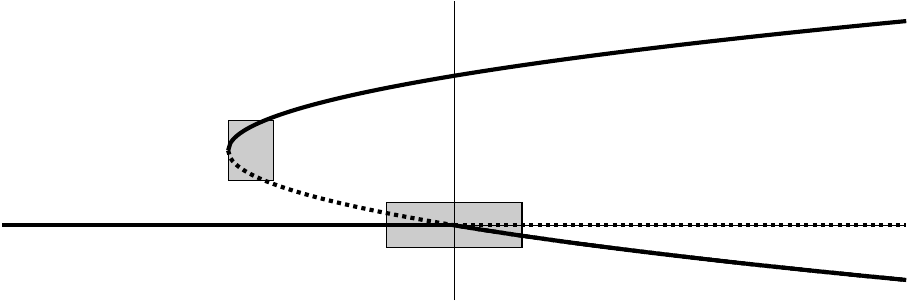}
\put(47,32){$x$}
\put(102,6){$\mu$}
\put(5,17){saddle node}
\put(52,13){transcritical}
\put(85,7){\Large {\bf U}}
\put(35,9.5){\Large {\bf U}}
\put(85,2){\Large {\bf S}}
\put(65,25.5){\Large {\bf S}}
\put(15,7){\Large {\bf S}}
\end{overpic}
\caption{The normal form bifurcation analysis of an example system.  The stable branches are denoted by ${\bf S}$ (solid line) and the unstable branches by ${\bf U}$ (dotted line).  Taylor expanding near the bifurcation points shows that the normal form around $(x,\mu)=(0,0)$ is a transcritical bifurcation and around $(x,\mu)=(1,-1)$ is a saddle node.}
 \label{fig:normal_form_example}
\end{figure}

Equilibrium solutions are given by solutions satisfying
\begin{equation}
  x_0 (x_0^2 - 2x_0 - \mu) =0 \,\,\,\, \rightarrow \,\,\,\,  x_0=0, 1\pm \sqrt{1+\mu} .
\end{equation}
Thus there are up to three branches of solutions total for $x_0(\mu)$. Stability of each branch is determined by a linear stability analysis where $x=x_0+\tilde{x}$ with $\tilde{x}\ll 1$.   This gives the linearized evolution
\begin{equation}
  \frac{d\tilde{x}}{dt} =  ( 4 x_0 - 3x_0^2 + \mu ) \tilde{x}.
\end{equation}
For the equilibrium solution $x_0=0$, the above yields $d\tilde{x}/dt=\mu\tilde{x}$ which has exponentially growing solutions (unstable) for $\mu>0$ and exponentially decaying solutions (stable) for $\mu<0$.  
For the equilibrium $x_0=1+\sqrt{1+\mu}$, the linearized evolution is given by $d\tilde{x}/dt=-2(1+\mu+\sqrt{1+\mu})\tilde{x}$ which produces exponentially decaying solutions (stable) for all value of $\mu>-1$ associated with this branch of solutions.  Finally, for the equilibrium $x_0=1-\sqrt{1+\mu}$, the linearized evolution is given by $d\tilde{x}/dt=-2(1+\mu-\sqrt{1+\mu})\tilde{x}$ which produces exponentially growing (unstable) solutions for $-1<\mu<0$ and exponentially decaying (stable) solutions for $\mu>0$.  Figure~\ref{fig:normal_form_example} shows graphically the full stability structure of the dynamics.

The transitions in the stability of the solution branches are of primary interest in normal form theory.  In particular, there is a transition in stability at $(x,\mu)=(0,0)$ and $(x,\mu)=(1,-1)$.  To determine the normal form at these locations, we can locally expand each bifurcation with a Taylor series.  Thus for the bifurcation at the origin, we expand
\begin{subeqnarray}
  && x=0 + \tilde{x} \\
  && \mu= 0 + \tilde{\mu}
\end{subeqnarray}
which when inserted to the governing equations gives
\begin{equation}
    \frac{d\tilde{x}}{dt} = \tilde{\mu}\tilde{x} + 2 \tilde{x}^2 .
\end{equation}
These are the leading order terms that emerge from the Taylor expansion, showing that the bifurcation at 
$(x,\mu)=(0,0)$ has a transcritical normal form.  We can also expand in a similar fashion around the bifurcation
point $(x,\mu)=(1,-1)$ so that
\begin{subeqnarray}
  && x=1+ \tilde{x} \\
  && \mu= -1 + \tilde{\mu}
\end{subeqnarray}
which when inserted to the governing equations gives
\begin{equation}
    \frac{d\tilde{x}}{dt} = \tilde{\mu} - \tilde{x}^2 .
\end{equation}
These are the leading order terms that emerge from the Taylor expansion, showing that the bifurcation at 
$(x,\mu)=(1,-1)$ has a saddle bifurcation normal form.  Thus as illustrated in the shaded regions of Fig.~\ref{fig:normal_form_example}, the two bifurcation locations have normal forms of the saddle and transcritical type. 

\subsection*{Normal form of the Lorenz system}
Another example to consider is the well known Lorenz system
\begin{equation}
\frac{d{\bf x}}{dt} = \left[ \begin{array}{ccc} \sigma (y-x) \\ rx-y-zx \\-bz+xy \end{array}   \right]
\end{equation}
which is the equation represented in vector form for the state space ${\bf x}=[x \,\,\, y \,\,\, z]^T$.  For $r<1$, the trivial solutions ${\bf x}=0$ is stable.  Once $r>1$, the solution is no longer stable and the solution bifurcates.  To characterize the bifurcation, one can make a coordinate change based upon perturbation theory ($\epsilon\\ 1$)
\begin{subeqnarray}
  && {\bf x} = {\bf 0} + \epsilon {\bf x}_1 + \epsilon^2 {\bf x}_2 + \cdots \\
  && \tau=\epsilon^2 t \\
  && \epsilon= \sqrt{r-1} \\
  &&  {A} = [1 \,\,\, 1 \,\,\, 0] \frac{{\bf x}}{2\epsilon}
\end{subeqnarray}
which gives the governing equations
\begin{equation}
   A_\tau =\frac{\sigma}{b(1+\sigma)} \left( b A - A^3 \right)
\end{equation}
which is a normal form for the pitchfork bifurcation.   This is the bifurcation structure of the Lorenz from the trivial solution.

\subsection{Imperfect Bifurcations}

Bifurcations are critically important for understanding many parametrically dependent dynamical systems.  Indeed, the onset of instability and switching of solution branches is a canonical process across the engineering and physical sciences. The normal forms  constructed are important for characterizing this process, but they also can be {\em structurally unstable} in their behavior.  Specifically, when normal forms are perturbed, does their basic instability structure survive.
To illustrate this, consider the perturbed form of the transcritical bifurcation
\begin{equation}
    \frac{dx}{dt}= \mu x- x^2 + \delta
\end{equation}
where $\delta\ll 1$ is a perturbation to the problem.  As before, we can consider the equilibrium solutions for which $dx/dt=0$.  This gives
\begin{equation}
    \mu x - x^2 + \delta=0 \,\,\,\,\, \rightarrow \,\,\,\,\, x_0 =\frac{\mu}{2}  \left[  1 \pm \sqrt{1 + 4 \frac{\delta}{\mu^2} } \right]
\end{equation}
with $\mu\neq 0$.

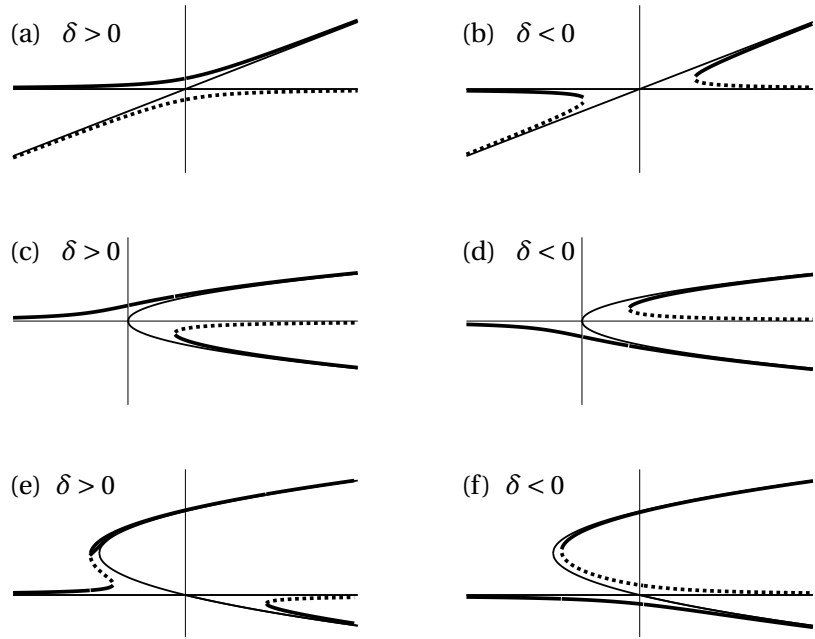
\begin{figure}[t]
\begin{overpic}[width=0.65\textwidth]{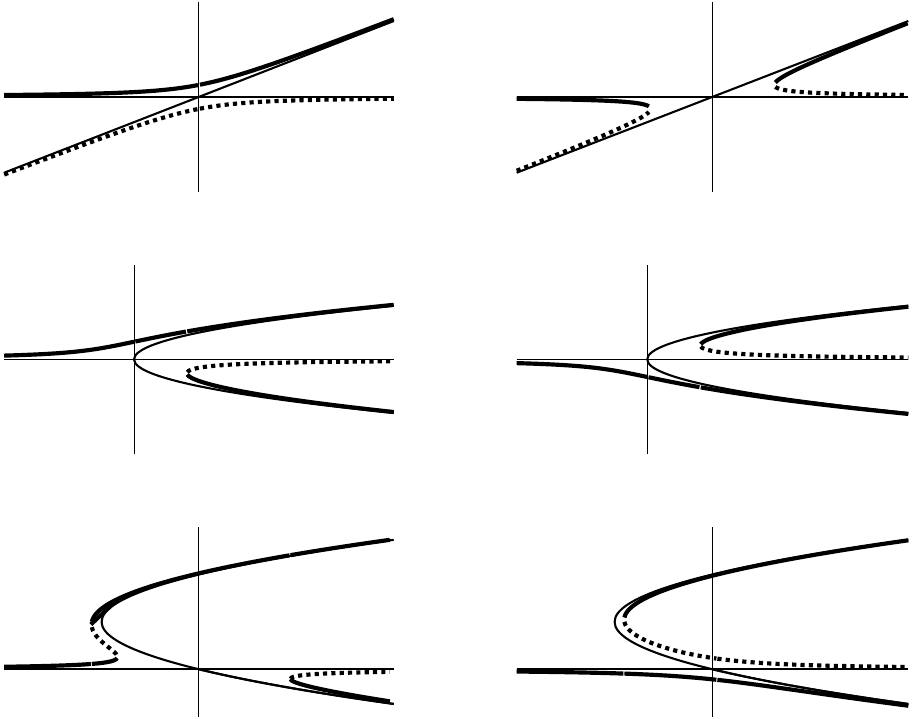}
\put(0,74){(a) \,\, $\delta>0$}
\put(56,74){(b) \,\, $\delta<0$}
\put(0,47){(c) \,\, $\delta>0$}
\put(56,47){(d) \,\, $\delta<0$}
\put(0,18){(e) \,\,$\delta>0$}
\put(56,18){(f)\,\, $\delta<0$}
\end{overpic}
\caption{The structural instabilities of normal form bifurcations, i.e. imperfect bifurcations.  Stable solutions are indicated by the bolded solid line and unstable solutions by the bolded dotted line.  In the top panel, the transcritical bifurcation is destroyed in one of two ways depending on wether $\delta$ is positive or negative.  The pitchfork bifurcation (middled panels) is also structurally unstable due to the considered perturbation.  Finally, the bifurcation structure of Fig.~\ref{fig:normal_form_example}  is also shown to be structurally unstable.  For every panel, either $\delta=\pm 0.1$.}
 \label{fig:bif_imperfect}
\end{figure}

Figure~\ref{fig:bif_imperfect} shows that the introduction of this perturbation destroys the normal form bifurcation for any $\delta\neq 0$.  Thus the transcritical bifurcation is structurally unstable to this perturbation.    The stability of each branch can once again be determined by linearization so that $x=x_0 + \tilde{x}$.  In this case, the top branch of solutions is stable while the bottom branch is unstable for $\delta>0$.  For $\delta<0$, there is now a forbidden region of $\delta$ around the original where solutions no longer exist.  The  top branch of solutions are stable for both the left and right solution branches.

The destruction of the bifurcation structure is also self-evident for the pitchfork normal form when perturbed
\begin{equation}
    \frac{dx}{dt}= \mu x- x^2 + \delta
\end{equation}
where $\delta\ll 1$ is a perturbation to the problem.
As before, we can consider the equilibrium solutions for which $dx/dt=0$.  This gives
\begin{equation}
    \mu x - x^3 + \delta=0 
\end{equation}
which is a cubic equation for determining the equilibrium solutions.
Figure~\ref{fig:bif_imperfect} shows that the introduction of this perturbation destroys the normal form bifurcation for any $\delta\neq 0$.   Note the resulting bifurcations are still important to understand even thought the structure itself is unstable.  Of course, not all bifurcations destroy the bifurcation structure, but it is important to understand that these canonical forms of instability, i.e. the normal forms, can be mitigated by appropriately applying a perturbation.  Also shown in Fig.~\ref{fig:bif_imperfect} is the example used in Fig.~\ref{fig:normal_form_example} of the normal form section where the bifurcation structure is shown to be structurally unstable.

\newpage
\section*{Lecture 21:  Pattern Forming Systems:  An Introduction}

Asymptotics and perturbation theory again play a critical role in stability theory for partial differential equations and pattern forming, spatio-temporal systems.  Generically, we can consider the nonlinear PDE
\begin{equation}
  u_t = N(u, u_x, u_{xx}, \cdots, x, t, \mu)
\end{equation}  
which characterizes the evolution of a spatio-temporal field $u(x,t)$ and where $\mu$ is a parameter that induces a bifurcation at some parameter value $\mu=\mu_c$.  If one can find a solution to this system, denoted by $u_0(x,t)$, then stability of the solution can be determined by a perturbation expansion of the form
\begin{equation}
  u = u_0(x,t) + \epsilon \tilde{u}(x,t)  .
\end{equation}
Plugging in this expansion into the governing equations and collecting orders in $\epsilon$ gives
\begin{eqnarray}
  &&  O(1):  \,\,\,\,\,  {u_0}_t = N(u_0, {u_0}_x, {u_0}_{xx}, \cdots, x, t,\mu)   \nonumber \\
  &&  O(\epsilon):  \,\,\,\,\,  {\tilde{u}}_t = L(u_0, {u_0}_x, {u_0}_{xx}, \cdots,x,t,\mu) \tilde{u}
\end{eqnarray}
where $O(\epsilon^2)$ terms are ignored.   Stability is determined by the dynamics of $\tilde{u}(x,t)$ which actually satisfied a {\em linear} evolution equation characterized by the linear operator $L$.  Exponential solutions are the canonical solutions to linear equations so that one assumes
\begin{equation}
   \tilde{u}(x,t)=v(x) \exp(\lambda t)
\end{equation}
which results in the eigenvalue problem
\begin{equation}
   Lv=\lambda v .
\end{equation}
The eigenvalues of the linear operator then determines the {\em spectral} linear stability of the solution.  Specifically, if any eigenvalue has $\Re \{\lambda\}>0$ the system is unstable.  In contrast if $\Re\{\lambda\} <0$ for all eigenvalues the system is stable.  If $\Re\{\lambda\} \leq 0$, then the system if typically classified as neutrally stable, although zero eigenvalues are often thought of as generating a center manifold for which higher-order perturbation theory can often determine stability.

It should be noted that the eigenvalues 
\begin{equation}
\lambda = \lambda(\mu),
\end{equation}
thus the behavior of the eigenvalues as a function of the bifurcation parameter $\mu$ is critical in determining the stability of the solutions.  In particular, often there are critical values of $\mu=\mu_c$ where an eigenvalue crosses from the left half plane (negative real part) to the right half plane (positive real part), leading to an instability of the solution $u_0(x,t)$.

\subsection*{The Fisher-Kolmogorov Equation}

Consider the Fisher-Kolmogorov nonlinear evolution equation
\begin{equation}
   u_t - u_{xx} + u^3  - \mu u = 0 .
   \label{eq:FK_stab}
\end{equation}
The evolution equation admits constant amplitude solutions when
\begin{equation}
  u^3-\mu u = u (u^2 - \mu) =0  \,\,\, \Longrightarrow \,\,\,  u_0=0, u_0=\pm \sqrt{\mu} .
\end{equation}
A linear stability analysis can be performed for each solution branch.  For $u=u_0 + \epsilon \tilde{u}$, the linearized dynamics gives
\begin{equation}
  \tilde{u}_t -\tilde{u}_{xx} + 3 u_0^2 \tilde{u} - \mu \tilde{u} = 0,
\end{equation}
Fourier transforming this equation in $x$ gives the differential equation
\begin{equation}
   \hat{u}_t + \left( k^2 + 3u_0^2  - \mu \right) \hat{u} =0  \,\,\, \Longrightarrow \,\,\, \hat{u} = \hat{u}_0 \exp \left[ -\left(k^2 +3u_0^2 -\mu \right) t \right] = \hat{u}_0 \exp \left[ -\omega t \right]
\end{equation}
For $\omega=k^2 +3u_0^2 -\mu>0$ the solution is stable, while for $\omega=k^2 +3u_0^2 -\mu<0$ exponential growth ensues making the solution unstable.

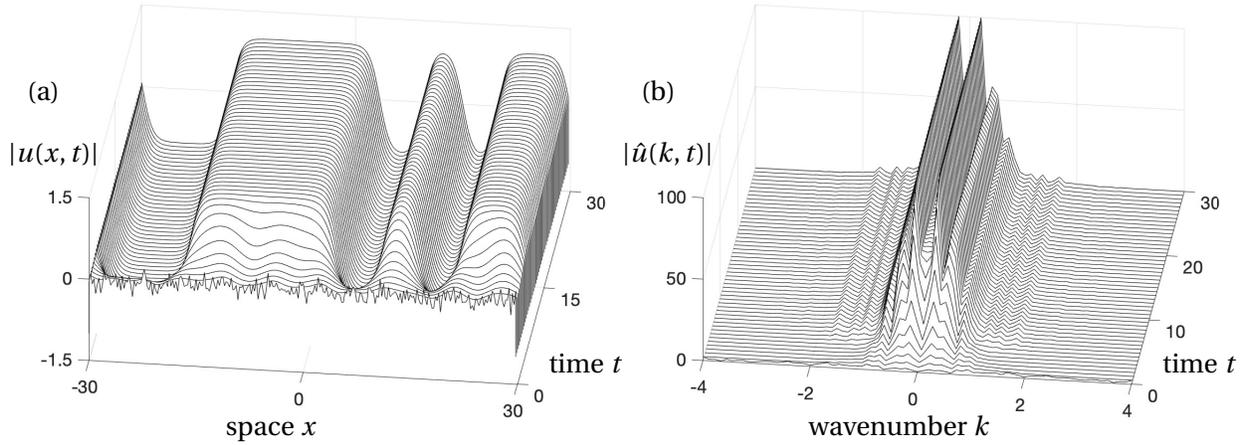
\begin{figure}[t]
\begin{overpic}[width=0.5\textwidth]{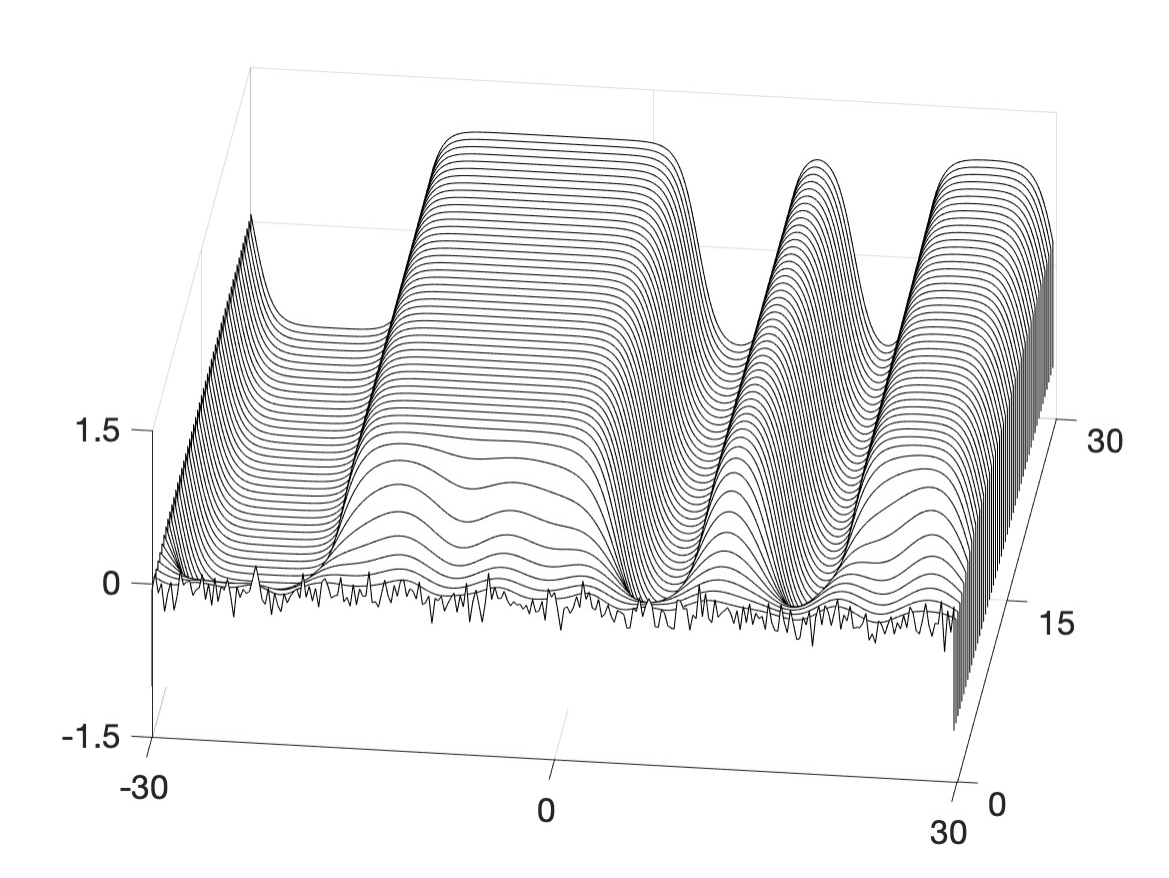}
\put(3,54){(a)}
\put(0,44){$|{u}(x,t)|$}
\put(35,0){space $x$}
\put(87,10){time $t$}
\end{overpic}
\hspace*{-.1in}
\begin{overpic}[width=0.5\textwidth]{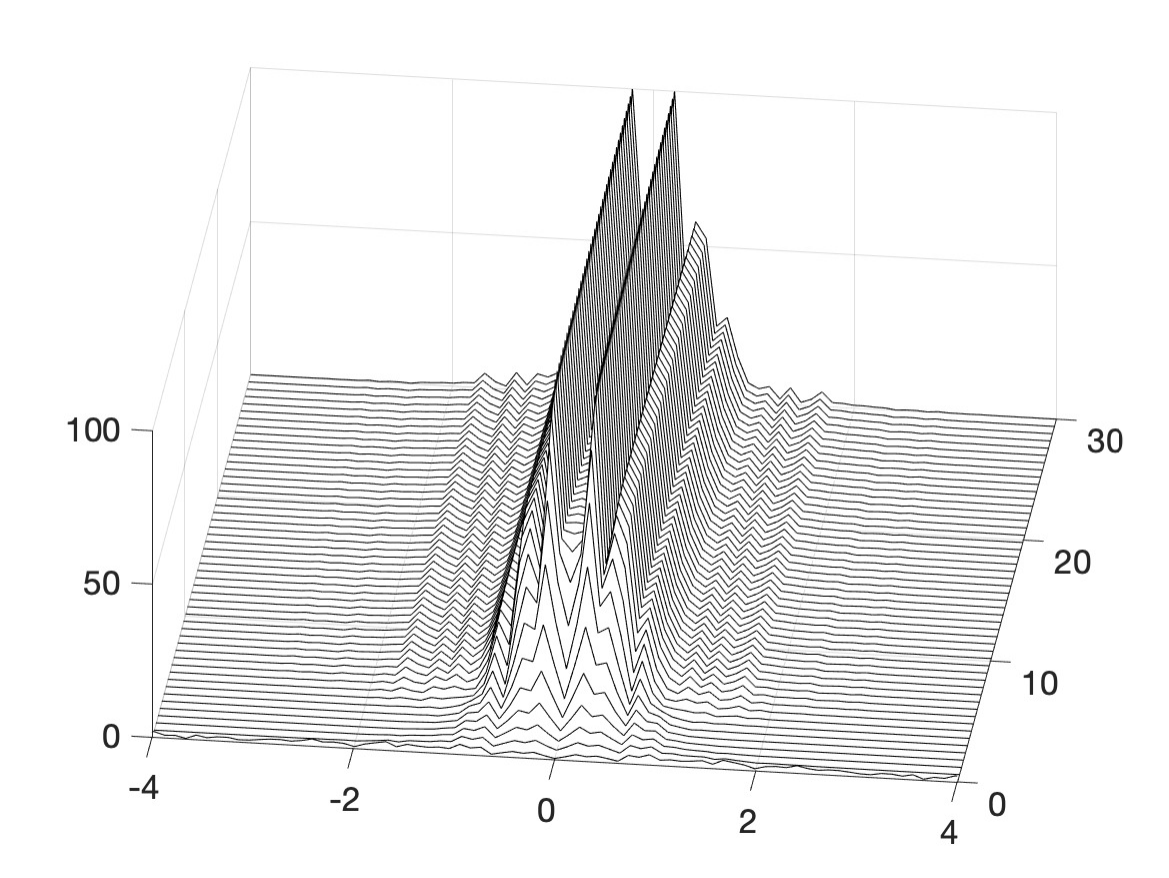}
\put(30,0){wavenumber $k$}
\put(87,10){time $t$}
\put(0,44){$|\hat{u}(k,t)|$}
\put(3,54){(b)}
\end{overpic}
\caption{Instability of plane wave solutions of the Fisher Kolmogorov equation (\ref{eq:FK_stab}) with $\mu=1$.  (a)  The spatio-temporal evolution of the field $u(x,t)$ shows that the initial conditions evolve to either the $u_0=\pm\sqrt{\mu}$ solution branches, generating a heteroclinic connection between these branches spatially.   (b)  Corresponding spectral evolution $\hat{u}(k,t)$ of the field.}
 \label{fig:mi_FK}
\end{figure}

For the trivial solutions $u_0=0$, when $\mu<0$ the dispersion relation gives $\omega=k^2+|\mu|>0$ which is always positive so that that trivial branch for $\mu<0$ is stable.  For $\mu>0$, the trivial branch of solutions is unstable since $\omega=k^2-\mu<0$ for wavenumbers less than $k<\sqrt{\mu}$.  Indeed, the most unstable wavenumber is for $d\omega/dk=0$ which occurs at $k_{max}=0$.  At $\mu=0$, there is a pitchfork bifurcation to the two branches of solutions $u_0=\pm \sqrt{u}$.  These solutions do not exist for $\mu<0$, but for $\mu>0$ they are stable since $\omega=k^2+3u_0^2-\mu = k^2 + 2\mu >0$ for all wavenumbers.  Thus the non-trivial plane wave solutions are stable, bifurcating in a pitchfork from $\mu=0$.  Figure~\ref{fig:mi_FK} shows the evolution of an initial noisy state near the unstable trivial solution for $\mu=1$.  Note that solution evolves toward either the $u_0=\pm\sqrt{\mu}=\pm 1$ branch of solutions for $\mu=1$.

\subsection*{The Kuramoto-Sivashinsky Equation}

Consider the nonlinear evolution equation known as the Kuramoto-Sivashinsky (KS) equation. 
\begin{equation}
   u_t +\mu u_{xxxx}  + u_{xx} + u u_x  = 0
   \label{eq:ks_mi}
\end{equation}
The KS admits the trivial solution $u_0=0$.  A linear stability analysis about the trivial solution yields
\begin{equation}
  \tilde{u}_t + \mu \tilde{u}_{xxxx} + \tilde{u}_{xx} = 0 .
\end{equation}
Fourier transforming this equation in $x$ gives the differential equation
\begin{equation}
   \hat{u}_t  + \left( \mu k^4 - k^2 \right) \hat{u} = 0 \,\,\, \Longrightarrow \,\,\, \hat{u} = \hat{u}_0 \exp \left[ \left(k^2-\mu k^4 \right) t \right] = \hat{u}_0 \exp \left[ \omega^2 t \right]
   \end{equation}
 If $\mu<0$, then the KS equation is ill-posed as the fourth-order diffusion leads to exponential blow-up of solutions.   For $\mu>0$, the fourth-order diffusion regularizes the evolution equation by providing higher-order diffusion.  In this case, there is a small band of unstable frequencies for which the exponential coefficient $\omega^2(k,\mu)$ is then positive and leads to growth.  Specifically, consider the argument of the exponential $\omega^2=\left(k^2-\mu k^4\right)$.  For $k>1/\sqrt{\mu}$, the argument is positive and such frequencies decay to zero.  However, all frequencies for which $k<1/\sqrt{\mu}$ are unstable.  In fact, the most unstable growth can be found to be at $k_{max}=1/\sqrt{2 \mu}$.  Thus a spatial modulational at wavenumber $k_{max}$ dominates the instability observed. 

\begin{figure}[t]
\begin{overpic}[width=0.45\textwidth]{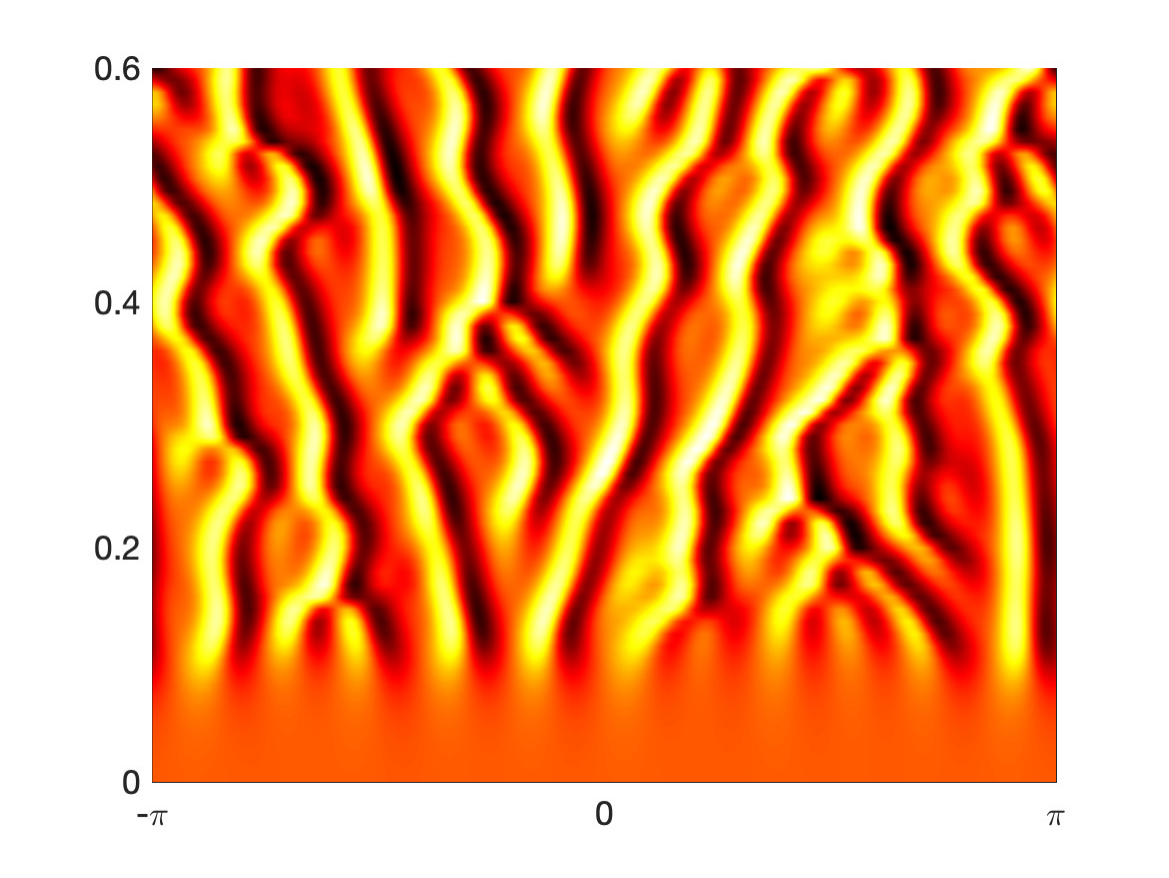}
\put(3,59){(a)}
\put(45,72){$|{u}(x,t)|$}
\put(45,0){space $x$}
\put(3,35){\rotatebox{90}{time $t$}}
\end{overpic}
\hspace*{-.1in}
\begin{overpic}[width=0.5\textwidth]{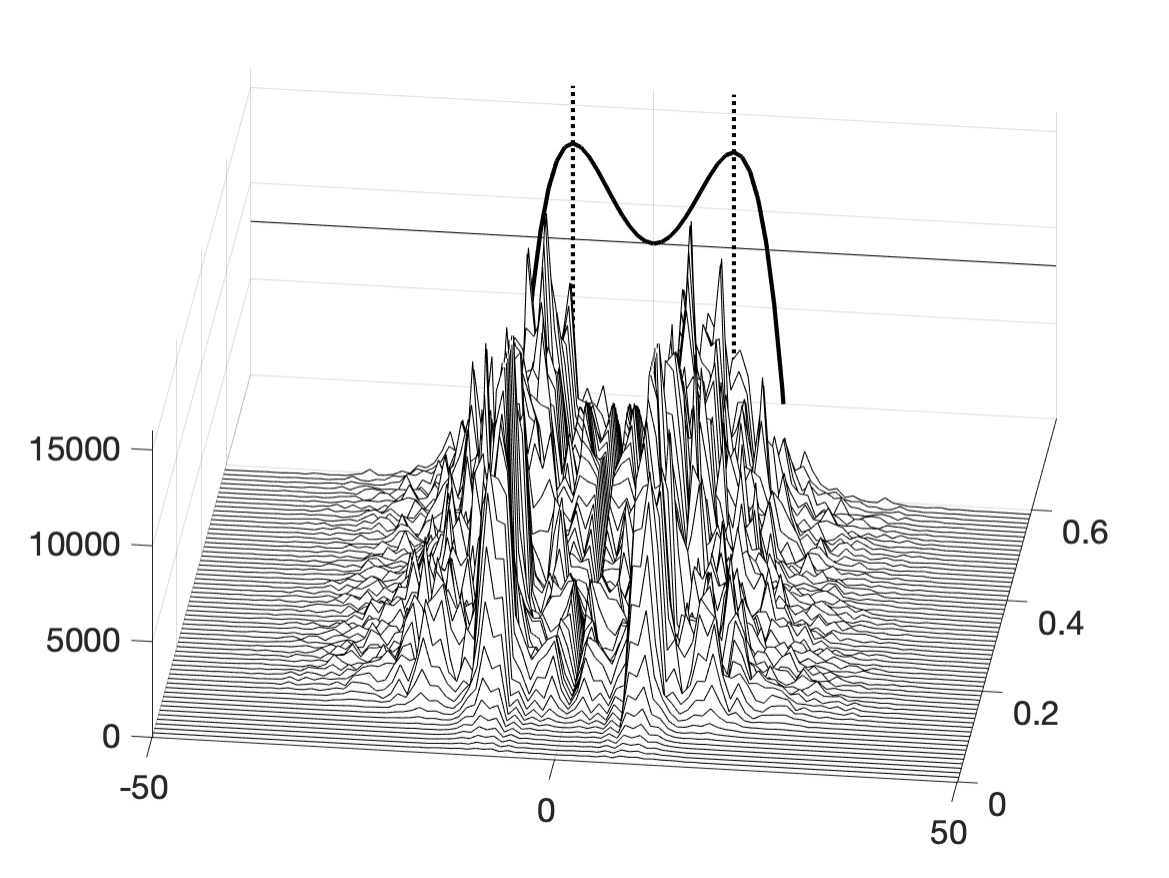}
\put(60,68){$k_{max}=1/\sqrt{2\mu}$}
\put(68,56){$\omega^2=k^2-\mu k^4$}
\put(89,48){$\omega^2=0$}
\put(30,0){wavenumber $k$}
\put(87,10){time $t$}
\put(0,44){$|\hat{u}(k,t)|$}
\put(3,54){(b)}
\end{overpic}
\begin{overpic}[width=0.45\textwidth]{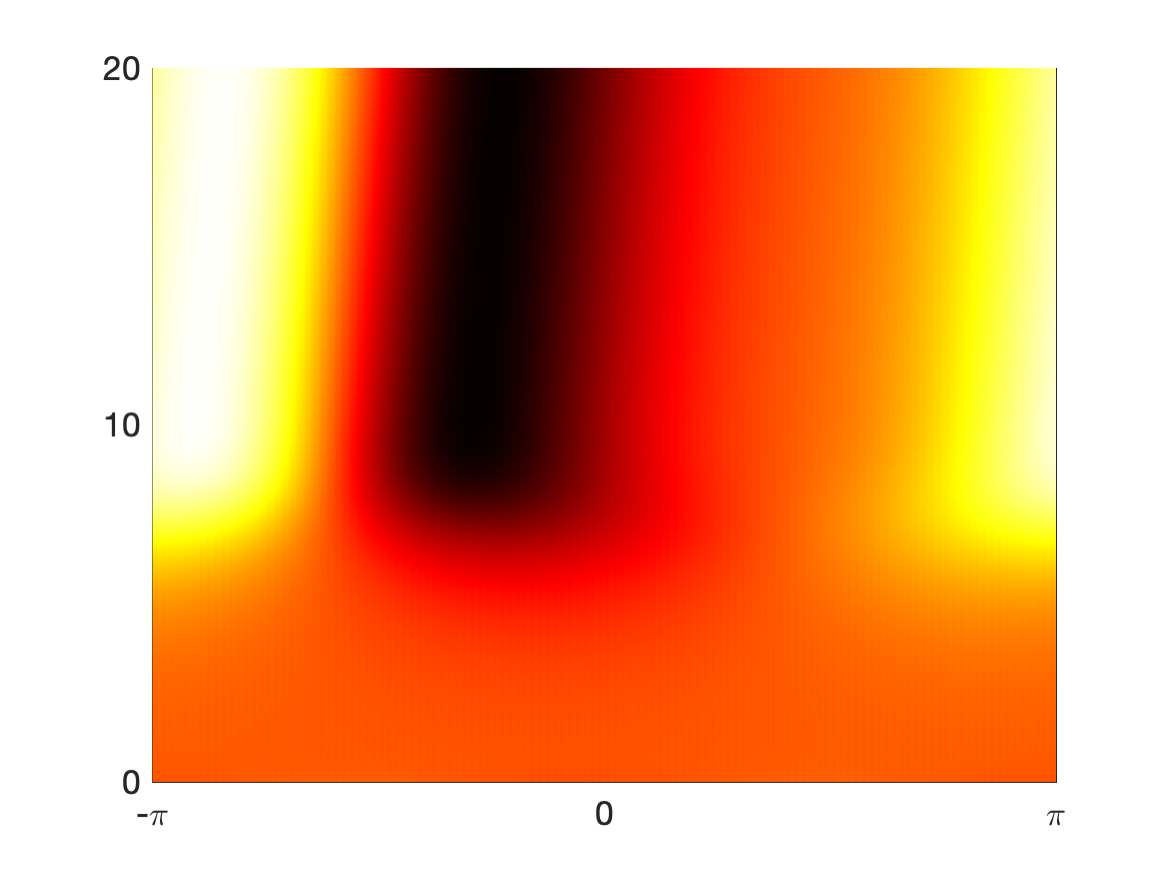}
\put(3,59){(c)}
\put(45,72){$|{u}(x,t)|$}
\put(45,0){space $x$}
\put(3,35){\rotatebox{90}{time $t$}}
\end{overpic}
\hspace*{-.1in}
\begin{overpic}[width=0.5\textwidth]{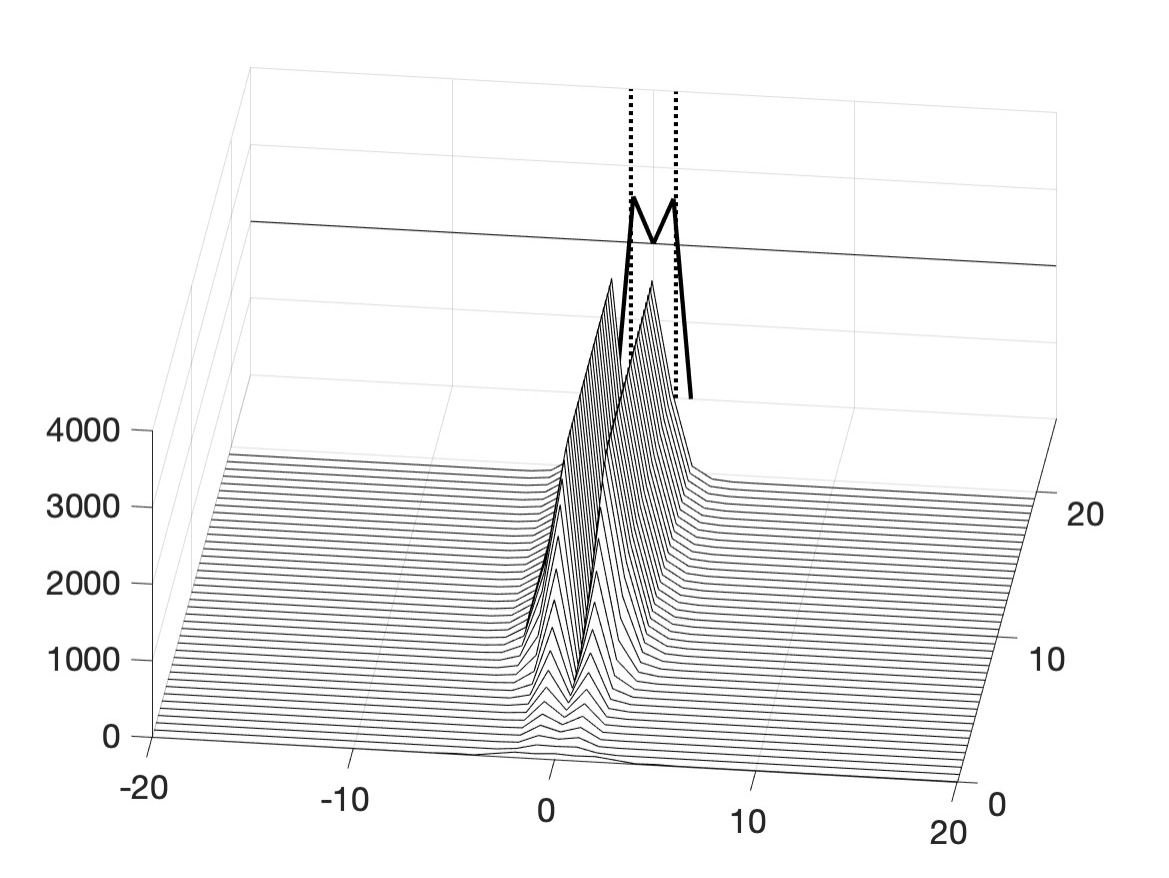}
\put(55,68){$k_{max}=1/\sqrt{2\mu}$}
\put(62,56){$\omega^2=k^2-\mu k^4$}
\put(89,48){$\omega^2=0$}
\put(30,0){wavenumber $k$}
\put(87,10){time $t$}
\put(0,44){$|\hat{u}(k,t)|$}
\put(3,54){(d)}
\end{overpic}
\caption{Modulational instability of trivial solution of the KS equation for $\mu=0.005$ (top panels) and $\mu=0.4$ (bottom panels).  (a)  The spatio-temporal evolution of the field $u(x,t)$ governed by (\ref{eq:ks_mi}) with $\mu=0.005$.  Note the broadband modulational structure of the instability as it develops in time.  (b)  Corresponding spectral evolution $\hat{u}(k,t)$ of the complex field for $\mu=0.005$.  Note the growth of side-bands at the predicted most unstable wavenumber $k_{max}=1/\sqrt{2\mu}$.  The growth rate of the entire unstable band of frequencies is shown by the curve $\omega^2 =k^2-\mu k^4$.  Note that the linear stability theory and numerical experiments are in strong agreement.  In (c) and (d), the numerical experiments are repeated with $\mu=0.4$.  The strong damping limits the unstable wavenumbers significantly, leading to a simple roll pattern.}
 \label{fig:mi_ks}
\end{figure}
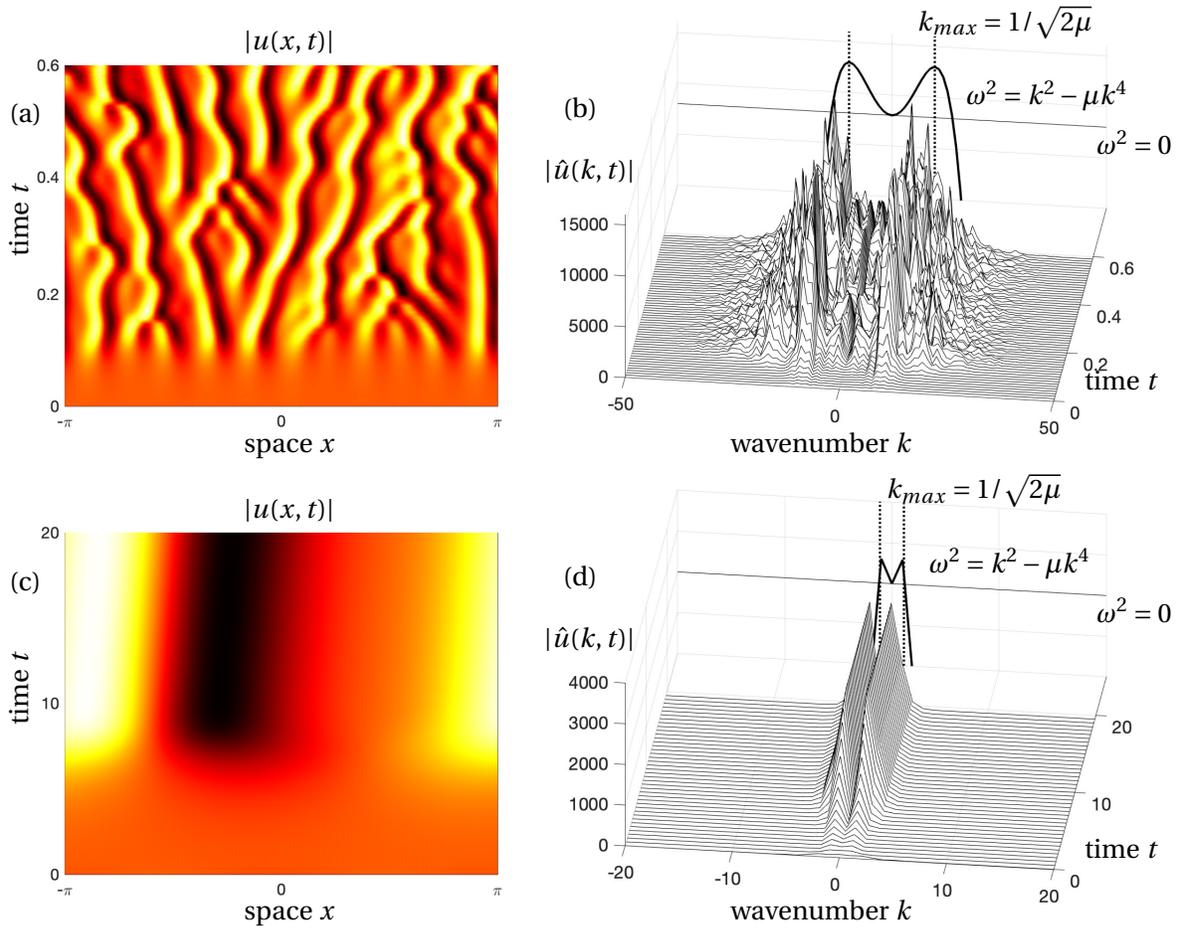

Figure~\ref{fig:mi_ks} compares the linear stability theory with direct numerical simulations of the governing KS equation for $\mu=0.4$ and $\mu=0.005$.   Specifically, for both values of $\mu$ a modulational instability is predicted to occur with a dominant growth spatial frequency of $k_{max}=1/\sqrt{2\mu}$.  As shown in Fig.~\ref{fig:mi_ks}, a modulational instability does indeed occur.  The corresponding frequency evolution $\hat{u}(k,t)$ is shown as well along with the predicted maximum grown rate.  The predicted dominant growth wavenumber $k_{max}$ is shown to be corroborated by the direct numerical simulation.  Note that for $\mu=0.005$, the fourth-order diffusion is small and a large band of frequencies goes unstable, leading to the chaotic behavior observed. For $\mu=0.4$, the fourth-order diffusion is much stronger and strongly damps the system, leading to only a single unstable wavenumber in the simulation.

\subsection*{Modulational Stability of Planes Waves of the Nonlinear Schr\"odinger Equation}

Consider the nonlinear Schr\"odinger equation (NLS)
\begin{equation}
   i u_t  + \frac{\mu}{2} u_{xx} + |u|^2 u = 0
   \label{eq:nls_mi}
\end{equation}
on the domain $x\in[-\infty,\infty]$.  Plane wave, or often called {\em continuous wave} (CW), solutions exist for the NLS in both the focusing ($\mu=1$, anomalous dispersion) and defocusing ($\mu=-1$, normal dispersion)  parameter regimes of this equation.  In particular, the solution is given by
\begin{equation}
   u_0(x,t)= A\exp(i|A|^2 t) 
   \label{eq:nls_cw}
\end{equation}
for an arbitrary constant $A$.  The stability of the CW solution is determined by linearizing about this solution
\begin{equation}
  u(x,t)= \left(  A + \tilde{u}(x,t) \right)  \exp( i |A|^2 t )
\end{equation}
where the nonlinear phase rotation has been explicitly removed from the linear perturbation $\tilde{u}(x,t)$.  The expansion yields at next order the linear evolution dynamics
\begin{equation}
  i \tilde{u}_t + \frac{\mu}{2} \tilde{u}_{xx} + A^2 \left(  \tilde{u} + \tilde{u}^* \right) =0 .
\end{equation}
This linear PDE can be reduced using a Fourier transform in $x$ so that
\begin{equation}
   i \hat{u}_t - \frac{\mu k^2}{2} \hat{u} + A^2 (\hat{u} + \hat{u}^*) = 0
\end{equation}
where $\hat{u}$ denotes the Fourier transform of $\tilde{u}$ and the parameter $k$ is the wavenumber associated with the transform.  
The resulting differential equations can be evaluated by separating $\hat{u}$ into its real ($R$) and imaginary ($I$) components so that 
\begin{equation}
  \hat{u}=R + iI . 
\end{equation}
Inserting this decomposition into the linear evolution equation for $\hat{u}$ yields the system
\begin{equation}
   \dot{\bf x} = \left[  \begin{array}{cc}  0 & \frac{\mu k^2}{2} \\ 2A^2 - \frac{\mu k^2}{2} & 0 \end{array} \right] {\bf x}
\end{equation}
where ${\bf x} = [R \,\,\, I]^T$.  Stability can be evaluated by considering exponential solutions of the form ${\bf x}={\bf v}\exp(\omega t)$ which yields an eigenvalue problem whose eigenvalues are then
\begin{equation}
  \omega^2 = - \frac{k^4}{4} +\mu k^2 A^2
\end{equation}
where $\mu^2 = 1$ for either $\mu=\pm 1$.
If $\mu=-1$, the defocusing case of NLS, then the eigenvalues $\omega = \pm i \sqrt{k^4/4 + k^2 A^2}$ are purely imaginary and the CW solutions are stable.  If $\mu=1$, the focusing case of NLS, then the eigenvalues $\omega^2 = - \frac{k^4}{4} + k^2 A^2$ have can have a real component which is positive for sufficiently small wavenumber values $k$.  In particular, for the range $0<k<2A$, the $\omega$ have a real, positive values which promote instability, whereas for $k>2A$ the real part of the eigenvalues are negative and do not promote instability.   The maximum growth rate is determined by $d\omega/dk=0$ which gives $k_{max} =\sqrt{2} A$ as the most unstable wavenumber.  Thus the CW solution is unstable and dominated by a spatially varying wave with wavenumber $k_{max}$.  This is a modulational instability since the growth is dominated by a non-zero spatial wavenumber.

\begin{figure}[t]
\begin{overpic}[width=0.5\textwidth]{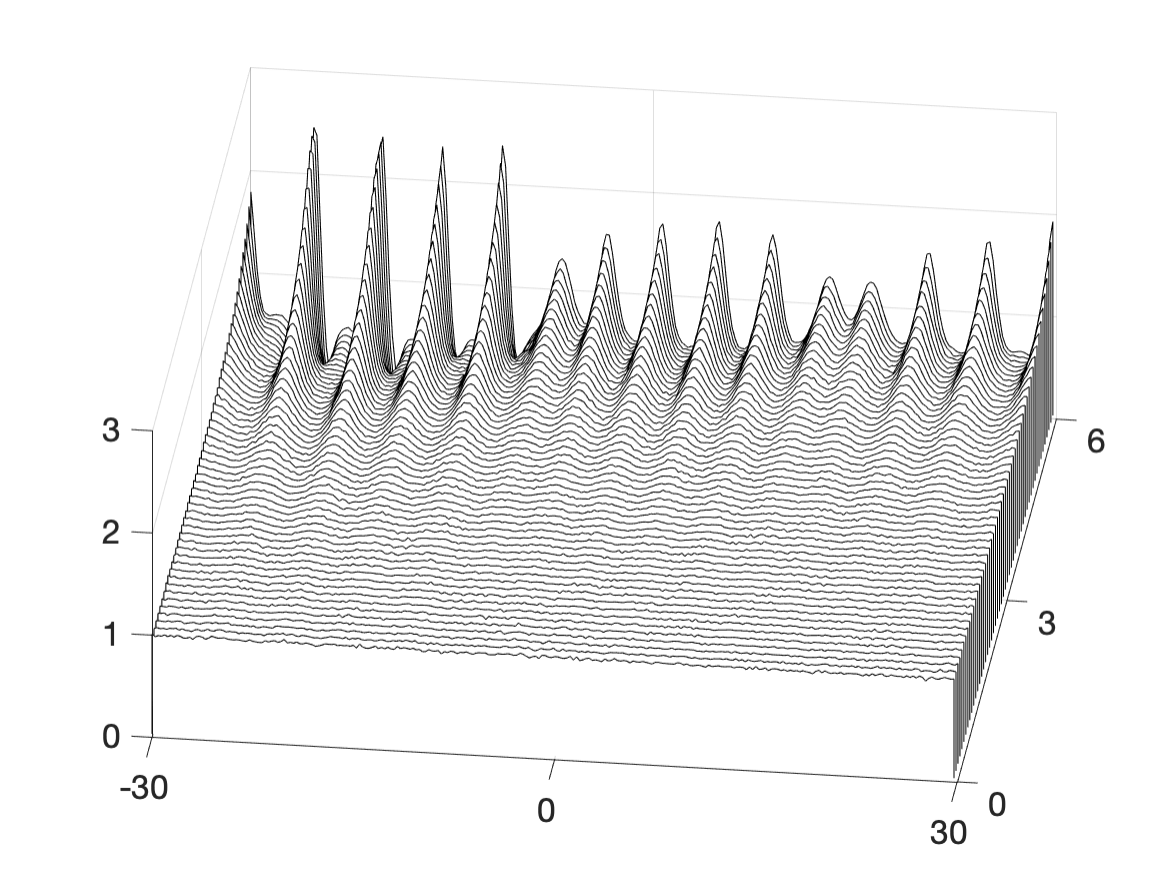}
\put(3,54){(a)}
\put(0,44){$|{u}(x,t)|$}
\put(35,0){space $x$}
\put(87,10){time $t$}
\end{overpic}
\hspace*{-.1in}
\begin{overpic}[width=0.5\textwidth]{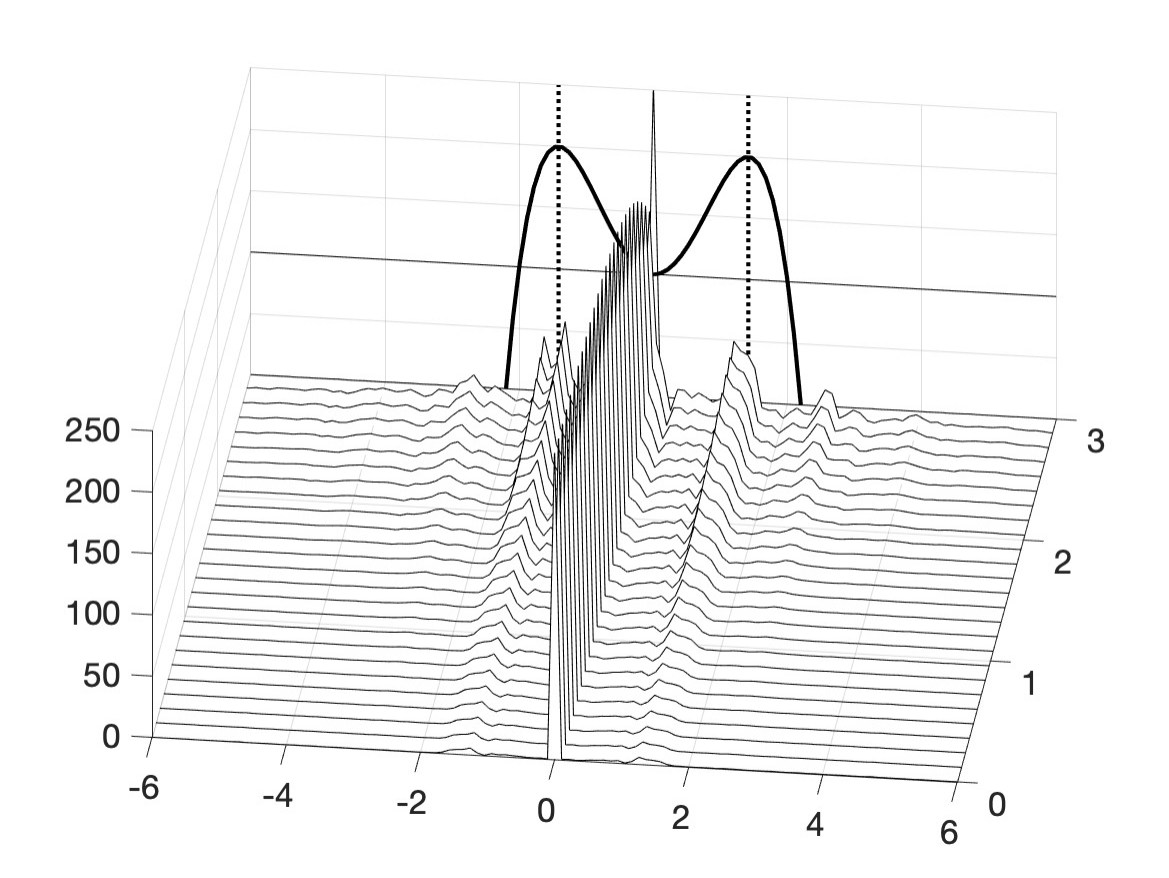}
\put(60,68){$k_{max}=\sqrt{2}$}
\put(68,56){$\omega^2=-k^4/4+k^2$}
\put(89,46){$\omega^2=0$}
\put(30,0){wavenumber $k$}
\put(87,10){time $t$}
\put(0,44){$|\hat{u}(k,t)|$}
\put(3,54){(b)}
\end{overpic}
\caption{Modulational instability of plane wave solution (\ref{eq:nls_cw}) with $A=1$ of the nonlinear Schr\"odinger equation.  (a)  The spatio-temporal evolution of the complex field $u(x,t)$ governed by (\ref{eq:nls_mi}) with $\mu=1$.  Note the modulational structure of the instability as it develops in time.  (b)  Corresponding spectral evolution $\hat{u}(k,t)$ of the complex field.  Note the growth of side-bands at the predicted most unstable wavenumber $k_{max}=\sqrt{2}A$ (note that $A=1$).  The growth rate of the entire unstable band of frequencies is shown by the curve $\omega^2 =-k^4/4+k^2$ with $\omega^2>0$.  Note that the linear stability theory and numerical experiments are in strong agreement.}
 \label{fig:mi_nls}
\end{figure}
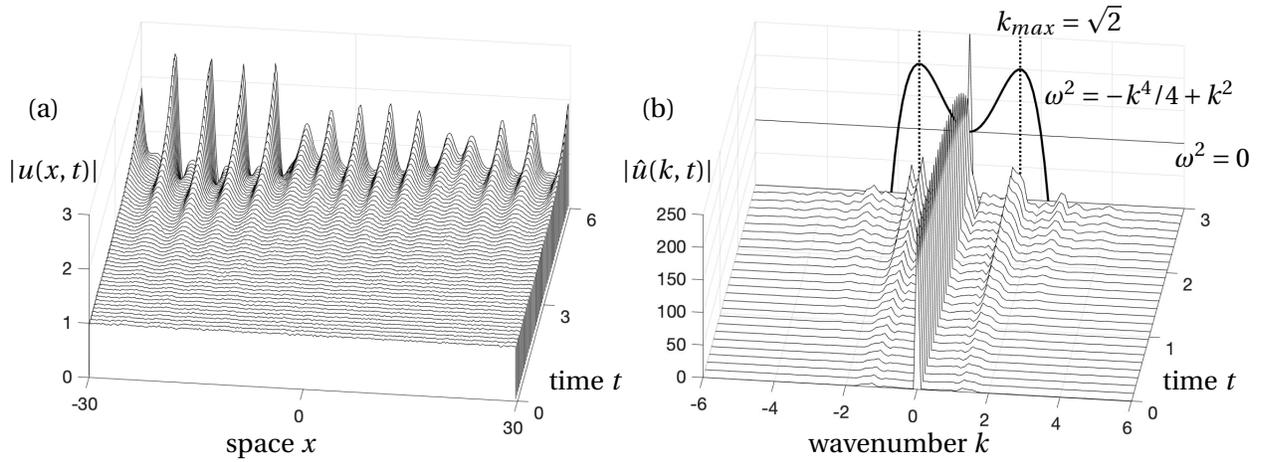

Figure~\ref{fig:mi_nls} compares the linear stability theory with direct numerical simulations of the governing NLS equation.   Specifically, for $\mu=1$, a modulational instability is predicted to occur with a dominant growth spatial frequency of $k_{max}=\sqrt{2}A$.  As shown in Fig.~\ref{fig:mi_nls}(a), a modulational instability does indeed occur.  The corresponding frequency evolution $\hat{u}(k,t)$ is shown in panel (b) along with the predicted maximum grown rate.  The predicted dominant growth wavenumber $k_{max}$ is shown to be corroborated by the direct numerical simulation.

\subsection*{Canonical Pattern Forming Instabilities}

Linear stability analysis plays a critical role in understanding the canonical instabilities that can occur in spatio-temporal PDE systems.  As was shown for both the nonlinear Schr\"odinger and Kuramoto-Sivashinsky equations, modulational instabilities can dominate the patterns that form in the PDE system.
To give a more general treatment of the instabilities which can occur in PDE systems we once again consider the  nonlinear PDE
\begin{equation}
  u_t = N(u, u_x, u_{xx}, \cdots, x, t, \mu)
  \label{eq:nonlinearPDE}
\end{equation}  
which characterizes the evolution of a spatio-temporal field $u(x,t)$ and where $\mu$ is a parameter that induces a bifurcation at some parameter value $\mu=\mu_c$.  Again, if one can find a solution to this system, denoted by $u_0(x,t)$, then stability of the solution can be determined by a perturbation expansion of the form $u = u_0(x,t) + \epsilon \tilde{u}(x,t)$ for which the perturbed field evolves according to the linear PDE
\begin{equation}
   {\tilde{u}}_t = L(u_0, {u_0}_x, {u_0}_{xx}, \cdots,x,t,\mu) \tilde{u} .
\end{equation}
Stability is determined by the dynamics of $\tilde{u}(x,t)$ which actually satisfied a {\em linear} evolution equation characterized by the linear operator $L$.  

A standard method for analyzing the behavior of the perturbed field is to evaluate the evolution dynamics of the exponential, traveling wave solution
\begin{equation}
  \tilde{u}(x,t)=A \exp \left[   i \left(  kx - \omega t  \right)  \right]
  \label{eq:TW}
\end{equation}
where $k$ is the wavenumber and $\omega$ is the frequency.  In multi-dimensional settings (2D and 3D, for instance), the spatial characterization $kx$ is replaced by ${\bf k}\cdot {\bf x}$.

Inserting (\ref{eq:TW}) into the linear evolution equation for $\tilde{u}(x,t)$ yields an equation relating the wavenumber $k$ to the frequency $\omega$
\begin{equation}
   G(k,\omega)=0
\end{equation}
which is known as the {\em dispersion relation}.  The dispersion relation provides a fundamental characterization of the linear dynamics.  Moreover, the dispersion relation will also determine the nature of the fundamental instabilities observed in the PDEs, and ultimately the patterns which form in the spatio-temporal domain.

To simplify the analysis, let $\lambda=-i\omega$ so that the exponential solution (\ref{eq:TW}) takes the form
\begin{equation}
  \tilde{u}(x,t)=A \exp \left[  \lambda t +  i  kx  \right] .
  \label{eq:TW2}
\end{equation}
Stability of solutions is determined by the real part $\Re \{\lambda\}$.  Specifically recall that if for any eigenvalue $\Re \{\lambda\}>0$, the solution grows and is therefore unstable.  In contrast, if all eigenvalues are such that  $\Re \{\lambda\}<0$, the solution is stable.  

Importantly, from the viewpoint of stability, one should amend the dispersion relation to include the bifurcation parameter $\mu$ so that
\begin{equation}
   G(k,\omega,\mu)=0 
\end{equation}
and the dependence on $\mu$ is explicitly accounted for.  This also implies
\begin{subeqnarray}
  &&  \lambda=\lambda(\mu) \\
  &&  k=k(\mu)
\end{subeqnarray}
so that the temporal and spatial parametrizations are explicitly made dependent upon the parameter $\mu$.   Of particular importance is when
\begin{equation}
   \Re \{\lambda \} = \Re \{\lambda (\mu_c) \} =0.
\end{equation}
where we assumed that for $\mu>\mu_c$ the $ \Re \{\lambda \} >0$ and the solution is unstable while for $\mu<\mu_c$ the $ \Re \{\lambda \} <0$ and the solution is stable.  Thus value of  $\mu=\mu_c$ is the bifurcation point leading to unstable behavior and one of the canonical patterns of instability in the spatio-temporal system.  At $\mu=\mu_c$, we denote
\begin{subeqnarray}
  &&  \Im \{\lambda (\mu_c)\}=\omega_0 \\
  &&  k (\mu_c) = k_0
  \end{subeqnarray}
which gives the oscillation frequency $\omega_0$ at the bifurcation point along with the corresponding spatial wavenumber $k_0$.

\begin{figure}[t]
\vspace*{-.6in}
\begin{overpic}[width=0.9\textwidth]{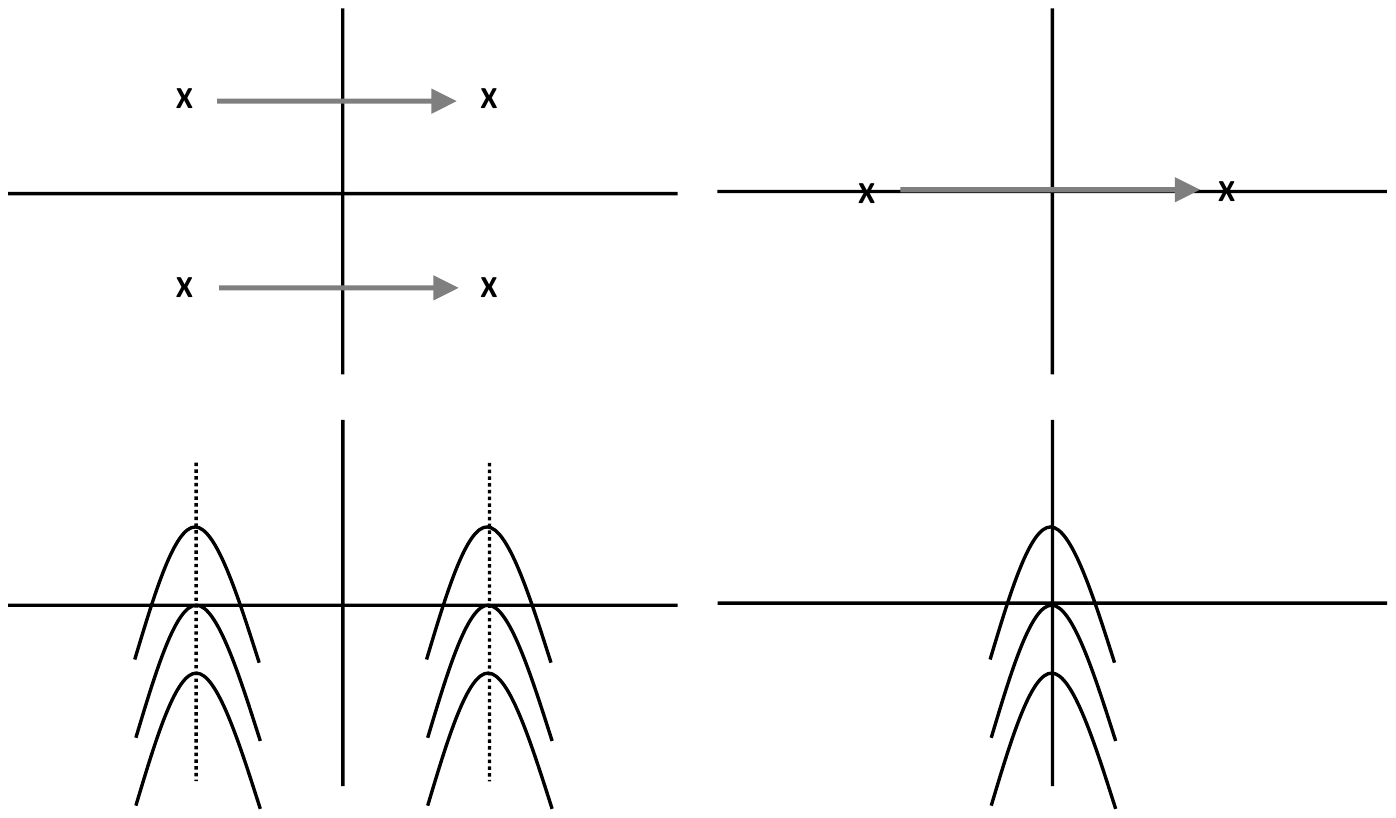}
\put(45,50){$\Re \{\lambda \}$}
\put(88,50){$\Re \{\lambda \}$}
\put(21,62){$\Im \{\lambda \}$}
\put(64,62){$\Im \{\lambda \}$}
\put(15,55){$\mu<\mu_c$}
\put(34,55){$\mu>\mu_c$}
\put(56,54.5){$\mu<\mu_c$}
\put(78,54.5){$\mu>\mu_c$}
\put(21,37){$\Re \{\lambda \}$}
\put(64,37){$\Re \{\lambda \}$}
\put(46,25){$k$}
\put(89,25){$k$}
\put(37,33){$k_0$}
\put(19,33){$-k_0$}
\put(71,33){$k_0=0$}
\put(7,15){$\mu<\mu_c$}
\put(7,19.5){$\mu=\mu_c$}
\put(7,24){$\mu>\mu_c$}
\put(59,15){$\mu<\mu_c$}
\put(59,19.5){$\mu=\mu_c$}
\put(59,24){$\mu>\mu_c$}
\put(5,62){(a)}
\put(55,62){(b)}
\put(5,37){(c)}
\put(55,37){(d)}
\end{overpic}
\vspace*{-.9in}
\caption{Canonical characterization of the instabilities in spatio-temporal systems as $\mu$ crosses from stable ($\mu<\mu_c$) to unstable ($\mu>\mu_c$).  (a)  A pair of eigenvalues (${\bf x}$) cross into the right-half plane with the imaginary part of the eigenvalues being $\pm i \omega_0$.  This is the canonical Hopf bifurcation.  (b). An eigenvalue (${\bf x}$) crosses into the right-half plane with $\omega_0=0$, resulting an a growth mode with with no temporal oscillations.  (c)  The corresponding unstable wavenumbers which has a band of spatial wavenumbers generating a modulational instability.  (d)  The spatial instability which grows with wavenumber $k_0=0$, thus generating the growth of a spatially homogenous plane wave.  }
 \label{fig:patterns}
\end{figure}

The values of $\omega_0$ and $k_0$ determine the nature of the pattern forming instability.  Figure~\ref{fig:patterns} depict the canonical instabilities in spatio-temporal systems as $\mu$ crosses from stable ($\mu<\mu_c$) to unstable ($\mu>\mu_c$).  At $\mu=\mu_c$, the real part of the eigenvalues lie on the imaginary axis.  In Fig.~\ref{fig:patterns}(a)  A pair of eigenvalues (${\bf x}$) cross into the right-half plane with the imaginary part of the eigenvalues being $\pm i \omega_0$.  This is the canonical Hopf bifurcation where temporal oscillations are observed at the onset of the instability.   In Fig.~\ref{fig:patterns}(b), the dominant instability is given by eigenvalue (${\bf x}$) crossing into the right-half plane with $\omega_0=0$, resulting an a growth mode with with no temporal oscillations.  Fig.~\ref{fig:patterns}(c) and (d) give the two potential corresponding spatio-temporal growth modes, in (c) the unstable wavenumbers are centered around a band at $k_0$ so that a modulational instability is observed.  For panel (d), the spatial instability grows with wavenumber $k_0=0$, thus generating the growth of a spatially homogenous plane wave.  There are four possible combinations of instabilities:  (i) $k_0=0$ and $\omega_0$=0, (ii) $k_0\neq 0$ and $\omega_0=0$, (iii) $k_0=0$ and $\omega_0\neq 0$, and (iv)  $k_0\neq$ and $\omega_0\neq 0$.

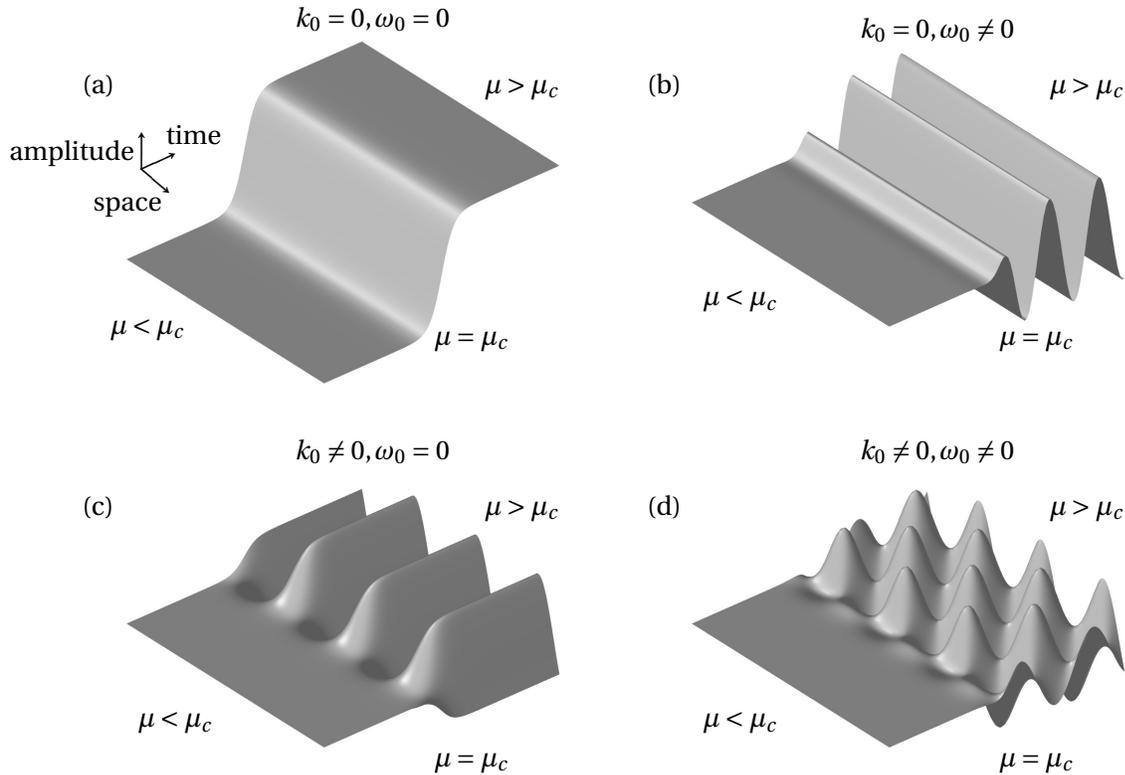
\begin{figure}[t]
\begin{overpic}[width=0.45\textwidth]{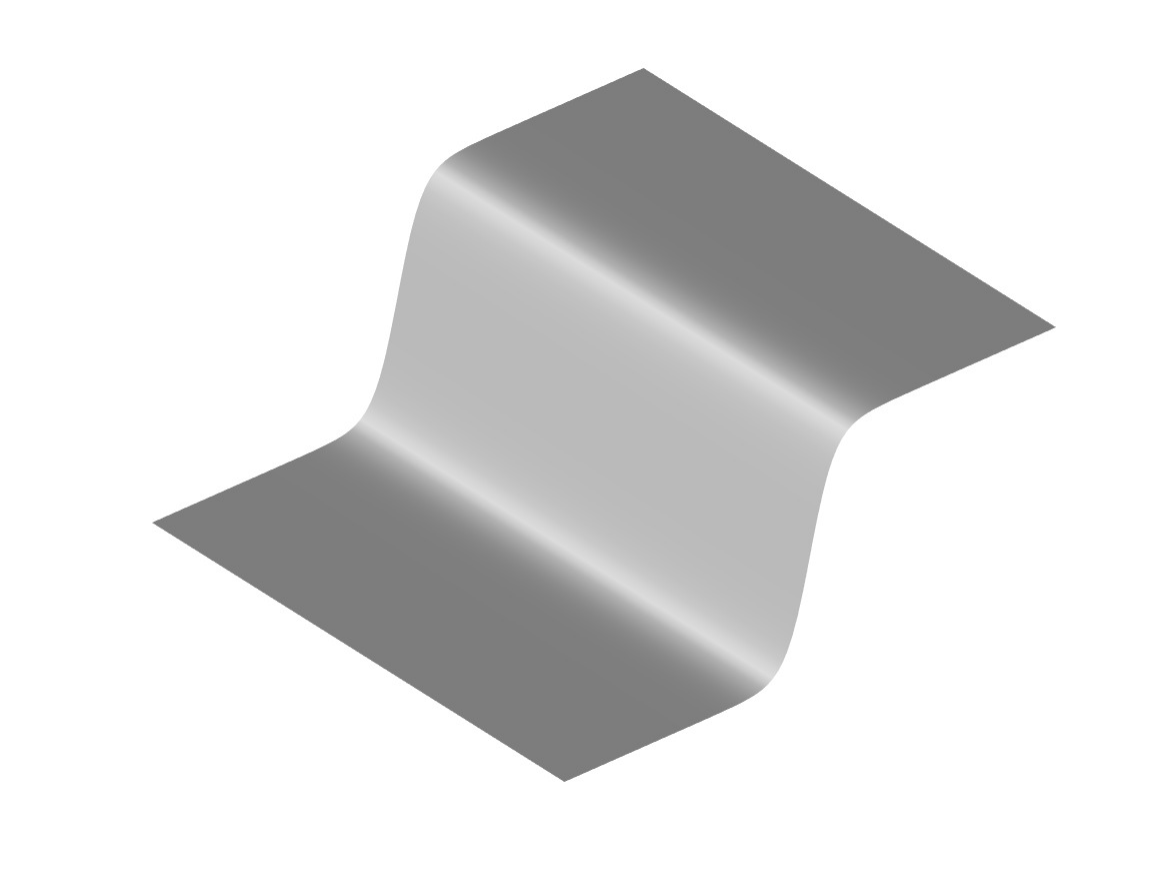}
\put(5,60){(a)}
\put(68,15){$\mu=\mu_c$}
\put(10,17){$\mu<\mu_c$}
\put(77,60){$\mu>\mu_c$}
\put(43,72){$k_0=0, \omega_0=0$}
\put(15,45){\rotatebox{25}{$\longrightarrow$}}
\put(14.3,45.7){\rotatebox{-40}{$\longrightarrow$}}
\put(15.0,45.9){\rotatebox{90}{$\longrightarrow$}}
\put(20,51){time}
\put(7,39){space}
\put(-8,48){amplitude}
\end{overpic}
\begin{overpic}[width=0.45\textwidth]{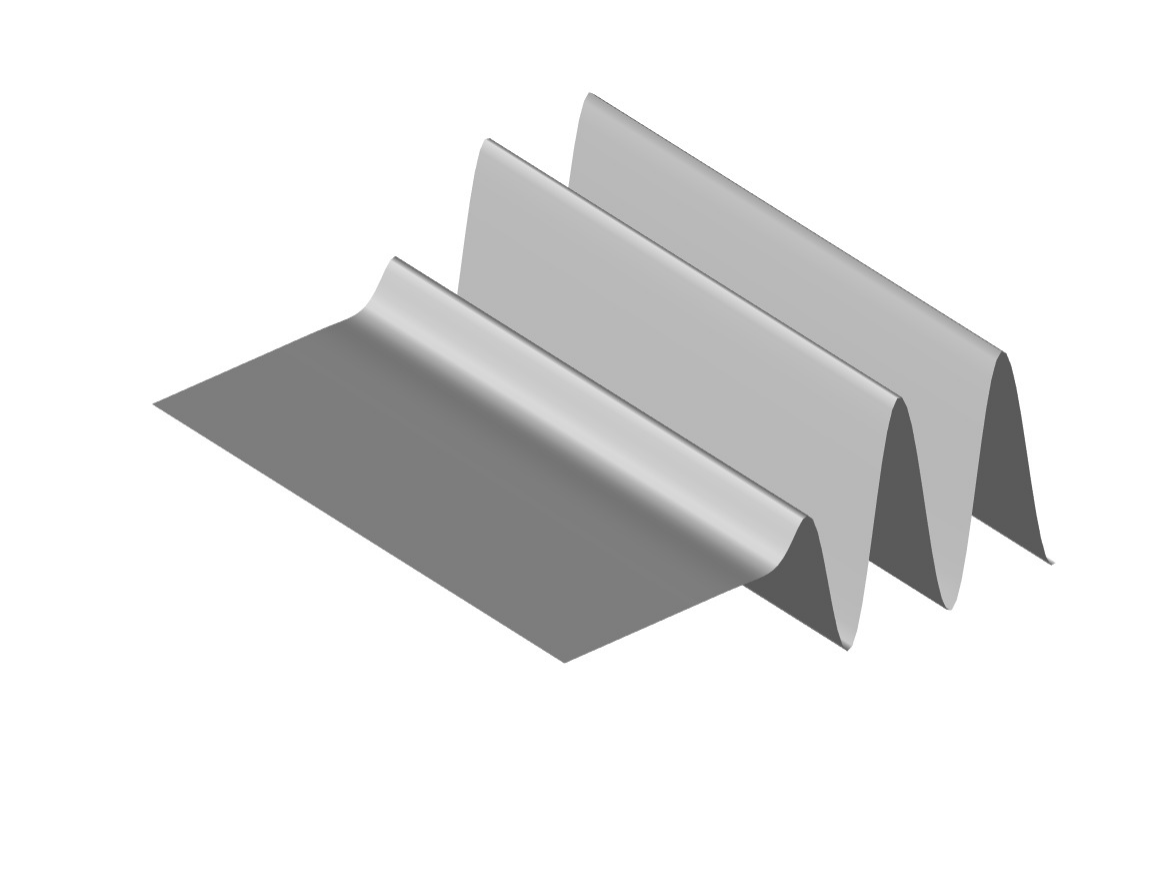}
\put(68,15){$\mu=\mu_c$}
\put(15,22){$\mu<\mu_c$}
\put(77,60){$\mu>\mu_c$}
\put(43,70){$k_0=0, \omega_0\neq0$}
\put(5,60){(b)}
\end{overpic}
\begin{overpic}[width=0.45\textwidth]{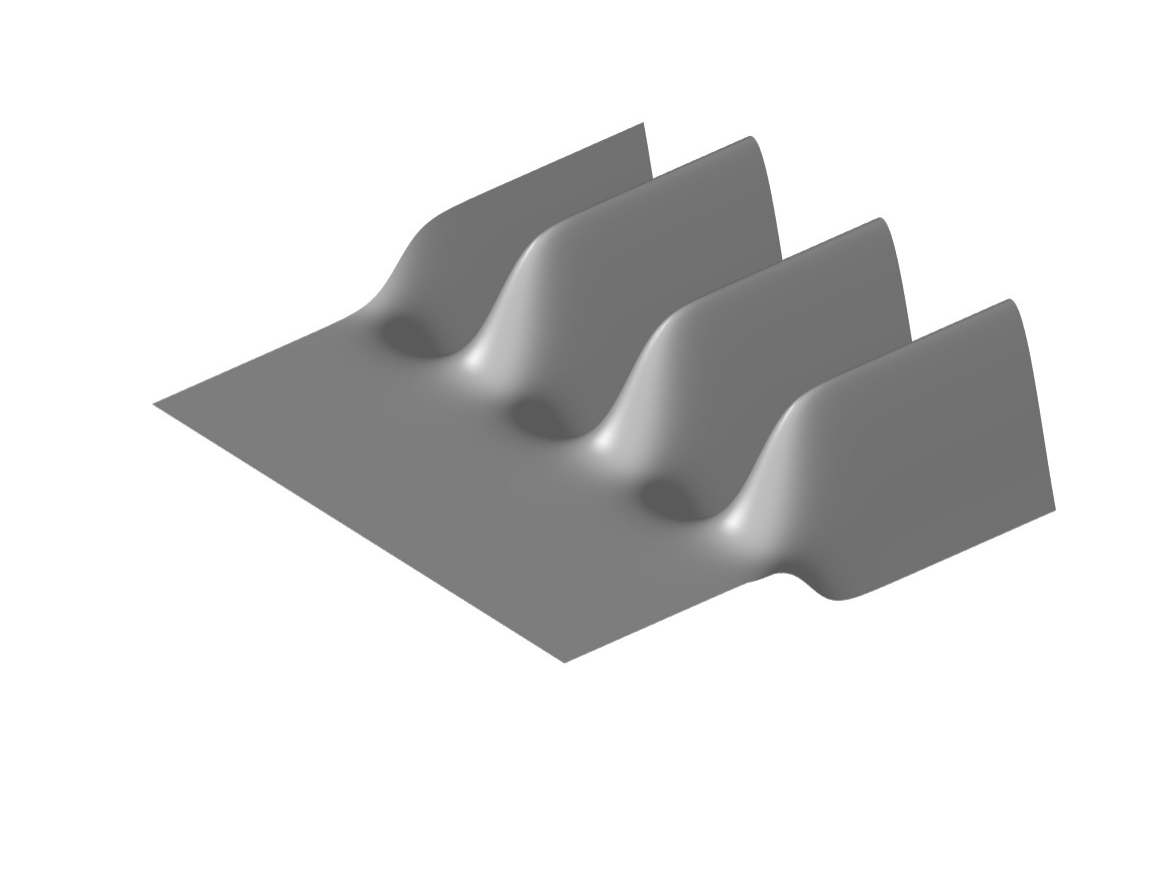}
\put(68,15){$\mu=\mu_c$}
\put(15,22){$\mu<\mu_c$}
\put(77,60){$\mu>\mu_c$}
\put(43,70){$k_0\neq0, \omega_0=0$}
\put(5,60){(c)}
\end{overpic}
\begin{overpic}[width=0.45\textwidth]{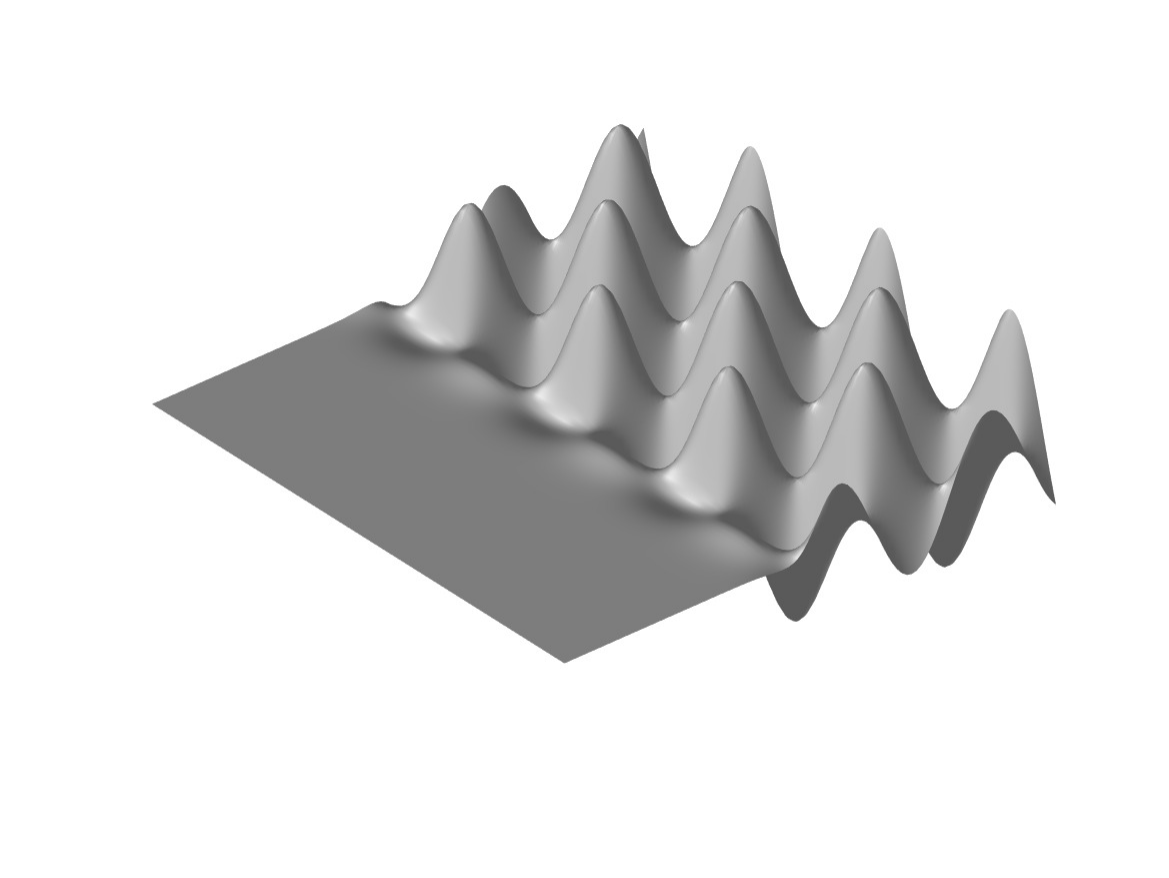}
\put(68,15){$\mu=\mu_c$}
\put(15,22){$\mu<\mu_c$}
\put(77,60){$\mu>\mu_c$}
\put(43,70){$k_0\neq0, \omega_0\neq0$}
\put(5,60){(d)}
\end{overpic}
\caption{Instability types for spatio-temporal systems as dictated by the most unstable wavenumber $k_0$ and most unstable temporal frequency $\omega_0$ at the instability points $\mu=\mu_c$.  These include four standard varieties:  (a) $k_0=0$ and $\omega_0$=0, (b) $k_0\neq 0$ and $\omega_0=0$, (c) $k_0=0$ and $\omega_0\neq 0$, and (d)  $k_0\neq$ and $\omega_0\neq 0$.  These panels show the canonical spatio-temporal patterns that form the basis of the overall pattern forming instabilities.  }
 \label{fig:patterns2}
\end{figure}

Figure~\ref{fig:patterns2} demonstrates the various spatio-temporal patterns that can emerge from linear instabilities.  At the bifurcation point $\mu=\mu_c$, the computation of the dominant unstable wavenumber and frequency, $k_0$ and $\omega_0$ respectively, determine the spatio-temporal patterns as the bifurcation parameter goes from $\mu<\mu_c$ to $\mu>\mu_c$.  The four instabilities shown in the subpanels of Fig.~\ref{fig:patterns2}  include 
(a) $k_0=0$ and $\omega_0=0$, (b) $k_0\neq 0$ and $\omega_0=0$, (c) $k_0=0$ and $\omega_0\neq 0$, and (d)  $k_0\neq$ and $\omega_0\neq 0$.
Cross and Hohenberg~\cite{cross1993pattern} classify these instabilities into three types:  (i) Type I$_{s}$:  stationary periodic ($k_0\neq 0$ and $\omega_0=0$ in Fig.~\ref{fig:patterns2}(b)), (ii) Type II$_{o}$:  oscillatory periodic ($k_0\neq$ and $\omega_0\neq 0$ in Fig.~\ref{fig:patterns2}(d)), and (iii) Type III$_{o}$ oscillatory uniform ($k_0=0$ and $\omega_0\neq 0$ in Fig.~\ref{fig:patterns2}(c)).  The case labeled Type III$_{s}$ ( $k_0=0$ and $\omega_0=0$ in 
Fig.~\ref{fig:patterns2}(a)) does not involve pattern formation in an essential way and is typically not considered part of the pattern forming analysis~\cite{cross1993pattern}.

\newpage
\section*{Lecture 22:  Linear Stability and Order Parameters}

Linear stability analysis gives an informative and interpretable diagnostic tool for understanding the pattern forming instabilities that arise in PDEs.  However, it is fundamentally limited.  Specifically, the onset of instability is predicted to grow exponentially at the rate of the largest $\Re \{ \lambda \}$.   Exponential growth is unbounded and leads to infinite amplitude solutions.  Observations of spatio-temporal systems, in experiment and computation, show that growth eventually saturates due to nonlinearities neglected in the perturbative theory of the linearization process.  What is then equally important is to understand the nature of the saturating nonlinearities which lead to new solutions $u_0(x,t)$ of the governing PDE system.

To amend the linear stability analysis and extract significant more information from the PDE evolution, a multiple-scale analysis can be used to enrich the analysis.   Again consider the  nonlinear PDE
\begin{equation}
  u_t = N(u, u_x, u_{xx}, \cdots, x, t, \mu) + \epsilon F(u,x,t)
  \label{eq:PDE_ms}
\end{equation}  
which now includes a potential (small) forcing term $F(\cdot)$ and 
which characterizes the evolution of a spatio-temporal field $u(x,t)$ and where $\mu$ is a parameter that induces a bifurcation at some parameter value $\mu=\mu_c$.   Suppose there exists a solution to this PDE which takes the form
\begin{equation}
    u(x,t)=u_0(x,t,\boldsymbol{\Theta})
\end{equation}
where $\boldsymbol{\Theta}$ is a vector of the constants that parametrize the solution.  The dependence on the parameters $\boldsymbol{\Theta}$ is made explicit since they will be leveraged in a multiple scale analysis.

A multiple scale analysis can be performed by introducing the slow time and space variables
\begin{subeqnarray}
  &&  \tau=\epsilon_1 t \\
  && \xi = \epsilon_2 x 
\end{subeqnarray}
where $\epsilon_1 \ll1 $ and $\epsilon_2 \ll1 $ are related to the small parameter $\epsilon$ used in a perturbation analysis.  It is typically  determined from dominant balance arguments, i.e. they are specific to the PDE system under consideration.  Importantly, the leading order solution $u_0(x,t,\boldsymbol{\Theta})$ is now modified by explicitly acknowledging that the parameters in   $\boldsymbol{\Theta}$ are no longer constant but rather a function of the slow variables
\begin{equation}
     \boldsymbol{\Theta} =  \boldsymbol{\Theta} (\tau,\xi) .
     \label{eq:slow_theta}
\end{equation}
This follows the standard theory of a multiple scale analysis. The governing PDE also requires modification in order to allow for slow time and space variation.  Thus derivatives in (\ref{eq:PDE_ms} produce the following
\begin{subeqnarray}
  &&  \frac{\partial}{\partial t}  \,\, \rightarrow \,\, \frac{\partial}{\partial t} + \epsilon_1 \frac{\partial}{\partial \tau}   \\
  &&  \frac{\partial}{\partial x}  \,\, \rightarrow \,\, \frac{\partial}{\partial x} + \epsilon_2 \frac{\partial}{\partial \xi}  \\
  &&  \frac{\partial^2}{\partial x^2}  \,\, \rightarrow \,\, \frac{\partial^2}{\partial x^2} + \epsilon_2 \frac{\partial^2}{\partial \xi \partial x} +
  \epsilon_2^2  \frac{\partial^2}{\partial \xi^2}   \\  
  && \hspace*{.5in} \vdots
\end{subeqnarray}
which completes the change of variables from $(x,t)$ to $(x,t,\xi,\tau)$.

We can now reconsider the linear stability analysis of the last section in order to understand the nature of the saturating instabilities produced by the full governing PDE (\ref{eq:PDE_ms}).  In this case, the parameter $\epsilon$ will be defined as follows
\begin{equation}
    \epsilon_3 =  \mu-\mu_c \ll 1
\end{equation}
where $\epsilon_3$ is a small parameter related to $\epsilon$ via a dominant balance of the specific problem considered.  
The perturbation analysis considers perturbations from the solution $u_0(x,t,\boldsymbol{\Theta})$ so that
\begin{equation}
  u(x,t,\xi,\tau) =  u_0(x,t,\boldsymbol{\Theta}(\tau,\xi)) + \epsilon \tilde{u} (x,t,\xi,\tau) .
\end{equation}
Inserting this modified perturbation expansion into (\ref{eq:PDE_ms}) yields the linearized evolution equation
\begin{equation}
   {\tilde{u}}_t = L(u_0, {u_0}_x, {u_0}_{xx}, \cdots,x,t,\boldsymbol{\Theta},\mu) \tilde{u} + f(u_0,{u_0}_\tau,{u_0}_\xi, \cdots ,x,t,\xi,\,\tau) 
   \label{eq:lin_slow}
\end{equation}
where the forcing $f(\cdot)$ accounts for both the perturbing force $F(\cdot)$ in (\ref{eq:PDE_ms}), but also the slow scale variations in the leading order solution through (\ref{eq:slow_theta}).  This leads to the consideration of the canonical equation $Lu=f$ for which a great deal is known, including the fact that the forcing term $f$ must be orthogonal to the null space of $L^\dag$.  This will allow us to derive a number of equations that parametrized the dynamics beyond the simple linear stability analysis.  

The homogenous solution to (\ref{eq:lin_slow}) would take the form now of
\begin{equation}
    \tilde{u}(x,t,\xi,\tau) = A(\xi,\tau)  \exp \left[  \lambda (\tau) t +  i  k(\tau) x  \right]  .
\end{equation}
where the constant amplitude $A(\xi,\tau)$ is now dependent upon the slow parameters, and the leading order wavenumber $k$ and temporal evolution $\lambda$ also now can be shifted in slow time.  Application of the Fredholm-Alternative theorem can give a number of important results, including equations of slow dynamics of the form~\cite{benney1967propagation,Newell,kodama,cross1993pattern}
\begin{subeqnarray}
  &&  A_\tau = N_1 (A,A_\xi, A_{\xi\xi}, \cdots , \tau, \xi) \\
  &&  \boldsymbol{\Theta}_\tau = N_2 (  \boldsymbol{\Theta},  \boldsymbol{\Theta}_\xi ,  \boldsymbol{\Theta}_{\xi\xi}, \cdots , \tau, \xi) \\
  &&  \boldsymbol{\Theta}_\tau = N_3 (  \boldsymbol{\Theta}) \\
  && \lambda_\tau = N_4 (\lambda)
  \label{eq:modulation}
\end{subeqnarray}
where $N_j(\cdot)$ is determined by application of the Fredholm-Alternative.  The slow spatio-temporal modulation of the amplitude $A(\cdot)$ equation or parameter dependency $  \boldsymbol{\Theta}$ in (\ref{eq:modulation})(a) and (b) is often called the {\em order parameter} or {\em envelope equation}.  These equations determine how the linear instability is saturated through a dominant balance nonlinearity.  The evolution given by $ \boldsymbol{\Theta}_\tau$ and/or $\lambda_\tau$ determines how the eigenvalues, or solution parametrization, move slowly in its position in time due to the forcing.  Examples of these different scenarios will be given in what follows.

\subsection*{Stability of Solitons of the Nonlinear Schr\"odinger Equation}

An insightful example of the application of the multiple scales method applied to perturbation theory comes from the nonlinear Schr\"odinger equation~\cite{weinstein,kodama,kaup,elgin1993perturbations}.  This equation has already been considered in context of its modulational instability.  Here, its soliton solutions are instead the focal point of the analysis.  Additionally, the effects of perturbations on the governing equations are also introduced.  The governing evolution equation is given by
\begin{equation}
   i u_t  + \frac{1}{2} u_{xx} + |u|^2 u = \epsilon F (u, u_t, u_x, \cdots)
   \label{eq:nls_solitons}
\end{equation}
where $u(x,t)$ is a complex amplitude and $F$ is a small perturbative forcing to the system since $\epsilon \ll 1$.  The unforced system (with $\epsilon=0$) admits localized soliton solutions of the form
\begin{equation}
     u(x,t) = \eta \sech \left[ \eta (x-x_0) \right]  
     \exp \left[  i \xi (x-x_0)  +  i (\phi-\phi_0)   \right]  = u_0 \exp (i\psi)
\end{equation}
where
\begin{subeqnarray}
  && \frac{dx_0}{dt} = \xi \\
  && \frac{d\phi}{dt} = \frac{1}{2} \left( \xi^2 + \eta^2 \right) .
\end{subeqnarray}
The solution thus admits four parameters which characterize the solution $\boldsymbol{\Theta}^T = [ \eta \,\,\,  \xi \,\,\, x_0 \,\,\, \phi_0 ]$.
Each of these parameters will be assumed to vary on the slow time $\tau=\epsilon t$.

A perturbation analysis expands about the soliton solution in the following form
\begin{equation}
     u(x,t) = \left(  u_0 + \epsilon \tilde{u} + \cdots \right) \exp(i \psi )
     \label{eq:NLS_expansion}
\end{equation}
where  for algebraic convenience the leading order phase factor $i\psi$ has been factored out in the expansion.  Inserting (\ref{eq:NLS_expansion}) into
(\ref{eq:nls_solitons}) yields the following linearized evolution at $O(\epsilon)$
\begin{equation}
   i \tilde{u}_t + \frac{\eta^2}{2} \tilde{u}_{\zeta\zeta} - \frac{\eta^2}{2} \tilde{u}  + 2 u_0^2 \tilde{u}  + u_0^2 \tilde{u}^* = \hat{F}
\end{equation}
where
\begin{equation}
  \hat{F}
   = \exp(-i\psi)   F ( u_0 \exp(i\psi) ) - i {u_0}_\tau - \xi_\tau (x-x_0)  u_0  + (\xi {x_0}_\tau + {\phi_0}_\tau )  u_0
   \label{eq:rhs_nls}
\end{equation}
where $\zeta=\eta(x-x_0)$.

\begin{figure}[t]
\begin{overpic}[width=0.8\textwidth]{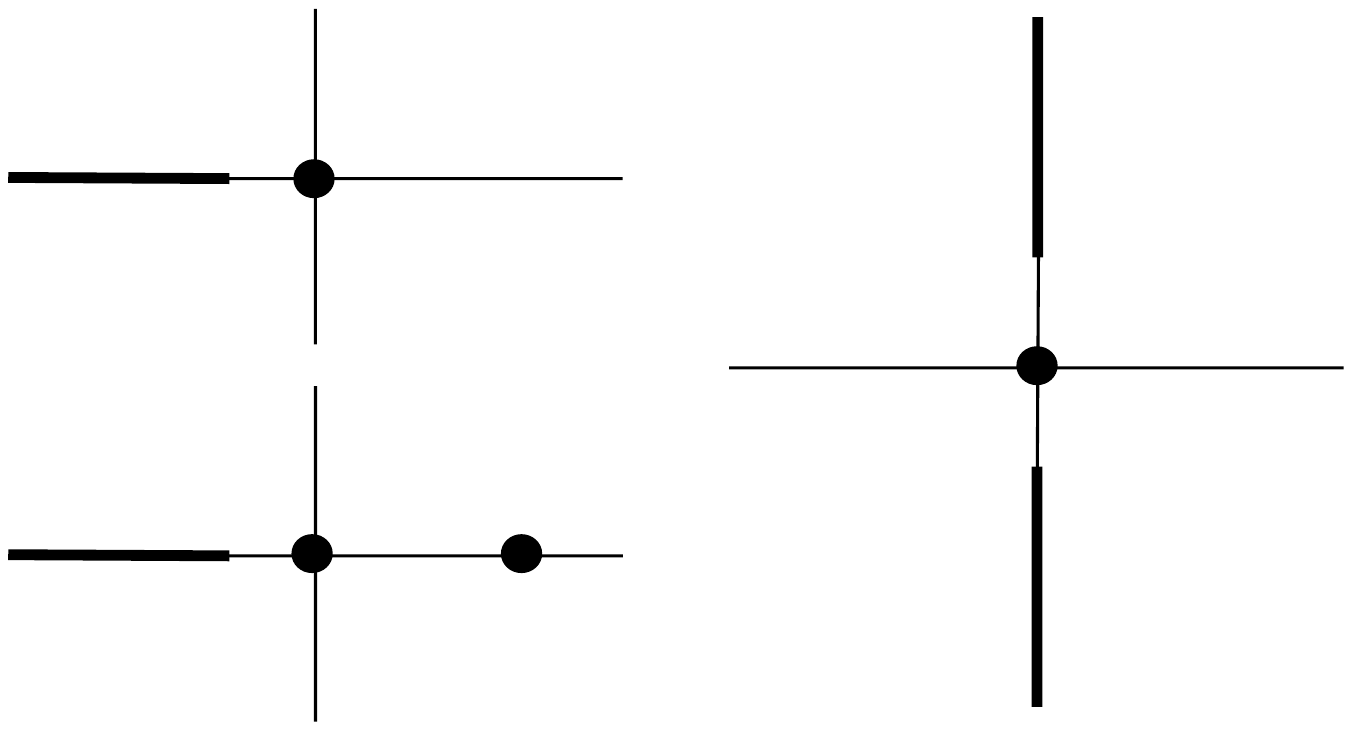}
\put(3,60){(a) \,\,\, $L_-$}
\put(3,35){(b) \,\,\, $L_+$}
\put(55,60){(c) \,\,\, ${\bf L}, {\bf L}^\dag$}
\put(41,24){$3$}
\put(30,24){$0$}
\put(30,47){$0$}
\put(22,24){$-1$}
\put(22,47){$-1$}
\put(75,36){$0$}
\put(75,45){$1$}
\put(75,32){$-1$}
\put(47,24){$\Re \{\lambda\}$}
\put(31,36){$\Im \{\lambda\}$}
\end{overpic}
\vspace*{-.7in}
\caption{Spectra of nonlinear Schr\"odinger operators $L_+$, $L_-$, ${\bf L}$ and ${\bf L}^\dag$.  The self-adjoint $L_-$ operator has a single zero eigenvalue and a continuous spectra starting at $\lambda=-1$.  The self-adjoint $L_+$ operator has two discrete eigenvalues at $\lambda=0$ and $\lambda=3$ and a continuous spectra starting at $\lambda=-1$.  The  ${\bf L}$ and ${\bf L}^\dag$ operators have a continuous spectral along the imaginary axis starting at $\lambda = \pm i$ with a degenerate eigenvalue of multiplicity four at the origin. }
 \label{fig:nls_spec}
\end{figure}

The perturbative field is decomposed  into its real and imaginary parts $\tilde{u}=R + iI$.  Defining the vector
 ${\bf w} = [R \,\,\, I]^T$ gives the linearized system
 \begin{equation}
  {\bf w}_t + {\bf L} {\bf w} = {\bf G}
\end{equation}
where 
\begin{subeqnarray}
 && {\bf L} = \frac{\eta^2}{2}   \left[  \begin{array}{cc} 0 & L_- \\ -L_+ & 0 \end{array}  \right] \\
 && {\bf G} = \left[  \begin{array}{c} \Im \{ \hat{F} \}  \\  \Re \{ \hat{F} \} \end{array}  \right]
\end{subeqnarray}
and the linear operators $L_\pm$ are given by
\begin{subeqnarray}
  && L_- = \partial^2_\zeta + 2 \sech^2  \zeta -1 \\
    && L_+ = \partial^2_\zeta + 6 \sech^2  \zeta -1 .
\end{subeqnarray}
Both $L_\pm$ are self-adjoint operators~\cite{weinstein,kodama,kaup,elgin1993perturbations}.  As shown in Fig.~\ref{fig:nls_spec}, the operator $L_+$ has a continuous spectrum starting at $\lambda=-1$ and two point spectra ($\lambda=0$ and $\lambda=3$) given by
\begin{subeqnarray}
  && L_+ (\sech \zeta \tanh \zeta ) = 0 \\
  && L_+ (\sech^2 \zeta ) = -3 \sech^2 \zeta .
\end{subeqnarray}
The operator $L_-$ has a continuous spectra starting at $\lambda=-1$ and a point spectra ($\lambda=0$) given by
\begin{equation}
  L_- (\sech \zeta) = 0  .
\end{equation}
The continuous spectra for the operators $L_\pm$ is found by considering when $x\rightarrow \pm \infty$ where each operator becomes $L_\pm \rightarrow \partial^2_\zeta  -1$.  The continous spectra is found by considering the solution of the of the form $u=\exp(ikx)$ in $L_\pm u = \lambda_\pm u$.  This gives $k^2=-(\lambda_\pm +1)$ which enforces real values of the wavenumber for $k\leq 1$.  Thus all wavenumbers $k\leq 1$ along the real axis compose the real spectra.
Two more operations are of note, especially for constructing the null space required for computing the Fredholm-Alternative theorem. 
\begin{subeqnarray}
  && L_+ ( \zeta \sech \zeta \tanh \zeta  - \sech \zeta) = 2 \sech \zeta \\
  && L_+ (\zeta \sech \zeta ) = 2 \sech \zeta \tanh \zeta .
\end{subeqnarray}
and
\begin{subeqnarray}
  && L_- L_+ ( \zeta \sech \zeta \tanh \zeta  - \sech \zeta) = 0 \\
  && L_+ L_- (\zeta \sech \zeta ) =0 .
\end{subeqnarray}
With these operations, the null space of the adjoint operator
\begin{equation}
   {\bf L}^\dag = \frac{\eta^2}{2}   \left[  \begin{array}{cc} 0 & -L_+ \\ L_- & 0 \end{array}  \right]
\end{equation}
can be constructed.  Specifically, the null space of ${\bf L}^\dag {\bf v}=0$ is spanned by the vectors
\begin{subeqnarray}
  && \left[ \begin{array}{c} \sech\zeta \\ 0 \end{array} \right] \\
  && \left[ \begin{array}{c} 0 \\ \sech\zeta \tanh\zeta  \end{array} \right] .
\end{subeqnarray}
For this problem, one must also consider the generalized null space satisfying ${\bf L}^2{\bf v}=0$.  The generalized null space is spanned by the generalized eigenvectors
\begin{subeqnarray}
  && \left[ \begin{array}{c}  \zeta \sech\zeta \\ 0\end{array} \right] \\
  && \left[ \begin{array}{c} 0 \\  \zeta \sech \zeta \tanh \zeta  - \sech \zeta \end{array} \right] .
\end{subeqnarray}

By enforcing the solvability condition, i.e. the forcing ${\bf G}$ must be orthogonal to the null space and generalized null space vectors of the adjoint, one can arrive at the following evolution equations for the slow time dynamics of the parameters $\boldsymbol{\Theta}$:
\begin{subeqnarray}
  &&  \left\langle  \Im \{\hat{F} \} , \sech\zeta \right\rangle =0 \\
  &&  \left\langle  \Im \{\hat{F} \} , \zeta \sech\zeta \right\rangle =0 \\
  &&  \left\langle  \Re \{\hat{F} \} , \sech\zeta \tanh\zeta \right\rangle =0 \\
 &&  \left\langle  \Re \{\hat{F} \} ,\zeta \sech \zeta \tanh \zeta  - \sech \zeta \right\rangle =0 .
\end{subeqnarray}
In general, the right hand side forcing has to be orthogonal to the null space and generalized null space of the adjoint operator.
These conditions are enforced on the $\hat{F}$ derived in (\ref{eq:rhs_nls}).  This results in the slow evolution dynamics of the parameters
\begin{subeqnarray}
  && \frac{d\eta}{d\tau} = \frac{1}{2}  \left\langle    \sech\zeta, \Im  \{  \exp(-i\psi) F (u_0\exp(i\psi)) \}  \right\rangle  \\
  &&   \frac{dx_0}{d\tau} = \frac{1}{\eta^2}  \left\langle \zeta    \sech\zeta, \Im  \{  \exp(-i\psi) F (u_0\exp(i\psi)) \}  \right\rangle  \\
  &&   \frac{d\xi}{d\tau} = \eta  \left\langle    \sech\zeta \tanh \zeta, \Re  \{  \exp(-i\psi) F (u_0\exp(i\psi)) \}  \right\rangle  \\
 &&  \xi \frac{dx_0}{d\tau} + \frac{d\phi_0}{d\tau} = \frac{1}{2}  \left\langle   \zeta  \sech\zeta \tanh \zeta, \Re  \{  \exp(-i\psi) F (u_0\exp(i\psi)) \}  \right\rangle  .
 \label{eq:nls_perturb}
\end{subeqnarray}
 These are the slow evolution equations $\boldsymbol{\Theta}_\tau=N_3(\boldsymbol{\Theta})$ which are induced by the perturbation to the governing NLS equations through $F(\cdot)$.

 This is clearly not the full story as we have neglected the effects of perturbations on the continuous spectrum of ${\bf L}$.  However, such computations are beyond the scope of the current work and one is referred to Kapitula and Promislow~\cite{kapitula2013spectral} for a more detailed treatment of this problem.  What is of importance here is that the original eigenvalues of the linearized operator has four zero modes which are neutrally stable and a continuous spectra which lives on the imaginary axis for $\lambda\geq 1$ and $\lambda\leq -1$.  The slow evolution equations (\ref{eq:nls_perturb}) characterize how the four zero eigenvalues are perturbed under a forcing $F(\cdot)$.  As an example, we can consider the perturbation
 \begin{equation}
  F  = - i \gamma u
\end{equation}
where the parameter $\gamma$ determines a linear damping applied to the governing NLS.  Application of this perturbation in (\ref{eq:nls_perturb})
gives the slow evolution dynamics
\begin{subeqnarray}
  && \frac{d\eta}{d\tau} = -2\gamma \eta \\
  &&   \frac{dx_0}{d\tau} = 0 \\
  &&   \frac{d\xi}{d\tau} = 0 \\
 &&   \frac{d\phi_0}{d\tau} = 0 .
\end{subeqnarray}
 In this case, the amplitude parameter slowly decays so that $\eta=\eta_0 \exp(-2\gamma \tau)$, which would be expected from the addition of linear damping.  The equations (\ref{eq:nls_perturb}) have been critical for theoretically understanding the optical pulse propagation fiber optic communications systems and mode-locked lasers.

\newpage
\section*{Lecture 23:  Order Parameters and Dominant Balance}

Pattern forming systems are often modeled by considering bifurcations from the trivial solution.  Thus as $\mu$ crosses the bifurcation point $\mu_c$, considerations is given to what grows from the trivial state $u_0(x,t)=0$.  In this case, the linear stability analysis would give an exponentially growing solution which would be well-known to saturate due to nonlinear dynamics in the governing equations (\ref{eq:nonlinearPDE}).  In this case, the expansion would be constructed as follows
\begin{equation}
  u(x,t,\xi,\tau) =  0 + \epsilon \tilde{u} (x,t,\xi,\tau) .
  \label{eq:order}
\end{equation}
where the evolution of the perturbation gives the order parameter, or envelope equation, dynamics determined from solvability conditions~\cite{benney1967propagation,Newell,kodama,cross1993pattern}.  

\subsection*{Order Parameter in a Nonlinear, 4th-Order Diffusion Problem}

To demonstrate how the multiple-scale analysis can be used to derive envelope equations, consider the governing evolution equation with bifurcation parameter $\mu$:
\begin{equation}
   u_t + \frac{1}{4} \left( {\partial_x}^2 - 1  \right) u - \mu u -  u^3 + u^5 + 3\beta u u_x^2 + (\beta+1) u^2 u_{xx} = 0
\end{equation}
which models the slowly-varying electric field envelop in an optical transmission line which uses parametric amplification~\cite{Kutz,Kutz1,Kutz2}.   The governing PDE has a trivial solution $u_0=0$ which can be linearized around to determine the onset of instability from the trivial state.  
One can consider the that the perturbative solution (\ref{eq:order}) takes the form $\tilde{u}(x,t)=\exp(\lambda t + i k x)$.  This gives the dispersion relation
\begin{equation}
   \lambda+\frac{1}{4} k^4 +\frac{1}{2} k^2 + \left( \frac{1}{4} - \mu \right) =0
\end{equation}
where $\lambda=0$ is achieved when the bifurcation parameter $\mu_c=(k^2+1)^2/4$.  Thus for $\mu>1/4$, the $k=0$ wavenumber is unstable and grows exponentially, while for $\mu<1/4$, the zero solution is stable.
Thus we consider (\ref{eq:order}) with the slow space and time scales
\begin{subeqnarray}
  &&  \tau=\epsilon^2 t \\
  &&  \xi = \epsilon x .
\end{subeqnarray}
The trivial solution is found to be stable for $\mu<1/4$ and unstable once $\mu>1/4$.  At the bifurcation point, the $\mu$ is expanded in the following form
\begin{equation}
   \mu=  \mu_c + \epsilon^2 \rho + \cdots=   \frac{1}{4} + \epsilon^2 \rho + \cdots .
\end{equation}
Inserting these asymptotic expansion in to the governing evolution equations yields after application of the solvability conditions
\begin{equation}
   \tilde{u}_\tau  = \frac{1}{2} \tilde{u}_{\xi\xi} + \tilde{u}^3 - \rho \tilde{u} .
\end{equation}
This is the order parameter description of the instability.  Note that instead of simple exponential growth or decay on either side of $\mu_c=1/4$, this description characterizes the simplified, dominant balance physics that occurs to produce a pattern forming system.  In particular, the envelope equation has the localized solution $\tilde{u}(\tau,\xi) = \sqrt{2\rho} \sech \sqrt{2\rho}\xi$.  This shows that at the instability, the nonlinearity interacts with the spatial derivatives to produce a localized solution.

\subsection*{Optical Parametric Oscillator }

As a specific example, we will consider the evolution dynamics in the optical parametric oscillator (OPO) since it yields a number of dominant balance physics, or order parameter descriptions, depending upon the parameter regimes considered~\cite{hewitt2005dynamics}.
The physical effects which influence the propagation of electromagnetic
energy in an OPO crystal arise from the
interaction of diffraction, parametric coupling between the signal, idler and
pump fields, attenuation, and the external pumping (driving) of the second
harmonic field.   The equations governing the leading order behavior for a continuous
wave can be derived from Maxwell's equations via a high frequency
expansion~\cite{derv1,derv2}.   The three key asymptotic reductions for this long
wavelength expansion come from applying the slowly-varying envelope
approximation, the paraxial approximation and rotating-wave approximation in
succession.    We consider the mean-field model of
a degenerate optical parametric oscillator (OPO) for which the signal and
idler fields coalesce.   The dimensionless
signal ($U$) and pump ($V$) field envelopes at the
fundamental and second harmonic are governed by the
nondimensionalized coupled equations~\cite{Lugiato,Trillo1,Trillo2}:
\begin{subeqnarray}
  && U_t = \frac{i}{2} U_{xx} + VU^* - (1+i\Delta_1) U \\
  && V_t = \frac{i}{2} \rho V_{xx} - U^2 - (\alpha + i\Delta_2) V + S 
\label{eq:opo}
\end{subeqnarray}
where $\Delta_1$ and $\Delta_2$ are the cavity detuning
parameters, $\rho$ is the diffraction ratio between signal
and pump fields, $\alpha$ is the pump-to-signal loss ratio,
and $S$ represents the external pumping term.  

Below the onset of any pattern generating instabilities, it is well known 
that the stable uniform steady-state response of the OPO is given by:
\begin{subeqnarray}
  && U=0 \\
  && V= \frac{S}{\alpha + i \Delta_2} \, .
\label{eq:opo-steady}
\end{subeqnarray}
This solution to (\ref{eq:opo}) becomes unstable once
a critical amount of pumping $|S|$ is applied.  A linear stability
analysis gives the critical value of pumping strength to be
\begin{equation}
  S_c= (\alpha + i\Delta_2)(1+i \Delta_1) \, .
\label{eq:critical}
\end{equation}
Thus once $|S|>|S_c|$, the trivial solution for the signal
field is unstable and a non-trivial, spatial structure
may arise.   Our aim is to investigate asymptotic reductions for
which nontrivial spatial structures in the OPO system can
be found and analyzed.   Note that solutions of
(\ref{eq:opo}) can also be found by modulating the external pumping
field~\cite{hewitt}.

Near $S=S_c$
it is appropriate to utilize an order parameter description
for the onset of instability~\cite{cross1993pattern,Newell,Segel}.
Thus we define the slow scales
\begin{subeqnarray}
   && \tau=\epsilon^2 t \\
   && {\sc T}=\epsilon^4 t \\
   && \xi =\epsilon x \\
   && {\sc X}=\epsilon^2 x
\end{subeqnarray}
where we define $\epsilon^2=|S-S_c|\ll 1$.  The fact that
we have defined two slow scales for both time and space
allows for the balance of different physical effects. 

Expand about the steady-state trivial solution of (\ref{eq:opo-steady}) by defining
\begin{subeqnarray}
  && U= 0 + \epsilon u (\tau,T,\xi,{\sc X})  \\
  && V= \frac{S}{\alpha + i\Delta_2} + \epsilon^2 v (\tau,T,\xi,{\sc X}) 
\label{eq:expand}
\end{subeqnarray}
with
\begin{equation}
  S- S_c = \epsilon^2 C +\epsilon^3 C_1 + \cdots ,
\label{eq:S}
\end{equation}
where $C$ and $C_j$ are constants.
After some manipulation it can be found that substituting (\ref{eq:expand})
into (\ref{eq:opo}) results in
\begin{subeqnarray}
  && \hspace*{-.15in} (1+ i \Delta_1) (u-u^*) = \epsilon^2 \left[
        \frac{i}{2} u_{\xi\xi} - u_\tau
     +vu^* + \frac{C}{\alpha+i\Delta_2} u^* 
     \right] + \epsilon^4 \left[ \frac{i}{2} 
        u_{{\sc X}{\sc X}} - u_T \right]  \\
  && \hspace*{-.15in} (\alpha+i\Delta_2) v = -u^2 + \epsilon^2 \left[
            \frac{i}{2} \rho v_{\xi\xi} - v_\tau \right] +
           \epsilon^4 \left[ \frac{i}{2}\rho 
         v_{{\sc X}{\sc X}} - v_T \right]  \, .
\label{eq:exact}
\end{subeqnarray}
which, upon manipulation, yields:
\begin{subeqnarray*}
 &&  v= - \frac{u^2}{\alpha+i\Delta_2} + \frac{\epsilon^2}{\alpha+i\Delta_2}
      \left[ \frac{i}{2}\rho v_{\xi\xi} - v_\tau \right] + O(\epsilon^4) \\
 &&  u^* =  u - \frac{\epsilon^2}{1+i\Delta_1} 
         \left[ \frac{i}{2} u_{\xi\xi} - u_\tau + vu^*
            + \frac{C}{\alpha + i\Delta_2} u^* \right] + O(\epsilon^4)
\end{subeqnarray*}
By applying an iterative proceedure, it can be found that
\[
 v u^* = - \frac{1}{\alpha+i\Delta_2} |u|^2 u + 
    \frac{\epsilon^2}{(\alpha+i\Delta_2)^2}
    \left( u (u^2)_\tau - \frac{i}{2}\rho u (u^2)_{\xi\xi} \right)
    + O(\epsilon^4) \, . 
\]
Making use of the above expressions, the right hand forcing of
(\ref{eq:exact}a) can then be found recursively to give
\begin{eqnarray}
 && R= \epsilon^2 \left[ \frac{i}{2} u_{\xi\xi} - u_\tau -
         \frac{ |u|^2 u}{\alpha + i\Delta_2} + \frac{C}{\alpha+i\Delta_2} u^*
      \right]   \nonumber \\
 && \hspace*{.5in} + \epsilon^4 \left[ \frac{i}{2} u_{XX} - u_T
         +\frac{1}{(\alpha+i\Delta_2)^2} \left( 
           u (u^2)_\tau -\frac{i}{2} \rho u (u^2)_{\xi\xi} \right) \right]
         + O(\epsilon^6)
\label{eq:R}
\end{eqnarray}
where we neglect contributions higher than $O(\epsilon^6)$.

The right hand forcing of (\ref{eq:exact}a) must satisfy certain orthogonality
conditions.  Fredholm's alternative theorem states that $R$ must be orthogonal
to the null space of the adjoint operator associated with the leading order
behavior of (\ref{eq:exact}a)~\cite{Friedman}.  For this case, the solvability
condition is given by
\begin{equation}
  (1-i\Delta_1) R + (1+i\Delta_1)R^*=0 \, .
\label{eq:solvabilityOPO}
\end{equation}
Applying this solvability condition, using the expression for R in
(\ref{eq:R}), gives the relevant order parameter equation to $O(\epsilon^6)$.

The most general nontrivial onset of instability arising from 
(\ref{eq:solvabilityOPO}) produces the real cubic, diffusive amplitude 
(Ginzburg-Landau or Fisher-Kolmogorov) equation:
\[
  \phi_\tau - \phi_{\zeta\zeta} \pm \phi^3 \mp \gamma \phi = 0
\]
where $u=((\alpha^2 + \Delta_2^2)/(\Delta_1 \Delta_2 - \alpha))^{1/2} \phi$,
$\xi=(\Delta_1/2)^{1/2} \zeta$, and $\gamma=-|C|(1+\Delta_1^2)/S_c>0$.  The
$\pm$ sign corresponds to the sign of the parameter $C$ respectively and
$\Delta_1>0$ in order for the equation to remain well-posed.  We
assume that both $\Delta_1=O(1)$ and $\alpha-\Delta_1\Delta_2=O(1)$ so that no
resonant detuning is implied.  Both solitary wave ($C<0$), which are unstable,
and front solutions~\cite{KE} ($C>0$) can be found for this case.  In order
for the case of $\Delta_1<0$ to make sense, the higher order (fourth-order)
diffusive corrections such as given by the Swift-Hohenberg equation must be
retained.

Another self-consistent characteristic balance arises in the evolution
dynamics when the slow-scales $\xi$ and $\tau$ are neglected.  The slow scales
are pushed to $O(\epsilon^2)$ and are captured by $T$ and $X$.  Application of
the solvability conditions (\ref{eq:solvabilityOPO}) yields the quintic, diffuse
equation
\[
  \phi_\tau - \phi_{\zeta\zeta}  
     + a \phi + b \phi^3 - \phi^5 =0
\]
where we have taken the near-resonance condition
$\alpha-\Delta_1\Delta_2 = \epsilon^2 \kappa$ and
rescaled so that:  $\tau=T/(2(\alpha^2+\Delta_2^2))$,
$\zeta=x/\sqrt{\Delta_1(\alpha^2+\Delta_2^2)}$.  Further $\phi=u$
with $a=\kappa C - C^2 K$ and $b=C-\kappa+ C K$
where $K= [(\alpha^2-\Delta^2)(1-\Delta_1^2)-4\Delta_1\Delta_2]/
(1+\Delta_1^2)(\alpha^2+\Delta_2^2)$.  Again $\Delta_1>0$
since $\Delta_1<0$ gives an ill-posed equation.
This equation has been studied extensively by
both geometric and analytic methods (see~\cite{Kap2} and references therein).

Although the above equations are of interest in their own
right, they have been extensively studied in a variety of contexts.
In this manuscript, the more interesting balance arises when balancing
the quintic term in the regime where the detuning parameter $\Delta_1$ is
near resonance.  We now only retain the slow scales $\xi$ and $T$
and apply the solvability condition~(\ref{eq:solvabilityOPO}).
After manipulation and upon taking $\Delta_1=\epsilon^2 \kappa$
and $\alpha-\Delta_1\Delta_2=\epsilon^2 \beta$, the following governing
envelope equation is found:
\begin{equation}
  \phi_t + \frac{1}{4} \left( \partial_\zeta^2 - \omega \right)^2 \phi
     -\gamma \phi - \sigma \phi^3 + \phi^5 
      \pm 3  \phi (\phi_\zeta)^2
      \pm 2  \phi^2 \phi_{\zeta\zeta} =0   
\label{eq:KK}
\end{equation}
where the $\pm$ correspond to the sign of $\Delta_2$ and
\begin{subeqnarray*}
    && \omega = 2\kappa |\Delta_2| \\
    && \sigma = -2\beta  \\
    && \gamma = \omega^2/4 + 2C\beta + C^2.
\end{subeqnarray*}
Note that $\phi=u$, $t=T/(2|\Delta_2|^2)$ and $\zeta=\xi/\sqrt{|\Delta_2|}$.
An equation of the same form but with coefficients different than those in
(\ref{eq:KK}) has been derived previously by Kutz and Kath in context
of a soliton transmission line with periodic, parametric
amplification~\cite{Kutz,Kutz1,Kutz2}.  Equation~(\ref{eq:KK}) and its three
solutions types (solitary waves, fronts, and periodic wave-trains) are the
focal point of the remainder of this manuscript.

\newpage
\section*{Lecture 24:  Modal Analysis and Mode Coupling}

As shown in this chapter, spatio-temporal systems have a diverse set of mathematical strategies that can be used to characterize the underlying dynamics of the system.  Asymptotics and perturbation methods play a critical role in helping determine the dominant physical processes responsible for the manifested patterns of spatio-temporal activity.   In this section, the low-dimensional structure of the  spatio-temporal dynamics are exploited to understand how the dynamics can be dominated by modes and their linear and nonlinear coupling dynamics.  

The underlying theoretical construction of the mode coupling dynamics will once again be derived from slow spatio-temporal scales along with enforcing solvability conditions.   As a starting  point to understanding the canonical dynamics of mode coupling, consider the following formulation
\begin{equation}
   i u_t + L u = \epsilon F(u,x,t)
   \label{eq:mode_couple}
\end{equation}
where $L$ is a linear, differential operator, $F(u,x,t)$ specifies either a forcing function acting on the spatio-temporal sytems and/or the nonlinear dynamics of the PDE,  and $i$ is the standard complex number.   For $\epsilon=0$, the dispersion relation can be determined by letting the leading order solution take the form
\begin{equation}
   u_0 (x,t) =v(x) \exp(i\lambda t)
\end{equation}
which gives the eigenvalue problem
\begin{equation}
      Lv = \lambda v .
\end{equation}
For the moment, let us assume that $L$ is a self-adjoint operator so that we can determine 
\begin{subeqnarray}
  &&  v_n  \,\,\, \mbox{orthonormal eigenfunctions} \\
  &&  \lambda_n \,\,\, \mbox{eigenvalues}
\end{subeqnarray}
that satisfy $Lv_n = \lambda_n v_n$ and with the property that
\begin{equation}
    \langle v_n, v_m \rangle = \delta_{nm}
\end{equation}
where $\delta_{nm}$ is the usual Kronecker delta for an orthonormal basis set $v_n$. 

The leading order solutions are typically constructed by linear superposition which gives a solution form
\begin{equation}
   u(x,t) = \sum_{n=1}^{\infty} a_n v_n(x) \exp(i\lambda_n t) .
   \label{eq:ef_mc}
\end{equation}
This solution specifies the evolution dynamics in the modal basis (eigenfunctions) given by the $v_n(x)$.  Given an initial condition $u(x,0)=U(x)$, the coefficients are determined by setting $t=0$ into the above solution expression and taking the inner product with respect to $v_m(x)$.  The orthonormality of the eigenfunctions gives
\begin{equation}
      a_n = \langle U, v_n  \rangle .
\end{equation}
This uniquely specifies the solution to the leading order equations when $\epsilon=0$. 

To characterize mode-coupling, the constants $a_n$ are modified so that they depend upon slow time $\tau=\epsilon t$.  Thus $a=a(\tau)$ and the leading-order eigenfunction solution is then modified to 
\begin{equation}
   u_0 (x,t,\tau) = \sum_{n=1}^{\infty} a_n (\tau) v_n(x) \exp(i\lambda_n t) .
\end{equation}
A perturbation expansion $u=u_0 + \tilde{u}$ gives the linearized evolution equations
\begin{equation}
    i \tilde{u}_t + L \tilde{u} = -i {u_0}_\tau + F(u_0,x,t) .
\end{equation}
The right-hand side terms must be orthogonal to the null space of the adjoint operator $L^\dag$.  Since initially we are assuming that $L$ is self-adjoint, then this is equivalent to assuming that the forcing must be orthogonal to the null space of $L$, which are the solutions (\ref{eq:ef_mc}). 

Initially, let us assume that the forcing and/or nonlinear terms $F=0$.  Applying the Fredholm-Alternative theorem in this case gives
\begin{equation}
   \frac{d{a_n}}{d \tau} =0
\end{equation}
which implies that no mode-coupling occurs when the system is unperturbed through a forcing or through nonlinearity.  Indeed, for a linear system that is unperturbed, the initial condition determines the initial distribution of energy in each mode $v_n(x)$ through the coefficient $a_n$.  The distribution of energy never changes in time without a perturbation or without nonlinear interactions.  

Mode-coupling, however, is a dominant paradigm of interactions in many physical and engineering systems.  This requires an understanding of how energy can be transferred across spatial modes in a given system.  There are at least three ways this can occur:  (i) by driving the system with a forcing, (ii) through nonlinear interactions, and (iii) through interaction in a non-orthogonal basis.  Each of these paradigms will be explored in what follows.

\subsection*{Mode-coupling through forcing}

Forcing is one of the most common ways to engineer the interaction between modes in a system.  It is the basis of many ideas from atomic and optical physics, including the foundational concept behind the laser and stimulated emission.  In this scenario, we will consider the forcing signal to be active so that $F\neq 0$.  
Moreover, consider the case of a two-mode interaction where only the $n$th and $m$th modes are active so that the leading order solution is
\begin{equation}
   u(x,t) =  a_n(\tau) v_n(x) \exp(i\lambda_n t) + a_m(\tau) v_m(x) \exp(i\lambda_m t) 
      \label{eq:two_mode}
\end{equation}
where all other mode amplitudes $a_j=0$ with $j\neq n,m$.

Linearizing about this solution gives the following evolution equations at $O(\epsilon)$:
\begin{eqnarray}
 &&  i \frac{d{a_n}}{d \tau} = \langle F (u_0,x,t), v_n \rangle  \exp(i\lambda_n t)  \\
 &&  i \frac{d{a_m}}{d \tau} =  \langle F (u_0,x,t), v_m \rangle  \exp(i\lambda_m t) .
\end{eqnarray}
The form of $F(\cdot)$ then determines the evolution dynamics.  In many physically relevant scenarios, including quantum mechanics and waveguide optics, the forcing term takes the form
\begin{equation}
    F= V(x,t) u
\end{equation}
where $V(x,t)$ is, for example, a correction to a given quantum mechanical potential well or a spatio-temporal index of refraction modulation in an optical waveguide.  With this form of forcing, and with a two-mode expansion for the leading-order field $u_0(x,t)$, this gives
\begin{eqnarray}
 &&  i \frac{d{a_n}}{d \tau} = \alpha_{nn} a_n  + \alpha_{mn}  a_m \exp(i\Delta\lambda t)  \\
 && i  \frac{d{a_m}}{d \tau} =  \alpha_{mm}  a_m  + \alpha_{nm} a_n  \exp(-i\Delta\lambda t)  
\end{eqnarray}
where  $\Delta=\lambda_n-\lambda_m$ and $\alpha_{nm} = \langle V(x,t) v_n, v_m \rangle$.  This gives a two-by-two system of differential equations that describe the coupling dynamics between the two modes.  

The parameter $\Delta$ is called the detuning, and it measures the difference in the two eigenvalues.  For a large detuning $\Delta \gg 1$, the mode-coupling is weak since the {\em rotating wave approximation} (i.e. the terms with $\exp(\pm i \Delta\lambda t)$ with $\Delta \gg 1$), which zeroes out rapidly oscillating mean-zero terms, gives  
\begin{eqnarray}
 &&  i \frac{d{a_n}}{d \tau} = \alpha_{nn} a_n    \\
 &&  i \frac{d{a_m}}{d \tau} =  \alpha_{mm}  a_m  .
\end{eqnarray}
However, one can force the system at a given resonant frequency so that
\begin{equation}
   V(x,t)=V(x)\exp(i\Delta t) 
\end{equation}
which gives a resonant interaction and maximal coupling between the two eigenstates:
\begin{eqnarray}
 &&   i \frac{d{a_n}}{d \tau} = \alpha_{nn} a_n \exp(-i\Delta\lambda t) + \alpha_{mn}  a_m   \\
 &&   i \frac{d{a_m}}{d \tau} =  \alpha_{mm}  a_m \exp(i\Delta\lambda t)  + \alpha_{nm} a_n   .
\end{eqnarray}
This is done in practice in optical systems through the use of a grating profile in the index of refraction along the propagation direction, or through optical stimulation at a resonant energy difference between quantum states in a laser cavity.  Examples of both these scenarios will be given in what follows.  The rotating wave approximation then give the resonant coupling dynamics
\begin{eqnarray}
 &&   i \frac{d{a_n}}{d \tau} =  \alpha  a_m   \\
 &&   i \frac{d{a_m}}{d \tau} =   \alpha a_n   .
\end{eqnarray}
where $\alpha=\alpha_{mn}=\alpha_{nm}$.   This resonant interaction produces oscillatory dynamics where energy is periodically transferred back and forth between the $n$th and $m$th modes by driving the system at a resonant frequency of $\Delta \lambda$.  

\subsection*{Mode-coupling through nonlinearity}

Nonlinearity is an alternative to producing interactions between modes in a system.  In this case, superposition of solutions no longer holds and a mode-mixing dynamics is produced by nonlinear dynamics of the system.  In this case, the forcing signal is assumed to be a nonlinear function of the state so that $F(u,x,t)=N(u)$ where $N(u)$ specifies the form of the nonlinearity.   We again consider the nonlinearity to be perturbatively small so that for a two-mode interaction where only the $n$th and $m$th modes are active at leading order,  the perturtabive solution is
\begin{eqnarray}
 &&  i \frac{d{a_n}}{d \tau} = \langle N (u_0), v_n \rangle  \exp(i\lambda_n t)  \\
 &&  i \frac{d{a_m}}{d \tau} =  \langle N (u_0), v_m \rangle  \exp(i\lambda_m t) .
\end{eqnarray}
For a given nonlinear interaction function $N(u_0)=N\left(a_n(\tau) v_n(x) \exp(i\lambda_n t) + a_m(\tau) v_m(x) \exp(i\lambda_m t) \right)$, the projection of the dynamics onto each eigenfunction direction produces a nonlinear dynamical system coupling the modal amplitudes $a_n$ and $a_m$.   A specific form of nonlinearity, which will be illustrated in the examples, must be considered in order to derive the nonlinear coupling dynamics.

\subsection*{Mode-coupling through non-orthogonal modes}

In certain scenarios, orthogonality of the modes is not guaranteed, thus the inner product and orthonormal conditions used to derive the above mode-coupling equations no longer holds.  Indeed, the inner product between modes in this scenario is generically not zero, forcing modes to interact.  Alternatively, one could decompose the dynamics into a basis set using methods such as the {\em dynamic mode decompositon} (DMD) which does not require modes to be orthogonal~\cite{Kutz2016book}.  A DMD decomposition represents the data in the form
\begin{equation}
   u (x,t) = \sum_{n=1}^{\infty} a_n(t) (\tau) v_n(x) \exp(i\lambda_n t) .
\end{equation}
where the $v_n(x)$ are the DMD modes, the $\lambda_n$ are the DMD eigenvalues, and the $a_n(t)$ are the time-dependent loadings for the DMD modes and eigenvalues.  Inserting this into the governing evolution equations, even with $F=0$, and taking the inner product with respect to $v_n$ and $v_m$ gives a coupled set of differential equations.  Thus the non-orthogonality creates mode-mixing since each direction is not independent of the other.   In many fluid flows where a DMD analysis is used, non-orthogonal modes are present and mode-mixing can occur simply due to such modes.

\subsection*{Mode-coupling in quantum mechanics}

As an example of the mode-coupling dynamics that can occur in both linear and nonlinear systems,
consider the Schr\"odinger equation modeling quantum mechanics
\begin{equation}
  i \hbar \psi_t + \frac{\hbar^2}{2m} \nabla^2 \psi - V({\bf x}) \psi =0
\end{equation}
where $\psi ({\bf x},t)$ is the wavefunction and $V({\bf x})$ is the potential.   The eigenfunctions and eigenvalues of this quantum system can be used to express the solution of this problem
\begin{equation}
   \psi= \sum_{n=1}^N A_n \phi_n({\bf x}) \exp \left( \frac{i E_n t}{\hbar} \right) 
\end{equation}
where the eigenfunctions $\phi_n$ and eigenvalues $E_n$ are determined from
\begin{equation}
    \frac{\hbar^2}{2m} \nabla^2 \phi_n - [ V({\bf x}) + E_n ] \phi_n = 0. 
    \label{eq:phi}
\end{equation}
This is the generic solution technique for quantum systems.

In our specific example, we would like to consider driving the system and mode-coupling.  In this case, the governing equations are modified to 
\begin{equation}
  i \hbar \psi_t + \frac{\hbar^2}{2m} \nabla^2 \psi - V({\bf x}) \psi = \epsilon  \cos(\omega t)
  \label{eq:forced}
\end{equation}
where $\epsilon\ll1$ shows this to be a perturbation to the quantum system which is driven with frequency $\omega$.  The solution form is then modified to
\begin{equation}
   \psi= \sum_{n=1}^N A_n (t) \phi_n({\bf x}) \exp \left( \frac{i E_n t}{\hbar} \right) 
   \label{eq:sol}
\end{equation}
where each mode amplitude $A_n(t)$ is now allowed to vary in time.  Note that because of $\hbar$, the rotating waves happen on a very fast time scale because of the rotation frequency of each mode ${E_n}/{\hbar}$.  Thus $A_n(t)$ is a slowly varying amplitude of each mode.  At leading order, the $\phi_n$ still satisfy (\ref{eq:phi}).  The derivation below follows from optical waveguide theory, which is analogous to quantum mechanical systems~\cite{kutz1997nonlinear}.  It can also be modified to handle wavepackets if necessary~\cite{kutz2014solitons}.

Mode-coupling is accomplished within the rotating wave approximation.  Thus the solution form (\ref{eq:sol}) is inserted into (\ref{eq:forced}) and the rapidly rotating components are discarded.  We also consider coupling between the $j$th and $k$th mode.  In particular we define the quantity
\begin{equation}
  \Delta = \omega - (E_j - E_k)/\hbar
\end{equation}
which is known as the detuning parameter.  If $\Delta=0$, then the forcing $\omega$ is resonantly tuned to the transition between the  $j$th and $k$th state.

To determine the dynamics of the coupling between the $j$th and $k$th state, we define two new variables
\begin{subeqnarray}
   && u= A_j \exp (-i \Delta t/2) \\
   && v= A_k \exp(+i \Delta t/2)
\end{subeqnarray}
which gives
\begin{subeqnarray}
   && i \frac{du}{dt} + c v - \frac{\Delta}{2} u = 0  \\
   && i \frac{dv}{dt} + c u + \frac{\Delta}{2} v = 0
   \label{eq:linear}
\end{subeqnarray}
where the coupling constant $c$ is proportional to the overlap integral between the $j$th and $k$th state, i.e. $c\propto ((V({\bf x}) \phi_j, \phi_k)$.

The linear mode-coupling dynamics admits a simple solution.  Given $u(0)=U_0$ and $v(0)=V_0$, one can show~\cite{kutz1997nonlinear,boyd2003nonlinear}
\begin{subeqnarray}
  &&u=U_0 \left[  \cos\left( \sqrt{c^2 + \frac{\Delta^2}{4} } t \right) - \frac{i\Delta}{ 2 \sqrt{c^2 + {\Delta^2}/{4} } }
      \sin \left( \sqrt{c^2 + \frac{\Delta^2}{4} } t \right)   \right] \\
   && v= V_0  \left[  \frac{ic}{ \sqrt{c^2 + {\Delta^2}/{4} } }
      \sin \left( \sqrt{c^2 + \frac{\Delta^2}{4} } t \right)   \right] .
      \label{eq:linsolution}
\end{subeqnarray}
This is the general solution.  

Consider the special case when $c=0$, but $\omega$ is used to drive transitions.  This gives
\begin{subeqnarray}
  &&u=U_0 [ \cos (\Delta/2) t - i \sin(\Delta/2)t \\
   && v= i V_0   \sin (\Delta/2) t  
\end{subeqnarray}
which shows the driving frequency can be used to induce transitions.

If one includes a cubic nonlinearity, which typically models an intensity-dependent response, then the governing mode-coupling equations are modified to 
\begin{subeqnarray}
   && i \frac{du}{dt} + c v - \frac{\Delta}{2} u +  \epsilon ( c_{jj} |u|^2 + c_{jk} |v|^2 )u = 0  \\
   && i \frac{dv}{dt} + c u + \frac{\Delta}{2} v +  \epsilon ( c_{kj} |u|^2 + c_{kk} |v|^2 )v = 0
   \label{eq:nonlinear}
\end{subeqnarray}
where the nonlinear coupling coefficients $c_{jk} =(|\phi_j |^2, |\phi_k|^2 )$ and $\epsilon\ll 1$ mades the nonlinear contribution perturbatively small.  We can then look for asymptotic solutions
\begin{subeqnarray}
  &&u=u_0 + \epsilon u_1 + \cdots \\
   && v= v_0 + \epsilon v_1 + \cdots   
\end{subeqnarray}
The leading order solution for $u_0$ and $v_0$ is already given by (\ref{eq:linsolution}), so we just need to find corrections to this approximation from the perturbation theory.  One can solve for the higher-order solutions, but also let $U_0$ and $V_0$ depend on slow-time $\tau=\epsilon t$ in order to get the slow phase shifts induced by the nonlinearity.

To start this, let us rewrite the above system in matrix form with ${\bf x} = ( u \,\,\, v)^T$ so that
\begin{equation}
   \frac{d{\bf x}}{dt} = i {\bf A} {\bf x}  - i \epsilon {\bf B} {\bf x}
\end{equation}
with
\begin{subeqnarray}
  &&{\bf A} \left[ \begin{array}{cc}  -\Delta/2 & c \\ c & \Delta/2 \end{array}  \right]  \\
   && {\bf B} = \left[  \begin{array}{cc} c_{jj} |u|^2  & c_{jk} |v|^2  \\ c_{kj} |u|^2 & c_{kk} |v|^2 \end{array} \right].
\end{subeqnarray}
We now use a perturbation expansion
\begin{equation}
   {\bf x} = {\bf x}_0 + \epsilon {\bf x}_1 + \cdots
\end{equation}
which has the leading order problem
\begin{equation}
   \frac{d{\bf x}_0}{dt} = i {\bf A} {\bf x}_0.
\end{equation}
The solution can be found by inserting
\begin{equation}
  {\bf x}_0 = {\bf v} \exp(i\lambda t)
\end{equation}
which gives
\begin{equation}
  {\bf A}{\bf v}= \lambda {\bf v}.
\end{equation}
The eigenvalue problem gives $\lambda_\pm =\pm \sqrt{c^2+\Delta^2/4}$.  

The leading order solution to this problem is then
\begin{equation}
  {\bf x}_0 = a_1(\tau) \exp (i \sqrt{c^2+\Delta^2/4} t) + a_2(\tau)  \exp (-i \sqrt{c^2+\Delta^2/4} t)
\end{equation}

The $O(\epsilon)$ problem is now given as follows
\begin{equation}
   \frac{d{\bf x}_1}{dt} - i {\bf A} {\bf x}_1  =  -\frac{d {\bf x}_0}{d\tau} + i {\bf B}({\bf x}_0) {\bf x}_0.
\end{equation}
The right hand side of this equation must be orthogonal to the null space of the adjoint operator, which is almost self-adjoint aside from complex conjugation, thus the solutions $\exp (\pm i \sqrt{c^2+\Delta^2/4} t)$ must be orthogonal to the adjoint solutions given by $\exp (\mp i \sqrt{c^2+\Delta^2/4} t)$.  Applying this at next order gives
\begin{subeqnarray}
  &&  \frac{da_1}{d\tau} = ( c_{jj} |a_1|^2 + c_{jk} |a_2|^2 ) a_1 \\
   && \frac{da_2}{d\tau} = ( c_{kj} |a_1|^2 + c_{kk} |a_2|^2 ) a_2
\end{subeqnarray}
which has the general solution
\begin{subeqnarray}
  &&   a_1 = a_1(0) \exp [ -i ( c_{jj} |a_1(0) |^2 + c_{jk} |a_2 (0) |^2 ) t] \\
   &&  a_2 = a_2(0) \exp [ -i ( c_{kj} |a_1(0) |^2 + c_{kk} |a_2 (0) |^2 ) t].
\end{subeqnarray}
This shows the explicit phase shift induced by the nonlinearity.  We can then put this all together to get the leading-order solution which accounts for nonlinear phase shifts:
\begin{eqnarray}
  {\bf x}_0 = a_1(0) \exp \left[ -i \left( c_{jj} |a_1(0) |^2 + c_{jk} |a_2 (0) |^2  - \sqrt{c^2+\Delta^2/4} \right) t \right]   \nonumber \\
                   + a_2(0) \exp \left[ -i \left( c_{kj} |a_1(0) |^2 + c_{kk} |a_2 (0) |^2  + \sqrt{c^2+\Delta^2/4} \right) t \right]
\end{eqnarray}
This gives an approximation to a forced, nonlinear mode coupling dynamics to leading order.
Nonlinearity and more exotic waveguides can also be considered~\cite{mahmud2002bose}. 

\subsection*{Mode-coupling in optical waveguides}

Electrodynamics, specifically the evolution of the electric field through optical waveguides, provides another common application area where mode-coupling dominates the dynamics.  In this case, the governing equations are given by
\begin{equation}
  i u_t + \frac{1}{2}{u}_{xx} - n(x,t) u = \epsilon F(u,x,t)
\end{equation}
where $u(x,t)$ is the normalized and rescaled electric-field envelope and $n(x,t)$ is a similarly normalized index of refraction profile that forms the waveguide.  This governing equation is identical to the Schr\"odinger equation of quantum mechanics.  The perturbation $F(\cdot)$ will be used to create mode-coupling dynamics.

Mode-coupling will be demonstrated via two of the paradigms highlighted:  forcing and nonlinearity.  For forcing, we will modify the index of refraction as follows:
\begin{equation}
    n(x,t)=n_0(x) \left[ 1 + \epsilon \cos (\omega t) \right]
    \label{eq:mc_force}
\end{equation}
where $\omega$ determines the frequency of the forcing.   This is done in practice using a optical grating~\cite{kutz1997nonlinear,boyd2003nonlinear}.  Nonlinearity can also be used to create mode-coupling dynamics.  In this case
\begin{equation}
   F(u,x,t)= | u|^2 u
   \label{eq:mc_non}
\end{equation}
which arises from the intensity dependent index of refraction.  This is known as a Kerr nonlinearity.  Such an intensity dependent response is critical for the generation of optical solitons.   This concept has been used in laser dynamics where waveguide coupling is controlled with nonlinearity~\cite{proctor2005passive,zhang2015semiconductor}.

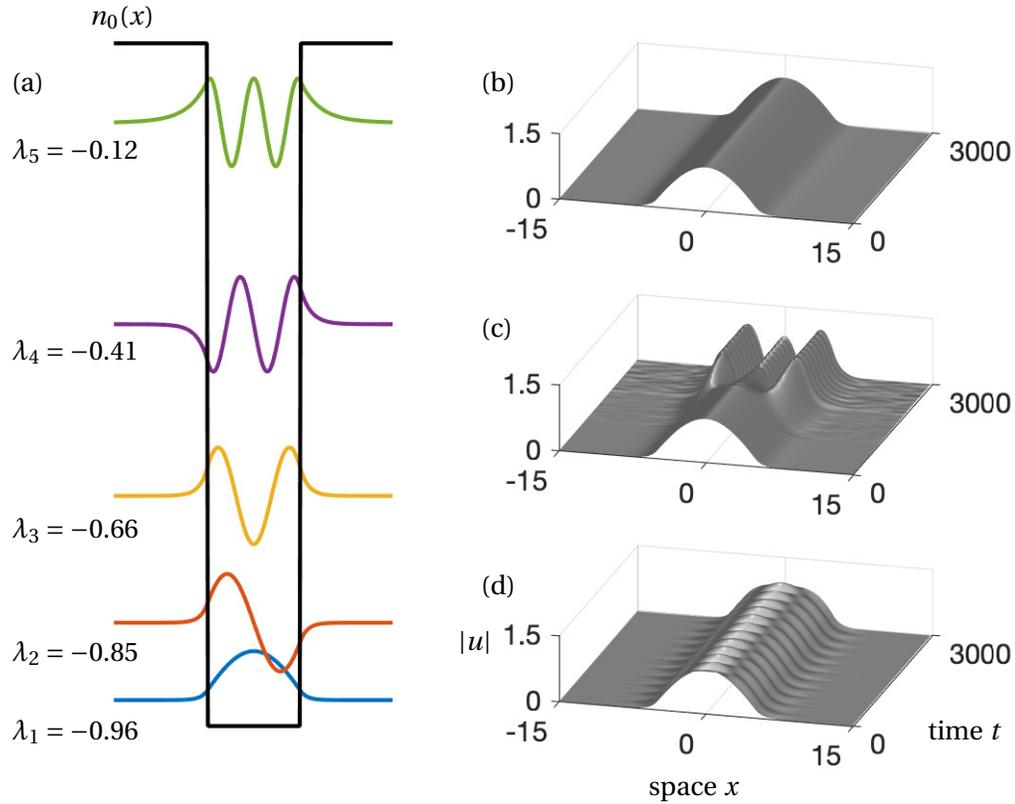
\begin{figure}[t]
\begin{overpic}[width=0.9\textwidth]{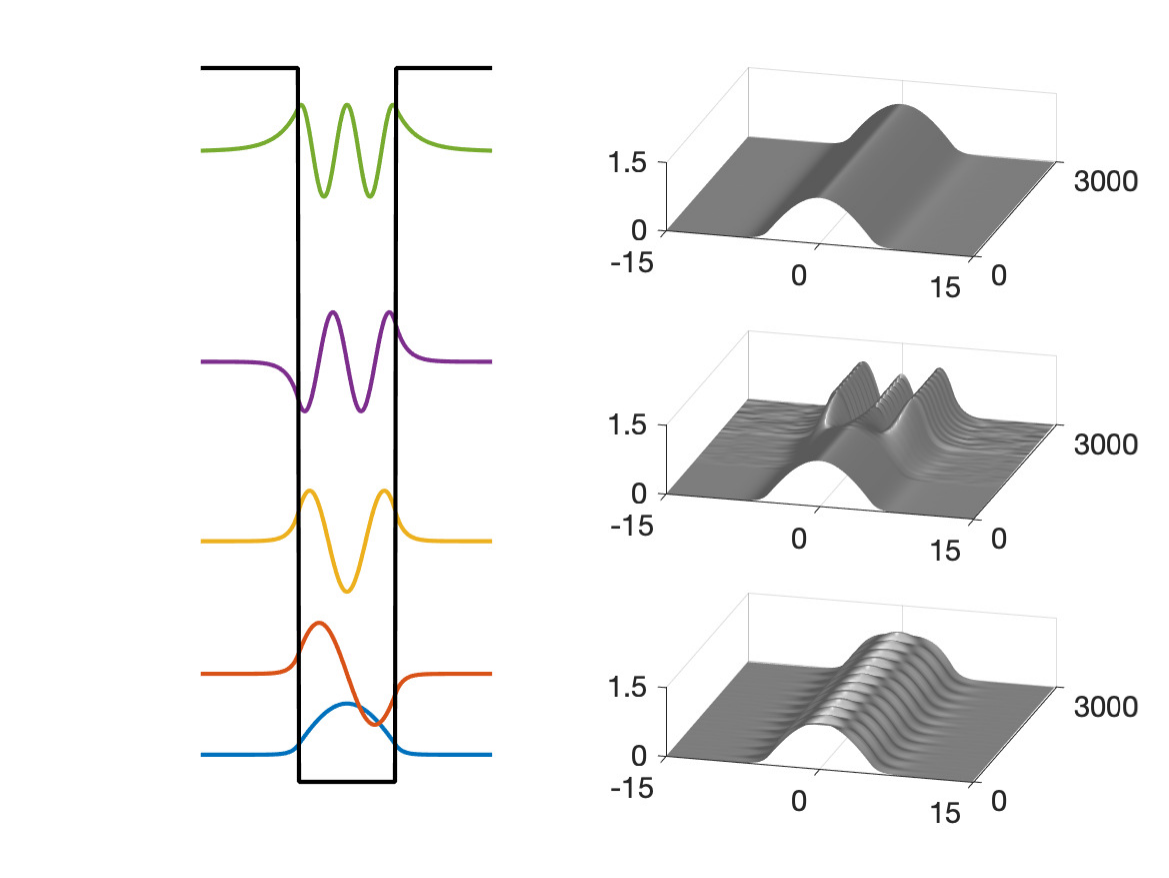}
\put(8,65){(a)}
\put(50,65){(b)}
\put(50,43){(c)}
\put(50,20){(d)}
\put(8,7){$\lambda_1=-0.96$}
\put(8,14){$\lambda_2=-0.85$}
\put(8,25){$\lambda_3=-0.66$}
\put(8,41){$\lambda_4=-0.41$}
\put(8,59){$\lambda_5=-0.12$}
\put(15,71){$n_0(x)$}
\put(48,15){$|u|$}
\put(65,2){space $x$}
\put(90,7){time $t$}
\end{overpic}
\vspace*{-.0in}
\caption{Mode coupling with eigenvalues (propagation constants):  $\lambda_n=-0.96, 0.85, 0.66, 0.41,$ and $ 0.12$.  (a) Index of refraction profile $n_0(x)$ and the mode (eigenfunction) profiles of the five discrete (bound) eigenvalues of the system.  (b)  Evolution dynamics for the linear system with no perturbation and initial conditions in the ground state $u(x,0)=v_1(x)$.  (b)  Evolution dynamics when a forcing is applied from $t\in [1000,2000]$ with $\omega=\lambda_1-\lambda_3$.  Note that nearly complete transfer of energy from mode one ($v_1(x)$) to mode three ($v_3(x)$).  (c) Mode-coupling dynamics induced by nonlinearity and no forcing.}
 \label{fig:mode_couple}
\end{figure}

For $\epsilon=0$, the leading order solution of the optical waveguide problem can be formulated as an eigenvalue problem where $u(x,t)=v(x) \exp(-i \lambda t)$ to that
\begin{equation}
  Lv_n= \frac{1}{2}{v_n}_{xx} - n_0(x) v_n = -\lambda_n v_n
\end{equation}
where $v_n$ are the eigenfunctions and $\lambda_n$ are the eigenvalues given an index profile $n_0(x)$.  
The eigenvalues are often denoted as $\beta_n$ and are called the propagation constants in optics.  The choice of a negative sign on the right hand side of the equation gives positive values of the index of refraction.  The index profile we consider takes the form
\begin{equation}
  n_0(x)= \left\{ \begin{array}{rl} -1 & |x|<5 \\ 0 & \mbox{elsewhere} \end{array} \right.
\end{equation}
Solutions of the eigenvalue problem are shown in Fig.~\ref{fig:mode_couple}(a).  Specifically, there are five modes localized in the potential well created by the index profile $n_0(x)$.  These are the discrete spectra of the problem and the modes of interest in mode-coupling.

Since the operator for the eigenvalue problem is a self-adjoint Sturm-Liouville operator, the eigenfunctions form an orthonormal set.  Thus mode-coupling cannot occur between them.
Mode coupling occurs by either forcing the system via (\ref{eq:mc_force}) or via nonlinearity (\ref{eq:mc_non}).  Figure~\ref{fig:mode_couple}(b) shows the evolution of the solution when the initial condition is given by 
\begin{equation}
   u(x,0)=v_1 (x) .
\end{equation}
Thus the initial condition is the ground state (lowest eigenvalue) of the system.  As the evolution shows, no mode coupling occurs between any of the states.  If the system is forced resonantly between mode one and three such that
\begin{equation}
   \omega= \lambda_1-\lambda_3
\end{equation}
in (\ref{eq:mc_force}), then mode coupling occurs between $v_1(x)$ and $v_3(x)$ as shown in Fig.~\ref{fig:mode_couple}(c).  Note that the forcing is on between $t\in[1000,2000]$ and is then turned off again.  For this simulation, $\epsilon=0.2$ gives a switching time between the two eigenfunctions of approximately $t\approx 1000$.  Once the forcing is turned off, the energy is now constrained to the third eigenstate.  This shows that principled, temporal forcing of a system can move the system from one spatio-temporal state to another.  For the nonlinear perturbatoin, $\epsilon=0.1$ in (\ref{eq:mc_non}) and the dynamics shows that partial coupling occurs between the ground state and other modes.  A detailed analysis of nonlinear mode-coupling in optical waveguides is given by Jensen~\cite{jensen1982nonlinear}.  In either case, mode-coupling can occur in a forced linear system or simply through nonlinear interactions, giving insight into important ways that spatio-temporal systems can switch or couple dominant spatio-temporal patterns of activity.

\newpage
\section*{Lecture 26: Floquet Theory}

Thus far, we have primarily been focused on considering the dynamics, stability and bifurcation structure of many steady-state solutions.  However, there are many nonlinear dynamical systems that produce periodic solutions that persistent and so are often considered stable in nature.  Such limit cycles are now the focus of our efforts.  In particular, our objective is to study the stability of such solutions using what is called Floquet theory.

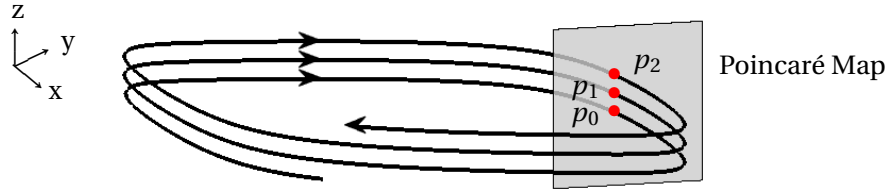
\begin{figure}[t]
\begin{overpic}[width=1\textwidth]{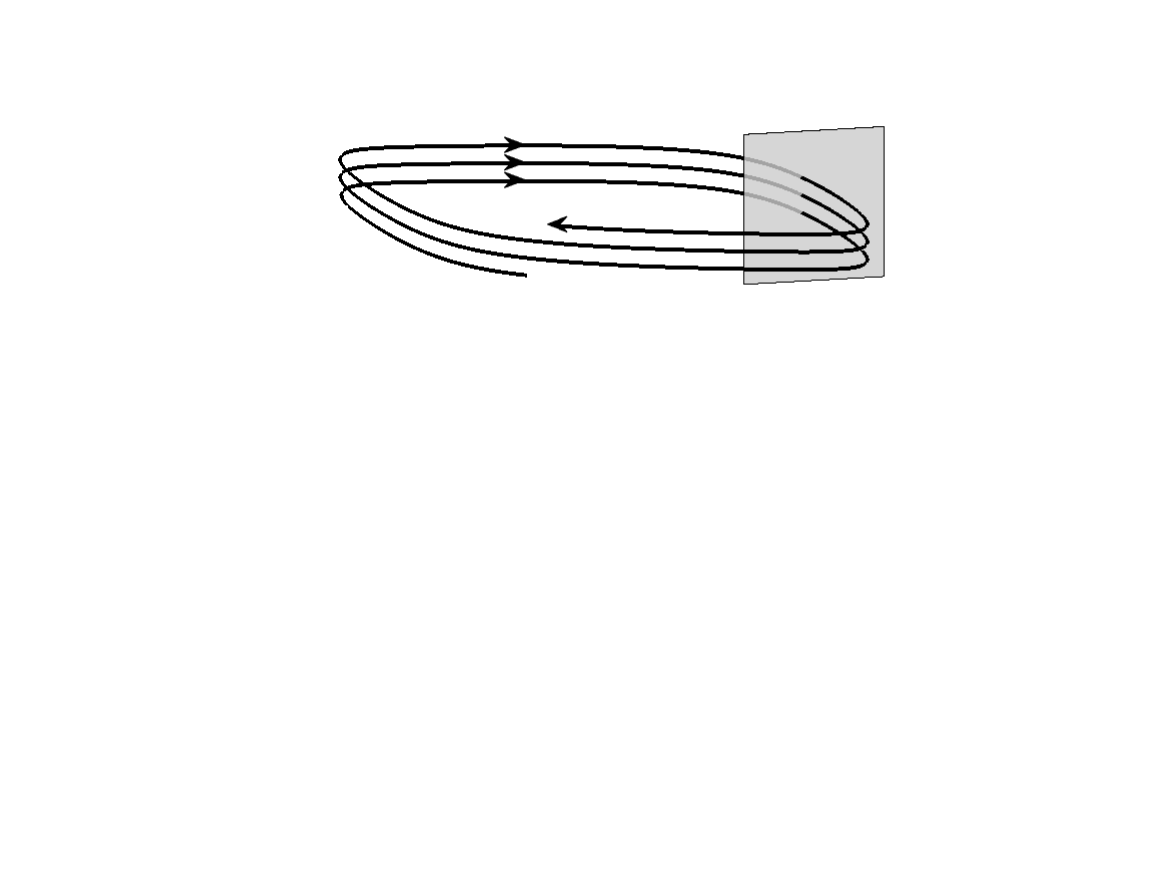}
\put(77,60){Poincar\'e Map}
\put(65,56.2){$p_0$}
\put(65,58.4){$p_1$}
\put(70,60.5){$p_2$}
\put(67.5,57.5){ {\color{red}{\Large \textbullet}} }
\put(67.5,56){ {\color{red}{\Large \textbullet}} }
\put(67.5,59){ {\color{red}{\Large \textbullet}} }
\put(20,60){\rotatebox{25}{$\longrightarrow$}}
\put(19.6,60.4){\rotatebox{-40}{$\longrightarrow$}}
\put(20.0,60.5){\rotatebox{90}{$\longrightarrow$}}
\put(23,58){x}
\put(24,62){y}
\put(20,64.5){z}
\end{overpic}
\vspace*{-3.2in}
\caption{Graphical representation of the Poincar\'e map (first recurrence map or return map).  The 3-dimensional dynamics crosses the 2-dimensional Poincar\'e section (shaded region) once per round trip of the dynamics (red dots).  The dynamics crosses at the points $p_n$, with the dynamics progressing from $p_n$ to $p_{n+1}$.}
 \label{fig:poincare_map}
\end{figure}

To start, we consider a graphical representation of a solution in an $n$-dimensional space which is nearly periodic. In a lower-dimensional space, we define what is called a Poincar\'e map, or a first recurrence map, which is the intersection of a periodic orbit in the state space of a continuous dynamical system with a prescribed lower-dimensional subspace, called the Poincar\'e section, transversal to the dynamical flow of the system.  Figure~\ref{fig:poincare_map} shows a representative example of a Poincar\'e map (return map) for a simple $3$-dimensional flow projected to a $2$-dimensional map.  
In modern applications, such mappings can be used to discover slow-time dynamics of periodic behaviors~\cite{bramburger2020poincare,bramburger2020sparse}. One can also average over the periodic fast times in order to study the slow-scale dynamics~\cite{ding2009operating}.

More formally, we can characterize the Poincar\'e dynamics as a map
\begin{equation}
  p_{n+1} = \phi (p_n) .
\end{equation}
This then allows us to define the concepts of stability in terms of the map $\phi(\cdot)$.  Specifically, we have the following two possibilities
\begin{subeqnarray}
 && \mbox{stable periodic orbit:} \,\,\,\,\,  p_{n+1}=\phi(p_n) \,\, \mbox{where} \,\, p_{n+1}=p_n \\
  && \mbox{unstable periodic orbit:} \,\,\,\,\,  p_{n+1}=\phi(p_n) \,\, \mbox{where} \,\, p_{n+1}\neq p_n, \,\, |p_{n}-p_{n+1}|>|p_{n-1}-p_n|.
\end{subeqnarray}
These are intuitive concepts for stability.  What is required is a principled method to compute the criteria.  This is exactly what Floquet theory allows us to do.

\subsection*{Floquet Theory}

The analysis of the stability of periodic solutions, either temporally or spatially periodic, is known as Floquet theory.   To develop the theory, we will consider the linear, second order differential equation
\begin{equation}
  \frac{d^2 x}{dt^2} + p_1(t)  \frac{d x}{dt}  + p_2 (t) x = 0 
  \label{eq:hill}
\end{equation}
whose coefficients are periodic on a time interval $T$ so that
\begin{subeqnarray}
  &&  p_1(t) = p_1(t+T) \\
  &&  p_2(t) = p_2(t+T) .
\end{subeqnarray}
As a linear, second order equation, there exits two linearly independent solutions.  This allows us to construct the general solution
\begin{equation}
  x= c_1 x_1(t) + c_2 x_2(t)
\end{equation}
where $x_1(t)$ and $x_2(t)$ are the two, linearly-independent fundamental solutions of the differential equations.  The governing equation with $p_1(t)=0$ is known as {\em Hill's equation}.

Periodicity of the coefficients of (\ref{eq:hill}) gives 
\begin{equation}
  \frac{d^2 x (t+T)}{dt^2} + p_1(t+T)  \frac{d x(t+T)}{dt}  + p_2 (t+T) x (t+T) = 0 
\end{equation}
which collapses to 
\begin{equation}
  \frac{d^2 x (t+T)}{dt^2} + p_1(t)  \frac{d x(t+T)}{dt}  + p_2 (t) x (t+T) = 0 
\end{equation}
due the periodicity of $p_1(t)$ and $p_2(t)$.  Thus the solutions
$x_1(t+T)$ and $x_2(t+T)$ must also be fundamental solutions of the governing equation (\ref{eq:hill}).   A linear transform allows us to express these fundamental solutions in terms of the original fundamental solutions so that
\begin{subeqnarray}
 &&  x_1(t+T) = a_{11} x_1(t) + a_{12} x_2(t) \\
&&  x_2(t+T) = a_{21} x_1(t) + a_{22} x_2(t) .
\label{eq:period_map}
\end{subeqnarray}
Or in matrix form with ${\bf x}=[x_1 \,\,\, x_2]^T$, this gives
\begin{equation}
  {\bf x}(t+T) = {\bf A} {\bf x}(t).
\end{equation}
We can now transform this via diagonalization using the eigenvectors of the matrix ${\bf A}$.  Specifically, if we have the eigenvalue decomposition of the matrix ${\bf A}{\bf V}={\bf V}{\bf \Lambda}$, then we can define
\begin{equation}
  {\bf x}={\bf V} {\bf v}
\end{equation}
so that
\begin{equation}
  {\bf v} (t+T) = {\bf \Lambda} {\bf v}(t)
\end{equation}
where ${\bf \Lambda}={\bf V}^{-1} {\bf A}{\bf V}$.  This is now a diagonalized coordinate system since the matrix ${\bf \Lambda}$ is diagonal with eigenvalues (diagonal components) $\lambda_1$ and $\lambda_2$.  This gives the Floquet solutions
\begin{equation}
  {\bf v}_j (t+T) = \lambda_j {\bf v}_j (t)
\end{equation}
for $j=1,2$.  Applying this mapping over $n$ periods gives the iterative mapping
\begin{equation}
    {\bf v}_j (t+nT) = \lambda^n_j {\bf v}_j (t) 
    \label{eq:floquet_exp}
\end{equation}
which upon taking their norms gives $|{\bf v}_j (t+nT)|^2 = |\lambda^n_j|^2 |{\bf v}_j (t)|^2 $.   Thus the stability of periodic solutions can easily be determined from the eigenvalues of the matrix ${\bf A}$.  Specifically, this would give 
\begin{equation}
   \lim_{n\rightarrow \infty}  {\bf v}_j (t) = \left\{ \begin{array}{cc} 0 & |\lambda_j|<1 \\ \infty & |\lambda_j|>1 \end{array}  \right. 
\end{equation}
with the solution actually not growing or decaying if eigenvalues are exactly on the unit circle with  $\lambda_j|=1$.  

The $\lambda_j$ are called the {\em Floquet multipliers}.  Once determined, stability can be evaluated.  The question is:  how do we determine the matrix ${\bf A}$ so that we can compute the eigenvalues.  To do this we again consider (\ref{eq:floquet_exp}) and multiply both sides by
$\exp[-\gamma_j (t+T)]$ which gives
\begin{equation}
    {\bf v}_j (t+nT) \exp[-\gamma_j (t+T)] = \lambda^n_j {\bf v}_j (t) 
 \exp[-\gamma_j (t+T)] = \lambda^n_j \exp[-\gamma_j T] {\bf v}_j (t) 
 \exp[-\gamma_j (t)] .
\end{equation}
Periodic solutions are constructed by choosing
\begin{subeqnarray}
  && \lambda^n_j \exp[-\gamma_j T] =1 \\
  && \phi(t) =  {\bf v}_j (t) \exp[-\gamma_j (t)]
\end{subeqnarray}
so that
\begin{equation}
  \phi(t+T)=\phi(t) .
\end{equation}
Thus $\phi(t)$ is periodic and we have the characteristic exponent $\gamma_j=\ln [\lambda_j]/T$ which is unique to within a multiple of $2 i n\pi/ T$.  Our periodic solutions are then
\begin{equation}
  {\bf v}_j(t) = \exp(\gamma_j t) \phi_j(t)
\end{equation}
where the $\phi_j$ satisfy $\phi_j(t+T)=\phi_j(t)$.

We can determine the characteristic exponents and Floquet multipliers in practice by recalling (\ref{eq:period_map}) and considering the solution at time $t=0$ where we impose
\begin{subeqnarray}
&&  x_1(0)=1, \,\, \frac{dx_1(0)}{dt} =0 \\
&&  x_2(0)=0, \,\, \frac{dx_2(0)}{dt} =1 .
\label{eq:fd_bc}
\end{subeqnarray}
The Wronskian in this case guarantees that these two fundamental solution are linearly independent.   Importantly, inserting these initial conditions into (\ref{eq:period_map}) gives
\begin{subeqnarray}
 &&  x_1(T) = a_{11} x_1(0) + a_{12} x_2(0) \\
&&  x_2(T) = a_{21} x_1(0) + a_{22} x_2(0) 
\end{subeqnarray}
with then $x_1(T)=a_{11}$ and $x_2(T)=a_{21}$.  Differentiating these equations and inserting the initial conditions gives also
$x_1'(T)=a_{12}$ and $x_2'(T)=a_{22}$.  Thus we have in total
\begin{equation}
   {\bf A} = \left[ \begin{array}{cc} a_{11} & a_{12} \\ a_{21} & a_{22} \end{array} \right] = \left[ \begin{array}{cc} x_1(T) & x_1'(T)  \\ x_2(T) & x_2'(T) \end{array} \right] .
\end{equation}
The Floquet multipliers (eigenvalues) of this matrix can be found to be determined from
\begin{equation}
   \lambda^2 - \left[  x_1(T)+x_2'(T) \right] \lambda
   + \left[  x_1(T) x_2'(T) - x_1'(T)x_2(T) \right] = 0
\end{equation}
where $\Delta(T)= \left[  x_1(T) x_2'(T) - x_1'(T)x_2(T) \right] $ is the Wronskian between the solutions $x_1(T)$ and $x_2(T)$.  Since the Wronskian is constant, we then have
\begin{equation}
  \Delta(T)=\Delta(0)= \left[  x_1(0) x_2'(0) - x_1'(0)x_2(0) \right] =1 
\end{equation}
which makes the eigenvalue equation
\begin{equation}
  \lambda^2 - 2\alpha \lambda +1 = 0
\end{equation}
with $\alpha= (x_1(T)+x_2'(T))/2$.  The eigenvalues can be trivially computed to be
\begin{equation}
  \lambda_\pm=\alpha\pm \sqrt{\alpha^2 -1 }
\end{equation}
with $\lambda_+ \lambda_- =1$.  We can conclude the following:\\

\noindent (i) $|\alpha|>1$:  There is one bounded and one unbounded root and the solutions are unstable.\\

\noindent (ii) $|\alpha|<1$:  The eigenvalue are complex conjugate pairs that reside on the unit circle since $\lambda_+ \lambda_- =1$.  This is a stable solution.\\

The requirement that $|\alpha|<1$ gives the  the condition
\begin{equation}
  |\Gamma|=| x_1(T)+x_2'(T) | < 2
\end{equation}
where $\Gamma$ is the Floquet discriminant.  Thus if we can find the Floquet discriminant to be less than the absolute value of two, the periodic solution is stable.  This is one of the fundamental and practical results of Floquet theory.

\newpage
\section*{Lecture 27:  The Pendulum and Floquet Theory}

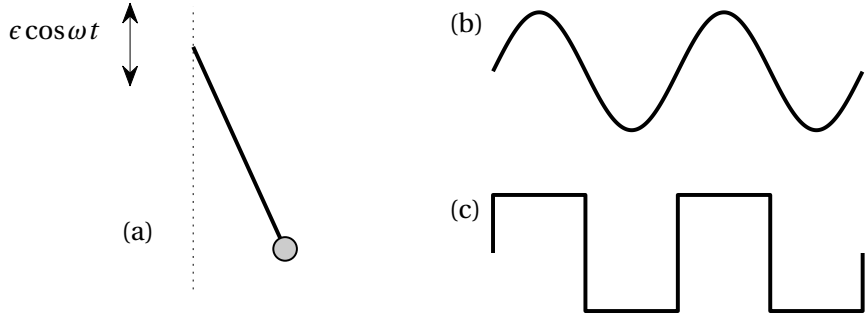
\begin{figure}[t]
\begin{overpic}[width=0.6\textwidth]{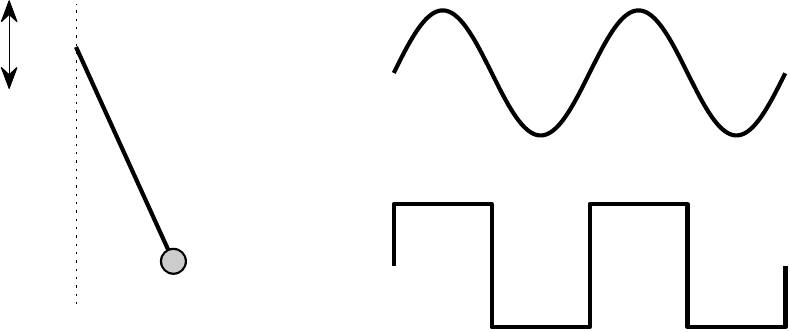}
\put(0,10){(a)}
\put(44,38){(b)}
\put(44,13){(c)}
\put(-15,37){$\epsilon \cos \omega t$}
\end{overpic}
\caption{(a) Schematic of a pendulum whose support is driven by a periodic forcing $\epsilon \cos \omega t$. We will show that for sufficiently high driving frequency $\omega$, the pendulum can be stabilized in the inverted position.  (b)  The driving is sinusoidal, but the analytic theory approximates this by a piecewise linear forcing (c).  This allows for an explicit computation of the Floquet discriminant.}
 \label{fig:invert}
\end{figure}

As an application of Floquet theory, we consider the classic problem of the periodically forced pendulum.  Specifically, the pendulum is considered to be forced by its support with a frequency $\omega$.  Figure~\ref{fig:invert} shows the periodically forced pendulum and its oscillating support.  The equations of motion can be derived by standard free-body diagrams and $F=ma$ physics principles to yield
\begin{equation}
  \frac{d^2 x}{dt^2}  + ( \delta + \epsilon \cos \omega t ) \sin x = 0
\end{equation}
where $x(t)$ is the angular position of the pendulum, $\delta$ measures the natural oscillation frequency of the pendulum, $\omega$ is the frequency of the oscillating support, and $\epsilon$ is the magnitude of the oscillations of the support.  Note that this is Hill's equation in the small amplitude oscillation limit where $\sin x \approx x$, $p_1(t)=0$ and $p_2(t)=\delta + \epsilon \cos\omega t$.

By approximating the pendulum in the down or inverted position, we can produce analytic solutions for the Floquet theory.  These two important positions are achieved when $x\approx 0$ (down position) and $x\approx \pi$ (inverted position).  The two cases produce the approximation
\begin{subeqnarray}
 &&  x\approx 0:  \,\,\,\,\, \sin x \approx x - \frac{x^3}{3!} + \cdots \\
 && x\approx \pi:  \,\,\,\,\, \sin (x+\pi) =-\sin x \approx -x + \frac{x^3}{3!} + \cdots
\end{subeqnarray}
which can be used to approximate the pendulum in the down and upright position
\begin{equation}
  \frac{d^2 x}{dt^2}  \pm ( \delta + \epsilon \cos \omega t ) x = 0
  \label{eq:matthieu}
\end{equation}
which is known as Matthieu's equation.\\

\noindent {\em Down Pendulum:}  Consideration of the pendulum in the down position gives (\ref{eq:matthieu}) with the positive sign.  In this case, a great deal can be said about the governing behavior.  Specifically, if the forcing frequency is resonant with the natural frequency then linear oscillations can grow unbounded.  A proper consideration in this case is given by an asymptotic treatment of the Duffing equation in the asymptotics chapter.  This shows that the nonlinearity and any damping suppresses the unbounded growth and also generates a frequency shift that can be well characterized with theory.\\

\noindent {\em Inverted Pendulum:}  This is the case we consider here for 
(\ref{eq:matthieu}) with the negative sign.  Specifically, we will show that you can stabilize the pendulum in the inverted position if you drive the forcing frequency $\omega$ at a high enough value.  This is a surprising, yet quite consequential result that arises directly from mathematical considerations.
Shifting time by $t\rightarrow t+\pi/2$ gives the governing equation in this case
\begin{equation}
  \frac{d^2 x}{dt^2}  - ( \delta + \epsilon \sin \omega t ) x = 0
\end{equation}
where the $\sin\omega t$ term take the place of the $\cos\omega t$ when the time is shifted by $\pi/2$.

For analytic purposes, one further assumption will be made concerning an approximation for the $\sin \omega t$ forcing.  Specifically the sinusoidal function will be replaced by a step-wise constant function with the same period and amplitude of modulation.  Figure~\ref{fig:invert} shows the approximation to be made.  This then gives the governing equaions
\begin{subeqnarray}
 && \frac{d^2 x}{dt^2}  - ( \delta + \epsilon ) x = 0 \hspace*{.5in} 0<t<T/2 \\
 && \frac{d^2 x}{dt^2}  - ( \delta - \epsilon ) x = 0  \hspace*{.5in}  T/2<t<T .
 \label{eq:piecewise_pendulum}
\end{subeqnarray}
This approximation will allow us to analytically compute the Floquet discriminant for the inverted pendulum.

Recall that the Floquet discriminant requires the calculation of the following quantity
\begin{equation}
  \Gamma= x_1(T) + x_2'(T)
  \label{eq:fd_gamma}
\end{equation}
where the solutions $x_1(t)$ and $x_2(t)$ satisfy the initial conditions (\ref{eq:fd_bc}).  For the piecewise constant approximation (\ref{eq:piecewise_pendulum}), the solutions are given by
\begin{subeqnarray}
  &&  x(t) = c_1 \exp \left[ \sqrt{\delta+\epsilon} t \right] 
  + c_2 \exp \left[ -\sqrt{\delta+\epsilon} t \right]   \hspace*{.5in} 0<t<T/2 \\
  && x(t) = c_3 \exp \left[ \sqrt{\delta-\epsilon} t \right] 
  + c_4 \exp \left[ -\sqrt{\delta-\epsilon} t \right]  \hspace*{.5in}  T/2<t<T .
  \label{eq:piece_wise_ic}
\end{subeqnarray}
With these solutions, there are two particular solutions that need to be computed:  $x_1(t)$ and $x_2(t)$.  These two different solutions satisfy 
(\ref{eq:piecewise_pendulum}).  To solve this, the initial conditions are first applied to the solutions (\ref{eq:piece_wise_ic}a).  This solution is then evaluated at $t=T/2$ in order to produce initial conditions for 
(\ref{eq:piece_wise_ic}b).  After quite a bit of algebra, the two following solutions are found
\begin{subeqnarray}
  &&  \hspace*{-.8in} x_1\!=\! \cosh \left[ \sqrt{\delta+\epsilon} \frac{T}{2} \right]
      \cosh \left[ \sqrt{\delta-\epsilon} \left( t- \frac{T}{2} \right) \right]
      \!+\! \frac{\sqrt{\delta+\epsilon}}{\sqrt{\delta-\epsilon}}
      \sinh \left[ \sqrt{\delta+\epsilon} \frac{T}{2} \right]
      \sinh \left[ \sqrt{\delta-\epsilon} \left( t- \frac{T}{2} \right) \right] \\
      &&  \hspace*{-.8in} x_2\!=\! 
      \frac{1}{\sqrt{\delta-\epsilon}}
      \cosh \left[ \sqrt{\delta+\epsilon} \frac{T}{2} \right]
      \sinh \left[ \sqrt{\delta-\epsilon} \left( t- \frac{T}{2} \right) \right]
      \!+\! \frac{1}{\sqrt{\delta+\epsilon}}
      \sinh \left[ \sqrt{\delta+\epsilon} \frac{T}{2} \right]
      \cosh \left[ \sqrt{\delta-\epsilon} \left( t- \frac{T}{2} \right) \right]
\end{subeqnarray}
which is sufficient to compute the Floquet discriminant (\ref{eq:fd_gamma})
\begin{equation}
  \Gamma=2 \cosh \left[ \sqrt{\delta+\epsilon} \frac{T}{2} \right]
      \cosh \left[ \sqrt{\delta-\epsilon} \frac{T}{2} \right]
      + \left(  \frac{\sqrt{\delta+\epsilon}}{\sqrt{\delta-\epsilon}}
      + \frac{\sqrt{\delta-\epsilon}}{\sqrt{\delta+\epsilon}}  \right)
      \sinh \left[ \sqrt{\delta+\epsilon} \frac{T}{2} \right]
      \sinh \left[ \sqrt{\delta-\epsilon} \frac{T}{2} \right] .
\end{equation}
Finally, we recall that the forcing frequency can be related to the period by $T=2\pi/\omega$.  This allows us to express the Floquet discriminant in terms of the driving frequency $\Gamma = \Gamma(\omega)$.  

If $\delta>\epsilon$, then as $\omega\rightarrow \infty$, then $\Gamma=2$ and the Floquet discriminant shows the system to be on the stability boundary (neural stability).  If $\delta<\epsilon$, then $|\Gamma|<2$ for windows in $\omega$ since
\begin{equation}
   \Gamma\rightarrow 2\cos \left[ \sqrt{\epsilon-\delta} \frac{\pi}{\omega} \right] \hspace*{.5in} \omega\rightarrow \infty .
\end{equation}
Thus $|\Gamma|\leq 2$ for all sufficiently large $\omega$ showing that as the driving frequency is increased, the pendulum can be stabilized in the inverted position.   

\subsection*{The Floquet discriminant for the pendulum}
What has been proposed is a method by which one can compute the stability of periodic solutions.  We will do this computationally now for the pendulum for both the down and inverted pendulum.  All that is required is to simulate the system from time $t\in[0,T]$ where $T=2\pi/\omega$.  In particular, we do two simulations do determine $x_1(t)$ and $x_2(t)$ with the initial conditions (\ref{eq:fd_bc}).

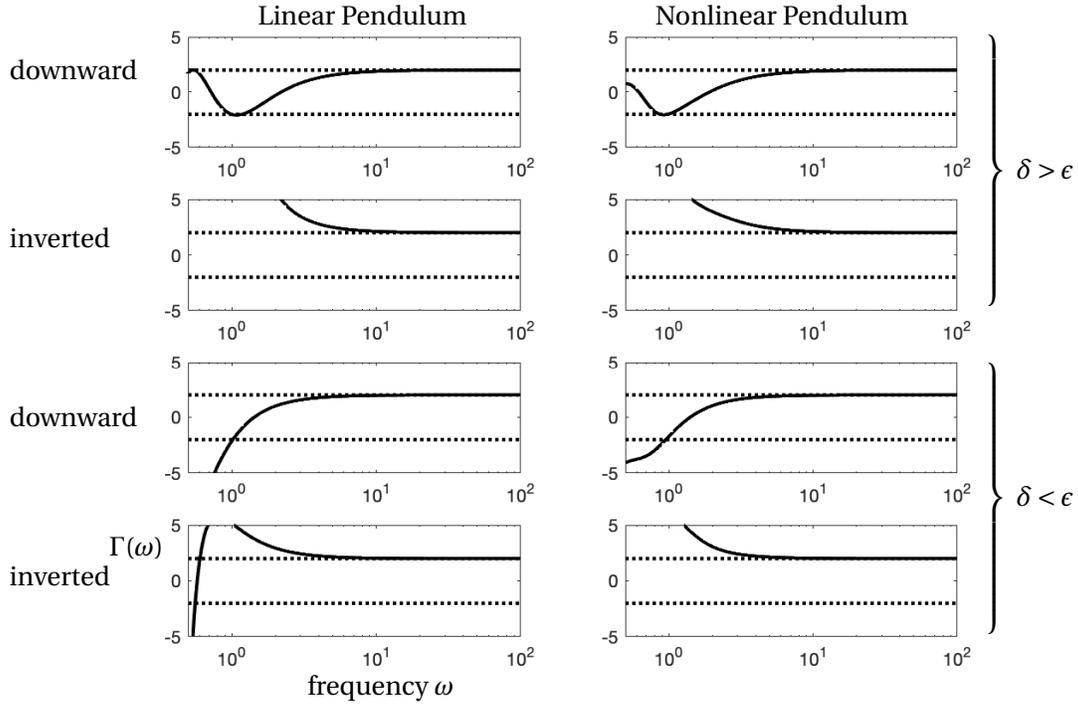
\begin{figure}[t]
\begin{overpic}[width=0.8\textwidth]{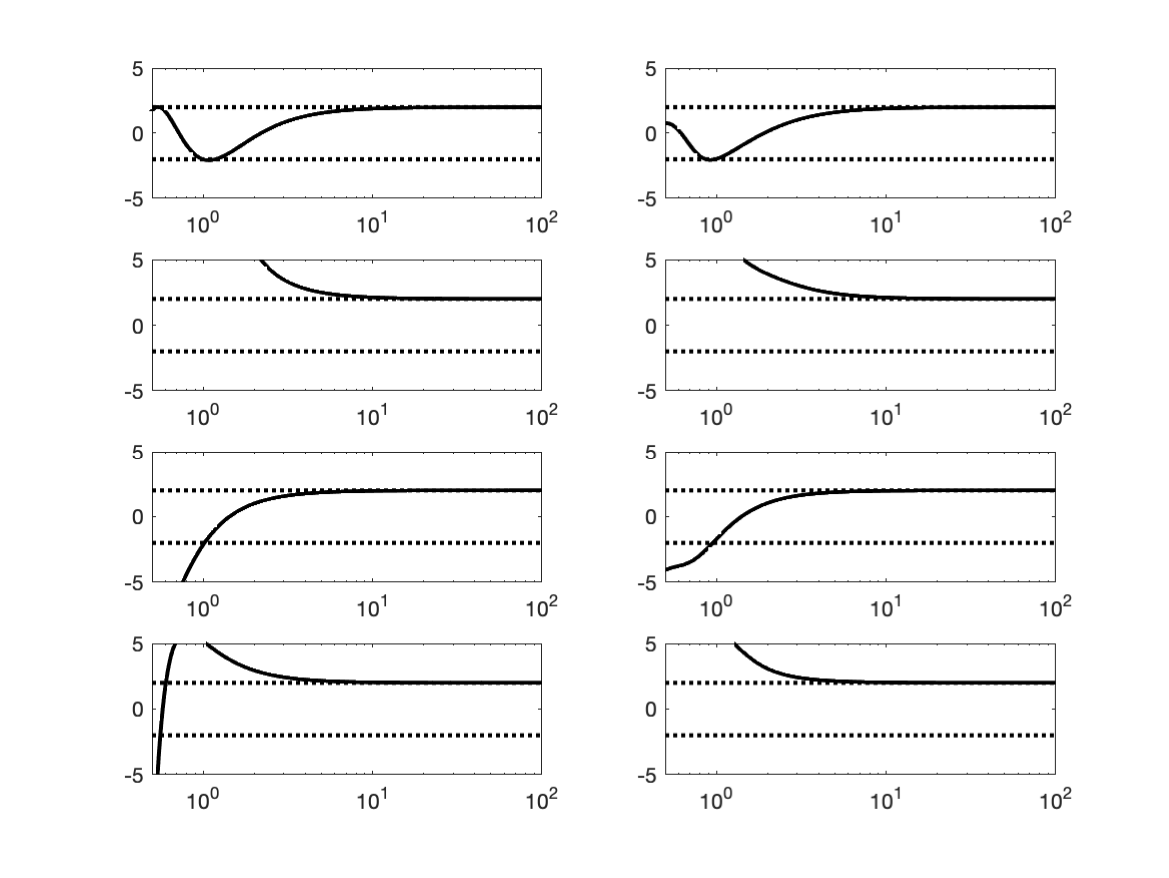}
\put(20,70.5){Linear Pendulum}
\put(60,70.5){Nonlinear Pendulum}
\put(5,17){$\Gamma(\omega)$}
\put(25,3){frequency $\omega$}
\put(-5,65){downward}
\put(-5,48){inverted}
\put(-5,30){downward}
\put(-5,14){inverted}
\put(90,55){$ \left. \begin{array}{c} \\ \\ \\ \\ \\ \\ \\ \\ \end{array} \right\} \,\,\, \delta>\epsilon$ }
\put(90,22){$ \left. \begin{array}{c} \\ \\ \\ \\ \\ \\ \\ \\ \end{array} \right\} \,\,\, \delta<\epsilon$}
\end{overpic}
\caption{Floquet discriminant for the linear (left column) and nonlinear (right column) pendulum for the case of $\delta>\epsilon$ and $\delta<\epsilon$ for both the downward and inverted pendulum. The dotted lines denote the stability boundaries where $\Gamma=\pm 2$. The downward pendulum (top row) is always stable when $\delta>\epsilon$ as expected.  However, the downward pendulum can be destabilized for $\delta<\epsilon$ and for sufficiently low frequency forcing.  The inverted pendulum can be stabilized for sufficiently high values of frequency $\omega$.  }
 \label{fig:discriminant_pendulum}
\end{figure}

We consider the Floquet analysis for both the downward pendulum and the inverted pendulum, in both the linear pendulum approximation and the full nonlinear pendulum dynamics.  For the pendulum near the down position, the two governing models are given by
\begin{subeqnarray}
 &&   \frac{d^2 x}{dt^2}  + ( \delta + \epsilon \sin \omega t ) x = 0 \\
 &&   \frac{d^2 x}{dt^2}  + ( \delta + \epsilon \sin \omega t ) \sin x = 0 .
\end{subeqnarray}
The Floquet discriminant is computed for both models.  It is expected that since the pendulum is in the downward position, it should be stable.  This is true when $\delta>\epsilon$.  Figure~\ref{fig:discriminant_pendulum} show the Floquet discriminant $\Gamma(\omega)$ for the downward pendulum in various scenarios.

The inverted pendulum has the two governing models
\begin{subeqnarray}
 &&   \frac{d^2 x}{dt^2}  - ( \delta + \epsilon \sin \omega t ) x = 0 \\
 &&   \frac{d^2 x}{dt^2}  - ( \delta + \epsilon \sin \omega t ) \sin x = 0 .
\end{subeqnarray}
As with the downward pendulum, 
the Floquet discriminant is computed for both models.  It is expected that since the pendulum is in the inverted position, it should be unstable.  However, as shown in   Fig.~\ref{fig:discriminant_pendulum}, the Floquet discriminant $\Gamma(\omega)$ computed drops to $|\Gamma|=2$ when $\omega\rightarrow \infty$, indicated that the inverted pendulum can be stabilized for a sufficiently high forcing frequency $\omega$.

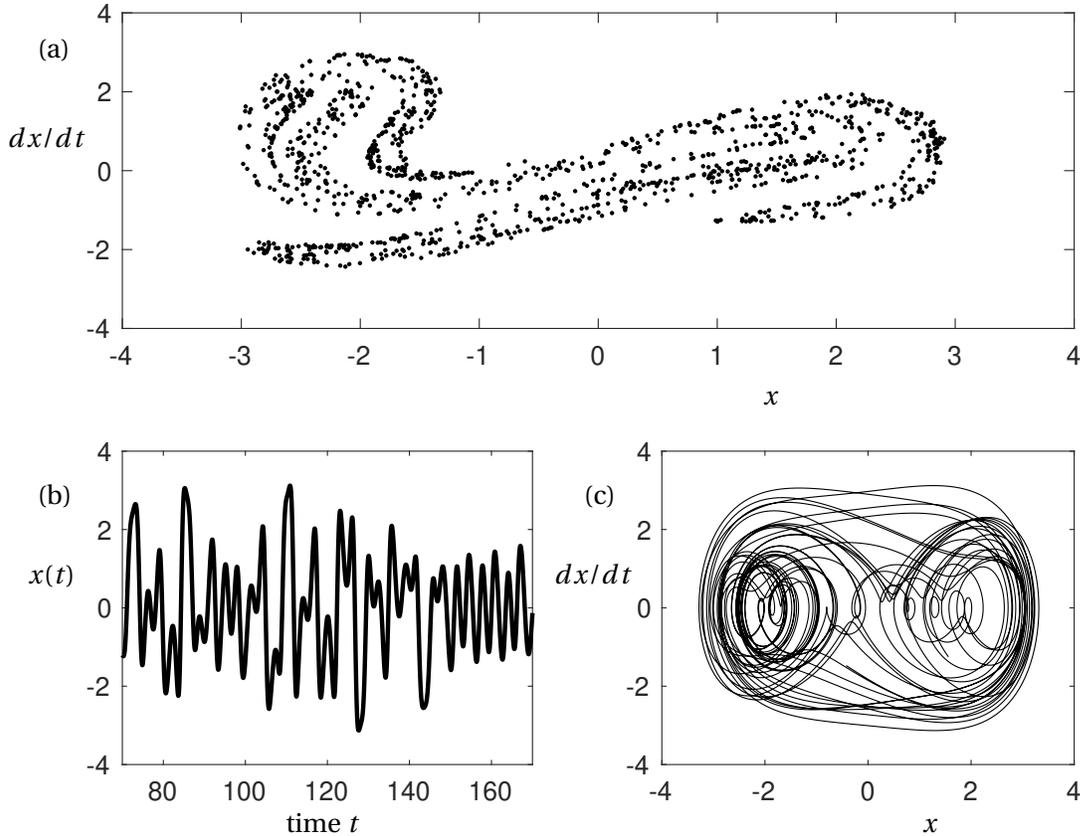
\begin{figure}[t]
\begin{overpic}[width=0.8\textwidth]{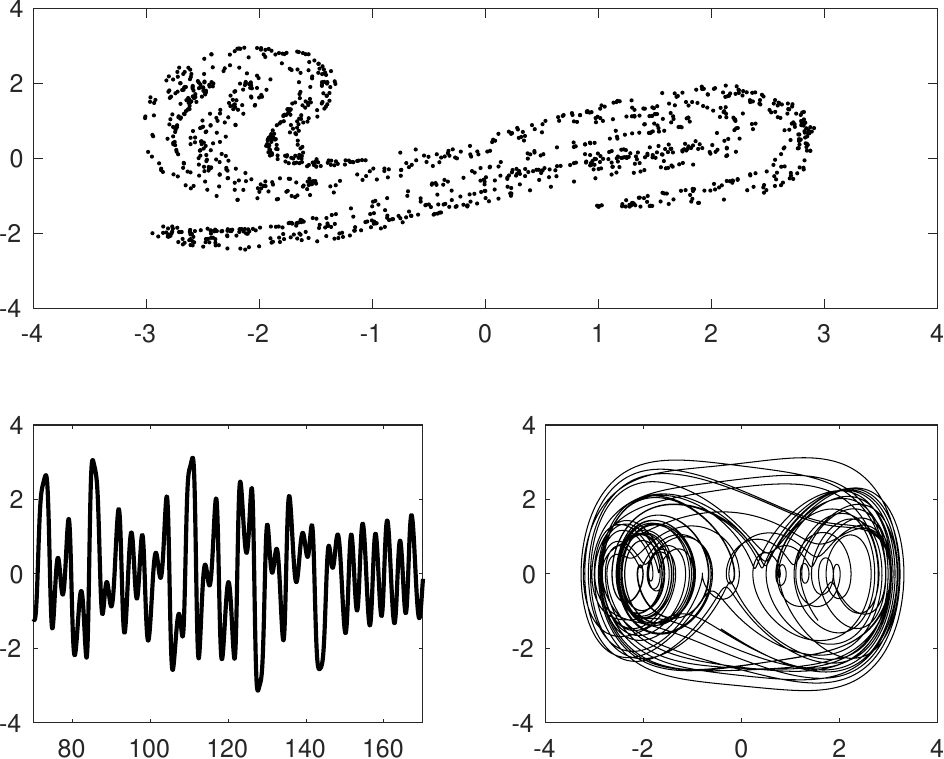}
\put(-5,75){(a)}
\put(-5,30){(b)}
\put(50,30){(c)}
\put(-8,66){${dx}/{dt}$}
\put(47,22){${dx}/{dt}$}
\put(68,40){$x$}
\put(-6,22){$x(t)$}
\put(20,-3){time $t$}
\put(84,-3){$x$}
\end{overpic}
\caption{(a) Poincar\'e map of the forced Duffing equation.  Note the chaotic behavior of the system that is illustrated in the map dynamics.  (b) A portion of the time-series $x(t)$ showing the irregular, but approximately quasi-periodic, behavior.  (c) The phase-plane portrait of the chaotic dynamics.  The parameters of the Duffing model (\ref{eq:duffing_chaos}) are given by $\omega=1$, $\delta=0.1$, $\kappa=1/4$ and $\gamma=1.5$.}
 \label{fig:duffing_chaos}
\end{figure}

\subsection*{Poincare Maps of the Forced Duffing Equation}

In addition to periodic phenomena, the Floquet theory analysis considers the dynamics on the Poincar\'e section as illustrated in Fig.~\ref{fig:poincare_map}.  In the following example, we will illustrate the chaotic dynamics that can emerge on such maps.  Specifically, we will consider the damped-driven Duffing equation
\begin{equation}
   \frac{d^2 x}{dt^2} +\delta \frac{dx}{dt}  -  u + \kappa u^3 = \gamma \cos \omega t .
\label{eq:duffing_chaos} 
\end{equation}
This system will be studied in detail in the asymptotics chapter for a different sign of the linear and nonlinear terms in $u$.  Here, we study the driving dynamics that occurs when forced with the driving term $\gamma \cos \omega t$.  The driving frequency $\omega$ will determine the Poincar\'e map sampling period $T=2\pi/\omega$. 
Damped-driven systems of this kind are common in many physical systems~\cite{akhmediev2008dissipative}, such as mode-locked lasers~\cite{spaulding2002nonlinear,kim2000pulse,kutz2008passive,kutz1998noise} and detonation engines~\cite{koch2020mode,koch2020multi,koch2020modeling}

Figure~\ref{fig:duffing_chaos} shows a number of features of the dynamics.  Most importantly, the Poincar\'e section is illustrated showing the chaotic dynamics induced by the forcing.  Indeed, many damped-driven systems can produce  remarkably complex, chaotic behavior when nonlinearity plays a role in the underlying dynamics.  The Poincar\'e map provides a powerful analysis tool for understanding the chaotic mapping that occurs on the cycle of a period of the forcing.

\newpage

\bibliographystyle{abbrv}
\bibliography{ALLBIB}

\begin{thebibliography}{10}

\bibitem{abramowitz1948handbook}
M.~Abramowitz and I.~A. Stegun.
\newblock {\em Handbook of mathematical functions with formulas, graphs, and
  mathematical tables}, volume~55.
\newblock US Government printing office, 1948.

\bibitem{akhmediev2008dissipative}
N.~Akhmediev and A.~Ankiewicz.
\newblock {\em Dissipative solitons: from optics to biology and medicine},
  volume 751.
\newblock Springer Science \& Business Media, 2008.

\bibitem{bender2013advanced}
C.~M. Bender and S.~A. Orszag.
\newblock {\em Advanced mathematical methods for scientists and engineers I:
  Asymptotic methods and perturbation theory}.
\newblock Springer Science \& Business Media, 2013.

\bibitem{benney1967propagation}
D.~Benney and A.~Newell.
\newblock The propagation of nonlinear wave envelopes.
\newblock {\em Journal of mathematics and Physics}, 46(1-4):133--139, 1967.

\bibitem{boyd2003nonlinear}
R.~W. Boyd.
\newblock {\em Nonlinear optics}.
\newblock Elsevier, 2003.

\bibitem{bramburger2020sparse}
J.~J. Bramburger, D.~Dylewsky, and J.~N. Kutz.
\newblock Sparse identification of slow timescale dynamics.
\newblock {\em arXiv preprint arXiv:2006.00940}, 2020.

\bibitem{bramburger2020poincare}
J.~J. Bramburger and J.~N. Kutz.
\newblock Poincar{\'e} maps for multiscale physics discovery and nonlinear
  floquet theory.
\newblock {\em Physica D: Nonlinear Phenomena}, page 132479, 2020.

\bibitem{derv1}
G.~Bukauskas, V.~Kabelka, A.~Piskarskas, and A.~Y. Stabinis.
\newblock Features of three-photon parametric interaction of ultrashort light
  packets in the nonlinear amplification regime.
\newblock {\em Soviet Journal of Quantum Electronics}, 4(3):290, 1974.

\bibitem{callaham2021learning}
J.~L. Callaham, J.~V. Koch, B.~W. Brunton, J.~N. Kutz, and S.~L. Brunton.
\newblock Learning dominant physical processes with data-driven balance models.
\newblock {\em Nature communications}, 12(1):1016, 2021.

\bibitem{cojocaru2009parabolic}
E.~Cojocaru.
\newblock Parabolic cylinder functions implemented in matlab.
\newblock {\em arXiv preprint arXiv:0901.2220}, 2009.

\bibitem{cross1993pattern}
M.~C. Cross and P.~C. Hohenberg.
\newblock Pattern formation outside of equilibrium.
\newblock {\em Reviews of modern physics}, 65(3):851, 1993.

\bibitem{derv2}
R.~Danelyus, G.~Dikchyus, V.~Kabelka, A.~Piskarskas, A.~Stabinis, and
  Y.~Yasevichyute.
\newblock Parametric excitation of light in the picosecond rang.
\newblock {\em Soviet Journal of Quantum Electronics}, 7(11):1360, 1977.

\bibitem{ding2009operating}
E.~Ding and J.~N. Kutz.
\newblock Operating regimes, split-step modeling, and the haus master
  mode-locking model.
\newblock {\em Journal of the Optical Society of America B}, 26(12):2290--2300,
  2009.

\bibitem{elgin1993perturbations}
J.~Elgin.
\newblock Perturbations of optical solitons.
\newblock {\em Physical Review A}, 47(5):4331, 1993.

\bibitem{Friedman}
B.~Friedman.
\newblock {\em Principles and techniques of applied mathematics}.
\newblock Courier Dover Publications, 1990.

\bibitem{hewitt}
S.~E. Hewitt, K.~Intrachat, and J.~N. Kutz.
\newblock Dynamics and stability of a new class of periodic solutions of the
  optical parametric oscillator.
\newblock {\em Optics communications}, 240(4-6):423--436, 2004.

\bibitem{hewitt2005dynamics}
S.~E. Hewitt and J.~N. Kutz.
\newblock Dynamics of the optical parametric oscillator near resonance
  detuning.
\newblock {\em SIAM Journal on Applied Dynamical Systems}, 4(4):808--831, 2005.

\bibitem{jensen1982nonlinear}
S.~Jensen.
\newblock The nonlinear coherent coupler.
\newblock {\em IEEE Journal of Quantum Electronics}, 18(10):1580--1583, 1982.

\bibitem{Kap2}
T.~Kapitula.
\newblock Bifurcating bright and dark solitary waves for the perturbed
  cubic-quintic nonlinear schr{\"o}dinger equation.
\newblock {\em Proceedings of the Royal Society of Edinburgh Section A:
  Mathematics}, 128(3):585--629, 1998.

\bibitem{kapitula2013spectral}
T.~Kapitula and K.~Promislow.
\newblock {\em Spectral and dynamical stability of nonlinear waves}, volume
  185.
\newblock Springer, 2013.

\bibitem{kaup}
D.~Kaup.
\newblock Perturbation theory for solitons in optical fibers.
\newblock {\em Physical Review A}, 42(9):5689, 1990.

\bibitem{kevorkian2013perturbation}
J.~Kevorkian and J.~D. Cole.
\newblock {\em Perturbation methods in applied mathematics}, volume~34.
\newblock Springer Science \& Business Media, 2013.

\bibitem{kim2000pulse}
A.~Kim, J.~Kutz, and D.~Muraki.
\newblock Pulse-train uniformity in optical fiber lasers passively mode-locked
  by nonlinear polarization rotation.
\newblock {\em IEEE Journal of Quantum electronics}, 36(4):465--471, 2000.

\bibitem{koch2020mode}
J.~Koch, M.~Kurosaka, C.~Knowlen, and J.~N. Kutz.
\newblock Mode-locked rotating detonation waves: Experiments and a model
  equation.
\newblock {\em Physical Review E}, 101(1):013106, 2020.

\bibitem{koch2020multi}
J.~Koch, M.~Kurosaka, C.~Knowlen, and J.~N. Kutz.
\newblock Multi-scale physics of rotating detonation engines: Autosolitons and
  modulational instabilities.
\newblock {\em arXiv preprint arXiv:2003.06655}, 2020.

\bibitem{koch2020modeling}
J.~Koch and J.~N. Kutz.
\newblock Modeling thermodynamic trends of rotating detonation engines.
\newblock {\em Physics of Fluids}, 32(12):126102, 2020.

\bibitem{kodama}
Y.~Kodama and M.~J. Ablowitz.
\newblock Perturbations of solitons and solitary waves.
\newblock {\em Studies in Applied Mathematics}, 64(3):225--245, 1981.

\bibitem{kutz2008passive}
J.~N. Kutz.
\newblock Passive mode-locking using phase-sensitive amplification.
\newblock {\em Physical Review A}, 78(1):013845, 2008.

\bibitem{Kutz2016book}
J.~N. Kutz, S.~L. Brunton, B.~W. Brunton, and J.~L. Proctor.
\newblock {\em Dynamic Mode Decomposition: Data-Driven Modeling of Complex
  Systems}.
\newblock SIAM, 2016.

\bibitem{kutz1997nonlinear}
J.~N. Kutz, B.~J. Eggleton, J.~B. Stark, and R.~E. Slusher.
\newblock Nonlinear pulse propagation in long-period fiber gratings: theory and
  experiment.
\newblock {\em IEEE Journal of Selected Topics in Quantum Electronics},
  3(5):1232--1245, 1997.

\bibitem{KE}
J.~N. Kutz, T.~Erneux, S.~Trillo, and M.~Haelterman.
\newblock Curvature dynamics and stability of topological solitons in the
  optical parametric oscillator.
\newblock {\em JOSA B}, 16(11):1936--1941, 1999.

\bibitem{kutz2014solitons}
J.~N. Kutz and E.~Farnum.
\newblock Solitons and ultra-short optical waves: the short-pulse equation
  versus the nonlinear schr{\"o}dinger equation.
\newblock {\em Edited by Hugo E. Hern{\'a}ndez-Figueroa, Erasmo Recami, and},
  page 148, 2014.

\bibitem{Kutz1}
J.~N. Kutz, C.~V. Hile, W.~L. Kath, R.-D. Li, and P.~Kumar.
\newblock Pulse propagation in nonlinear optical fiber lines that employ
  phase-sensitive parametric amplifiers.
\newblock {\em JOSA B}, 11(10):2112--2123, 1994.

\bibitem{Kutz}
J.~N. Kutz and W.~L. Kath.
\newblock Stability of pulses in nonlinear optical fibers using phase-sensitive
  amplifiers.
\newblock {\em SIAM Journal on Applied Mathematics}, 56(2):611--626, 1996.

\bibitem{Kutz2}
J.~N. Kutz, W.~L. Kath, R.-D. Li, and P.~Kumar.
\newblock Long-distance pulse propagation in nonlinear optical fibers by using
  periodically spaced parametric amplifiers.
\newblock {\em Optics letters}, 18(10):802--804, 1993.

\bibitem{kutz1998noise}
J.~N. Kutz and P.~Wai.
\newblock Noise-induced amplitude and chirp jitter in dispersion-managed
  soliton communications.
\newblock {\em Optics letters}, 23(13):1022--1024, 1998.

\bibitem{mahmud2002bose}
K.~Mahmud, J.~Kutz, and W.~Reinhardt.
\newblock Bose-einstein condensates in a one-dimensional double square well:
  Analytical solutions of the nonlinear schr{\"o}dinger equation.
\newblock {\em Physical Review A}, 66(6):063607, 2002.

\bibitem{murdock1999perturbations}
J.~A. Murdock.
\newblock {\em Perturbations: theory and methods}.
\newblock SIAM, 1999.

\bibitem{Newell}
A.~C. Newell and J.~A. Whitehead.
\newblock Finite bandwidth, finite amplitude convection.
\newblock {\em Journal of Fluid Mechanics}, 38(2):279--303, 1969.

\bibitem{Lugiato}
G.-L. Oppo, M.~Brambilla, and L.~A. Lugiato.
\newblock Formation and evolution of roll patterns in optical parametric
  oscillators.
\newblock {\em Physical Review A}, 49(3):2028, 1994.

\bibitem{proctor2005passive}
J.~L. Proctor and J.~N. Kutz.
\newblock Passive mode-locking by use of waveguide arrays.
\newblock {\em Optics letters}, 30(15):2013--2015, 2005.

\bibitem{Segel}
L.~A. Segel.
\newblock Distant side-walls cause slow amplitude modulation of cellular
  convection.
\newblock {\em Journal of Fluid Mechanics}, 38(1):203--224, 1969.

\bibitem{spaulding2002nonlinear}
K.~M. Spaulding, D.~H. Yong, A.~D. Kim, and J.~N. Kutz.
\newblock Nonlinear dynamics of mode-locking optical fiber ring lasers.
\newblock {\em JOSA B}, 19(5):1045--1054, 2002.

\bibitem{Trillo2}
S.~Trillo and M.~Haelterman.
\newblock Excitation and bistability of self-trapped signal beams in optical
  parametric oscillators.
\newblock {\em Optics letters}, 23(19):1514--1516, 1998.

\bibitem{Trillo1}
S.~Trillo, M.~Haelterman, and A.~Sheppard.
\newblock Stable topological spatial solitons in optical parametric
  oscillators.
\newblock {\em Optics letters}, 22(13):970--972, 1997.

\bibitem{weinstein}
M.~I. Weinstein.
\newblock Modulational stability of ground states of nonlinear schr{\"o}dinger
  equations.
\newblock {\em SIAM journal on mathematical analysis}, 16(3):472--491, 1985.

\bibitem{zhang2015semiconductor}
X.~Zhang, M.~Williams, S.~T. Cundiff, and J.~N. Kutz.
\newblock Semiconductor diode laser mode-locking by a waveguide array.
\newblock {\em IEEE Journal of Selected Topics in Quantum Electronics},
  22(2):34--39, 2015.

\end{thebibliography}

\end{document}